\documentclass[10pt]{preprint}

\usepackage[margin=3 cm]{geometry}
\usepackage[utf8x]{inputenc}
\usepackage[english]{babel}
\usepackage[full]{textcomp}
\usepackage[osf]{newtxtext}
\usepackage{amssymb}
\usepackage{amsmath}
\usepackage{amsthm}
\usepackage{breakurl}
\usepackage{mhequ}
\usepackage{booktabs}
\usepackage{tikz}
\usepackage{mathrsfs}
\usepackage{microtype}
\usepackage{slashed}
\usepackage{mathtools}
\usepackage{centernot}
\usepackage{footnote}
\usepackage{enumerate}
\usepackage[shortlabels]{enumitem}
\usepackage{stackrel}
\usepackage{longtable}
\usepackage{cprotect}
\usepackage{xstring}
\usepackage{calc}
\usepackage{bbm}
\usepackage[colorlinks=true, pdfstartview=FitV, linkcolor=colorLink, citecolor=colorCite, urlcolor=colorLink, linktocpage=true]{hyperref}
\usepackage{orcidlink}
\usepackage{esint}
\usepackage{stmaryrd}

\DeclareSymbolFont{timesoperators}{T1}{ptm}{m}{n}
\SetSymbolFont{timesoperators}{bold}{T1}{ptm}{b}{n}
\DeclareMathAlphabet{\mathbb}{U}{jkpsyb}{m}{n}
\SetMathAlphabet{\mathbb}{bold}{U}{jkpsyb}{bx}{n}

\allowdisplaybreaks
\definecolor{colorLink}{RGB}{0,100,162}
\definecolor{colorCite}{RGB}{8,124,100}

\newcommand{\opP}{\operatorname{{P}}}
\newcommand{\opbfP}{\pmb{\operatorname{{P}}}}
\newcommand{\opF}{{\operatorname{{F}}}}
\newcommand{\bfPi}{\operatorname{\pmb{\Pi}}}
\newcommand{\mfK}{\mathfrak{K}}

\def\R{\mathbb{R}}

\def\E{\mathbb{E}}

\def\CS{\mathcal{S}}

\def\tgamma{\tilde{\gamma}}

\newcommand{\bfK}{\pmb{\mathcal{K}}}
\newcommand{\mcK}{\mathcal{K}}

\def\one{\mathbf{1}}

\let\eps\varepsilon
\let\epsilon\varepsilon
\def\fraks{\mathfrak{s}}
\let\s\fraks

\let\f\frac
\let\d\partial

\DeclareMathOperator{\supp}{supp}

\DeclareMathOperator{\avL}{\slashed{L}}
\def\avL{\rlap{\kern0.05em-}L}

\def\eqdef{\coloneq}

\newcommand{\vertii}[1]{{\left\Vert #1 
		\right\Vert}}
\newcommand{\vertiii}[1]{{\left\vert\kern-0.35ex\left\vert\kern-0.35ex\left\vert #1 
		\right\vert\kern-0.35ex\right\vert\kern-0.35ex\right\vert}}

\def\dash{\leavevmode\unskip\kern0.18em--\penalty\exhyphenpenalty\kern0.18em}
\def\slash{\leavevmode\unskip\kern0.15em/\penalty\exhyphenpenalty\kern0.15em}

\makeatletter
\renewcommand{\operator@font}{\mathgroup\symtimesoperators}
\makeatother

\makeatletter
\DeclareRobustCommand{\TitleEquation}[2]{\texorpdfstring{\StrLeft{\f@series}{1}[\@firstchar]$\if%
b\@firstchar\boldsymbol{#1}\else#1\fi$}{#2}}
\makeatother

\renewcommand{\hat}{\widehat}
\renewcommand{\tilde}{\widetilde}

\newcommand{\overbar}[1]{\mkern 1.5mu\overline{\mkern-1.5mu#1\mkern-1.5mu}\mkern 1.5mu}
\def\bsquare{\overbar{\square}}

\colorlet{symbols}{blue!90!black}
\colorlet{symbols2}{red!90!black}
\colorlet{testcolor}{green!60!black}

\usetikzlibrary{calc}
\usetikzlibrary{shapes.misc}
\usetikzlibrary{shapes.symbols}
\usetikzlibrary{shapes.geometric}
\usetikzlibrary{snakes}
\usetikzlibrary{decorations}
\usetikzlibrary{decorations.markings}
\usetikzlibrary{external}

\tikzset{
	root/.style={circle,fill=testcolor,inner sep=0pt, minimum size=2mm},
	broot/.style={circle,fill=gray,inner sep=0pt, minimum size=2mm},
	dot/.style={circle,fill=black,inner sep=0pt, minimum size=1mm},
		reddot/.style={circle,fill=red,inner sep=0pt, minimum size=1mm},
			bluedot/.style={circle,fill=blue,inner sep=0pt, minimum size=1mm},
	eps/.style={circle,fill=white,draw=symbols,inner sep=0pt,minimum size=0.8mm},
	reps/.style={circle,fill=white,draw=symbols2,inner sep=0pt,minimum size=0.8mm},
	int/.style={circle,fill=black,draw=black,inner sep=0pt,minimum size=0.7mm},
	var/.style={circle,fill=black!10,draw=black,inner sep=0pt, minimum size=2mm},
		bluevar/.style={circle,fill=blue,draw=blue,inner sep=0pt, minimum size=1.5mm},
		redvar/.style={circle,fill=red,draw=red,inner sep=0pt, minimum size=1.5mm},
		mixvar
	dotred/.style={circle,fill=black!50,inner sep=0pt, minimum size=2mm},
	generic/.style={semithick,shorten >=1pt,shorten <=1pt},
	dist/.style={ultra thick,draw=testcolor,shorten >=1pt,shorten <=1pt},
	testfcn/.style={ultra thick,testcolor,shorten >=1pt,shorten <=1pt,<-},
	testfcnx/.style={ultra thick,testcolor,shorten >=1pt,shorten <=1pt,<-,
		postaction={decorate,decoration={markings,mark=at position 0.6 with {\drawx}}}},
	keps/.style={semithick,shorten >=1pt,shorten <=1pt,densely dashed,->},
	kprimex/.style={semithick,shorten >=1pt,shorten <=1pt,densely dashed,->,
		postaction={decorate,decoration={markings,mark=at position 0.4 with {\drawx}}}},
	kernel/.style={semithick,shorten >=1pt,shorten <=1pt,->},
	multx/.style={shorten >=1pt,shorten <=1pt,
		postaction={decorate,decoration={markings,mark=at position 0.5 with {\drawx}}}},
	kernelx/.style={semithick,shorten >=1pt,shorten <=1pt,->,
		postaction={decorate,decoration={markings,mark=at position 0.4 with {\drawx}}}},
	kernel1/.style={->,semithick,shorten >=1pt,shorten <=1pt,postaction={decorate,decoration={markings,mark=at position 0.45 with {\draw[-] (0,-0.1) -- (0,0.1);}}}},
	kernel2/.style={->,semithick,shorten >=1pt,shorten <=1pt,postaction={decorate,decoration={markings,mark=at position 0.45 with {\draw[-] (0.05,-0.1) -- (0.05,0.1);\draw[-] (-0.05,-0.1) -- (-0.05,0.1);}}}},
	kernelBig/.style={semithick,shorten >=1pt,shorten <=1pt,decorate, decoration={zigzag,amplitude=1.5pt,segment length = 3pt,pre length=2pt,post length=2pt}},
	rho/.style={dotted,semithick,shorten >=1pt,shorten <=1pt},
	renorm/.style={shape=circle,fill=white,inner sep=1pt},
	labl/.style={shape=rectangle,fill=white,inner sep=1pt},
	xi/.style={circle,fill=symbols!10,draw=symbols,inner sep=0pt,minimum size=1.2mm},
	xix/.style={crosscircle,fill=symbols!10,draw=symbols,inner sep=0pt,minimum size=1.2mm},
	xib/.style={circle,fill=symbols!10,draw=symbols,inner sep=0pt,minimum size=1.6mm},
	xibx/.style={crosscircle,fill=symbols!10,draw=symbols,inner sep=0pt,minimum size=1.6mm},
	not/.style={circle,fill=symbols,draw=symbols,inner sep=0pt,minimum size=0.5mm},
cumu2n/.style={inner sep=3pt},
cumu2/.style={draw=red!80,fill=red!40},
cumu2b/.style={draw=blue!80,fill=blue!40},
cumu2nv/.style={inner sep=3pt},
cumu2v/.style={draw=red!80,fill=white,very thick},
cumu3/.style={regular polygon, regular polygon sides=3,draw=red!80,rounded corners=3pt,fill=red!40,minimum size=5mm},
cumu4/.style={regular polygon, regular polygon sides=4,draw=red!80,rounded corners=3pt,fill=red!40,minimum size=7mm},
cumu5/.style={regular polygon, regular polygon sides=5,draw=red!80,rounded corners=3pt,fill=red!40,minimum size=7mm},
	>=stealth,
	not/.style={circle,fill=symbols,draw=symbols,inner sep=0pt,minimum size=0.5mm},
kernels2/.style={very thick,segment length=12pt},
arrho/.style={dotted,semithick,shorten >=1pt,shorten <=1pt, ->},
	}
	
\makeatletter
\def\DeclareSymbol#1#2#3{%
	\expandafter\gdef\csname MH@symb@#1\endcsname{\tikzsetnextfilename{symbol#1}%
		\tikz[baseline=#2,scale=0.15,draw=symbols,line join=round,line cap=round]{#3}}%
	\expandafter\gdef\csname MH@symb@#1s\endcsname{\scalebox{0.75}{\tikzsetnextfilename{symbol#1}%
			\tikz[baseline=#2,scale=0.15,draw=symbols,line join=round,line cap=round]{#3}}}%
	\expandafter\gdef\csname MH@symb@#1ss\endcsname{\scalebox{0.65}{\tikzsetnextfilename{symbol#1}%
			\tikz[baseline=#2,scale=0.15,draw=symbols,line join=round,line cap=round]{#3}}}%
}

\def\<#1>{\ifmmode\mathchoice{\csname MH@symb@#1\endcsname}{\csname MH@symb@#1\endcsname}{\csname MH@symb@#1s\endcsname}{\csname MH@symb@#1ss\endcsname}\else\csname MH@symb@#1\endcsname\fi}
\makeatother

\DeclareSymbol{X}{-2.4}{\node[eps] {};}
\DeclareSymbol{I}{0}{\draw[white] (-.5,0) -- (.5,0); \draw (0,0)  -- (0,1.5) {};}
\DeclareSymbol{bI}{0}{
\draw[white] (-.5,0) -- (.5,0); \draw[kernels2] (0,0)  -- (0,1.5);
}

\DeclareSymbol{1r}{0}{\draw[white] (-.5,0) -- (.5,0); \draw[red] (0,0)  -- (0,1.5) node[reps] {};}

\DeclareSymbol{1}{0}{\draw[white] (-.5,0) -- (.5,0); \draw (0,0)  -- (0,1.5) node[eps] {};}
\DeclareSymbol{2}{0}{\draw (-0.5,1.5) node[eps] {} -- (0,0) -- (0.5,1.5) node[eps] {};}
\DeclareSymbol{3}{0}{\draw (0,0) -- (0,1.5) node[eps] {}; \draw (-.7,1.3) node[eps] {} -- (0,0) -- (.7,1.3) node[eps] {};}
\DeclareSymbol{E3}{0}{\draw (0,0) -- (0,1.5) node[not] {}; \draw (-.7,1.3) node[not] {} -- (0,0) -- (.7,1.3) node[not] {}; \node[eps] (0,0) {};}

\DeclareSymbol{E5E4}{1}{
\draw (.5,2) node[not] {} -- (1.5,1.5) -- (2.5,2) node[not] {};
\draw (1.5,1.5) -- (1.5,2.7) node[not] {}; \draw (0.8,2.5) node[not] {} -- (1.5,1.5) -- (2.2,2.5) node[not] {};
\draw (0,0) -- (0,1.5) node[not] {}; \draw (-.7,1.3) node[not] {} -- (0,0) -- (.7,1.3) node[not] {};
\draw (-1.2,0.7) node[not] {} -- (0,0);
\draw (0,0) to[bend right=50] (1.5,1.5); 
\node[eps] at (0,0) {};\node[eps] at (1.5,1.5) {};
}
\DeclareSymbol{E4E4}{1}{
\draw (.5,2) node[not] {} -- (1.5,1.5) -- (2.5,2) node[not] {};
\draw (1,2.6) node[not] {} -- (1.5,1.5) -- (2,2.6) node[not] {};
\draw (0,0) -- (0,1.5) node[not] {}; \draw (-.7,1.3) node[not] {} -- (0,0) -- (.7,1.3) node[not] {};
\draw (-1.2,0.7) node[not] {} -- (0,0);
\draw (0,0) to[bend right=50] (1.5,1.5); 
\node[eps] at (0,0) {};\node[eps] at (1.5,1.5) {};
}
\DeclareSymbol{3E4}{1}{
\draw (1.5,1.5) -- (1.5,2.7) node[not] {}; \draw (0.8,2.5) node[not] {} -- (1.5,1.5) -- (2.2,2.5) node[not] {};
\draw (0,0) -- (0,1.5) node[not] {}; \draw (-.7,1.3) node[not] {} -- (0,0) -- (.7,1.3) node[not] {};
\draw (-1.2,0.7) node[not] {} -- (0,0);
\draw (0,0) to[bend right=50] (1.5,1.5); 
\node[eps] at (0,0) {};
}

\DeclareSymbol{3E3}{1}{
\draw (1.5,1.5) -- (1.5,2.7) node[not] {}; \draw (0.8,2.5) node[not] {} -- (1.5,1.5) -- (2.2,2.5) node[not] {};
\draw (0,0) -- (0,1.5) node[not] {}; \draw (-.7,1.3) node[not] {} -- (0,0) -- (.7,1.3) node[not] {};
\draw (0,0) to[bend right=50] (1.5,1.5); 
\node[eps] at (0,0) {};
}

\DeclareSymbol{E5}{0}{\draw (0,0) -- (0,1.5) node[not] {}; \draw (-.7,1.3) node[not] {} -- (0,0) -- (.7,1.3) node[not] {};
\draw (-1.2,0.7) node[not] {} -- (0,0) -- (1.2,0.7) node[not] {}; \node[eps] (0,0) {};}
\DeclareSymbol{E52}{-3}{\draw (0,0) -- (0,-1); \draw (1.5,-0.4) node[not] {} -- (0,-1) -- (-1.5,-0.4) node[not] {}; 
\draw (-1,0.5) node[not] {} -- (0,0) -- (1,0.5) node[not] {};
\draw (0,0) -- (0,1.2) node[not] {}; \draw (-.7,1) node[not] {} -- (0,0) -- (.7,1) node[not] {};\node[eps] (0,0) {};}
\DeclareSymbol{E32}{-3}{\draw (0,0) -- (0,-1); \draw (1.5,-0.4) node[not] {} -- (0,-1) -- (-1.5,-0.4) node[not] {}; 
\draw (0,0) -- (0,1.2) node[not] {}; \draw (-.7,1) node[not] {} -- (0,0) -- (.7,1) node[not] {};\node[eps] (0,0) {};}
\DeclareSymbol{E50}{-3}{\draw (0,0) -- (0,-1); 
\draw (-1,0.5) node[not] {} -- (0,0) -- (1,0.5) node[not] {};
\draw (0,0) -- (0,1.2) node[not] {}; \draw (-.7,1) node[not] {} -- (0,0) -- (.7,1) node[not] {};\node[eps] (0,0) {};}

\DeclareSymbol{3E2}{-3}{\draw (0,0) -- (0,-1); \draw (1,0) node[not] {} -- (0,-1) -- (-1,0) node[not] {}; \draw (0,0) -- (0,1.2) node[not] {}; \draw (-.7,1) node[not] {} -- (0,0) -- (.7,1) node[not] {};\node[eps] at (0,-1) {};}

\DeclareSymbol{31}{-3}{\draw (0,0) -- (0,-1) -- (1,0) node[eps] {}; \draw (0,0) -- (0,1.2) node[eps] {}; \draw (-.7,1) node[eps] {} -- (0,0) -- (.7,1) node[eps] {};}
\DeclareSymbol{30}{-3}{\draw (0,0) -- (0,-1); \draw (0,0) -- (0,1.2) node[eps] {}; \draw (-.7,1) node[eps] {} -- (0,0) -- (.7,1) node[eps] {};}
\DeclareSymbol{32}{-3}{\draw (0,0) -- (0,-1); \draw (1,0) node[eps] {} -- (0,-1) -- (-1,0) node[eps] {}; \draw (0,0) -- (0,1.2) node[eps] {}; \draw (-.7,1) node[eps] {} -- (0,0) -- (.7,1) node[eps] {};}
\DeclareSymbol{32b}{-4.3}{\draw (0,0.5) -- (0,-1); \draw (1,0) node[not] {} -- (0,-1) -- (-1,0) node[not] {};}
\DeclareSymbol{22}{-2}{\draw (0,0.3) -- (0,-1); \draw (1,0) node[eps] {} -- (0,-1) -- (-1,0) node[eps] {};\draw (-.7,1) node[eps] {} -- (0,0.3) -- (.7,1) node[eps] {};}
\DeclareSymbol{20}{-3}{\draw (0,0) -- (0,-1);\draw (-.7,1) node[not] {} -- (0,0) -- (.7,1) node[not] {};}
\DeclareSymbol{12}{-3}{\draw (0,0.3) -- (0,-1); \draw (1,0) node[not] {} -- (0,-1) -- (-1,0) node[not] {};\draw (-.7,1) node[not] {} -- (0,0.3);}
\DeclareSymbol{10}{-3}{\draw (0,0.3) -- (0,-1);\draw (-.7,1) node[not] {} -- (0,0.3);}
\DeclareSymbol{21}{-3}{\draw (0,0.3) -- (0,-1) -- (1,0) node[not] {};\draw (-.7,1) node[not] {} -- (0,0.3) -- (.7,1) node[not] {};}

\DeclareSymbol{1'}{0}{\draw[white] (-.5,0) -- (.5,0); \draw (0,0)  -- (0,1.5) node[not] {}; \draw ( -0.2,1)-- ( 0.2, 1);}

\DeclareSymbol{3'}{0}{\draw (0,0) -- (0,1.5) node[not] {};
\draw  (-0.25,1.1) -- ( 0.25, 1.1);
 \draw (-.7,1.3) node[eps] {} -- (0,0) -- (.7,1.3) node[eps] {};}

\DeclareSymbol{3''}{0}{\draw (0,0) -- (0,1.5) node[eps] {};
 \draw (-.7,1.3) node[not] {} -- (0,0) -- (.7,1.3) node[not] {};
 \draw  (-0.8,0.7) -- ( -0.25, 1);
  \draw  (0.8,0.7) -- ( 0.25, 1);
 }
 
 \DeclareSymbol{3'''}{0}{\draw (0,0) -- (0,1.5) node[not] {};
\draw  (-0.3,1.05) -- ( 0.3, 1.05);
 \draw (-.7,1.3) node[not] {} -- (0,0) -- (.7,1.3) node[not] {};
 \draw  (-0.8,0.7) -- ( -0.25, 1);
  \draw  (0.8,0.7) -- ( 0.25, 1);
 }
 
 \DeclareSymbol{3'2}{-3}{\draw (0,0) -- (0,-1); \draw (1,0) node[eps] {} -- (0,-1) -- (-1,0) node[eps] {}; \draw (0,0) -- (0,1.4) node[not] {}; \draw (-.7,1) node[eps] {} -- (0,0) -- (.7,1) node[eps] {};
\draw  (-0.25,1) -- ( 0.25, 1); 
 }
 
 \DeclareSymbol{3''2}{-3}{\draw (0,0) -- (0,-1); \draw (1,0) node[eps] {} -- (0,-1) -- (-1,0) node[eps] {}; \draw (0,0) -- (0,1.3) node[eps] {}; \draw (-.7,1) node[not] {} -- (0,0) -- (.7,1) node[not] {};
 \draw  (-0.8,0.5) -- ( -0.25, 0.8);
  \draw  (0.8,0.5) -- ( 0.25, 0.8); 
 }
 
  \DeclareSymbol{3'''2}{-3}{\draw (0,0) -- (0,-1); \draw (1,0) node[eps] {} -- (0,-1) -- (-1,0) node[eps] {}; \draw (0,0) -- (0,1.4) node[not] {}; \draw (-.7,1) node[not] {} -- (0,0) -- (.7,1) node[not] {};
 \draw  (-0.8,0.5) -- ( -0.25, 0.8);
  \draw  (0.8,0.5) -- ( 0.25, 0.8); 
  \draw ( -0.25,1)-- ( 0.25, 1);
 }
 
  \DeclareSymbol{3'2'}{-3}{\draw (0,0) -- (0,-1); \draw (1,0) node[eps] {} -- (0,-1) -- (-1,0) node[not] {}; \draw (0,0) -- (0,1.4) node[not] {}; \draw (-.7,1) node[eps] {} -- (0,0) -- (.7,1) node[eps] {};
\draw  (-0.25,1) -- ( 0.25, 1); 
\draw  (-0.8,-0.6) -- ( -0.4, -0.2);
 }
 
  \DeclareSymbol{3'2''}{-3}{\draw (0,0) -- (0,-1); \draw (1,0) node[not] {} -- (0,-1) -- (-1,0) node[not] {}; \draw (0,0) -- (0,1.4) node[not] {}; \draw (-.7,1) node[eps] {} -- (0,0) -- (.7,1) node[eps] {};
\draw  (-0.25,1) -- ( 0.25, 1); 
\draw  (-0.8,-0.6) -- ( -0.4, -0.2);
\draw  (0.8,-0.6) -- ( 0.4, -0.2);
 }

   \DeclareSymbol{3''2'}{-3}{\draw (0,0) -- (0,-1); \draw (1,0) node[eps] {} -- (0,-1) -- (-1,0) node[not] {}; \draw (0,0) -- (0,1.3) node[eps] {}; \draw (-.7,1) node[not] {} -- (0,0) -- (.7,1) node[not] {};
 \draw  (-0.8,0.5) -- ( -0.25, 0.8);
  \draw  (0.8,0.5) -- ( 0.25, 0.8); 
\draw  (-0.8,-0.6) -- ( -0.4, -0.2);
 }
  \DeclareSymbol{32'}{-3}{\draw (0,0) -- (0,-1); \draw (1,0) node[eps] {} -- (0,-1) -- (-1,0) node[not] {}; \draw (0,0) -- (0,1.3) node[eps] {}; \draw (-.7,1) node[eps] {} -- (0,0) -- (.7,1) node[eps] {};
\draw  (-0.8,-0.6) -- ( -0.4, -0.2);
 }
  \DeclareSymbol{3'0}{-3}{\draw (0,0) -- (0,-1); 
   \draw (0,0) -- (0,1.4) node[not] {}; \draw (-.7,1) node[eps] {} -- (0,0) -- (.7,1) node[eps] {};
\draw  (-0.25,1) -- ( 0.25, 1); 
 }

\DeclareSymbol{1}{0}{\draw[white] (-.4,0) -- (.4,0); \draw (0,0)  -- (0,1.2) node[eps] {};}
\DeclareSymbol{2}{0}{\draw (-0.5,1.2) node[eps] {} -- (0,0) -- (0.5,1.2) node[eps] {};}
\DeclareSymbol{11}{0}{\draw (0,1.8) node[eps] {} -- (-0.7,0.9) -- (0,0) -- (0.7,1) node[eps] {};}
\DeclareSymbol{10}{0}{\draw (0,1.8) node[eps] {} -- (-0.8,0.9) -- (0,0);}
\DeclareSymbol{21}{0.7}{\draw (-1,1.8) node[eps] {} -- (-0.5,0.9); \draw (0,1.8) node[eps] {} -- (-0.5,0.9) -- (0,0) -- (0.5,0.9) node[eps] {};}
\DeclareSymbol{20}{0.7}{\draw (-1,1.8) node[eps] {} -- (-0.5,0.9); \draw (0,1.8) node[eps] {} -- (-0.5,0.9) -- (-0.5,0);}
\DeclareSymbol{210}{1.1}{\draw (-0.5,2.4) node[eps] {} -- (-1,1.6);\draw (-1.5,2.4) node[eps] {} -- (-0.5,0.8); \draw (0,1.6) node[eps] {} -- (-0.5,0.8) -- (-0.5,0);}
\DeclareSymbol{211}{1.1}{\draw (-0.5,2.4) node[eps] {} -- (-1,1.6);\draw (-1.5,2.4) node[eps] {} -- (-0.5,0.8); \draw (0,1.6) node[eps] {} -- (-0.5,0.8) -- (0,0) -- (0.5,0.8) node[eps] {};}
\DeclareSymbol{22j}{0.7}{\draw (-1.5,1.8) node[eps] {} -- (-1,0.9) -- (0,0) -- (1,0.9) -- (1.5,1.8) node[eps] {};
\draw (-0.5,1.8) node[eps] {} -- (-1,0.9);\draw (0.5,1.8) node[eps] {} -- (1,0.9);}
\DeclareSymbol{2'}{0}{\draw (-0.65,1.3) node[not] {} -- (0,0) -- (0.6,1.3) node[eps] {};
\draw  (-0.8,0.7) -- ( -0.15, 1.1);} 

\DeclareSymbol{2''}{0}{\draw (-0.65,1.3) node[not] {} -- (0,0) -- (0.6,1.3) node[not] {};
\draw  (-0.8,0.7) -- ( -0.15, 1.1);
\draw  (0.8,0.7) -- ( 0.15, 1.1);
} 
 
\DeclareSymbol{2'1}{0.7}{\draw (-1,1.8) node[not] {} -- (-0.5,0.9); \draw (0,1.8) node[eps] {} -- (-0.5,0.9) -- (0,0) -- (0.5,0.9) node[eps] {};
\draw  (-1.1,1.2) -- ( -0.5, 1.6);} 

\DeclareSymbol{2''1}{0.7}{\draw (-1,1.8) node[not] {} -- (-0.5,0.9); \draw (0,1.8) node[not] {} -- (-0.5,0.9) -- (0,0) -- (0.5,0.9) node[eps] {};
\draw  (-1.1,1.2) -- ( -0.5, 1.6);
\draw  (.15,1.2) -- ( -0.45, 1.6);} 
 
\DeclareSymbol{2''1'}{0.7}{\draw (-1,1.8) node[not] {} -- (-0.5,0.9); \draw (0,1.8) node[not] {} -- (-0.5,0.9) -- (0,0) -- (0.55,1) node[not] {};
\draw  (-1.1,1.2) -- ( -0.5, 1.6);
\draw  (.15,1.2) -- ( -0.45, 1.6);
\draw  (0.8,0.4) -- ( 0, 0.8);
} 
\DeclareSymbol{2'1'}{0.7}{\draw (-1,1.8) node[not] {} -- (-0.5,0.9); \draw (0,1.8) node[eps] {} -- (-0.5,0.9) -- (0,0) -- (0.55,1) node[not] {};
\draw  (-1.1,1.2) -- ( -0.5, 1.6);
\draw  (0.8,0.4) -- ( 0, 0.8);} 

\DeclareSymbol{21'}{0.7}{\draw (-1,1.8) node[eps] {} -- (-0.5,0.9); \draw (0,1.8) node[eps] {} -- (-0.5,0.9) -- (0,0) -- (0.55,1) node[not] {};
\draw  (0.8,0.4) -- ( 0, 0.8);} 

\DeclareSymbol{2'0}{0.7}{\draw (-1,1.8) node[not] {} -- (-0.5,0.9); \draw (0,1.8) node[eps] {} -- (-0.5,0.9) -- (-0.5,0);
\draw  (-1.1,1.2) -- ( -0.5, 1.6);}

\DeclareSymbol{b1}{0}{
\draw[white] (-.5,0) -- (.5,0); \draw[kernels2] (0,0)  -- (0,1.5);
 \draw  (0,1.5) node[eps] {};
}

\DeclareSymbol{b2}{0}{
\draw[kernels2] (0,0) -- (-0.5,1.4);
\draw[kernels2] (0,0) -- (+0.5,1.4); 
 \draw  (-0.5,1.4) node[eps] {};
 \draw  (0.5,1.4) node[eps] {};
 }
 
 \DeclareSymbol{20b}{0.7}{\draw (-1,1.8) node[eps] {} -- (-0.5,0.9);
  \draw (0,1.8) node[eps] {} -- (-0.5,0.9);
  \draw[kernels2](-0.5,0.9) -- (-0.5,0);}
  
\DeclareSymbol{210b}{1.1}{\draw (-0.5,2.4) node[eps] {} -- (-1,1.6);\draw (-1.5,2.4) node[eps] {} -- (-0.5,0.8); \draw (0,1.6) node[eps] {} -- (-0.5,0.8);
\draw[kernels2] (-0.5,0.8) -- (-0.5,0);}
 
\DeclareSymbol{bI}{0}{\draw[white] (-.5,0) -- (.5,0); \draw[kernels2] (0,0)  -- (0,1.5) {};} 

\DeclareSymbol{Xi2}{0}{\draw[white] (-.5,0) -- (.5,0); \draw (0,0)  -- (0,1.5) node[eps] {};
\draw (0,0) node[eps] {};}

 \DeclareSymbol{1'Xi}{0}{\draw[white] (-.5,0) -- (.5,0); \draw (0,0)  -- (0,1.5) node[not] {}; \draw ( -0.2,1)-- ( 0.2, 1);
 \draw (0,0) node[eps] {};
 }

\theoremstyle{plain}
\newtheorem{theorem}{Theorem}[section]
\newtheorem{corollary}[theorem]{Corollary}
\newtheorem{lemma}[theorem]{Lemma}
\newtheorem{proposition}[theorem]{Proposition}
\newtheorem{prop}[theorem]{Proposition}

\theoremstyle{definition}
\newtheorem{definition}[theorem]{Definition}

\newtheorem{assumption}[theorem]{Assumption}
\newtheorem{remark}[theorem]{Remark}

\numberwithin{equation}{section}
\renewcommand{\theequation}{\thesection.\arabic{equation}}

\colorlet{darkblue}{blue!90!black}
\colorlet{darkgreen}{green!50!black}
\colorlet{darkred}{red!90!black}
\let\comm\margincomment

\DeclareMathOperator{\id}{id}
\DeclarePairedDelimiter{\abs}{\lvert}{\rvert}

\begin{document}

\title{Homogenisation of singular SPDEs}
\author{Martin Hairer$^{1}$\orcidlink{0000-0002-2141-6561}, Harprit Singh$^{2}$\orcidlink{0000-0002-9991-8393}}

\institute{
EPFL, Switzerland and Imperial, UK. \email{martin.hairer@epfl.ch}
\and University of Vienna, AT. \email{singhh92@univie.ac.at}}

\maketitle
\begin{abstract}

We introduce an approach to study homogenisation of a large class of singular SPDEs of the form 
$$
\partial_t u_\eps -  \nabla\cdot  {A}(x/\eps,t/\eps^2) \nabla u_\eps 
= F(x/\eps , t/\eps^2, u_\eps , 
\nabla u_\eps , \xi )
$$
which is based on the idea of importing (classical) homogenisation results into the framework of regularity structures and the insight that one can rewrite the SPDE under consideration in terms of a model, where
the correctors (from homogenisation theory) are seen as further `abstract noises'. 

As applications, we establish periodic space-time homogenisation results for oscillatory generalisations of the 2d g-PAM and $\Phi^4_3$ equation proving that when the noise is regularised at scale $\delta\ll 1$ solutions to the equations with coefficient field ${A}(x/\eps,t/\eps^2)$, when appropriately renormalised, converge to solutions to the corresponding homogenised equation along any sequence $(\eps,\delta)\to 0$. 
We make the observation that the unbounded divergences can be written as sums of two types of terms: `small scale' terms, the spatial dependence of which is an explicit local function of finitely many derivatives of the coefficient field and `large scale' terms, which for logarithmic divergences are explicit involving the homogenised matrix and correctors. Furthermore, in order to recover the same solution to the corresponding homogenised equation along any joint limit $(\eps, \delta)\to 0$ one has to subtract additional bounded renormalisation constants $c_{\eps,\delta}$, which appear due to oscillations at mesoscopic scales, as well as resonances between the coefficient field and the oscillations in the nonlinearity and which in general have the property that $\lim_{\eps\to 0} \lim_{\delta\to 0} c_{\eps,\delta} \neq \lim_{\delta\to 0} \lim_{\eps\to 0} c_{\eps,\delta}$.

%
%
\end{abstract}

\setcounter{tocdepth}{2}
\tableofcontents


\section{Introduction}\label{sec:intro}
The area of stochastic partial differential equations has seen rapid progress
over the past decade, spurred by the introduction of the theory of
regularity structures~\cite{Hai14} and of para-controlled calculus~\cite{GIP15}.
Despite the 
close connections of singular SPDEs to physical phenomena---for instance, $\Phi^4$ as a model of ferromagnetism and the parabolic Anderson model
relating to branching processes---the
theory of singular SPDEs has to date been developed primarily in homogeneous settings,
involving constant-coefficient operators. In many situations, however, such an
assumption is not justified and physical models often exhibit heterogeneities on
small scales. In the case of classical (non-singular) PDEs
it is by now very well understood 
when, and up to which scales, the behaviour of solutions to heterogeneous equations is still governed by an
effective (constant coefficient) homogenised equation.


In this article we develop a general framework for homogenisation of singular SPDEs that combines quantitative homogenisation with the analytic machinery of regularity structures. The key idea is to lift the correctors from homogenisation as new abstract noises into an enlarged regularity structure. This allows to encode the two scale expansions into the abstract formulation of the equation and thereby to treat
homogenisation and renormalisation within the same fixed-point problem.
We illustrate the approach by establishing periodic space-time homogenisation results for oscillatory variants of the g-PAM equation on $\mathbb{T}^2\times \mathbb{R}$, formally given by
\begin{equ}\label{eq:g-PAM}
\partial_t u - \nabla \cdot A(x/\eps,t/\eps^2)\nabla u=  \sum_{i,j=1}^2 f_{i,j}(x/\eps, t/\eps^2) \partial_i u \partial_j u  +{\sigma}(u)\xi \;,
\end{equ}
where $\xi$ is a spatial white noise and $f_{i,j}$ are bounded measurable functions as well as 
for (oscillatory) $\Phi^4_3$ equations on $\mathbb{T}^3\times \mathbb{R}$ formally given by
\begin{equ}\label{eq:phi^4}
\partial_t u - \nabla \cdot A(x/\eps,t/\eps^2) \nabla u= -f(x/\eps, t/\eps^2) u^3 + \xi \;,
\end{equ}
where $\xi$ is a space-time white noise.

When the matrix $A$ and the functions $f,f_{i,j}$ are constant, Equations \eqref{eq:g-PAM}\&\eqref{eq:phi^4} were first solved in \cite{Hai14}, 
see also \cite{GIP15, kup16, BB16, cat18, jag23}, which initiated substantial progress in the field. On the one hand, alternative approaches to regularity structures and para-controlled calculus have been developed such as the flow approach \cite{kup16,Duc21} and a multi-index formulations of models \cite{OW19,LOT, BOS25}.
On the other hand the theory of regularity structures \cite{Hai14} now provides a fully automated `black-box'-solution framework for constant coefficient equations when combined with \cite{BHZ19, BCCH20, CH16}, and if further combined with \cite{BB21, Sin23, BSS25} this black box extends to equations with (regular) variable coefficients.


These advances in the area of singular SPDEs have also provided a rigorous framework for `weak universality' results, previously conjectured in the physics literature, see \cite{GP16univ, HQ18, furlan2019weak, MP19, EX22,yang23hairer,  KWX24}. 
Roughly speaking, these works show that for each of the 
SPDEs considered therein, there exists a large class of perturbed equations\slash models —
obtained by modifying the differential operator and\slash or the nonlinearity—
whose solutions converge under suitable rescaling to those of the original
equation.

The homogenisation problems studied here may be viewed as natural instances of weak universality. 
However, in contrast to the previously analysed settings,
where perturbations of the differential operator are typically spatially homogeneous and can
be understood directly at the level of the symbol, for the oscillatory
space-time coefficients considered here this not possible and, as we shall see, substantial changes in the renormalisation 
 procedure have to be made.

Let us also mention that many singular SPDEs, in particular also the $\Phi^4$ equation, arise as stochastic quantisation equations. The rigorous understanding of these SPDEs has led to new constructions and understanding about the respective invariant measures \cite{gubinelli2021pde, HS22varphi, BDFT23a, CS23invariant, klose2024large, duch24, BDW25, DHYZ25}. While this is not the focus of this article, its content does have clear implications in that direction. 


Returning to \eqref{eq:g-PAM} and \eqref{eq:phi^4}, it is a priori not even clear how a homogenisation result for singular SPDEs should be formulated since \eqref{eq:g-PAM} and \eqref{eq:phi^4} as written are only formal and, when the matrix $A$ is 
not constant, the required renormalisation may have to be chosen in an inhomogeneous way.\footnote{Note that a similar situation arises when considering singular SPDEs in geometric settings, c.f.\ \cite{BB16, DDD19, mouzard2022weyl, HS23m, BDFT23a, hao2024singular, MS23}, but in contrast to the setting here, renormalisation constants are often still sufficient for covariant equations.}
A notion of solution was put forth in \cite{Sin23}, which formally corresponds for \eqref{eq:g-PAM} to
\begin{equs}
 \partial_t u - \nabla \cdot A(x,t) \nabla u &=\sum_{i,j=1}^2 f_{i,j}(x,t)  \left( \partial_i u \partial_j u - \infty \cdot \frac{(A_s^{-1})_{i,j}(x,t)}{\sqrt{\det(A_s(x,t))}}  \sigma^2(u)\right)\\
 &\qquad \qquad + \sigma(u)\left( \xi  -\infty \cdot  \frac{\sigma'(u)}{\sqrt{\det(A_s(x,t))}}\right) \ ,
\end{equs}
respectively, for \eqref{eq:phi^4} to
$$\partial_t u - \nabla  \cdot A(x,t) \nabla u= -u_\eps^3 +  \infty \cdot   \frac{u_\epsilon}{\sqrt{\det(A_s(x,t))}} - \infty  \cdot \frac{u_\epsilon}{\det(A_s(x,t))} +\xi_\epsilon \ ,$$
where $A_s$ denotes the symmetric part of $A$.
First homogenisation results for the $\Phi^4_2$, resp.\ $P(\phi)_2$, equation were then established in \cite{HS23per} (for $f=1$), resp. \cite{CX23} (for $f=1$ and symmetric $A$ not depending on $t$).
Analytically, both \cite{HS23per} and \cite{CX23} work directly with the remainder equation following Da Prato--Debussche \cite{DD03}, together with the observation that
classical homogenisation results provide convergence of heat kernels as $\eps\to 0$ in a sufficiently strong topology to conclude by continuity of the solution map as a function of the heat kernel. 
\subsection*{Homogenisation by means of Regularity Structures}
For more singular equations such as \eqref{eq:g-PAM}--\eqref{eq:phi^4}, besides the fact that we use the theory of regularity structures since a more sophisticated solution theory for singular SPDEs is required, if one tries to conclude by convergence of the heat kernel, the following facts pose a direct obstacle:
\begin{itemize}
\item 
Denoting $A^{\eps}(x,t)= A( x/\eps, t/\eps^2)$ the rescaled coefficient field,
 even for smooth functions $f$ the solutions to the equation 
\begin{equ}\label{eq:smooth homogenisation}
(\partial_t- \nabla \cdot A^\eps \nabla )u_\eps= f 
\end{equ}
 converge as $\eps\to 0$ in $C^\gamma$ to the solution $\bar{u}$ of a homogenised equation 
$
(\partial_t- \nabla \cdot \bar{A} \nabla )\bar{u}= f ,
$
 only for $\gamma<1$ but not for $\gamma>1$.
 \item Solutions to the singular SPDEs \eqref{eq:g-PAM} and \eqref{eq:phi^4} are both necessarily obtained as the reconstructions of a modelled distributions belonging 
 to $\mathcal{D}^\gamma$ (which is a generalisation of the space $C^\gamma$) for $\gamma>1$.\footnote{This is required in order to be able to close the fixed point problem. While the heuristic is explained here in terms of regularity structures, one encounters an analogous requirement when using for example para-controlled calculus.}
\end{itemize}
On the other hand, it is well understood that if we denote by $\psi_i^{\eps}:= \psi_i(x/\eps, t/\eps^2)$ the rescaled (parabolic) correctors of homogenisation theory, see \eqref{eq:correctors} for the precise definition, then
$$u_{\eps} - \bar{u} - \eps \sum_{i=1}^d\psi_i^{\eps} \partial_i \bar{u}  \to 0, \qquad \text{ in $C^{\gamma}$ for $\gamma\in (1,2)$ } \ .$$
This suggests to add correctors as new basis elements to the regularity structure. 
Thus, roughly speaking, we shall study homogenisation of singular SPDEs by rewriting the fixed point problem 
$u_{\eps}= \Gamma_\eps F_\eps(u_{\eps}, \nabla u_{\eps}, \xi)$ (say for vanishing initial conditions)
as
\begin{equ}\label{eq:formal rewriting}
 u_{\eps}=\bar{\Gamma} F_\eps(u_{\eps}, \nabla u_{\eps}, \xi)  + \eps \psi^{\eps}\nabla \bar{\Gamma}  F_\eps(u_{\eps}, \nabla u_{\eps}, \xi)+ (\Gamma_\eps-\bar{\Gamma}  - \eps \psi^{\eps} \nabla\bar{\Gamma} ) F_\eps(u_{\eps}, \nabla u_{\eps}, \xi) \;,
\end{equ}
and lift this to an abstract fixed point problem at the level of modelled distributions.
The kernel $\bar{\Gamma}$ in the first term is the (constant coefficient) heat kernel associated to the homogenised operator $\nabla \cdot \bar A \nabla$ and can be treated as in \cite{Hai14}. To treat the second term we add new abstract noises $\Psi_i$ to the regularity structure which correspond to the correctors, i.e.\ $ \Pi_x^{\eps}\Psi_i= \psi_i^{\eps}$ and also lift multiplication with $\eps$ to an abstract operator $\mathcal{E}$, somewhat reminiscent of \cite{HQ18}.
Lastly,
 we use that the kernel $ (\Gamma_\eps-\bar{\Gamma}  - \eps \psi^{\eps} \nabla\bar{\Gamma} ) $ 
converges to $0$ in a sufficiently strong topology to argue similarly to \cite{GH19, Sin23}.\footnote{To see that this strategy is plausible recall that for an abstract noise symbol $\Xi$ of homogeneity $|\Xi|<0$  
the map on modelled distributions $f\mapsto \Xi \cdot f$ is a map $ \mathcal{D}^{\gamma}\mapsto \mathcal{D}^{\gamma+|\Xi|}$,
in contrast to the map on functions $f\mapsto (\Pi\Xi) \cdot f$, which maps  $C^\gamma \to C^{|\Xi|}$ for $\gamma>-| \Xi | $.}
\subsection{Application to concrete equations}
We apply the strategy outlined above to the (oscillatory) g-PAM and $\Phi^4_3$ equations. For both equations we consider noise regularised at length scale $\delta>0$ and prove that in order to recover solutions to the corresponding homogenised singular SPDEs 
as $(\eps,\delta)\to 0$, further `large scale' (namely varying at scales larger than $\eps$) renormalisation is necessary in addition to the (local) `small scale' renormalisation functions varying at scale $\eps$ appearing in \cite{Sin23}. For the g-PAM equation \eqref{eq:g-PAM}, as one might expect due to the presence of derivatives on the right hand side, these `large scale' counterterms 
involve the correctors $\psi_i$.
More precisely, Theorem~\ref{thm:g-pam} states that there exist unbounded\footnote{Their precise asymptotic behaviour is given in Theorem~\ref{thm:g-pam}.} renormalisation constants 
$\alpha_{\eps,\delta}^{\<Xi2>}, \ \alpha_{\eps,\delta}^{\<b2>}, \ \bar{\alpha}_{\eps,\delta}^{\<Xi2>}, \ \bar{\alpha}_{\eps,\delta}^{\<b2>}$
and bounded constants $c_{\eps,\delta}^{\<Xi2>}$, $c_{\eps,\delta}^{{\<b2>{\mu,\nu}}}$, $\gamma_{\eps,\delta}$ for $\mu,\nu=1,2$
 such that if we write $u_{\eps,\delta}$ for the solution to 
\begin{align*}
\partial_t u - \nabla \cdot A^\eps \nabla u_{\eps,\delta}&=\sum_{\mu,\nu=1}^2 f_{\mu,\nu}^\eps \Bigg( \partial_\mu u \partial_\nu u_{\eps,\delta} - \Big( \frac{ (A_s^\eps)_{\mu,\nu}^{{-1}}}{\sqrt{\det(A_s^\eps)}}\alpha^{\<b2>}_{\eps,\delta}\\
&\qquad + \sum_{i,j} (\mathbf{1}_{\mu=i}+(\partial_{\mu} \phi_i)^\eps)(\mathbf{1}_{\nu=j}+(\partial_{\nu} \phi_j)^\eps) 
\frac{ (\bar{A}^{{-1}})_{i,j}    }{\sqrt{\det(\bar{A})}} \bar{\alpha}^{\<b2>}_{\eps,\delta} 
 + c^{\<b2>\mu,\nu}_{\eps,\delta}  \Big) \cdot  \sigma^2(u_{\eps,\delta}) \Bigg)
  \\& \qquad
+ \sigma(u_{\eps,\delta})\Big( \xi_{\eps, \delta}  -\big(\frac{\alpha_{\eps,\delta}^{\<Xi2>}}{\sqrt{\det(A_s^\eps)}} +\frac{\bar{\alpha}_{\eps,\delta}^{\<Xi2>}}{\sqrt{\det(\bar{A})}}  + c_{\eps,\delta}^{\<Xi2>}  \big)  \sigma'(u_{\eps,\delta})\Big)  
- \gamma_{\eps, \delta} \cdot \sigma^2(u_{\eps,\delta})
\ ,
\end{align*}
any joint limit $(\eps,\delta)\to 0$ recovers the same solution to the homogenised SPDE. 
Notably, the constants 
$c_{\eps,\delta}^{\<Xi2>}$, $c_{\eps,\delta}^{{\<b2>\mu,\nu}}$, $\gamma_{\eps,\delta}$ 
cannot be extended continuously to all of $(\eps,\delta)\in [0,1]^2$ in general. The former two constants are 
required for the same reason as the analogous bounded constant $c_{\eps,\delta}$ in \cite{HS23per} for the $\Phi^4_2$ 
equation. Roughly speaking, they capture interactions between homogenisation and 
renormalisation at mesoscopic scales. The constant $\gamma_{\eps,\delta}$ can be interpreted as occurring due to resonances between the oscillations in the nonlinearity and the coefficient field.
Importantly, for each $\eps>0$ the process $u_{\eps, 0}$ belongs to the class of solutions constructed in \cite{Sin23} and for $\eps=0$ in \cite{Hai14}.
Let us mention that while, here, $\xi_{\eps, \delta}$ is the rather specific regularisation by the heat semigroup of the operator $\nabla \cdot A^\eps \nabla$ itself, Theorem~\ref{thm:main_translation invariant} provides the analogous result for usual homogeneous regularisation.

\begin{remark}
A special case of \eqref{eq:g-PAM} where the nonlinearity does not depend on the gradient of the solution, i.e.\ $f_{\mu,\nu}=0$ was previously studied in \cite{CFX23}.
The main difference\footnote{
Another difference is that \cite{CFX23} considers coefficient fields $A\in C^\alpha$ for $\alpha>0$ which do not depend on time, while we consider space-time homogenisation assuming the slightly stronger regularity assumption $A\in C^{1,\alpha}$ since we allow for the additional gradient term in the nonlinearity.
}
between the corresponding special case of Theorem~\ref{thm:g-pam} and Theorem~\ref{thm:main_translation invariant} and the results therein is that for each $\eps>0$ we work with the notion of solution introduced in \cite{Sin23} to the singular (inhomogeneous) SPDE (which is independent of the homogenisation problem) and track precisely how homogenisation effects deform the (required) renormalisation functions while \cite{CFX23} treats renormalisation more qualitatively (working with an implicitly defined notion of solution). 
\end{remark}

Theorem~\ref{thm:phi4} establishes a homogenisation result for the (oscillatory) $\Phi^4_3$ equation.  The main difference, compared the case of g-PAM and $\phi^{4}_2$ in \cite{HS23per} is that this time the large scale part of the non-logarithmically diverging diagram does not seem to be cleanly expressible in terms of only the homogenised matrix for all regimes of $(\eps,\delta)$, but only for $\eps\lesssim \delta^2$, see Theorem~\ref{thm:phi4 restricted}.

\begin{remark}
The proofs of Theorems~\ref{thm:g-pam}, \ref{thm:main_translation invariant}, and \ref{thm:main_translation invariant2} on the g-PAM and Theorem~\ref{thm:phi4} on the $\Phi^4_3$ equation actually provide (up to a modification) a.s.\ convergence of solutions together with a quantitative (very bad polynomial) convergence rate/modulus of continuity in $(\eps,\delta)$. Since these results are already lengthy to state we prefer to formulate them simply in terms of convergence in probability as is more common in the SPDE literature.
\end{remark}

\begin{remark}
Given that \cite{Sin23} provides a finite dimensional solution family for the variable coefficient g-PAM and $\Phi^4_3$ equation for any $\eps>0$ as well as continuity in $\eps \in (0,1]$, on might ask whether that solution map extends continuously to $\eps=0$. Note that the results 
of Section~\ref{sec:Homogenisation of (oscillatory)...} imply that the answer is no, see also \cite[Rem.~3.1]{CX23} for a discussion of this point for the $\Phi^4_2$ equation.
\end{remark}

\paragraph{Functionality of our approach}
Let us mention that the approach to homogenisation pursued here has the advantage that it is functional (in the programming sense) and 
separates the homogenisation problem into distinct mathematical subtasks which fit into the rough analysis\slash pathwise approach to singular SPDEs.
\begin{enumerate}
\item 
It takes as input appropriate estimates on the fundamental solutions $\Gamma_\eps$. At small scales these estimates are well known and follow from classical parabolic regularity theory assuming appropriate regularity on the coefficient field. At large scales one needs estimates on the difference between the kernel  $\Gamma_\eps$ and its two scale expansion to an equation dependent order. For the equations considered here, first order two scale expansions are sufficient and we only need a (slight) upgrade from the results already available in the literature. 
\item 
We provide a general abstract fixed point theorem for modelled distributions on a regularity structure where the equation involves lifts of integral operators of the type appearing in the two scale expansion (to arbitrary order),
which is continuous in the homogenisation length scale $\eps\in [0,1]$, the point being that it is also continuous at $\eps = 0$. In applications to specific equations this theorem takes as an input the above mentioned estimates on kernels and 
a model on the regularity structure.
\item The construction of an equation dependent regularity structure and renormalisation group as well as the identification of the renormalised equation then follows mostly along understood lines, c.f.\ \cite{BHZ19, BCCH20, BB21, HS23m, BSS25}, though see Remark~\ref{rem:polynomial-nonlinearities}.
This is carried out only for the specific equations considered, mainly in order to keep the presentation brief and accessible to readers not familiar with these works. 
\item Convergence of models, which we establish for the two example equations by hand, again takes as input appropriate estimates on the fundamental solution.
Here the work is split between a part which is rather regularisation agnostic and follows along similar arguments as the corresponding bounds in \cite{Hai14}, as well as
the identification of the (regularisation dependent) counterterms. The small scale counter terms are identified by `freezing coefficients' following \cite{Sin23}.
In order to identify the explicit `large scale' counterterms for logarithmic divergences we require sharp\footnote{In contrast, to the previous uses of such bounds where there is always some `wiggle room'. 
} error bounds on the homogenisation error for the fundamental solution. Notably, for g-PAM involving derivatives in the nonlinearity we use $L^{p}$-bounds and (optimal) $L^\infty$-estimates do not seem to be sufficient.
\end{enumerate}
This separation into distinct subtasks allows for a streamlined exposition by leveraging established understanding in the area of singular SPDEs, in particular results from \cite{Hai14}, instead of rederiving essentially known SPDEs estimates. This in particular results in Section~\ref{sec:RS} being rather short despite the main concepts of the theory of regularity structures being recalled. 
It furthermore makes it clear that our approach is not limited to the specific setting considered here (for instance the boundedness assumption on correctors can be relaxed). 
However, while the approach itself is quite general, there are several places in its implementation where it would be desirable to develop more precise or more broadly applicable results:
\begin{itemize}
\item It is well understood that the two scale expansion can be taken to arbitrary order in periodic homogenisation, c.f.\ \cite{KMS07}. It would be nice to have corresponding heat kernel estimates, which could then be used as input to the machinery developed here, see Remark~\ref{rem:expanding to different orders}. Alternatively, it would also be interesting to develop alternative kernel free approaches to homogenisation of singular SPDEs (c.f.\ \cite{OW19,BOS25} where such kernel free approaches are implemented in a different setting).
\item It would be desirable to allow for more general initial conditions, in particular in the presence of gradient terms in the nonlinearity, see Remark~\ref{rem:initial condition} and Section~\ref{sec:convergence initial conditions}. 
\item Finally, it would be desirable to have a more explicit descriptions of non-logarithmic divergences, such as for the $\Phi^4_3$ equation in Theorem~\ref{thm:phi4}, see also Remark~\ref{rem:measo div}. As a first step it would already be interesting to have such a description in larger ranges of $(\eps,\delta)\in (0,1]^2$
than $\eps\lesssim \delta^{2}$ in Theorem~\ref{thm:phi4 restricted}. 
\end{itemize}

During the preparation of this work an alternative approach to tackle periodic homogenisation problems for singular SPDEs closer to para-controlled calculus was developed in \cite{CX_paracontrolled}. While that approach at the moment seems to be restricted to spatial homogenisation of equations with symmetric coefficient fields, it otherwise seems to compare to ours analogously to the relationship between the theory of regularity structures \cite{Hai14} and para-controlled calculus \cite{GIP15} for constant coefficient SPDEs.
Of course periodic homogenisation of classical PDEs is a highly developed theory c.f.~\cite{sanchez1980non, BLP78, jikov2012homogenization, bakhvalov2012homogenisation} and we build on already established understanding, in particular, we use estimates from \cite{GS15,GS20}. For the more recently developed theory of quantitative stochastic homogenisation which is related but quite distinct, we refer to \cite{armstrong2019quantitative, josien2022annealed,armstrong2022elliptic, gloria2020regularity}.

%
%
\subsection{Structure of article}
In Section~\ref{sec:Parabolic Operators, Heat Kernels and Homogenisation} we recall mostly known results, starting with estimates on the fundamental solution of uniformly parabolic operators. In Section~\ref{sec:periodic_hom} we recall elements of the theory of periodic homogenisation, in particular the definition of the homogenised matrix $\bar{A}$ and the correctors $\psi_i$ as well as uniform (parabolic) regularity estimates. In Section~\ref{sec:kernel est} we state the homogenisation estimates on the fundamental solution
which will be used as an input when applying our framework, with Theorem~\ref{thm:further kernel estimate} and Lemma~\ref{lem:lp kernel bound} being the only results not directly taken from the existing literature. The auxiliary estimates for the proof of this theorem are outsourced to Section~\ref{sec:aux_est}.

Section~\ref{sec:Homogenisation of (oscillatory)...} states the main results on the (oscillatory) g-PAM and $\Phi^4$ equations, which are obtained by applying the machinery developed in the rest of the article.

Section~\ref{sec:RS} starts by recalling elements of the theory of regularity structures \cite{Hai14}. Section~\ref{sec:A slightly weakened topology on kernels} introduces topologies on spaces of singular kernels $\bfK^\beta_{L, R}$ which are strong enough to be compatible with the theory of regularity structures, but weak enough so that the (post-processed) two scale expansion error of the fundamental solution converges to zero in this topology. Section~\ref{sec:lifting corrector terms} introduces the necessary infrastructure on the regularity structure to lift the corrector terms to the abstract formulations of the SPDE and establishes appropriate Schauder estimates. Section~\ref{sec:An abstract fixed point theorem} provides a general abstract fixed point theorem.
In Section~\ref{sec_post_process} 
we prove that the kernel estimates of Section~\ref{sec:Parabolic Operators, Heat Kernels and Homogenisation} indeed imply convergence 
of the appropriately post-processed two scale expansion error of the fundamental solution in $\bfK^\beta_{L, R}$ for $L,R$ slightly larger than $1$ and $\beta$ slightly less than $2$ (which are the exponents needed for the equations considered here).
Finally, in Section~\ref{sec_post_process} we implement the rewriting sketched in \eqref{eq:formal rewriting} at the level of modelled distributions, which then fits exactly into the setting of the abstract fixed point Theorem~\ref{thm:fixed point}.

In Section~\ref{sec:application gpam} we apply the developed machinery to the g-PAM equation and prove Theorems~\ref{thm:g-pam},\ref{thm:main_translation invariant}, and~\ref{thm:main_translation invariant2}. In Section~\ref{sec:abstract form gpam} we perform the rewriting of the equation as alluded to in \eqref{eq:formal rewriting} and explained in Section~\ref{sec:Homogenisation in RS} and construct the associated regularity structure and model. Identification of the renormalised equation for general renormalisation functions is performed along the usual lines in Section~\ref{sec:renorm eq}. In Section~\ref{eq: convergence gpam} we provide the stochastic estimates on the model in a way agnostic to the specific regularisation and the parts specific to each regularisation are performed in subsections of Section~\ref{sec:identification} which then concludes with the proof of the main results in Subsection~\ref{sec:proof gpam}.

In Section~\ref{sec:application phi} we prove results on the $\Phi^4_3$ equation, Theorems~\ref{thm:phi4} and \ref{thm:phi4 restricted}. While we here only work with one regularisation, we structure the section in the same way as 
Section~\ref{sec:application gpam}, where several regularisation are considered, in order for the interested reader to
be able to adapt the proof to other regularisations.

Appendix~\ref{ap:A} contains the (standard) Lemma~\ref{lem:appendix} about oscillatory functions as well as the less standard Lemma~\ref{lem:ap_osc_kernel} on how such functions interact with renormalised singular kernels.
%

\paragraph{Acknowledgements}
\textit{This article is dedicated to the memory of Kenneth R.\ Zuckerman.}
HS would like to thank Antoine Gloria and David Lee for their hospitality during a visit to Sorbonne University and Scott Armstrong and Nicolas Clozeau for valuable discussions. HS gratefully acknowledges financial support from the Swiss National Science
Foundation (SNSF), grant number 225606, and previously from Ilya Chevyrev’s New Investigator Award EP/X015688/1.
Both authors thank Weijun Xu for discussions on the topic.


\subsection{Notation and function spaces}
We equip $\mathbb{R}^{d+1}$ with the parabolic scaling $(1,\ldots,1,2)$ in the sense of \cite{Hai14} and write $|z|_\fraks= |x|+ |t|^{1/2} $ for $z=(x,t)\in \mathbb{R}^{d+1}$. 
For $\lambda>0$ set
\begin{equ}
\mathcal{S}_{\fraks}^\lambda: \mathbb{R}^{d+1}\to \mathbb{R}^{d+1}, \qquad z=(x,t) \mapsto \mathcal{S}_{\fraks}^\lambda(z)= (\lambda^{-1}x, \lambda^{-2}t)\  .
\end{equ}
For a function $\phi: \mathbb{R}^{d+1}\to \mathbb{R}$ and $z_0\in \mathbb{R}^{d+1}$ we write $\phi_{z_0}^\lambda(z)= \frac{1}{\lambda^{|\fraks|}}\phi\left(\mathcal{S}_{\fraks}^\lambda (z-z_0) \right)$. 
We introduce the backward and forward parabolic cylinders
\begin{equ}\label{eq:backwards parabolic cylinder}
Q_{r}(x,t)=B_r(x)\times (t-r^2,t) \qquad \text{and} \qquad \tilde{Q}_{r}(x,t)=B_r(x)\times (t, t+r^2)\subset \mathbb{R}^{d+1}\;
\end{equ}
at $(x,t)\in \mathbb{R}^d\times \mathbb{R}$ of radius $r>0$.
We denote 
$\triangle:=\left\{(z,\bar{z})\in (\mathbb{R}^{d+1})^{\times 2} \ : \ z= \bar{z} \right\}$.
We shall mostly consider function spaces on $\mathbb{R}^{d+1}$. We shall often (and often freely) identify functions and distributions on $\mathbb{T}^d\times \mathbb{R}$ with their counterpart on $\mathbb{R}^{d+1}$ by pullback under the projection
\begin{equ}
\pi_{d,1}: \mathbb{R}^{d+1}\to \mathbb{T}^d\times \mathbb{R}\, \qquad
(x,t) \mapsto (\pi_d x, t)\ .
\end{equ}

\subsubsection{Functions of one variable}
For a measurable function $f\in L^{0}(Q)$ we shall write
$$\| F\|_{L^p(Q)}:=\left( \int_{Q} |f|^p\right)^{1/p} , \qquad \| F\|_{\avL^p(Q)}:=\left( \fint_{Q} |f|^p\right)^{1/p} \ ,$$
where $\fint_{Q} f :=\frac{1}{|Q|} \int_{Q} f$. 
For $f: \mathbb{R}^{d+1}\to \mathbb{R}$ and a multi-index $k\in \mathbb{N}^{d+1}$, we write $D^{k}f(x,t)= \partial^{k_1}_{x_1}\dots\partial^{k_d}_{x_d} \partial^{k_{d+1}}_{t} f(x,t)$ (whenever this makes sense) while reserving the gradient notation 
$\nabla$ for only spatial derivatives.
For $\gamma>0$ such that $D^{k}f(z_0)$ exists whenever $|k|_\fraks<\gamma$, let
$$ \opP_{z_{0}}^\gamma [f](z):= \sum_{|k|_\fraks<\gamma} \frac{D^{k}f(z_0)}{k!} (z-z_0)^k\ .$$
For $\gamma\in (0,\infty)\setminus \mathbb{N}$, let
$$\|f\|_{C_\fraks^\gamma(Q)}= \sup_{z,\bar{z}\in Q} \frac{|f(z)- \opP_{\bar{z} }^\gamma[f](z)| }{|z-\bar{z}|^\gamma } $$ 
and for $\gamma\in \mathbb{N}$, write $\|f\|_{C_\fraks^\gamma(Q)}= \sup_{|k|_\fraks=\gamma}\sup_{z,\bar{z}\in Q}  |D^kf(z)|$.
Define
$$ \vertiii{f}_{C_\fraks^\gamma(Q)} = \max_{|k|_\fraks<\gamma} \|D^k f\|_{L^\infty(Q)} + \|f\|_{C_\fraks^\gamma(Q)}\ .$$

\begin{lemma}\label{lem:check_higher_holder}
Let $\gamma\in (1,2)$, then for any parabolic cylinder $Q_r$ of radius $r\leq R$ the following estimate holds 
$$ \|f\|_{C_\fraks^\gamma(Q)} \lesssim_R \sup_{(x,t),(x,s) \in Q} \frac{|f(x,t)-f(x,s)|}{|t-s|^{\gamma/2}} + \|\nabla f\|_{C_\fraks^{\gamma-1}(Q)}\ .$$
\end{lemma}
\begin{proof}
Writing
\begin{equ}\label{eq:loc_rewrite holder}
f(x,t)-\opP^{\gamma}_{(y,s)}[f](x,s)= \big(f(x,t)-f(x,s)\big)+ \big(f(x,s)-\opP^{\gamma}_{(y,s)}[f](x,s)\big)
\end{equ}
the claim follows from
$$ f(x,s)-\opP^{\gamma}_{(y,s)}[f](x,s) = f(x,s) - f(y,s) - \langle \nabla f(y,s), x-y \rangle =  \int_{0}^1 \langle \nabla f(y + r (x-y),s  )   -  \nabla f(y,s), x-y \rangle dr \ .$$
\end{proof}

\begin{remark}\label{rem:useful notation norm}
Let $\mathfrak{B}_{R}:= \{ \phi\in C(B_1)\ : \ \vertiii{\phi}_{C_\fraks^R(Q)}\leq 1\}$. For  notational convenience we also let
$$\| \psi\|_{\mathfrak{B}_{R}^\lambda}=\begin{cases}
 \max_{|k|_\fraks<\gamma} \left(\lambda^{|\fraks|+|k|_\fraks} \|D^k \psi \|_{L^\infty} \right) +   \lambda^{|\fraks|+R}\|\psi\|_{C_\fraks^R} & \text{if } \supp \psi \subset B_{\lambda}\\
+\infty & \text{else. }
\end{cases} 
 $$
Note that 
$ \| \psi\|_{\mathfrak{B}_{R}^\lambda} <C $ if and only if $\psi= C \phi^{\lambda}$ for some $\phi\in \mathfrak{B}_{R}$.
\end{remark}

%

\subsubsection{Functions of two variables}

We introduce the following H\"older semi-norms for bounded subsets $B,B'\subset \mathbb{R}^{d+1}$
\begin{equs}
\|F\|_{C^{0,\gamma}_\fraks(B\times B')}&= \sup_{z\in B} \sup_{z',\bar{z}'\in B'}  \frac{|F(z,z') - \opP^\gamma_{\bar{z}'}[F({z}, \,\cdot\,  )](z') |}{ |z'-\bar{z}'|_\fraks^{\gamma} } \ ,\\
\|F\|_{C^{\gamma,0}_\fraks(B\times B')}&=\sup_{z,\bar{z}\in B} \sup_{z'\in B'}  \frac{|F(z,z') - \opP^\gamma_{\bar{z}} [F(\,\cdot\, , {z}')](z) | }{|z-\bar{z}|_\fraks^\gamma} \ ,
\end{equs}
Defining 
$$\opP^{(\gamma,\gamma')}_{({z}_0, \bar{z}_0)}[F](z, \bar{z}) = \sum_{|k|_\fraks<\gamma, |l|_\fraks<\gamma' } \frac{D_1^{k} D_2^{l} F({z}_0, \bar{z}_0) }{k!\ l!} (z-z_0)^k(\bar{z}-\bar{z}_0)^l \ ,$$
we set
\begin{equ}
\|F\|_{C^{\gamma,\gamma'}_\fraks(B\times B')}= \sup_{z,\bar{z}\in B} \sup_{z',\bar{z}'\in B'}  \frac{|F(z,z') -\opP^\gamma_{\bar{z}} [F(\,\cdot\, , {z}')](z) -\opP^{\gamma'}_{\bar{z}'}[F({z}, \,\cdot\,  )](z')+  \opP^{(\gamma,\gamma')}_{(\bar{z}, \bar{z}')}[F](z, {z}')| }{|z-\bar{z}|_\fraks^\gamma  |z'-\bar{z}'|_\fraks^{\gamma'} } \ 
\end{equ}
and
\begin{align*}
 \vertiii{F}_{C^{\gamma,\gamma'}_\fraks(B\times B')} &= \max_{|k|_\fraks<\gamma, |l|_\fraks<\gamma' }\|D_1^{k} D_2^{l} F\|_{L^\infty(B\times B')} + 
\max_{|k|_\fraks<\gamma}\|D^k_1F\|_{C^{0,\gamma'}_\fraks(B\times B')} \\
&\qquad +\max_{|l|_\fraks<\gamma' }\|D^l_2F\|_{C^{\gamma,0}_\fraks(B\times B')} + \|F\|_{C^{\gamma,\gamma'}_\fraks(B\times B')}\ .
\end{align*}

\begin{lemma}\label{lem:higehr hölder criterion 2}
Let $Q, Q'$ be parabolic cylinders. For $\gamma \in (1,2), \gamma ' \in (0,1)$ the 
following estimate holds 
\begin{align*}
\|F\|_{C^{\gamma,\gamma'}_\fraks(Q\times Q')} \lesssim 
\| \nabla_1 F\|_{C^{\gamma-1,\gamma'}_\fraks(Q\times Q')}
&+ 
\sup_{(x,t),(x,s) \in Q} \frac{\|  F(x,t, \,\cdot\, )-  F(x,s, \,\cdot\, )\|_{C_{\fraks}^\gamma(Q')}}{|t-s|^{\gamma'/2}} \ .
\end{align*}
Similarly, for $\gamma, \gamma ' \in (1,2)$  the following estimate holds 
\begin{align*}
\|F\|_{C^{\gamma,\gamma'}_\fraks(Q\times Q')} \lesssim 
&  
 \sup_{(x,t),(x,s) \in Q}\sup_{(x',t'),(x',s') \in Q'}   \frac{|F(x,t,x',t')-F(x,s,x',t') -F(x,t,x',s') +F(x,s,x',s')|}{|t-s|^{\gamma/2}|t'-s'|^{\gamma'/2}}   \\
&+\|\nabla_1 \nabla_2 F\|_{C^{\gamma-1,\gamma'-1}_\fraks(Q\times Q')} 
+
\sup_{(x',t'),(x',s') \in Q'}  \frac{ \| \nabla_1 F(\,\cdot\, ,x',t')- \nabla_1 F(\,\cdot\, ,x',s')\|_{C_{\fraks}^{\gamma-1}(Q)}}{|t'-s'|^{\gamma'/2}} \\
&+ 
\sup_{(x,t),(x,s) \in Q} \frac{ \| \nabla_2 F(x,t, \,\cdot\, )- \nabla_2 F(x,s, \,\cdot\, )\|_{C_{\fraks}^{\gamma'-1}(Q')}}{|t-s|^{\gamma/2}} 
\ .
\end{align*}
\end{lemma}
\begin{proof}
The former inequality follows similarly to the proof of Lemma~\ref{lem:check_higher_holder}.
For the second inequality, introduce $F^{z',\bar{z}'}(z):= F(z,z')- \opP^{\gamma'}_{\bar{z}'} [F(z, \,\cdot\,  )](z')$,
then one finds that
$$F(z,z') -\opP^\gamma_{\bar{z}} [F(\,\cdot\, , {z}')](z) -\opP^{\gamma'}_{\bar{z}'}[F({z}, \,\cdot\,  )](z')+  \opP^{(\gamma,\gamma')}_{(\bar{z}, \bar{z}')}[F](z, {z}')
= F^{z',\bar{z}'}(z)- \opP^\gamma_{\bar z}[F^{z',\bar{z}'}](z)\ .$$
One concludes by rewriting both Taylor remainders as in \eqref{eq:loc_rewrite holder}.
\end{proof}

%


\section{Parabolic Operators, Heat Kernels and Homogenisation}\label{sec:Parabolic Operators, Heat Kernels and Homogenisation}
Throughout this article we make the following assumption.
\begin{assumption}\label{assumption}
 Fix $A: \mathbb{R}^{d+1}\to \mathbb{R}^{d\times d}$, uniformly elliptic, $\mathbb{Z}^{d+1}$-periodic and belonging to $C^{\gamma}$ for some $\gamma>1$.
\end{assumption} 
This regularity assumption\footnote{
which is more commonly written as $C^{1,\eta}$ for some $0<\eta<(\gamma-1)\wedge1 $ in the PDE literature, 
} is not optimal for the estimates in this section, but will be convenient. 
We first recall several well known results on non-constant coefficient heat kernels.
We shall denote by $\Gamma_\eps$ the heat kernel of the differential operator
\begin{equ}
\partial_t +\mathcal{L}_\eps\ \qquad \text{with }\qquad  \mathcal{L}_\eps= - \nabla \cdot (A(x/\eps, t/\eps^2)  \nabla)  \ .
\end{equ}
Writing $\tilde{A}:  \mathbb{R}^{d+1}\to \mathbb{R}^{d\times d}$ for the matrix with entries $\tilde{a}_{i,j}(y,s)={a}_{j,i}(y,-s)$, we set 
$\tilde{\mathcal{L}}_\eps= -\nabla \cdot (\tilde A(x/\eps, t/\eps^2) \nabla)$. One finds that the fundamental solution $\tilde{\Gamma}_\eps(x,t,\zeta,\tau)$ of $\partial_t +\tilde{\mathcal{L}}_\eps$ satisfies 
\begin{equation}\label{eq:adjoint}
\tilde{\Gamma}_\eps (x,t;\zeta,\tau)= \Gamma_\eps (\zeta, -\tau;x, -t)\ .
\end{equation}


\begin{proposition}\label{prop:small scale heat kernel bounds}
There exists $\mu= \mu(A)>0$ such that for $T>0$,
$a,a' \in \mathbb{N}^d,\ b,b'\in \mathbb{N}$ satisfying $|a|+2b\leq 2$ , $|a'|+2b'\leq 2$ ,
\begin{equ}\label{eq:smalls_scale_kernel_bound}
|\partial_x^{a}\partial_\zeta^{a'} \partial_{t}^b \partial_\tau^{b'} \Gamma_1(x,t; \zeta, \tau)|\lesssim_{a,a',b, b',T}(t-\tau)^{-(|a|+|a'| + 2b+2b' + d)/2}  \exp\left(-\frac{\mu |x-\zeta|^2}{t-\tau} \right) \
\end{equ}
uniformly over $x,\zeta\in \mathbb{R}^d$ and $-\infty <\tau<t<\infty$ satisfying $|t-\tau|<T$.
\end{proposition}
\begin{proof}
The estimate for the case $a=b=a'=b'=0$ is well known, c.f.\  \cite[Thm~1.1]{HK04}. 
Next, recall the classical\footnote{The first \eqref{eq:a priori estimate0} is just Corollary~\ref{cor:technical} below with $\eps=1$, while \eqref{eq:a priori estimate} for example follows by combining
Theorem~\ref{thm:technical} and Corollary~\ref{cor:technical} again with $\eps=1$.} interior regularity estimates, stating that for any $p>d+2$ there exists $C>0$ such that for any $r<1$ and any weak solution $u$ of 
\begin{equ}\label{local:eq}
\partial_t u +\mathcal{L}_1 u = \nabla f , \qquad \text{ on } Q_{2r}
\end{equ}
the following inequalities holds 
\begin{equs}
\| u\|_{L^\infty(Q_r)}&\leq C \left(  \| u_\eps \|_{\avL^2(  Q_{2r}(z_0))} + r \|  f \|_{\avL^p(  Q_{2r}(z_0))}\right) \label{eq:a priori estimate0} \\
\|\nabla u\|_{L^\infty(Q_r)} &\leq C \left(\frac{1}{r} \|u\|_{\avL^2(Q_2r)} + \| f\|_{\avL^q(Q_2r)}+  r\|\nabla  f\|_{\avL^q(Q_2r)} \right)\ . \label{eq:a priori estimate}
\end{equs}
Thus, choosing $r= \frac{|t-s|}{8}$ in \eqref{eq:a priori estimate} one concludes that \eqref{eq:smalls_scale_kernel_bound} holds for $a=b=b'=0$ and $a=1$.
By differentiating \eqref{local:eq} one then concludes the remaining cases where all derivatives fall on the first variable.
Looking at the adjoint equation and using the Identity~\eqref{eq:adjoint} one concludes the remaining cases similarly.
\end{proof}
One notes (by uniqueness of fundamental solutions) the scaling property
\begin{equ}\label{eq:formula for Gamma_epsilon}
\Gamma_\eps (z;\bar{z})= \frac{1}{\eps^d} \Gamma_1 (\mathcal{S}^{\eps}_{\fraks} (z); \mathcal{S}_\fraks^\eps\bar{z}) \ .
\end{equ}
which together with Proposition~\ref{prop:small scale heat kernel bounds} implies the following.
\begin{corollary}\label{cor:small scale}
There exists $\mu= \mu(A)>0$ such that for $T>0$,
$a,a' \in \mathbb{N}^d,\ b,b'\in \mathbb{N}$ satisfying $|a|+2b\leq 2$, $|a'|+2b'\leq 2$ ,
$$|\partial_x^{a}\partial_\zeta^{a'} \partial_{t}^b \partial_\tau^{b'} \Gamma_\eps(x,t; \zeta, \tau)|\lesssim_{a,a',b, b',T} (t-\tau)^{-(|a|+|a'| + 2b+2b' + d)/2}  \exp\left( -\frac{\mu |x-\zeta|^2}{t-\tau} \right) \ $$
uniformly over $x,\zeta\in \mathbb{R}^d$, $-\infty <\tau<t<\infty$ satisfying $|t-\tau|<T\epsilon^2$ and $\epsilon\in (0,1]$.
\end{corollary}
%
%

\subsection{Behaviour close to the diagonal}
Let $C_\eps(x,t)= (4\pi)^{-d/2} \det\big(A_s(x/\eps, t/\eps^2)\big)^{-1/2}$ where $A_s$ denotes the symmetric part of $A$  and
\begin{equ}[e:deftheta]
\vartheta_\eps^{(x,t)}(\zeta)= \sum_{i,j} a^{i,j}(x/\eps, t/\eps^2) \zeta_i \zeta _j, \qquad
 w_\eps^{z}(\zeta, \tau)= \frac{\mathbf{1}_{\{\tau > 0 \}}}{\tau^{d/2}} \exp \left(-\frac{\vartheta_\eps^{z}(\zeta)}{4\tau}\right)\ ,
\end{equ}
where $a^{i,j}:= (A^{-1})_{i,j}$.
The fundamental solution of the differential operator with coefficients ``frozen'' at $z=(\zeta, \tau)$ is 
\begin{equation}\label{eq:exlicit Z_0}
Z_{ \eps;0 }(x,t; \zeta,\tau): =  C_\eps(\zeta,\tau) w_\eps^{(\zeta,\tau)}(x-\zeta,t-\tau)\ .
\end{equation}
We further define
$${Z}_{\eps;1}(x,t; \zeta, \tau)\eqdef \sum_{k,i,j}  \eps^{-1}\partial_k a_{i,j}(\zeta/\eps,\tau/\eps^2) C_\eps(\zeta, \tau) \int_{\tau}^t \int_{\mathbb{R}^d} w_\eps^{\zeta,\tau}(x-\eta,t-\sigma)
  (\eta-\zeta)_k \partial_i \partial_j Z_{\eps;0}(\eta, \sigma; \zeta,\tau) \, d\eta\, d\sigma\ . $$
Writing
\begin{equation}\label{eq:reflection}
\mathfrak{R}_\zeta: \mathbb{R}^d\to \mathbb{R}^d,\qquad  x\mapsto -x+2\zeta\ .
\end{equation}
for the reflection map at the point $\zeta\in \mathbb{R}^d$, 
one notes that
$$Z_{\eps;0}(\mathfrak{R}_\zeta(x),t; \zeta, \tau)=Z_{\eps;0}(x,t; \zeta, \tau)\ . $$
and
 \begin{equation}\label{eq:reflection_symmetry}
{Z_{\eps;1}}(\mathfrak{R}_\zeta(x),t; \zeta, \tau)=-{Z}_{\eps;1}(x,t; \zeta, \tau). 
 \end{equation}

\begin{lemma}\label{lem:close to diagonal}
There exists $\mu= \mu(A)>0$ such that for $T>0$,
$a, a' \in \mathbb{N}^d,\ b,b'\in \{0,1\}$ satisfying $a+a' +2b+ 2b'\leq 2$ ,
\begin{equ}
|\nabla_x^a\nabla_x^{a'} \partial_t^b\partial_t^{b'} \left(\Gamma_\eps- Z_{\eps;0}-Z_{\eps;1}\right)(x,t; \zeta, \tau) |\lesssim_{T}(t-\tau)^{-(d-2+|a|+|a'|+2b+2b' )/2}  \exp\left(-\frac{\mu |x-\zeta|^2}{t-\tau} \right) \
\end{equ}
uniformly over $x,\zeta\in \mathbb{R}^d$ and $-\infty <\tau<t<\infty$ satisfying $|t-\tau|<\eps^2 T$.
\end{lemma}

\begin{proof}
The proof follows using the classical parametrix construction of the heat kernel, c.f.\ \cite[Ch.~1\& Ch.~6]{Fri08}.
 Recall that one can explicitly write
\begin{equs}\label{eq:explicit expansion}
\Gamma_\eps(x,t; \zeta, \tau)
&= \sum_{\nu=0}^\infty Z_{\eps;\nu}(x,t; \zeta, \tau) \ ,
\end{equs}
see \cite[Eq.~2.11]{Sin23}, where the $Z_{\eps;\nu}$ for $\nu\notin \{1,2\}$ are explicitly given by \cite[Eq.~2.5]{Sin23} and for $\nu\in \{1,2\}$ in \cite[Lem.~2.5]{Sin23}. Then 
the estimates for the case $\eps=1$ are standard. For $a=a'=b=b'=0$ it follows for example directly from \cite[Lem.~2.7 \& Cor.~2.4.1]{Sin23}. 
In order to prove the case $\eps\in (0,1)$ one argues again by scaling. One directly checks that 

$$C_1\circ \mathcal{S}^\eps=C_\eps, \qquad  \vartheta_1^{ \mathcal{S}^\eps(z)}= \vartheta_\eps^z, \qquad 
 w_1^{\mathcal{S}^\eps(z)}\circ \mathcal{S}^\eps=  \eps^d w_\eps^{z}$$
Thus, it follows that $Z_{\eps;0}= \eps^d Z_{1,0}\circ(\mathcal{S}^\eps\otimes \mathcal{S}^\eps)$. 

Next, let
$$\tilde{w}^{i,j,k}_\eps(x,t;\zeta,\tau )\eqdef \int_{\tau}^t \int_{\mathbb{R}^d} w_\eps^{\zeta,\tau}(x-\eta,t-\sigma)
  (\eta-\zeta)_k \partial_i \partial_j Z_{\eps;0}(\eta, \sigma; \zeta,\tau) \, d\eta\, d\sigma \ ,$$
  then a (slightly tedious) computation shows that $\tilde{w}^{i,j,k}_\eps= \eps^{d+1} \tilde{w}^{i,j,k}_1\circ (\mathcal{S}^\eps\otimes \mathcal{S}^\eps)$, which, since 
${Z}_{\eps;1}(x,t; \zeta, \tau)= \sum_{k,i,j}  \eps^{-1}\partial_k a_{i,j}(\zeta/\eps,\tau/\eps^2) C_\eps(\zeta, \tau)\tilde{w}^{i,j,k}_\eps (x,t;\zeta,\tau )$,
 implies that 
$${Z}_{\eps;1}= {Z}_{1;1}\circ (\mathcal{S}^\eps\otimes \mathcal{S}^\eps)\ .$$
Thus, one concludes the proof by rescaling the case $\eps=1$.

The cases involving derivatives follow by for example differentiating the summands in \eqref{eq:explicit expansion} term by term, using that $\eps^{-1}\lesssim |t-\tau|^{-1/2}$ and the inequalities in \cite[Rem~2.6]{Sin23}.
\end{proof}

Next, recall that the map $\Gamma_{\eps}^{*}(z,z')\eqdef \Gamma_{\eps}(z'; z)$ is the fundamental solution of the formal (space-time) adjoint of $\partial_t - \nabla \cdot A \nabla$. Setting 
\begin{equ}\label{eq:Z^*}
Z^*_{ \eps;0 }(x,t; \zeta,\tau): =  C_\eps(x,t) w_\eps^{(x,t)}(x-\zeta,t-\tau)\ ,
\end{equ}
we observe that it satisfies the same scaling property as $Z_{ \eps;0 }$ above. 
The following is then a direct consequence of Lemma~\ref{lem:close to diagonal}.
\begin{corollary}\label{cor:Z*}
There exists $\mu= \mu(A)>0$ such that for $T>0$,
\begin{equ}
\abs[\Big]{\left(\Gamma_1- Z^*_{1;0}\right)(x,t; \zeta, \tau)} \lesssim_{T}(t-\tau)^{-(d-1)/2}  \exp\left(-\frac{\mu |x-\zeta|^2}{t-\tau} \right) \
\end{equ}
uniformly over $x,\zeta\in \mathbb{R}^d$ and $-\infty <\tau<t<\infty$ satisfying $|t-\tau|< T$.
\end{corollary}
We define 
$Z^\eps_{0;\mu} ( \eta, \tau; x,t ) =
 -\frac{\sum_l \big( A^\eps_{s}(x,t) \big)^{-1}_{\mu,l} (\eta-x)_l}{2(t-\tau )} Z_{\eps,0}^*( \eta, \tau; x,t )  $ for $\mu\in \{1,...,d\}$. 
 Note that 
 \begin{equ}\label{eq:scaling linearisation of derivative}
 Z^\eps_{0;\mu} = \eps^{d+1} Z^1_{0;\mu} \circ (\mathcal{S}^\eps\otimes \mathcal{S}^\eps)\ .
 \end{equ}
 We recall the following bound, which is a combination of Lemma~\ref{lem:close to diagonal} and \cite[Lem.~2.8]{Sin23}.
\begin{lemma}\label{lem:bound on Z_{0;1}}
There exists $\mu= \mu(A)>0$ such that for $T>0$,
\begin{equ}
|\partial_{x_i}\Gamma_1- Z^1_{0;\mu}(x,t; \zeta, \tau) |\lesssim_{T}(t-\tau)^{-d/2}  \exp\left(-\frac{\mu |x-\zeta|^2}{t-\tau} \right) \
\end{equ}
uniformly over $x,\zeta\in \mathbb{R}^d$ and $-\infty <\tau<t<\infty$ satisfying $|t-\tau|< T$.

\end{lemma}

\subsection{Periodic homogenisation}\label{sec:periodic_hom}
For $i=1,\ldots,d$ one defines the \textit{correctors} $\phi_i: \mathbb{T}^d\times \mathbb{R}\to \mathbb{R}$ as the unique (weak) $1$-periodic solution to the cell problem
\begin{equation}\label{eq:correctors}
(\partial_t+\mathcal{L}_1)\phi_i= \nabla \cdot (A(x,t) e_j)\; , \qquad 
\int_{[0,1]^{d+1}} \phi_i(z) dz =0\;,
\end{equation}
where $e_i$ denotes the $i$-th basis vector of $\mathbb{R}^d$, c.f.\ \cite[Sec.~2]{GS17}. Similarly denote by ${\phi}'_i$
the correctors to the adjoint homogenisation problem associated to $\tilde{A}$ and set $\tilde{\phi}_i(y,s)\eqdef {\phi}'_i(y,-s)$. Under Assumption~\ref{assumption} it is a standard consequence of parabolic regularity theory that $\phi_i, \tilde{\phi}_i\in C_\fraks^{\gamma}$ for some $\gamma>2$. 
The \textit{homogenised matrix} $\bar A= \big(\bar a_{i,j}\big)_{i,j=1}^d$, characterised by 
\begin{equation}\label{eq:hom_matrix}
\bar A e_j = \int_{ [0,1]^{d+1}} A(z) \left( e_j+ \nabla \phi_j(z) \right)   dz \;,
\end{equation}
is then seen to be strictly positive definite.
For $j,i=1,\ldots,d$, let
\begin{equ}
{b}_{i,j}\eqdef a_{i,j} +\sum_{k=1}^d a_{i,k}(x,t) \partial_k \phi_j(x,t) - \bar{a}_{i,j} \ .
\end{equ}
as well as $b_{d+1,j}= - \phi_j\ .$
The functions $\{\psi_{k,i,j}\}_{k,i=1; j=1}^{d+1;d},$ characterised by the next Lemma, \cite[Lem.~2.1]{GS17}, are called \textit{flux} or \textit{dual correctors}.
\begin{lemma}
Under Assumption~\ref{assumption}, there exists $\gamma>2$ and $\mathbb{Z}^{d+1}$ periodic functions $\psi_{k,i,j}\in C^\gamma(\mathbb{R}^{d+1})$ satisfying
\begin{equ}
b_{i,j}=   \partial_{t} \psi_{{d+1},i,j} + \sum_{k=1}^{d} \partial_k \psi_{k,i,j}, \qquad \psi_{k,i,j}=-\psi_{k,j,i}\ .
\end{equ}
\end{lemma}

One finds that whenever  
$
(\partial_t- \mathcal{L}_\eps) u_\eps= (\partial_t- \bar{\mathcal{L}}) \bar u
$
on $Q$, the \textit{two scale expansion}
\begin{equ}\label{eq:two_scale}
w_\eps= u_\eps-\bar{u}-\eps\sum_{i=1}^d \phi_i^\eps \partial_i \bar{u} - \eps^2 \sum_{i,j=1^d}  \psi^\eps_{d+1,i,j} \partial_i \partial_j \bar{u} \ , 
\end{equ}
satisfies 
$\partial_{t}+\mathcal{L}_\eps w_\eps= \eps \nabla\cdot F_\eps\ $ on $Q$,
where
\begin{equs}\label{eq:two_scale_error}
F_{\eps, i}(x,t)& = ( a_{i,j}^\eps \phi_k^\eps +\psi_{i,k,j}^\eps )\partial_j \partial_k \bar{u}
+ \eps\psi_{i,d+1,j}^\eps \partial_t \partial_j \bar{u}\\
&\quad +a_{i,j}^\eps (\eps\partial_j \psi_{d+1,l,k}^\eps) \partial_l \partial_k u_0
+\eps a_{i,j}^\eps \psi_{d+1,l,k}^\eps \partial_j \partial_l \partial_k \bar{u} \ ,
\end{equs}
see \cite[Thm~2.2]{GS17}.
We note in particular  that under Assumption~\ref{assumption} one has the bounds
\begin{equs}
\|F\|_{\avL^{p}(Q)} &\lesssim \|\nabla^{2}\bar{u}\|_{\avL^{p}(Q)} + \eps(\|\nabla^{3}\bar{u}\|_{\avL^{p}(Q)}  +\|\nabla \partial_t \bar{u}\|_{\avL^{p}(Q)})\label{Item,bound} \\
\| \nabla F\|_{\avL^{p}(Q)}  &\lesssim \eps^{-1} \|\nabla^{2}\bar{u}\|_{\avL^{p}(Q)}
+\|\nabla^{3}\bar{u}\|_{\avL^{p}(Q)}  +\|\nabla \partial_t \bar{u}\|_{\avL^{p}(Q)} \\
&\qquad +\eps( \|\nabla^{4}u_0\|_{\avL^{p}(Q)}  +\|\nabla^2 \partial_t \bar{u}\|_{\avL^{p}(Q)}) \label{Item,derivative bound}\\
 \| \partial_t F\|_{\avL^{p}(Q)}  &\lesssim \eps^{-2} \|\nabla^{2}\bar{u}\|_{\avL^{p}(Q)} 
+ \eps^{-1}(\|\nabla^{3}\bar{u}\|_{\avL^{p}(Q)}  +\|\nabla \partial_t \bar{u}\|_{\avL^{p}(Q)})
+ \|\nabla^{2}\partial_t \bar{u}\|_{\avL^{p}(Q)}  \\
&\qquad
+ \eps (\|\nabla^{3}\partial_t \bar{u}\|_{\avL^{p}(Q)}   + \|\nabla \partial_t^2 \bar{u}\|_{\avL^{p}(Q)})\;. \label{Item,time derivative bound}
\end{equs}
%
%
%
%
%
%
The following is \cite[Thm.~1.1]{GS15}.
\begin{theorem}\label{thm:technical}
Let $R>0$, $2<p<+\infty$ and assume $A$ is Hölder continuous, uniformly elliptic, bounded and $\mathbb{Z}^{d+1}$ periodic. There exists a constant $C_p=C(R,d,p,A)$ such that for all $\eps\in (0,1]$ and $f=(f_1,\ldots,f_d)\in L^p(Q_{2r}(z_0))$
any solution
 $u_\eps$ to $(\partial_t+\mathcal{L}_\eps)u_\eps= \nabla\cdot f$ on $Q_{2r}(z_0)$ satisfies the bound
 \begin{equ}\label{eq:technical_thm_1}
\| \nabla u_\eps \|_{\avL^p(  Q_{r}(z_0))} \leq C_p \left( \frac{1}{r}  \| u_\eps \|_{\avL^2(  Q_{2r}(z_0))} + \|  f \|_{\avL^p(  Q_{2r}(z_0))}\right)
 \end{equ}
uniformly over $r<R$.
If, furthermore, $p>d+2$ let $\alpha=1-\frac{d+2}{p}$. Then exists a constant $C=C(R,d,p,A)$ such 
\begin{equ}\label{eq:technical_thm_2}
\|u_\eps\|_{C_\fraks^{\alpha}(Q_{r}(z_0))} \leq C \left( \frac{1}{r^\alpha}  \| u_\eps \|_{\avL^2(  Q_{2r}(z_0))} + r^{1-\alpha} \|  f \|_{\avL^p(  Q_{2r}(z_0))}\right)\
\end{equ}
uniformly over $r<R$.
\end{theorem}

\begin{corollary}\label{cor:technical}
Let $R>0$ and $p>d+2$, 
then exists a constant $C=C(R,d,p,A)$ such that for $f=(f_1,\ldots,f_d)\in L^p(Q_{2r}(z_0))$ and
any 
 $u_\eps$ satisfying $(\partial_t+\mathcal{L}_\eps)u_\eps= \nabla\cdot f$ on $Q_{2r}(z_0)$
\begin{equ}
\|u_\eps\|_{L^\infty(Q_{r}(z_0))} \leq C \left(  \| u_\eps \|_{\avL^2(  Q_{2r}(z_0))} + r \|  f \|_{\avL^p(  Q_{2r}(z_0))}\right)\
\end{equ}
uniformly over $r<R$.
\end{corollary}
\begin{proof}
Since
\begin{align*}
|u(z)|&\lesssim \left(\fint_{Q_{r/2}(z)} |u(z)-u(w)|^2 dw\right)^{1/2}  + \|u\|_{\avL^{2}(Q_{r/2}(z))} dy\\
&\lesssim \left(\fint_{Q_{r/2}(z)} |z-w|^{2\alpha} dw\right)^{1/2} \|u\|_{C^\alpha(Q_{r/2}(z) )}  + \|u\|_{\avL^{2}(Q_{r/2}(z))} dy\\
&\lesssim r^\alpha\|u\|_{C^\alpha(Q_{r/2}(z) )}  + \|u\|_{\avL^{2}(Q_{r/2}(z))}, 
\end{align*}
we conclude by Theorem~\ref{thm:technical} that
$$\sup_{z\in Q_{r}(z_0) } |u(z)| \lesssim \sup_{z\in Q_{r}(z_0) }  (\| u_\eps \|_{\avL^2(  Q_{r}(z))} + r \|  f \|_{\avL^p(  Q_{r}(z))})  \lesssim \| u_\eps \|_{\avL^2(  Q_{r}(z_0))} + r \|  f \|_{\avL^p(  Q_{r}(z_0))}\ ,$$
as claimed.
\end{proof}

\subsection{Kernel estimates}\label{sec:kernel est}

The next theorem is well known, c.f.\ \cite[Thm~1.1]{HK04} and \cite[Thm~4.1]{GS15}.
\begin{theorem}\label{thm:uniform}
There exists $\mu>0$ such that 
\begin{equ}
| \Gamma_\eps(x,t;y,s)|\lesssim \frac{1}{|t-s|^{\frac{d}{2}}} \exp\left(-\frac{\mu|x-y|^2}{t-s} \right)
\end{equ}
and $$|\nabla_x \Gamma_\eps(x,t;y,s)|+|\nabla_y \Gamma_\eps(x,t;y,s)|\lesssim \frac{1}{|t-s|^{\frac{d}{2}}} \exp\left(-\frac{\mu|x-y|^2}{t-s} \right)
$$
uniformly over $\eps\in (0,1]$, $x,y\in \mathbb{R}^d$ and $-\infty <s<t<\infty$.
\end{theorem}

We shall write $\bar{\Gamma}$ for the homogenised heat kernel and for $I,J\in \{0,1\}$
\begin{equ}
\bar{\Gamma}^{(I,J)}_\eps(x,t;y,s)\eqdef \Bigl(1+\eps  \mathbf{1}_{\{I>0\}} \sum_{i=1}^d \phi_i^\eps (x,t) \partial_{x_i}\Bigr)
\Bigl(1+ \eps \mathbf{1}_{\{J>0\}}\sum_{j=1}^d \tilde{\phi}_j^\eps (y,s) \partial_{y_j}\Bigr)
 \bar{\Gamma} (x,t; y, s) \ ,
\end{equ}
where ${\phi}_j^\eps(x,t)\eqdef {\phi}_j(x/\eps,t/\eps^2)$ and similarly for $\tilde{\phi}_j$. (Note that ${\bar{\Gamma}^{(0,0)}}= \bar{ \Gamma}$.)

The following is a direct consequence of \cite[Thm.~1.1, Thm.~1.2, Thm.~1.3]{GS20}.
\begin{corollary}\label{cor:pointwise}
There exists $C,\mu >0$ such that for $N,N'\in \{0,1\}$
$$
|\nabla^{N}_{x} \nabla^{N'}_y(\Gamma_\eps-\bar \Gamma^{N,N'})  (x,t; y, s)|\lesssim
\begin{cases}
\frac{\eps}{(t-s)^{\frac{d+1}{2}}} \exp\left( -\frac{\mu|x-y|^2}{t-\tau} \right) & \text{if } N= N'=0, \\
 \frac{\eps\log(2+\eps^{-1} (t-s)^{1/2}) }{(t-s)^{\frac{d+N+N' +1}{2}}} \exp\left( - \frac{\mu|x-y|^2}{t-\tau} \right) & \text{else.}
\end{cases}
$$
 uniformly over $x,y \in \mathbb{R}^d$ and $-\infty< s<t< +\infty$ satisfying 
$|t-s|>\eps^2/2$.
\end{corollary}

We shall need the following further kernel estimate. 
\begin{theorem}\label{thm:further kernel estimate}
There exists
$\mu>0$ such that for $L, R\in [0,2)$ with $R+L<3$ and any $T,\kappa>0$ the following holds.
Uniformly over $x,y \in \mathbb{R}^d$ and $-\infty< s<t< +\infty$ satisfying 
$\eps^2/2<|t-s|<T$ and setting $r = \sqrt{t-s}/8$,
\begin{align*}
\| \Gamma_{\eps}- &\Gamma^{\lfloor L \rfloor ,\lfloor R \rfloor}_\eps\|_{C_\fraks^{L, R}(Q_{r}(x,t)\times\tilde{Q}_{r}(y,s))} 
\lesssim_{L,R,\kappa} \frac{\eps^{A(L,R)}  }{|t-s|^{\frac{d+1+L+R}{2}}} \exp\left(-\frac{\mu|x-y|^2}{t-s} \right) 
\end{align*}
where 
$$
A(L,R)\eqdef
\begin{cases}
1-(L \vee R) & \text{ if } L,R<1, \\
(1-L) \wedge(2-R) -\kappa &\text{ if } L<1, R\geq 1, \\
(1-R) \wedge(2-L) -\kappa & \text{ if } R<1, L\geq 1,\\
3- L - R- \kappa    & \text{ if } L, R \geq 1.
\end{cases}
$$
\end{theorem}
Let us mention that while these estimates are not optimal when $L>1$ or $R>1$ they are sufficient to treat the equations considered here, c.f.\ Remark~\ref{rem:expanding to different orders}. 
\begin{proof}
We consider the different regimes of $R,L$ separately and set $r = \sqrt{t-s}/8$ as in the statement. Since the cases $L,R\in \{0,1\}$ follows directly from Corollary~\ref{cor:pointwise}, we consider the remaining cases. If $R\in(0,1)$ and $L=0$ it follows from 
\cite[Cor.~2.4.1]{HS23per} that 
\begin{align*}
\| \Gamma_{\eps}- \bar \Gamma\|_{C_\fraks^{0, R}(Q_{r}(x,t)\times\tilde{Q}_{r}(y,s))} &= 
\sup_{(\zeta,\tau)\in Q_{r}(x,t) } \| \big(\Gamma_{\eps}- \bar{\Gamma}\big)(\zeta,\tau; \cdot )\|_{C_\fraks^{R}(\tilde{Q}_{r}(y,s))}\\
&
\lesssim\sup_{(\zeta,\tau)\in Q_{r}(x,t) }  \frac{\eps^{1- R}  }{|\tau-s|^{\frac{d+1+R}{2}}} \exp\left(-\frac{\mu|\zeta-y|^2}{\tau-s} \right)\\
&
\lesssim \frac{\eps^{1- R}  }{|t-s|^{\frac{d+1+R}{2}}}\exp\left(-\frac{\mu|x-y|^2}{t-s} \right) \ .
\end{align*}
The case $R=0$, $L\in(0,1)$ follows analogously.  If $L,R\in (0,1)$ the estimate is a direct consequence of \cite[Prop.~2.5]{HS23per}.

We turn to the case $R\in [1,2)$, $L=0$ for which by Lemma~\ref{lem:check_higher_holder} it suffices to bound 
\begin{equ}
\| \nabla_2 \big(\Gamma_{\eps}- \Gamma^{0 ,1}_\eps\big)\|_{C_\fraks^{0,R-1}(Q_{r}(x,t)\times\tilde{Q}_{r}(y,s))}\;,
\end{equ}
as well as
$$ \sup_{(x',t')\in Q_{r}(x,t)}
\sup_{(\zeta, \tau),(\zeta, \tau')\in \tilde{Q}_{r}(y,s)  } 
\frac{| \big(\Gamma_{\eps}- \Gamma^{0 ,1}_\eps\big)(x',t'; \zeta,\tau ) - \big(\Gamma_{\eps}- \Gamma^{0 ,1}_\eps\big)(x',t'; \zeta, \tau' )  |}{|\tau-\tau'|^{R/2}} \;.
 $$
 The required estimate on the former term follows from Lemma~\ref{prop:holder convergence kernel}, while the estimate on the second term follows from Proposition~\ref{prop:time holder convergence kernel}.
 
For the case  $R\in [1,2)$, $L\in (0,1)$, we use Lemma~\ref{lem:higehr hölder criterion 2} which results in having to estimate two terms. The bound on the term involving a spatial gradient follows from
 Proposition~\ref{prop_loc} and the bound on the remaining term from Proposition~\ref{prop:time holder convergence kernel2}.

The case $L\in [1,2)$, $R\in (0,1)$ follows analogously to the case above.

Lastly, we turn to $L,R\in [1,2)$, in which case we use Lemma~\ref{lem:higehr hölder criterion 2} and estimate each resulting term separately. The bound on the term involving two spatial gradients is the content of
Proposition~\ref{prop:double_holder convergence kernel}. The estimate on the term involving no gradients follows by an interpolation the first inequality of Proposition~\ref{prop:time holder convergence kernel2} with Proposition~\ref{prop:double time increments}. Finally, the estimate on the terms involving exactly one spatial gradient is the content of the second inequality of Proposition~\ref{prop:time holder convergence kernel2}.
\end{proof}

\begin{lemma}\label{lem:lp kernel bound}
There exists
$\mu>0$ such that for $\mu\in \{1,...,d\}$, $2<p<\infty$,
\begin{equ}
\Big\|\partial_{\mu} \Gamma_\eps(\cdot;y,s) 
		- \sum_{j}  \big(\mathbf{1}_{\mu=j}+(\partial_{\mu} \phi_j)^\eps \big) 
	\partial_j \bar{\Gamma}(\cdot;y,s) \Big\|_{\avL^p(  Q_{r}(x,t) )}  
\lesssim 
\frac{\eps}{r |t-s|^{d/2+1/2}} \exp \Big( -\mu\frac{|x-y|^2}{|t-s|} \Big)\ ,
\end{equ}
uniformly over $x,y \in \mathbb{R}^d$, $-\infty< s<t< +\infty$ and $r,\eps \in (0,1]$ such that
$ \eps \leq r \leq \sqrt{|t-s|}/4 $.
%
\end{lemma}
\begin{proof}

Let
\begin{equ}
u_{\eps} (\cdot)\eqdef\Gamma_\eps (\cdot \ ;y,s ) \qquad \text{and} \qquad \bar{u}(\cdot)\eqdef \bar{\Gamma}(\cdot \ ;y,s )\ .
\end{equ}
Then, the two scale expansion $w_\eps$ defined in\eqref{eq:two_scale}
satisfies $
\partial_{t}+\mathcal{L}_\eps w^{z}_\eps= \eps \nabla\cdot F_\eps\ ,
$
for $F_\eps$ explicitly given in \eqref{eq:two_scale_error}. Thus, it follows from Lemma~\ref{thm:technical} that
 \begin{equ}
\| \nabla w_\eps \|_{\avL^p(  Q_{r}(x,t))} \leq C_p \left( \frac{1}{r}  \| u_\eps \|_{\avL^2(  Q_{2r}(x,t))} + \eps \|  F_\eps \|_{\avL^p(  Q_{2r}(x,t))}\right)
 \end{equ}
Since
$$
\partial_\mu w_\eps - \Big(\partial_{\mu} \Gamma_\eps(\cdot;y,s) 
		- \sum_{j}  \big(\mathbf{1}_{\mu=j}+(\partial_{\mu} \phi_j)^\eps \big) 
	\partial_j \bar{\Gamma}(\cdot;y,s) \Big) 
	=
\eps \phi^\eps  \partial_\mu \nabla \bar{u} - \eps^{2} \partial_\mu \left( \psi^\eps\nabla^2 \bar{ u} \right)$$
and by \eqref{Item,bound} it remains to bound
\begin{align*}
 \frac{1}{r}  \| u_\eps \|_{\avL^2(  Q_{2r}(x,t))} + 
 \max_{r'\in \{r,2r \}} \Big(
\eps\|\nabla^2 \bar{u}\|_{\avL^p(  Q_{r}(x,t))} +\eps^2 \|\nabla^3 \bar{u}\|_{\avL^p(  Q_{r}(x,t))} + \eps^2 \|\nabla \partial_t \bar{u}\|_{\avL^p(  Q_{r}(x,t))}\Big) \ .
\end{align*}
To estimate the first term use Theorem~\ref{thm:uniform}, for the remaining terms use that 
 that $ \eps \leq r \leq \sqrt{|t-s|}/4 $ together with the fact that 
\begin{itemize}
\item $\|\nabla^2\bar{u}\|_{\avL^p(  Q_{r'}(x,t))} \lesssim  \frac{1}{ |t-s|^{d/2 +1}} \exp \Big( -\mu\frac{|x-y|^2}{|t-s|} \Big)$
\item $\|\nabla^3 \bar{u}\|_{\avL^p(  Q_{'r}(x,t))}  +\|\nabla \partial_t \bar{u}\|_{\avL^p(  Q_{r'}(x,t))} \lesssim   \frac{1}{ |t-s|^{d/2 +3/2}} \exp \Big( -\mu\frac{|x-y|^2}{|t-s|} \Big)$
\end{itemize}
for $r'\in \{r,2r \}$.
\end{proof}

%

\subsection{Auxiliary estimates for the proof of Theorem~\ref{thm:further kernel estimate}}\label{sec:aux_est}
\begin{lemma}\label{prop:holder convergence kernel}
There exists 
$\mu>0$ such that for  $\alpha\in [0,1) $  and setting $r = \sqrt{t-s}/8$,
\begin{align*}
&\|  \nabla_2 \left(\Gamma_{\eps}- \Gamma^{1,1}_\eps\right)( \ \cdot\ , y,s)\|_{C_\fraks^{\alpha}(Q_{r}(x,t))}
+\| \nabla_1  \left(\Gamma_{\eps}- \Gamma^{1,1}_\eps\right)(x,t, \ \cdot \ )\|_{C_\fraks^{\alpha}(\tilde{Q}_{r}(y,s))}\\
&\qquad \qquad\qquad \qquad\qquad \qquad
 \lesssim_{\alpha} \frac{\eps\log(2+\eps^{-1} (t-s)^{1/2}) }{|t-s|^{\frac{d+2+\alpha}{2}}} \exp\left(-\frac{\mu|x-y|^2}{t-s} \right)\;,
\end{align*}
as well as for  $0<\alpha <1 $ and any $\kappa>0$
\begin{align*}
&\| \nabla_1 \nabla_2 \left(\Gamma_{\eps}- \Gamma^{1,1}_\eps\right)( \ \cdot\ , y,s)\|_{C_\fraks^{\alpha}(Q_{r}(x,t))}
+\| \nabla_1 \nabla_2 \left(\Gamma_{\eps}- \Gamma^{1,1}_\eps\right)(x,t, \ \cdot \ )\|_{C_\fraks^{\alpha}(\tilde{Q}_{r}(y,s))}\\
&\qquad \qquad\qquad \qquad\qquad \qquad
 \lesssim_{{\alpha},\kappa,T} \frac{\eps^{1-\alpha -\kappa} }{|t-s|^{\frac{d+3+\alpha}{2}}} \exp\left(-\frac{\mu|x-y|^2}{t-s} \right)\;,
\end{align*}
uniformly over $x,y\in \mathbb{R}^d$, $\eps\in (0,1]$ and $-\infty <s<t<\infty$ with $\eps^2<|t-s|<T$.
\end{lemma}

\begin{proof}
We shall only prove the estimates on the first term of each inequality, the latter ones follows similarly by considering the adjoint problem.
Let $r\eqdef \sqrt{ |t-s|}/8$, $z=(y,s)$ and $z_{0}=(x,t)$ as well as  
\begin{equ}\label{eq:ansatz}
u_{\eps}^{z}(\cdot)\eqdef \partial_{y_i}\Gamma_\eps (\cdot \ ;z ) \qquad \text{and} \qquad \bar{u}^{z}(\cdot)\eqdef \partial_{y_i} \bar{\Gamma}^{(0,1)}(\cdot \ ;z )\ .
\end{equ}
Then, the two scale expansion 
\begin{equ}\label{eq:two_scale_specific}
w^{z}_\eps= u^{z}_\eps-\bar{u}^{z}-\eps\sum_{i=1}^d \phi_i^\eps \partial_i u^{z}_0 - \eps^2 \sum_{i,j=1^d}  \psi^\eps_{d+1,i,j} \partial_i \partial_j u^{z}_0 \ , 
\end{equ}
satisfies
\begin{equ}\label{eq:parabolic}
\partial_{t}+\mathcal{L}_\eps w^{z}_\eps= \eps \nabla\cdot F^{z}_\eps\ ,
\end{equ}
for $F^{z}$ explicitly given in \eqref{eq:two_scale_error}. For the remainder of the proof we shall suppress the superscript $z$.

It follows from Theorem~\ref{thm:technical} that for $p>1$ such that $\alpha= 1-\frac{d+1}{p}$
\begin{equs}
\|w_\eps\|_{C_\fraks^{{\alpha}}(Q_{r}(z_0))} &\lesssim_{\bar{\alpha}} \frac{1}{r^{\bar{\alpha}}}  \| w_\eps \|_{\avL^2(  Q_{2r}(z_0))} + r^{1-\bar{\alpha}} 
\eps \|  F_\eps  \|_{\avL^p(  Q_{2r}(z_0))} \\
&\lesssim   \frac{\eps\log(2+\eps^{-1} (t-s)^{1/2})}{(t-s)^{\frac{d+2+\alpha}{2}}} \exp\left( \frac{\mu|x-y|}{t-\tau} \right) \ , \label{eq:local_wbound}
\end{equs}
where we used Corollary~\ref{cor:pointwise} and \eqref{Item,derivative bound} in the last line.
In view of \eqref{eq:two_scale_specific}
the first inequality follows.

Differentiating \eqref{eq:parabolic} one finds for $v^{z}_{\eps}\eqdef \nabla w^z$
\begin{equ}
\partial_{t}+\mathcal{L}_\eps v_{j;\eps}= \eps \nabla\cdot  \partial_j F_\eps + \nabla\cdot \left(\partial_j A^{\eps} v \right)  \ ,
\end{equ}
and thus by Theorem~\ref{thm:technical} for any $0<\bar{\alpha}<1$ and $p=\frac{d+2}{1-\alpha}$
\begin{equs}
\|v_\eps\|_{C_\fraks^{\bar{\alpha}}(Q_{r}(z_0))} &\lesssim_{\bar{\alpha}} \frac{1}{r^{\bar{\alpha}}}  \| v_\eps \|_{\avL^2(  Q_{2r}(z_0))} + r^{1-\bar{\alpha}} 
\left(
\eps \| \nabla  F_\eps  \|_{\avL^p(  Q_{2r}(z_0))}  +  \eps^{-1} \| v \|_{\avL^p(  Q_{2r}(z_0))}  \right)\\
&\lesssim  r^{1-\bar{\alpha}}\left( \eps^{-1}  \| v_\eps \|_{\avL^p(  Q_{2r}(z_0))} + 
\eps \| \nabla F_\eps  \|_{\avL^p(  Q_{2r}(z_0))} \right)\;, \label{eq:loc_estimate}
\end{equs}
where we used that $r\geq\eps/2$.
Since, on the one hand by \eqref{Item,derivative bound}
$$\eps \| \partial_j F_\eps  \|_{\avL^p(  Q_{2r}(z_0))} \lesssim  \sum_{n=2}^4  \|\nabla^{n}\bar{u}\|_{\avL^p( Q_{2r}(z_0))}  + \sum_{n=1}^2 \|\nabla^{n} \partial_t \bar{u}\|_{\avL^p( Q_{2r}(z_0))}\;, $$
and on the other hand 
\begin{equ}\label{eq:local_v_def}
v_\eps= \nabla w_\eps = \nabla_1 \nabla_2 \left(\Gamma_{\eps}- \Gamma_\eps^{1,1}\right)(\cdot, z) - \eps \phi^\eps  \nabla^{2} \bar{u} - \eps^{2} \nabla \left( \psi^\eps\nabla^2 \bar{ u} \right)
\end{equ}
implies
$$ \| v_\eps \|_{\avL^p(  Q_{2r}(z_0))} \lesssim \| \nabla_1 \nabla_2 \left(\Gamma_{\eps}- \Gamma^{1,1}\right) \|_{\avL^p(  Q_{2r}(z_0))}
+ \eps \sum_{n=2}^3 \| \nabla^{n} \bar{u} \|_{\avL^p(  Q_{2r}(z_0))} \;,$$
we conclude that the inequality~\eqref{eq:loc_estimate} with $r= \frac{\sqrt{t-s}}{8}$ together with Corollary~\ref{cor:pointwise} imply that
\begin{equs}
\|v_\eps\|_{C_\fraks^{\bar{\alpha}}(Q_{r}(z_0))}\lesssim_{\bar{\alpha}} & |t-s|^{\frac{1-\bar \alpha}{2}}
\left(\frac{\log(2+\eps^{-1} (t-s)^{1/2}) }{(t-s)^{\frac{d+3 }{2}}} + \frac{1 }{(t-s)^{\frac{d+4}{2}}} \right)\exp\left( \frac{\mu|x-y|}{t-\tau} \right)\\
&\lesssim 
 \frac{\log(2+\eps^{-1} (t-s)^{1/2})}{(t-s)^{\frac{d+3+\bar{\alpha}}{2}}} \exp\left( \frac{\mu|x-y|}{t-\tau} \right)\ . \label{eq:loc_eq2}
\end{equs}
Thus, by \eqref{eq:local_v_def} 
\begin{align*}
\|\nabla_1 \nabla_2 \left(\Gamma_{\eps}- \Gamma^{1,1}_\eps\right)( \ \cdot\ , y,s)\|_{C_\fraks^{\bar{\alpha}}(Q_{r}(z_0))}\lesssim_{\bar{\alpha}}  
 \frac{\log(2+\eps^{-1} (t-s)^{1/2})}{(t-s)^{\frac{d+3+\bar{\alpha}}{2}}} \exp\left( \frac{\mu|x-y|}{t-\tau} \right)\ .
\end{align*}
which by interpolation with Corollary~\ref{cor:pointwise} we conclude the proof.
\end{proof}

\begin{prop}\label{prop_loc}
For any $T>0$ there exists 
$\mu>0$ such that for $\alpha, \alpha'\in [0,1)$ and any $\kappa,T>0$, setting
$r = \sqrt{t-s}/8$,
\begin{align*}
&\| \nabla_1  \left(\Gamma_{\eps}- \Gamma^{1,1}_\eps\right)\|_{C_\fraks^{\alpha, \alpha'}(Q_{r}(x,t)\times\tilde{Q}_{r}(y,s))}
+\|  \nabla_2 \left(\Gamma_{\eps}- \Gamma^{1,1}_\eps\right)\|_{C_\fraks^{\alpha', \alpha}(Q_{r}(x,t)\times\tilde{Q}_{r}(y,s))} \\
&\lesssim_{{\alpha},\alpha',\kappa,T}\frac{\eps^{1-\alpha-\kappa} }{|t-s|^{\frac{d+2+\alpha+ \alpha'}{2}}} \exp\left(-\frac{\mu|x-y|^2}{t-s} \right)\;,
\end{align*}
uniformly over $x,y\in \mathbb{R}^d$, $\eps\in (0,1]$ and $-\infty <s<t<\infty$ with $\eps^2<|t-s| <T$.
\end{prop}
%

\begin{proof}
We shall again only prove the estimates on the first term.
Let $r\eqdef \sqrt{ |t-s|}/8$, $z=(y,s)$ and $z_{0}=(x,t)$ as well as  
\begin{equ}\label{eq:ansatz_primitive}
u_{\eps}^{z}(\cdot)\eqdef\Gamma_\eps (\cdot \ ;z ) \qquad \text{and} \qquad \bar{u}^{z}(\cdot)\eqdef \bar{\Gamma}^{(0,1)}(\cdot \ ;z )\ .
\end{equ}
(Note the distinction to \eqref{eq:ansatz})
For $z', \bar{z}'\in \tilde{Q}_{r}(y,s)$, using the notation of \eqref{eq:ansatz}, \eqref{eq:two_scale} set 
$$ u_\eps^{z', \bar{z}'}= u_\eps^{z'}-u_\eps^{\bar{z}'}, 
\qquad \bar{u}_{\eps}^{z', \bar{z}'} =\bar{u}_{\eps}^{z'}-\bar{u}_{\eps}^{\bar{z}'},
\qquad
w_{\eps}^{z', \bar{z}'}= w_{\eps}^{z'} -w_{\eps}^{ \bar{z}'},
\qquad
F^{z', \bar{z}'}=F^{z'} -F^{ \bar{z}'} 
$$
as well as $v_{\eps}^{z', \bar{z}'}=\nabla w_{\eps}^{z', \bar{z}'}$. It follows by Theorem~\ref{thm:technical} that for any $ \alpha' <1$
\begin{equs}\label{eq:local_difference}
\|v_\eps^{z', \bar{z}'}\|_{C_\fraks^{\bar{\alpha}}(Q_{r}(z_0))} 
&\lesssim  r^{1-\bar{\alpha}}\left( \eps^{-1}  \| v^{z', \bar{z}'}_\eps \|_{\avL^p(  Q_{2r}(z_0))} + 
\eps \| \partial_j F^{z', \bar{z}'}_\eps  \|_{\avL^p(  Q_{2r}(z_0))} \right). 
\end{equs}
Note that in view of \eqref{Item,derivative bound} together with an interpolation of \eqref{Item,bound} with \eqref{Item,time derivative bound} one finds that
\begin{equ}\label{eq:loc_to_insert}
\sup_{z', \bar{z}'\in \tilde{Q}_{r}(y,s)} \frac{\| \partial_j F^{z', \bar{z}'}_\eps  \|_{\avL^p(  Q_{2r}(z_0))}}{|z'-\bar{z}'|^{{\alpha}'}} \lesssim
\frac{\eps^{-1}}{(t-s)^{\frac{d+3+{\alpha}'}{2}}} \exp\left( \frac{\mu|x-y|}{t-\tau} \right)\;,
\end{equ}
for $\eps^2<|t-\tau|$.
Similarly, recalling the definition of $v_\eps^{z', \bar{z}'}$, it follows from 
Proposition~\ref{prop:holder convergence kernel} that for any $ {\alpha}'<1$
\begin{equ}\label{eq:loc_to_insert2}
\sup_{z', \bar{z}'\in \tilde{Q}_{r}(y,s)} \frac{\| v^{z', \bar{z}'}_\eps \|_{\avL^p(  Q_{2r}(z_0))}}{|z'-\bar{z}'|_\fraks^{{\alpha}'}} 
\lesssim \frac{\eps\log(2+\eps^{-1} (t-s)^{1/2})  }{|t-s|^{\frac{d+3+{\alpha}' }{2}}} \exp\left(-\frac{\mu|x-y|^2}{t-s}\right) \ .
\end{equ}
Finally, inserting \eqref{eq:loc_to_insert} and \eqref{eq:loc_to_insert2} into \eqref{eq:local_difference} gives
$$\sup_{z', \bar{z}'\in \tilde{Q}_{r}(y,s)}\frac{\|v_\eps^{z', \bar{z}'}\|_{C_\fraks^{\bar{\alpha}}(Q_{r}(z_0))} }{|z'-\bar{z}'|^{{\alpha}'}} \lesssim 
\frac{\log(2+\eps^{-1} (t-s)^{1/2})  }{|t-s|^{\frac{d+2+\bar{\alpha}+ {\alpha}' }{2}}} \exp\left(-\frac{\mu|x-y|^2}{t-s}\right) 
$$
which in turn, unraveling the definition of $v_\eps$ implies 
\begin{align*}
\| \nabla_1  \left(\Gamma_{\eps}- \Gamma^{1,1}_\eps\right)\|_{C_\fraks^{\bar{\alpha}, \alpha'}(Q_{r}(x,t)\times\tilde{Q}_{r}(y,s))}
&\lesssim_{\bar{\alpha}}\frac{\log(2+\eps^{-1} (t-s)^{1/2}) }{|t-s|^{\frac{d+2+\bar{\alpha}+ \alpha'}{2}}} \exp\left(-\frac{\mu|x-y|^2}{t-s} \right)\;.
\end{align*}
Finally, interpolation with the first inequality of Prop.~\ref{prop:holder convergence kernel} completes the proof.
\end{proof}

\begin{prop}\label{prop:double_holder convergence kernel}
There exists 
$\mu>0$ such that for $\alpha, \alpha'\in [0,1)$ and $T,\kappa>0$, setting $r = \sqrt{t-s}/8$,
\begin{align*}
\| \nabla_1 \nabla_2 &\Gamma_{\eps}-\nabla_1 \nabla_2 \Gamma^{1,1}_\eps\|_{C_\fraks^{\alpha, \alpha'}(Q_{r}(x,t)\times\tilde{Q}_{r}(y,s))}\\
&
 \lesssim_{\alpha,\alpha',\kappa} \frac{\eps^{1- \alpha-\alpha'-\kappa }} {|t-s|^{\frac{d+3+\alpha+ \alpha'}{2}}} \exp\left(-\frac{\mu|x-y|^2}{t-s} \right)\;,
\end{align*}
uniformly over $x,y\in \mathbb{R}^d$, $\eps\in (0,1]$ and $-\infty <s<t<\infty$ with $\eps^2<|t-s|<T $.
\end{prop}
\begin{proof}
Let $r\eqdef \sqrt{ |t-s|}/8$ and $z', \bar{z}'\in \tilde{Q}_{r}(y,s)$ and set again as in \eqref{eq:ansatz}
\begin{equ}\label{eq:ansatznew}
u_{\eps}^{z}(\cdot)\eqdef \partial_{y_i}\Gamma_\eps (\cdot \ ;z ) \qquad \text{and} \qquad \bar{u}^{z}(\cdot)\eqdef \partial_{y_i} \bar{\Gamma}^{(0,1)}(\cdot \ ;z )
\end{equ}
 as well as
$$ u_\eps^{z', \bar{z}'}= u_\eps^{z'}-u_\eps^{\bar{z}'}, 
\qquad \bar{u}_{\eps}^{z', \bar{z}'} =\bar{u}_{\eps}^{z'}-\bar{u}_{\eps}^{\bar{z}'},
\qquad
w_{\eps}^{z', \bar{z}'}= w_{\eps}^{z'} -w_{\eps}^{ \bar{z}'},
\qquad
F^{z', \bar{z}'}=F^{z'} -F^{ \bar{z}'} 
$$
as well as $v_{\eps}^{z', \bar{z}'}=\nabla w_{\eps}^{z', \bar{z}'}$. 
%
%
%
As in the proof of Proposition~\ref{prop_loc} for $\alpha'<\bar{\alpha}' <1$
\begin{equs}\label{eq:local_difference2}
\|v_\eps^{z', \bar{z}'}\|_{C_\fraks^{\bar{\alpha}}(Q_{r}(z_0))} 
&\lesssim  r^{1-\bar{\alpha}}\left( \eps^{-1}  \| v^{z', \bar{z}'}_\eps \|_{\avL^p(  Q_{2r}(z_0))} + 
\eps \| \partial_j F^{z', \bar{z}'}_\eps  \|_{\avL^p(  Q_{2r}(z_0))} \right). 
\end{equs}
but this time due to the extra derivative in the second variable in \eqref{eq:ansatznew}
$$\sup_{z', \bar{z}'\in \tilde{Q}_{r}(y,s)} \frac{\| \partial_j F^{z', \bar{z}'}_\eps  \|_{\avL^p(  Q_{2r}(z_0))}}{|z'-\bar{z}'|^{\bar{\alpha}'}} \lesssim
\frac{\eps^{-1}}{(t-s)^{\frac{d+4+\bar{\alpha}'}{2}}} \exp\left( \frac{\mu|x-y|}{t-\tau} \right)\ .
 $$
It follows from the second inequality of Proposition~\ref{prop:holder convergence kernel} that for any $ \kappa'>0$
$$\sup_{z', \bar{z}'\in \tilde{Q}_{r}(y,s)} \frac{\| v^{z', \bar{z}'}_\eps \|_{\avL^p(  Q_{2r}(z_0))}}{|z'-\bar{z}'|_\fraks^{{\alpha}'}} 
\lesssim \frac{\eps^{1-{\alpha'}-\kappa'}}{|t-s|^{\frac{d+3+{\alpha}' }{2}}} \exp\left(-\frac{\mu|x-y|^2}{t-s}\right) \ .
$$
Thus, by \eqref{eq:local_difference2}
\begin{align*}
\sup_{z', \bar{z}'\in \tilde{Q}_{r}(y,s)}  \frac{\|v_\eps^{z', \bar{z}'}\|_{C_\fraks^{\bar{\alpha}}(Q_{r}(z_0))} }{|z'-\bar{z}'|^{\bar{\alpha}'}} 
&\lesssim  r^{1-\bar{\alpha}}\left( 
\frac{\eps^{-\bar{\alpha}'-\kappa' } }{|t-s|^{\frac{d+3+\bar{\alpha}' }{2}}} 
 + 
\frac{1}{(t-s)^{\frac{d+4+\bar{\alpha}'}{2}}} 
\right)
\exp\left(-\frac{\mu|x-y|^2}{t-s}\right)\\
&\lesssim  \left( 
\frac{\eps^{-\bar{\alpha}'-\kappa' } }{|t-s|^{\frac{d+2+\bar{\alpha}+\bar{\alpha}' }{2}}} 
 + 
\frac{1}{(t-s)^{\frac{d+3+ \bar{\alpha}+\bar{\alpha}'}{2}}} 
\right)
\exp\left(-\frac{\mu|x-y|^2}{t-s}\right)\\
&\lesssim 
\frac{\eps^{-\bar{\alpha}'-\kappa' } }{|t-s|^{\frac{d+3+\bar{\alpha}+\bar{\alpha}' }{2}}} 
\exp\left(-\frac{\mu|x-y|^2}{t-s}\right)\ .
\end{align*}
Since
$$v^{z'}(z)= \nabla_{x} \nabla_{y} \left(\Gamma_{\eps}- \Gamma^{1,1}\right)(z,z') - \eps \phi^\eps(z)  \nabla^{2} \bar{u}^{z'} - \eps^{2} \nabla \left( \psi^\eps(z)\nabla^2 \bar{ u}^{z'}(z) \right)$$
we find that
\begin{align*}
\| \nabla_1 \nabla_2 \Gamma_{\eps}-\nabla_1 \nabla_2 \Gamma^{1,1}_\eps &\|_{C_\fraks^{\bar \alpha, \bar{\alpha}'}(Q_{r}(x,t)\times\tilde{Q}_{r}(y,s))}\\
&
\lesssim  
\frac{\eps^{-\bar{\alpha}'-\kappa' } }{|t-s|^{\frac{d+3+\bar{\alpha}+\bar{\alpha}' }{2}}} 
\exp\left(-\frac{\mu|x-y|^2}{t-s}\right). 
\end{align*}
On the other hand it follows from the second inequality of Proposition~\ref{prop:holder convergence kernel} that 
\begin{align*}
\| \nabla_1 \nabla_2 \Gamma_{\eps}-\nabla_1 \nabla_2 \Gamma^{1,1}_\eps &\|_{C_\fraks^{0, \mathring{\alpha}'}(Q_{r}(x,t)\times\tilde{Q}_{r}(y,s))}\\
&
\lesssim  
\frac{\eps^{1-\mathring{\alpha}'-\kappa' } }{|t-s|^{\frac{d+3+\mathring{\alpha}' }{2}}} 
\exp\left(-\frac{\mu|x-y|^2}{t-s}\right). 
\end{align*}
Choosing $\bar \alpha = 1-\kappa'$ for $\kappa'>0$ small enough and writing $\alpha,\alpha'\in (0,1)$ as
$$(\alpha,\alpha')=\lambda\cdot (1-\kappa', \bar{\alpha}') + (1-\lambda)\cdot (0,\mathring{\alpha}')$$ 
it follows by interpolation that 
\begin{align*}
\| \nabla_1 \nabla_2 \Gamma_{\eps}-\nabla_1 \nabla_2 \Gamma^{1,1}_\eps &\|_{C_\fraks^{\bar \alpha, \bar{\alpha}'}(Q_{r}(x,t)\times\tilde{Q}_{r}(y,s))}\\
&
\lesssim  
\frac{\eps^{(-\bar{\alpha}'-\kappa')\lambda } \eps^{(1-\mathring{\alpha}-\kappa')(1-\lambda)} }{|t-s|^{\frac{d+3+\bar{\alpha}+\bar{\alpha}' }{2}}} 
\exp\left(-\frac{\mu|x-y|^2}{t-s}\right). 
\end{align*}
which since 
$ (-\bar{\alpha}'-\kappa')\lambda +(1-\mathring{\alpha}-\kappa')(1-\lambda) = 1-\lambda-(\lambda \bar{\alpha}' +(1-\lambda) \mathring{\alpha} ) -\kappa'
= 1 - \frac{\alpha}{1-\kappa'} -\alpha' -\kappa'$
concludes the proof by choosing $\kappa'$ small enough.
\end{proof}
%
%

\begin{lemma}\label{lem:second_derivative_kernel estimate}
There exists 
$\mu>0$ such that for $p,T<\infty$
\begin{align*}
\|& \nabla_1^2 (\Gamma_{\eps}- \Gamma^{1,0}) (\,\cdot\, ,y,s) \|_{\avL^{p}(Q_{\sqrt{ |t-s|}/8}(x,t)} +\| \nabla_2^2 (\Gamma_{\eps}- \Gamma^{0,1}) (y,s, \,\cdot\, ) \|_{\avL^{p}(Q_{\sqrt{ |t-s|}/8}(x,t)} \\
&\lesssim_{p,T} \frac{\log(2+\eps^{-1} (t-s)^{1/2}) }{|t-s|^{\frac{d+3}{2}}} \exp\left(-\frac{\mu|x-y|^2}{t-s} \right)\;,
\end{align*}
uniformly over $x,y\in \mathbb{R}^d$, $\eps\in (0,1]$, and $-\infty <s<t<\infty$ with $\eps^2<|t-s|<T$.
\end{lemma}
\begin{proof}
Let $u_\eps= \Gamma_{\eps}( \cdot , y,s )$ and $\bar{u}= \bar{\Gamma}( \cdot , y,s ) $ as well as
$ w_\eps $ as in \eqref{eq:two_scale} and $F$ as in \eqref{eq:two_scale_error}.
Then $v=\nabla w_\eps$ satisfies
$$\partial_{t}+\mathcal{L}_\eps v_{j;\eps}= \eps \nabla\cdot  \partial_j F_\eps + \nabla\cdot \left(\partial_j A^{\eps} v \right)  \ ,$$
and the remainder of the proof follows the argument of the proof of the second inequality of Lemma~\ref{prop:holder convergence kernel}, but using \eqref{eq:technical_thm_1} instead of \eqref{eq:technical_thm_2}.
\end{proof}

\begin{prop}\label{prop:time holder convergence kernel}
There exists 
$\mu>0$ such that for any $1<\alpha<2$ and $T,\kappa>0$
\begin{align*}
&\sup_{ t - |t-s|/8<t'<t } \frac{| \left(\Gamma_{\eps}- \Gamma^{1,0}_\eps\right)( x,t , y,s) -\left(\Gamma_{\eps}- \Gamma^{1,0}_\eps\right)( x,t' , y,s) |}
{|t-t'|^{\alpha/2}}\\
&+\sup_{ s < s'<s+ |t-s|/8} \frac{| \left(\Gamma_{\eps}- \Gamma^{0,1}_\eps\right)( x,t , y,s) -\left(\Gamma_{\eps}- \Gamma^{0,1}_\eps\right)( x,t , y,s') |}
{|s-s'|^{\alpha/2}}\\
&\qquad \qquad\qquad \qquad\qquad \qquad
 \lesssim_{\alpha,\kappa, T} \frac{\eps^{2-\alpha-\kappa}  }{|t-s|^{\frac{d+1+2\alpha}{2}}} \exp\left(-\frac{\mu|x-y|^2}{t-s} \right)\;,
\end{align*}
uniformly over $x,y\in \mathbb{R}^d$, $\eps\in (0,1]$, and $-\infty <s<t<\infty$ with $\eps^2<|t-s|<T$.
\end{prop}

\begin{proof}
Let $u_\eps= \Gamma_{\eps}( \cdot , y,s )$ and $\bar{u}= \bar{\Gamma}( \cdot , y,s ) $ as well as
$ w_\eps $ as in \eqref{eq:two_scale} and $F$ as in \eqref{eq:two_scale_error}.
Then $\sigma_\eps (x,t) = \partial_t w(x,t)$ satisfies
$$
(\partial_t + \mathcal{L}_\eps) \sigma_\eps =  \nabla \cdot \left(\eps \partial_t F + (\partial_t A^\eps) \nabla w_\eps  \right) \ .
$$
Thus, by Corollary~\ref{cor:technical} choosing $\eps=r$
\begin{equs}
\|\sigma\|_{L^\infty(Q_r)} &\lesssim \|\sigma \|_{\avL^2(Q_{2r})} + \eps r\left( \| \partial_t F\|_{\avL^p(Q_{2r})} + \eps^{-2} \| \nabla w_\eps \|_{\avL^p(Q_{2r})}\ \right) \\
&\lesssim \|\sigma \|_{\avL^2(Q_{2r})} + \eps^2 \| \partial_t F\|_{\avL^p(Q_{2r})} + \| \nabla w_\eps \|_{\avL^p(Q_{2r})}\\
&\lesssim \eps^{-1}\| \nabla w\|_{\avL^2(Q_{2r}) } + \| \nabla w_\eps \|_{\avL^p(Q_{2r})} 
+ \| \nabla^2 w\|_{\avL^2(Q_{2r}) } \label{eq:local_time_est_1} \\
&\qquad + \eps \| \nabla \cdot F_\eps\|_{\avL^2(Q_{2r}) }  + \eps^2 \| \partial_t F\|_{\avL^p(Q_{2r})} \;. \label{eq:local_time_est_2}
\end{equs}
where in the last line we used that 
$\|\sigma \|_{\avL^2(Q_{2r}) } = \|\partial_t w\|_{\avL^2(Q_{2r}) }\leq \| \mathcal{L}_\eps w\|_{\avL^2(Q_{2r}) }+ \eps \| \nabla \cdot F_\eps\|\lesssim \eps^{-1}\| \nabla w\|_{\avL^2(Q_{2r}) }
+ \| \nabla^2 w\|_{\avL^2(Q_{2r}) } + \eps \| \nabla \cdot F_\eps\|_{\avL^2(Q_{2r}) } \ .
$
Thus the terms in \eqref{eq:local_time_est_1} are bounded by 
$$\frac{\log(2+\eps^{-1} (t-s)^{1/2}) }{|t-s|^{\frac{d+3}{2}}} \exp\left(-\frac{\mu|x-y|^2}{t-s} \right)\;,$$
where for the term  $\| \nabla^2 w\|_{\avL^2(Q_{2r}) } $ we used Lemma~\ref{lem:second_derivative_kernel estimate} and the estimate on $\eps^{-1}\| \nabla w\|_{\avL^2(Q_{2r}) } + \| \nabla w_\eps \|_{\avL^p(Q_{2r})} $ follows from Corollary~\ref{cor:pointwise}. 
The terms in \eqref{eq:local_time_est_2} can be estimated using the explicit expression in \eqref{Item,derivative bound} and \eqref{Item,time derivative bound}.
Thus, one finds
\begin{align}\label{eq:loc_eq}
&| \partial_{t_2} \left(\Gamma_{\eps}- \Gamma^{0,1}_\eps\right)( x,t , y,s) | 
 \lesssim \frac{\log(2+\eps^{-1} (t-s)^{1/2}) }{|t-s|^{\frac{d+3}{2}}} \exp\left(-\frac{\mu|x-y|^2}{t-s} \right) \ , 
\end{align}
uniformly over $x,y\in \mathbb{R}^d$, $\eps\in (0,1]$, and $-\infty <s<t<\infty$ with $\eps^2<|t-s|<T$.
Interpolating with \cite[Prop~2.4]{HS23per}, which states that for $\alpha\in (0,1)$ there exists $\mu>0$ such that
\begin{align*}
&\|\Gamma_{\eps}(\ \cdot\ ; y, s )- {\Gamma}^{1,0}_\eps (\ \cdot\ ; y, s )  \|_{C_\fraks^{\alpha}(Q_{\sqrt{ |t-s|}/8}(x,t))} \lesssim \frac{\eps}{|t-s|^{\frac{d+1+\alpha}{2}}} \exp\left(-\frac{\mu|x-y|^2}{t-s} \right)\ ,
\end{align*} 
uniformly over $x,y\in \mathbb{R}^d$, $-\infty <s<t<\infty$ and $\eps\in (0,1]$, concludes the proof.
\end{proof}

%

\begin{prop}\label{prop:time holder convergence kernel2}
There exists 
$\mu>0$ such that for  $\alpha\in [0,1), \alpha'\in(1,2) $ and $\kappa, T>0$
\begin{align*}
&\sup_{ s < s'<s+ |t-s|/8} \frac{\| \left(\Gamma_{\eps}- \Gamma^{1,1}_\eps\right)( \ \cdot\ , y,s)- \left(\Gamma_{\eps}- \Gamma^{1,1}_\eps\right)( \ \cdot\ , y,s) \|_{C_\fraks^{\alpha}(Q_{\sqrt{ |t-s|}/8}(x,t))}}{|s-s'|^{\alpha'/2}}\\
&+\sup_{ t - |t-s|/8<t'<t } \frac{\|  \left(\Gamma_{\eps}- \Gamma^{1,1}_\eps\right)(x,t, \ \cdot \ )- \left(\Gamma_{\eps}- \Gamma^{1,1}_\eps\right)(x,t', \ \cdot \ )\|_{C_\fraks^{\alpha}(\tilde{Q}_{\sqrt{ |t-t'|}/8}(y,s))}}{
|t-s|^{\alpha'/2}
}
\\
&\qquad \qquad\qquad \qquad\qquad \qquad
 \lesssim_{\alpha,\alpha',\kappa, T}\frac{\eps^{2-\alpha' -\kappa} }{|t-s|^{\frac{d+\alpha +\alpha'}{2}}} \exp\left(-\frac{\mu|x-y|^2}{t-s} \right)\;,
\end{align*}
as well as
\begin{align*}
&\sup_{ s < s'<s+ |t-s|/8} \frac{\| \nabla_1\left(\Gamma_{\eps}- \Gamma^{1,1}_\eps\right)( \ \cdot\ , y,s)- \nabla_1\left(\Gamma_{\eps}- \Gamma^{1,1}_\eps\right)( \ \cdot\ , y,s) \|_{C_\fraks^{\alpha}(Q_{\sqrt{ |s-s'|}/8}(x,t))}}{|t-s|^{\alpha'/2}}\\
&+\sup_{ t - |t-s|/8<t'<t } \frac{\| \nabla_2 \left(\Gamma_{\eps}- \Gamma^{1,1}_\eps\right)(x,t, \ \cdot \ )- \nabla_2 \left(\Gamma_{\eps}- \Gamma^{1,1}_\eps\right)(x,t', \ \cdot \ )\|_{C_\fraks^{\alpha}(\tilde{Q}_{\sqrt{ |t-s|}/8}(y,s))}}{
|t-t'|^{\alpha'/2}
}
\\
&\qquad \qquad\qquad \qquad\qquad \qquad
 \lesssim_{\alpha,\alpha',\kappa, T}\frac{\eps^{1-\alpha'-\kappa} }{|t-s|^{\frac{d+2+\alpha+\alpha'}{2}}} \exp\left(-\frac{\mu|x-y|^2}{t-s} \right)\;,
\end{align*}
uniformly over $x,y\in \mathbb{R}^d$, $\eps\in (0,1]$, and $-\infty <s<t<\infty$ with $\eps^2<|t-s|<T$.
\end{prop}

\begin{proof}
Since the proof is rather analogous to the proof of Lemma~\ref{prop:holder convergence kernel}, we shall stay brief.
For $r\eqdef \sqrt{ |t-s|}/8$, $z=(y,s)$ and $z_{0}=(x,t)$ as well as  
\begin{equ}
u_{\eps}^{z}(\cdot)\eqdef \partial_{s}\Gamma_\eps (\cdot \ ;y,s ) \qquad \text{and} \qquad \bar{u}^{z}(\cdot)\eqdef \partial_{s} \bar{\Gamma}^{(0,1)}(\cdot \ ;y,s )\ 
\end{equ}
the two scale expansion $w^{z}_\eps$ as in \eqref{eq:two_scale}
satisfies
$
\partial_{t}+\mathcal{L}_\eps w^{z}_\eps= \eps \nabla\cdot F^{z}_\eps\ ,
$
for $F^{z}$ as in \eqref{eq:two_scale_error}. 
It follows from Theorem~\ref{thm:technical} that for $p>1$ such that $\alpha= 1-\frac{d+1}{p}$
\begin{equs}
\|w_\eps\|_{C_\fraks^{{\alpha}}(Q_{r}(z_0))} &\lesssim_{{\alpha}} \frac{1}{r^{{\alpha}}}  \| w_\eps \|_{\avL^2(  Q_{2r}(z_0))} + r^{1-{\alpha}} 
\eps \|  F_\eps  \|_{\avL^p(  Q_{2r}(z_0))} 
\end{equs}
Thus one concludes using \eqref{eq:loc_eq} and \eqref{Item,bound} (analogously to the argument to obtain \eqref{eq:local_wbound}) that 
\begin{equ}
\|  \partial_{s} \left(\Gamma_{\eps}- \Gamma^{1,1}_\eps\right)( \ \cdot\ , y,s)\|_{C_\fraks^{\alpha}(Q_{r}(x,t))}
 \lesssim_{\alpha} \frac{\log(2+\eps^{-1} (t-s)^{1/2}) }{|t-s|^{\frac{d+3+\alpha}{2}}} \exp\left(-\frac{\mu|x-y|^2}{t-s} \right)\;.
\end{equ}
Similarly, one finds the same upper bound on $\| \partial_{t}  \left(\Gamma_{\eps}- \Gamma^{1,1}_\eps\right)(x,t, \ \cdot \ )\|_{C_\fraks^{\alpha}(\tilde{Q}_{r}(y,s))}$.
Finally, an interpolation with \cite[Prop~2.5]{HS23per} yields the first inequality.

For the latter inequality one finds that $v^{z}_{\eps}\eqdef \nabla w^z$ as in the proof of Lemma~\ref{prop:holder convergence kernel}  satisfies
\begin{equ}
\partial_{t}+\mathcal{L}_\eps v_{\eps;j}= \eps \nabla\cdot \left( \partial_j F_\eps + (\partial_j A^{\eps}) v \right)  \ ,
\end{equ}
and thus by Theorem~\ref{thm:technical} for any $0<\bar{\alpha}<1$ and $p=\frac{d+2}{1-\alpha}$ (respectively by Corollary~\ref{cor:technical} for $\bar{\alpha}=0$)
\begin{equs}
\|v_\eps\|_{C_\fraks^{\bar{\alpha}}(Q_{r}(z_0))} 
&\lesssim  r^{1-\bar{\alpha}}\left( \eps^{-1}  \| v_\eps \|_{\avL^p(  Q_{2r}(z_0))} + 
\eps \| \nabla F_\eps  \|_{\avL^p(  Q_{2r}(z_0))} \right)\;, 
\end{equs}
where we used that $r\geq\eps/2$.
In order to bound 
$\| v_\eps \|_{\avL^p(  Q_{2r}(z_0))}= \| \nabla w_\eps \|_{\avL^p(  Q_{2r}(z_0))}$ one notes that by \eqref{eq:technical_thm_1}
$$
\| \nabla w_\eps \|_{\avL^p(  Q_{r}(z_0))} \leq C_p \left( \frac{1}{r}  \| w_\eps \|_{\avL^2(  Q_{2r}(z_0))} + \eps \|  F_\eps  \|_{\avL^p(  Q_{2r}(z_0))} \right),
$$
the estimates on $\|  F_\eps  \|_{\avL^p(  Q_{2r}(z_0))}$ and $\| \nabla F_\eps  \|_{\avL^p(  Q_{2r}(z_0))}$ 
follow from \eqref{Item,bound} and \eqref{Item,derivative bound} and one concludes analogously to the derivation of \eqref{eq:loc_eq2} that 
\begin{align*}
&\| \nabla_1 \partial_{s}\left(\Gamma_{\eps}- \Gamma^{1,1}_\eps\right)( \ \cdot\ , y,s)\|_{C_\fraks^{\alpha}(Q_{r}(x,t))}
+\| \nabla_2 \partial_{t} \left(\Gamma_{\eps}- \Gamma^{1,1}_\eps\right)(x,t, \ \cdot \ )\|_{C_\fraks^{\alpha}(\tilde{Q}_{r}(y,s))}\\
&\qquad \qquad\qquad \qquad\qquad \qquad
 \lesssim_{{\alpha},T} \frac{\eps^{-1}\log(2+\eps^{-1} (t-s)^{1/2}) }{|t-s|^{\frac{d+4+\alpha}{2}}} \exp\left(-\frac{\mu|x-y|^2}{t-s} \right)\;. 
\end{align*}
An analogous interpolation to the one at the end of the proof of Proposition~\ref{prop:double_holder convergence kernel} of this inequality with the one in Proposition~\ref{prop:double_holder convergence kernel} completes the proof.
\end{proof}

\begin{prop}\label{prop:double time increments}
There exists 
$\mu>0$ such that for $T>0$
\begin{align*}
&| \partial_{t_1} \partial_{t_2}\left(\Gamma_{\eps}- \Gamma^{1,1}_\eps\right)( x,t , y,s)|
 \lesssim_{T} \frac{\eps^{-2} \log(2+\eps^{-1} (t-s)^{1/2})}{|t-s|^{\frac{d+5}{2}}} \exp\left(-\frac{\mu|x-y|^2}{t-s} \right)\;,
\end{align*}
uniformly over $x,y\in \mathbb{R}^d$, $\eps\in (0,1]$, and $-\infty <s<t<\infty$ with $\eps^2<|t-s|<T$.
\end{prop}
\begin{proof}
Let $u_\eps= \partial_s\Gamma_{\eps}( \cdot , y,s )$ and $\bar{u}= \partial_s \bar{\Gamma}^{0,1}( \cdot , y,s ) $ as well as
$ w_\eps $ as in \eqref{eq:two_scale} and $F$ as in \eqref{eq:two_scale_error}.
Then $\sigma_\eps (x,t) = \partial_t w(x,t)$ satisfies
$$
(\partial_t + \mathcal{L}_\eps) \sigma_\eps =  \nabla \cdot \left(\eps \partial_t F + (\partial_t A^\eps) \nabla w_\eps  \right) \ .
$$
Thus, as in the proof of Proposition~\ref{prop:time holder convergence kernel}
\begin{equs}
\|\sigma\|_{L^\infty(Q_r)} 
&\lesssim \eps^{-1}\| \nabla w\|_{\avL^2(Q_{2r}) } + \| \nabla w_\eps \|_{\avL^p(Q_{2r})} 
+ \| \nabla^2 w\|_{\avL^2(Q_{2r}) } \label{eq:local_time_est_11} \\
&\qquad + \eps \| \nabla \cdot F_\eps\|_{\avL^2(Q_{2r}) }  + \eps^2 \| \partial_t F\|_{\avL^p(Q_{2r})} \;. \label{eq:local_time_est_22}
\end{equs}
The terms $ \eps^{-1}\| \nabla w\|_{\avL^2(Q_{2r}) } + \| \nabla w_\eps \|_{\avL^p(Q_{2r})} $ are bounded  using the second inequality of Proposition~\ref{prop:time holder convergence kernel2}, the 
terms in \eqref{eq:local_time_est_22} can be estimated using the explicit expression in \eqref{Item,derivative bound} and \eqref{Item,time derivative bound} and lastly, to estimate 
$\| \nabla^2 w\|_{\avL^2(Q_{2r}) } $ one argues as in the proof of the second inequality in Lemma~\ref{prop:holder convergence kernel}, 
but using \eqref{eq:technical_thm_1} instead of \eqref{eq:technical_thm_2}.
\end{proof}

%


\section{Main Results on the g-PAM and \TitleEquation{\Phi^4_3}{Phi43} Equation}\label{sec:Homogenisation of (oscillatory)...}

In order to cleanly formulate\footnote{
As in \cite{HS23per} we shall often restrict to $(\eps,\delta)\in \square$ instead of $(\eps,\delta)\in (0,1]^2$ so that the differential operator  $\nabla \cdot A(x/\eps, t/\eps^2)\nabla$ can be pushed forward to the torus. We could 
just as well have formulated the results on the full plane but with periodic noise instead, in which case the statements 
holds for $\square$ replaced by $(0,1]^2$. 
}
 the results of this section, set $
\square\eqdef \{ (\eps, \delta) \in (0,1]^2 \ : \ \eps^{-1}\in \mathbb{N} \} \;, 
$ which
we will always view as a subspace of $[0,1]^2$. In particular its closure
$\bsquare$ equals $\bsquare = \{ (\eps, \delta) \in [0,1]^2 \ : \ \eps\in \{0\} \cup \mathbb{N}^{-1} \}$.
For functions $f, g: (0,1]^2 \to \mathbb{R}$ we shall write 
 $f\sim g$ 
 to mean that 
$f-g$ extends continuously to $[0,1]^2$.
\subsection{The oscillatory generalised Parabolic Anderson Model}
In this section we state our main theorems about the g-PAM \eqref{eq:g-PAM}. 
We first consider the following regularisation of spatial white noise $\xi\in \mathcal{D}'(\mathbb{T}^2)$
\begin{equation}\label{eq:heat kernel regularisation spatial white noise}
\xi_{\eps,\delta}(x,t)= \int_{\mathbb{R}^2} \Gamma_\eps(x,t, \zeta, t-\delta^2) \xi(\zeta) \,d\zeta\ ,
\end{equation}
where we implicitly identified $\xi$ with its pullback $\pi_{2}^* \xi$ to $\mathbb{R}^{2}$.
Note in particular that $\xi_{\eps,\delta}\in {C}(\mathbb{R}^{2+1})$ is a function of space and time, despite the original noise being constant in time. The following is a straightforward adaptation of \cite[Lem.~1.9]{HS23per}.
\begin{lemma}\label{lem:regularisation_for_homogenisation}
For every $\alpha<-1$ there exists a modification of \eqref{eq:heat kernel regularisation spatial white noise} which extends to a continuous map
\begin{equ}{}
[0,1]^2\to C_{\fraks}^{\alpha}(\mathbb{R}^{2+1}) \ , \qquad (\eps,\delta)\mapsto\xi_{\eps,\delta}\ .
\end{equ}
It has the property that for any $\eps\in [0,1]$ it holds that $\xi_{\eps,0}= \xi$ and that for any $\delta>0$
\begin{equ}\label{eq:homogenised noise}
\xi_{0,\delta} (x,t)= \int_{\mathbb{R}^{2+1}} \bar{\Gamma}(x,t, \zeta, t-\delta^2) \xi(d\zeta)\  .
\end{equ}
\end{lemma}

When the nonlinearity of the equation involves the derivative of the solution, i.e.\ $f_{i,j}\neq 0$, we shall only consider `prepared' initial conditions
for the reason explained in Remark~\ref{rem:initial condition} below.  
That is, for $v_0\in C^\eta(\mathbb{T}^2)$ with $\eta\geq 0$ we set $v_\eps:[-\eps,\infty)\times \mathbb{R}^2 \to \mathbb{R}$ to be the periodic  solution to the initial value problem 
\begin{equ}\label{eq:initial condition}
\partial_t v_\eps =\nabla \cdot A^\eps \nabla v_\eps \ , \qquad 
v_\eps(-\eps)= v_0\ .
\end{equ}

\begin{theorem}\label{thm:g-pam}
Let $A$ be as in Assumption~\ref{assumption} and denote by $\bar{A}$ the homogenised matrix defined in \eqref{eq:hom_matrix} and by $\{\phi_j\}_{j=1}^2$ the correctors defined in \eqref{eq:correctors}.
For $\mu,\nu=1,2$ let $f_{\mu,\nu}: \mathbb{T}^2\to \mathbb{R}$ be bounded and measurable,
let $\xi$ be white noise on $\mathbb{T}^2$ and $\xi_{\eps,\delta}$ as in \eqref{eq:heat kernel regularisation spatial white noise} and let $v\in \mathcal{C}^\alpha(\mathbb{T}^2)$ for $\alpha\geq 1/2$.
There exist constants
\begin{align*}
\alpha_{\eps,\delta}^{\<Xi2>}&\sim\frac{\abs{\log(\delta/\eps)\wedge 0}}{{ 2\pi}} ,
&\alpha_{\eps,\delta}^{\<b2>}\sim \frac{\abs{\log(\delta/\eps)\wedge 0}}{{ 4\pi}} \\
\bar{\alpha}_{\eps,\delta}^{\<Xi2>}&\sim \frac{\abs{\log(\delta)} \wedge \abs{\log(\eps)}}{{2\pi} } , 
&\bar{\alpha}_{\eps,\delta}^{\<b2>}\sim \frac{\abs{\log(\delta)} \wedge \abs{\log(\eps)}}{{4\pi} }
\end{align*}
and uniformly bounded constants $c_{\eps,\delta}^{\<Xi2>}$, $c_{\eps,\delta}^{{\<b2>{\mu,\nu}}}$, $\gamma_{\eps,\delta}$ for $\mu,\nu=1,2$
 such that if we write $u_{\eps,\delta}$ for the solution to 
\begin{align*}
\big(\partial_t - \nabla \cdot A^\eps \nabla\big) u_{\eps,\delta}
&=\sum_{\mu,\nu=1}^2 f_{\mu,\nu}^\eps \Bigg( \partial_\mu u_{\eps,\delta} \partial_\nu u_{\eps,\delta} - \Big( \frac{ (A_s^\eps)_{\mu,\nu}^{{-1}}}{\sqrt{\det(A_s^\eps)}}\alpha^{\<b2>}_{\eps,\delta}
\\
&\qquad + \sum_{i,j} (\mathbf{1}_{\mu=i}+(\partial_{\mu} \phi_i)^\eps)(\mathbf{1}_{\nu=j}+(\partial_{\nu} \phi_j)^\eps) 
\frac{ (\bar{A}^{{-1}})_{i,j}    }{\sqrt{\det(\bar{A})}} \bar{\alpha}^{\<b2>}_{\eps,\delta} 
 + c^{\<b2>{\mu,\nu}}_{\eps,\delta}  \Big) \cdot  \sigma^2(u_{\eps,\delta}) \Bigg)
  \\& 
+ \sigma(u_{\eps,\delta})\Big( \xi_{\eps, \delta}  -\big(\frac{\alpha_{\eps,\delta}^{\<Xi2>}}{\sqrt{\det(A_s^\eps)}} +\frac{\bar{\alpha}_{\eps,\delta}^{\<Xi2>}}{\sqrt{\det(\bar{A})}}  + c_{\eps,\delta}^{\<Xi2>}  \big)  \sigma'(u_{\eps,\delta})\Big)  
- \gamma_{\eps, \delta} \cdot \sigma^2(u_{\eps,\delta})
\ ,
\end{align*}

with initial condition $u_{\eps, \delta}(0)=v_\eps(0)$ as in \eqref{eq:initial condition}, there exists a random $T>0$ such that the solution map 
\begin{equ}
\square\to \L^{0}\left[C([0,T], C(\mathbb{T}^2))\right]\;, \qquad (\eps, \delta) \mapsto u_{\eps, \delta}\;,
\end{equ}
is well defined, has a unique continuous\footnote{Recall that the $L^0$-topology is characterised by convergence in probability.} extension to $\bsquare$ and, 
furthermore, the following hold.
\begin{enumerate}
\item\label{thm:main_item1} For $\delta>0$ one has $ \lim_{\eps\to 0} \beta_{\eps,\delta}=0$ for each $\beta_{\eps,\delta}\in \{ \alpha_{\eps,\delta}^{\<Xi2>}, \alpha_{\eps,\delta}^{\<b2>}, c_{\eps,\delta}^{\<Xi2>}, c_{\eps,\delta}^{{\<b2>{\mu,\nu}}}, \gamma_{\eps,\delta}  \}$, 
the limits $ \bar{\alpha}_{0,\delta}^{\<Xi2>}\eqdef\lim_{\eps\to 0} \bar{\alpha}_{\eps,\delta}^{\<Xi2>}$, $ \bar{\alpha}_{0,\delta}^{{\<b2>}}\eqdef\lim_{\eps\to 0} \bar{\alpha}_{\eps,\delta}^{{\<b2>}}   
$
 exist and
$u_{0, \delta}$ agrees with the classical solution to
\begin{equs}
\partial_t u - \nabla \cdot \bar A \nabla u &=\sum_{\mu,\nu=1}^2 \bar{f}_{\mu,\nu} \Bigg( \partial_\mu u \partial_\nu u - \frac{ (\bar{A}^{{-1}})_{\mu,\nu}    }{\sqrt{\det(\bar{A})}} \bar{\alpha}^{\<b2>}_{0,\delta}    \sigma^2(u) \Bigg)\\
&
+ \sigma(u)\Big( \xi_{0, \delta}  -\frac{\bar{\alpha}_{0,\delta}^{\<Xi2>}}{\sqrt{\det(\bar{A})}}    \sigma'(u)\Big) \ ,
\end{equs}
where $\bar{f}_{\mu,\nu}\eqdef \int_{[0,1]^3} f_{\mu,\nu}$ .
\item\label{thm:main_item2}
 For $\eps> 0$, the limits $\bar{\alpha}^{\tau}_{\eps,0}=\lim_{\delta\to 0} \bar{\alpha}^{\tau}_{\eps,\delta}$ for $\tau\in \{\<Xi2>, \<b2>\}$ and
  $c^{\tau}_{\eps, 0}=\lim_{\delta\to 0} c^{\tau}_{\eps,\delta}$ for  $\tau\in \{\<Xi2>, \<b2>_{\mu,\nu}\}$ exist. Furthermore, there exist constants 
$\hat{\alpha}_\eps^{\<Xi2>} ,\ \hat{\alpha}_\eps^{\<b2>}$ 
such that $u_{\eps,0}$ agrees with the\footnote{Assuming we choose the same way to affinely parametrise the solution family, e.g.\ by choosing the same cutoff function $\kappa(t)$.
}
 solution, in the sense of \cite{Sin23},
to the equation formally (omitting renormalisation) given by 
\begin{align*}
\partial_t u - \nabla \cdot A^\eps \nabla u
&=\sum_{\mu,\nu=1}^2 f_{\mu \nu}^\eps \Bigg( \partial_\mu u \partial_\nu u + \Big( \frac{ (A_s^\eps)_{\mu,\nu}^{{-1}}}{\sqrt{\det(A_s^\eps)}}\hat{\alpha}^{\<b2>}_{\eps}\\
&- \sum_{i,j=1}^2 (\mathbf{1}_{\mu=i}+(\partial_{\mu} \phi_i)^\eps)(\mathbf{1}_{\nu=j}+(\partial_{\nu} \phi_j)^\eps) \frac{ (\bar{A}^{{-1}})_{i,j}    }{\sqrt{\det(\bar{A})}} \bar{\alpha}^{\<b2>}_{\eps,0}    - c^{\<b2>{\mu,\nu}}_{\eps,0}  \Big) \cdot  \sigma^2(u) \Bigg)
  \\& 
+ \sigma(u)\Big( \xi  +\big(\frac{\hat{\alpha}_{\eps}^{\<Xi2>}}{\sqrt{\det(A_s^\eps)}} -\frac{\bar{\alpha}_{\eps,0}^{\<Xi2>}}{\sqrt{\det(\bar{A})}}  - c_{\eps,0}^{\<Xi2>}  \big)  \sigma'(u)\Big) 
- \gamma_{\eps, 0} \cdot \sigma^2(u)
\ .
\end{align*}
\item If $f_{ij}(z)$ does not depend on $z$, then $\gamma_{\eps,\delta}=0$. 
\end{enumerate}

\end{theorem}

\begin{remark}
Since $T$ is itself random, the notation $\L^{0}\left[C([0,T], C(\mathbb{T}^2))\right]$ is somewhat ambiguous. One
way of interpreting it is that one has embeddings $C([0,T], C(\mathbb{T}^2)) \subset C(\R_+, C(\mathbb{T}^2))$
(extending by a constant for times greater than $T$), so this is just the subset of 
$\L^{0}\left[C(\R_+, C(\mathbb{T}^2))\right]$ consisting of random variables $u$ 
such that $u(\omega) \in C([0,T(\omega)], C(\mathbb{T}^2))$ for almost every $\omega$.
\end{remark}

\begin{remark}\label{rem:explanation constants}
Note that the constant $\gamma_{\eps,\delta}$ can be interpreted as cancelling resonances between the oscillations of the functions $f^\eps_{i,j}$ and the oscillations appearing in the renormalised Faynman diagrams stemming from the variable coefficients $A^\eps$.
The two constants $c_{\eps,\delta}^{\<Xi2>}, c_{\eps,\delta}^{{\<b2>{\mu,\nu}}}$ appear by the same mechanism as the analogue constant in \cite[Thm.~1.10]{HS23per} for the $\Phi^4_2$ equation, namely, due to the remaining error when approximating the differential operator at small scales by the frozen coefficient operator and on large scales by the homogenised operator in the renormalisation counterterm.
All three constants $\beta_{\eps,\delta}\in \{ c_{\eps,\delta}^{\<Xi2>}, c_{\eps,\delta}^{{\<b2>{\mu,\nu}}}, \gamma_{\eps,\delta} \}$ have the property that 
 in general
\begin{equ}
\lim_{\delta\downarrow 0} \lim_{\eps\downarrow 0}  \beta_{\eps,\delta}= 0\neq  \lim_{\eps\downarrow 0}\lim_{\delta\downarrow 0}\beta_{\eps,\delta}\  ,
\end{equ}
see Remarks~\ref{rem:non-comm1},~\ref{rem:non-comm2}, and~\ref{rem:non-comm3}.
\end{remark}

\begin{remark}\label{rem:initial condition}
When $f_{i,j}=0$ we can, as in \cite{CFX23}, simply work with $L^\infty$ initial conditions and there is no need to work with `prepared' initial conditions in Theorem~\ref{thm:g-pam}.
In general, the need to prepare the initial condition stems from the fact that, as in \cite{Hai14}, we shall lift the solution to the linear equation to a modelled distribution in\footnote{The fact that we additionally assume $\eta>1/2$ stems from the choice of only working with first order corrector. } 
$\mathcal{D}^{\gamma, \eta}$ for $\eta>0$ and that at small times $t\ll\eps$ neither the correctors nor polynomials provide a precise description of said solution uniformly in $\eps$.
 We expect that further adding `initial layer correctors', c.f. \cite{KL20, shen2020regularity} to our regularity structure would allow to circumvent this, but 
optimising the class of initial conditions is left for future exploration.
\end{remark}

\subsection{The oscillatory g-PAM with generic homogeneous regularisation}

\begin{theorem}\label{thm:main_translation invariant}
Let $A$, $\bar{A}$, $\{\phi_j\}_{j=1}^2$ , $\{f_{\mu,\nu}\}_{\mu,\nu=1}^2$, $u_0\in \mathcal{C}^{\alpha}(\mathbb{T}^2)$ and $\xi$ be as in 
Theorem~\ref{thm:g-pam}. 
For $\rho\in C_c^\infty(B_1)$ even, non-negative with $\int_{\mathbb{R}^2} \phi=1$, consider for $\delta>0$ the regularisation $\xi_{\delta}(z)=\xi(\phi_z^\delta)$.
There exist constants
\begin{align*}
\bar{\alpha}_{\eps,\delta}^{\<Xi2>,\flat}&\sim \frac{\abs{\log(\delta)} \wedge \abs{\log(\eps)}}{{2\pi} } , 
&\bar{\alpha}_{\eps,\delta}^{\<b2>,\flat}\sim \frac{\abs{\log(\delta)} \wedge \abs{\log(\eps)}}{{4\pi} } ,
\end{align*}
 bounded constants $c_{\eps,\delta}^{\<Xi2>,\flat}$, $c_{\eps,\delta}^{{\<b2>{\mu,\nu},\flat }}$, $\gamma^\flat_{\eps,\delta}$
as well as for $\lambda\in (0,\infty)$ functions $D^{\<Xi2>}_{\lambda}: \mathbb{R}^{2\times 2} \to \mathbb{R}$ and 
$D^{\<b2>}_{\lambda}:\mathbb{R}^{2\times 2} \to \mathbb{R}^{2\times 2}$
satisfying\footnote{The second inequality of \eqref{eq:pam func_assym} is to be read with respect to the usual partial order on positive definite matrices.} for each positive definite $M\in\mathbb{R}^{2\times 2} $
\begin{equ}\label{eq:pam func_assym}
1+ \abs{\log (\lambda)}\lesssim_M D^{\<Xi2>}_{\lambda} (M) \lesssim_M 1+ \abs{\log (\lambda)}, \qquad 
 (1+ \abs{\log (\lambda)}) \id \lesssim_M D^{\<b2>}_{\lambda}(M) \lesssim_M (1+ \abs{\log (\lambda)}) \id 
\end{equ}
uniformly over $\lambda\in (0,1]$, such that 
if we denote by $u^{\flat}_{\eps, \delta}$ the solution to
\begin{align*}
\partial_t u^{\flat}_{\eps, \delta} - \nabla \cdot A^\eps \nabla u^{\flat}_{\eps, \delta}
&=\sum_{\mu,\nu=1}^2 f_{\mu,\nu}^\eps \Bigg(  \partial_\mu u^{\flat}_{\eps, \delta} \partial_\nu u^{\flat}_{\eps, \delta} - \Big(D^{\<b2>{\mu,\nu}}_{\eps/\delta}(A_s^\eps)\\
&- \sum_{i,j=1}^2 (\mathbf{1}_{\mu=i}+(\partial_{\mu} \phi_i)^\eps)(\mathbf{1}_{\nu=j}+(\partial_{\nu} \phi_j)^\eps) 
\frac{ (\bar{A}^{{-1}})_{i,j}    }{\sqrt{\det(\bar{A})}} \bar{\alpha}^{\<b2>,\flat}_{\eps,\delta} 
 + c^{\<b2>{\mu,\nu},\flat }_{\eps,\delta}  \Big) \cdot  \sigma^2(u^{\flat}_{\eps, \delta}) \Bigg)
  \\& 
+ \sigma(u^{\flat}_{\eps, \delta})\Big( \xi_{\delta}  -\big(D^{\<Xi2>}_{\delta/\eps}(A_s^\eps)+ \frac{\bar{\alpha}_{\eps,\delta}^{\<Xi2>,\flat}}{\sqrt{\det(\bar{A})}}  + c_{\eps,\delta}^{\<Xi2>,\flat}  \big)  \sigma'(u^{\flat}_{\eps, \delta})\Big)  
- \gamma^\flat_{\eps, \delta} \cdot \sigma^2(u^{\flat}_{\eps, \delta})
\ ,
\end{align*}
with initial condition initial condition $u^\flat_{\eps, \delta}(0)=v_\eps(0)$ as in \eqref{eq:initial condition}, there exists random $T>0$ such that the solution map
\begin{equ}
 \square
\to \L^{0}\left[C([0,T], \mathcal{D}'(\mathbb{T}^2))\right], \qquad (\eps, \delta) \mapsto u^{\flat}_{\eps, \delta}
\end{equ}
is well defined, has a unique continuous extension to $ \bsquare$ and the following hold. 
\begin{enumerate}
\item\label{thm:main_flat_item1} 
 For any $M$ positive definite, it holds that $\lim_{\lambda\to \infty} D^{\<Xi2>}_{\lambda}(M)=0$ and $\lim_{\lambda\to \infty} D^{\<b2>{\mu,\nu}}_{\lambda}(M)=0$.
For $\delta>0$  it holds that $\lim_{\eps\to 0} c^{\tau,\flat}_{\eps,\delta}=0$ for $\tau\in\big\{ \<Xi2>, \<b2>{\mu,\nu}\big\}$, that the limits $ \bar{\alpha}^{\tau,\flat}_{0,\delta}\eqdef\lim_{\eps\to 0} \bar{\alpha}^\flat_{\eps,\delta}$ exists for $\tau\in\big\{ \<Xi2>, \<b2>\big\}$ and that the process
$u^{\flat}_{0, \delta}$ agrees with the classical solution to 
\begin{equs}\label{eq:classical g pam,usual reg}
\partial_t u^{\flat}_{0, \delta} - \nabla \cdot \bar{A} \nabla u^{\flat}_{0, \delta}  &=\sum_{\mu,\nu=1}^2 \bar{f}_{\mu,\nu} \Bigg( \partial_\mu u^{\flat}_{0, \delta}  \partial_\nu u^{\flat}_{0, \delta}  - \frac{ (\bar{A}^{{-1}})_{\mu,\nu}    }{\sqrt{\det(\bar{A})}} \bar{\alpha}^{\<b2>,\flat}_{0,\delta}    \sigma^2(u^{\flat}_{0, \delta} ) \Bigg)\\
&\qquad 
+ \sigma(u^{\flat}_{0, \delta} )\Big( \xi_{\delta}  -\frac{\bar{\alpha}_{0,\delta}^{\<Xi2>,\flat} }{\sqrt{\det(\bar{A})}}   \sigma'(u^{\flat}_{0, \delta} )\Big) \ .
\end{equs}
\item\label{thm:main_flat_item2} For all $\eps\in [0,1]$ 
the process $u^{\flat}_{\eps,0}$ agrees with the process $u_{\eps,0}$ in Theorem~\ref{thm:g-pam} and is in particular independent of the choice of mollifier $\rho$.
For each $\eps>0$ and $\beta^{\tau,\flat}_{\eps,\delta} \in
\big\{\bar{\alpha}_{\eps,\delta}^{\<Xi2>,\flat},
\bar{\alpha}_{\eps,\delta}^{\<b2>,\flat},
 c_{\eps,\delta}^{\<Xi2>,\flat}, c_{\eps,\delta}^{{\<b2>{\mu,\nu},\flat }}, \gamma^\flat_{\eps,\delta} \big\}
$ the value $\lim_{\delta\to 0} \beta^{\tau,\flat}_{\eps,\delta}$ agrees with the corresponding value $\beta^\tau_{\eps,0}$ in Theorem~\ref{thm:g-pam}.
\item If $f_{ij}(z)$ does not depend on $z$, then $\gamma^\flat_{\eps,\delta}=0$. 
\end{enumerate}
\end{theorem}
\begin{remark}
Let us observe that the functions $D^{\<Xi2>}_{\delta}(A_s)$ and $D^{\<b2>{\mu,\nu}}_{\delta}(A_s)$ are exactly the renormalisation functions for the (inhomogeneous) g-PAM equation when working with the regularisation $\xi_{\delta}(x)= \xi(\rho^{\delta}_x)$ in \cite[Sec~3.3]{Sin23}.
Furthermore, an analogous remark to \cite[Rem.~1.18]{HS23per} applies to the asymptotic behaviour of these functions.
Lastly, Remark~\ref{rem:explanation constants} is also applicable to Theorem~\ref{thm:main_translation invariant}. 
\end{remark}

On the restricted set of parameters  $\triangle_C^{<}\eqdef \{(\eps,\delta)\in [0,1] \ : \ \eps\leq C \delta \}$, one obtains the following result.
\begin{theorem}\label{thm:main_translation invariant2}
Let $A$, $\{\phi_j\}_{j=1}^2$ , $\{f_{\mu,\nu}\}_{\mu,\nu=1}^2$, $u_0\in \mathcal{C}^{\alpha}(\mathbb{T}^2)$, $\xi_{\delta}$ and
$
\bar{\alpha}_{\eps,\delta}^{\<Xi2>,\flat}$, $\bar{\alpha}_{\eps,\delta}^{\<b2>,\flat}$, $c_{\eps,\delta}^{\<Xi2>,\flat}$, $c_{\eps,\delta}^{{\<b2>{\mu,\nu},\flat }}$, $\gamma^\flat_{\eps,\delta} $
as well as $D^{\<Xi2>}_{\lambda}$ and 
$D^{\<b2>}_{\lambda}$ 
be as in Theorem~\ref{thm:main_translation invariant}.
Then, the constants
$${c}_{\eps,\delta}^{\<Xi2>,\flat\flat}\eqdef c^{\<Xi2>, \flat}_{\eps,\delta} + \int_{[0,1]^{3}} D^{\<Xi2>}_{\delta/\eps} , \qquad
  {c}_{\eps,\delta}^{{\<b2>{\mu,\nu}, \flat\flat }}=c^{\<b2>{\mu,\nu}, \flat}_{\eps,\delta} + \int_{[0,1]^{3}} D^{\<b2>{\mu,\nu}}_{\delta/\eps}$$
are bounded on $\triangle^{<}_{C}$ for any $C>0$ and for $\delta>0$ vanish as $\eps\to 0$.
If we denote by $u^{\flat\flat}_{\eps, \delta}$ the solution to
\begin{align*}
&\partial_t u^{\flat\flat}_{\eps, \delta} - \nabla \cdot A^\eps \nabla u^{\flat\flat}_{\eps, \delta}\\
&=\sum_{\mu,\nu=1}^2 f_{\mu,\nu}^\eps \Bigg( \partial_\mu u^{\flat\flat}_{\eps, \delta} \partial_\nu u^{\flat\flat}_{\eps, \delta} - \Big( \sum_{i,j} (\mathbf{1}_{\mu=i}+(\partial_{\mu} \phi_i)^\eps)(\mathbf{1}_{\nu=j}+(\partial_{\nu} \phi_j)^\eps) 
\frac{ (\bar{A}^{{-1}})_{i,j}    }{\sqrt{\det(\bar{A})}} \bar{\alpha}^{\<b2>,\flat}_{\eps,\delta} \\
&\qquad
 + c^{\<b2>{\mu,\nu},\flat\flat}_{\eps,\delta}  \Big) \cdot  \sigma^2(u^{\flat\flat}_{\eps, \delta}) \Bigg)
+ \sigma(u^{\flat\flat}_{\eps, \delta})\Big( \xi_{ \delta}  -
\big(\frac{\bar{\alpha}_{\eps,\delta}^{\<Xi2>,\flat}}{\sqrt{\det(\bar{A})}} 
  + c_{\eps,\delta}^{\<Xi2>,\flat\flat}  \big)  \sigma'(u^{\flat}_{\eps, \delta})\Big)  
- \gamma^\flat_{\eps, \delta} \cdot \sigma^2(u^{\flat\flat}_{\eps, \delta})
\ ,
\end{align*}
with initial condition $u^{\flat\flat}_{\eps, \delta}(0)=v_\eps(0)$ as in \eqref{eq:initial condition}, there exists a random $T>0$ such that for every $C>0$ the solution map
\begin{equ}
\triangle^{<}_{C}\cap \square
\to \L^{0}\left[C([0,T], \mathcal{D}'(\mathbb{T}^2))\right], \qquad (\eps, \delta) \mapsto u_{\eps, \delta}
\end{equ}
is well defined and has a unique continuous extension to $\triangle^{<}_{C}\cap \bsquare$. Furthermore,
$u^{\flat\flat}_{0, \delta}$ agrees with $u^{\flat}_{0, \delta}$ in Theorem~\ref{thm:main_translation invariant} for all $\delta\in [0,1]$.
%
%
\end{theorem}

\subsection{The oscillatory \TitleEquation{\Phi^4_3}{Phi43} equation}
For conciseness, we shall only formulate the main result in the case of regularisation based on the heat kernel.\footnote{We believe that the interested reader will be able to formulate the applicable result for translation invariant regularisations, as the modifications necessary follow along similar lines as in the case of the parabolic Anderson model, i.e.\ Theorems~\ref{thm:main_translation invariant} and \ref{thm:main_translation invariant2}, and the proof in Section~\ref{sec:application phi} is structured in such a way to cleanly accommodate such an adaptation.} 
Let $\xi$ denote space time white noise on $\mathbb{T}^d\times \mathbb{R}$.
%
We work with the regularisation formally given by
\begin{equation}\label{eq:noise}
\xi_{\eps,\delta} (x, t)= \int_{\mathbb{R}^{d}}
 \Gamma_{\eps}(x,t, \zeta, t-\delta^2) \xi(d\zeta, t)\ ,
\end{equation}
where we implicitly identified $\xi$ with its pullback $\pi_{d,1}^* \xi$ to $\mathbb{R}^{d+1}$. The following is a direct variant of \cite[Lem.~1.9]{HS23per}.

\begin{lemma}\label{lem:regularisation_for_homogenisation2}
For every $\alpha< -\frac{d+2}{2}$ there exists a modification of \eqref{eq:noise} which extends to a continuous map
\begin{equ}{}
[0,1]^2\to C_{\fraks}^{\alpha}(\mathbb{R}^{d+1}) \ , \qquad (\eps,\delta)\mapsto\xi_{\eps,\delta}\ .
\end{equ}
It has the property that for any $\eps\in [0,1]$ it holds that $\xi_{\eps,0}= \xi$ and that for any $\delta>0$
\begin{equ}\label{eq:homogenised noise2}
\xi_{0,\delta} (x,t)= \int_{\mathbb{R}^{d}} \bar{\Gamma}(x,t, \zeta, t-\delta^2) \xi(d\zeta, t)\  .
\end{equ}
\end{lemma}
%

For similar reasons as explained in Remark~\ref{rem:initial condition} we shall only consider initial conditions which are regular perturbations of the solution to the linear equation.\footnote{
Having to choose initial conditions in this way for singular SPDEs is quite common, see also \cite[Sec.~5]{BCCH20}.
} For this we define
$$\Pi^{\eps,\delta}\<1>(0)\eqdef \int_{\mathbb{R}^3\times \mathbb{R}} \kappa(-s) \Gamma_\eps(\, \cdot \, ,0,y,s-\delta^2) \xi(dy,ds)\ ,
$$
which
extends to a continuous map
$[0,1]^2\to C_{\fraks}^{-\f12-\kappa}(\mathbb{T}^{3}) $
for any $\kappa>0$ by Lemma~\ref{lem:stochastic estimates_linear}.

\begin{theorem}\label{thm:phi4}
Let $\xi$ denote space-time white noise on $\mathbb{T}^3\times \mathbb{R}$ and let $v_0\in L^\infty (\mathbb{T}^3)$.  Consider for $\delta>0$ the  regularisation $\xi_{\eps,\delta}(x,t)$ as in \eqref{eq:noise}. For $\lambda\geq 1$, let $R^{\<2>}_\lambda$ be the 
 bounded $\mathbb{Z}^{3+1}$-periodic function
$$R^{\<2>}_\lambda : \mathbb{R}^{3+1}\to \mathbb{R}, \qquad 
z \mapsto \int_{\mathbb{R}^3\times [\lambda^2, \infty)} \Gamma_1^2(z; y,s ) dyds \ .$$
There exist constants $\alpha^{\<2>}, \alpha^{\<22>}\in \mathbb{R}$,
\begin{equ}\label{eq:assymptotics renormalsation constants}
\alpha_{\eps,\delta}^{\<2>} \sim
 \alpha^{\<2>}\big( \frac{1}{\delta}-\frac{1}{\eps})\vee 0 
, \quad 
%
%
\alpha_{\eps,\delta}^{\<22>} \sim {\alpha^{\<22>}} \abs{\log(\delta/\eps)\wedge 0}
, \quad 
\bar{\alpha}^{\<22>} _{\eps,\delta}\sim {\alpha^{\<22>}} \big(\abs{\log(\delta)} \wedge \abs{\log(\eps)}\big) \ ,
\end{equ}
and bounded constants ${c_{\eps,\delta}}$, $\gamma_{\eps,\delta}$ 
  such that if 
 we denote by $u_{\eps, \delta}$ the solution to
\begin{align*}
\partial_t u_{\eps, \delta} -\nabla \cdot A^\eps \nabla u_{\eps, \delta} &= 
f^\eps \Big[ -u_{\eps,\delta}^2 + \frac{3\alpha_{\eps, \delta}^{\<2>}}{ \sqrt{\det(A_s^\eps)}}+ \frac{3}{\eps} (R^{\<2>}_{(\delta/\eps)\vee 1}\circ \mathcal{S}^{\eps})  -
 \frac{9 f^\eps \alpha_{\epsilon, \delta }^{\<22>}}{ \det(A_s^\eps)}  -  \frac{ 9 \bar{f}\bar{\alpha}_{\epsilon, \delta }^{\<22>} }{ \det(\bar {A})} - c_{\eps,\delta}\Big] u_{\eps,\delta} \\
&\qquad + \gamma_{\eps,\delta} u_{\eps,\delta}
   +\xi_{\eps, \delta}\ 
\end{align*}
with initial condition $u_{\eps, \delta}(0)=\Pi^{\eps,\delta}\<1> (0)+ v$ with $v\in \L^\infty$, then there exists a (random) $T>0$ such that the solution map 
\begin{equ}
\square\to \L^{0}\left[C([0,T], \mathcal{D}'(\mathbb{T}^3))\right]\;, \qquad (\eps, \delta) \mapsto u_{\eps, \delta}\;,
\end{equ}
is well defined, has a unique continuous
 extension to $\bsquare$ and the following hold. 

\begin{enumerate}
\item\label{thm:main_item1a} For $\delta>0$ one has $\lim_{\eps\to 0} \beta_{\eps,\delta}=0$ for $\beta_{\eps,\delta}\in \big\{
\alpha^{\<2>}_{\eps,\delta}, \ \alpha^{\<22>}_{\eps,\delta}, \ c_{\eps,\delta},\ \gamma_{\eps,\delta}
 \big\}
 $
  and the limits 
$\bar{\alpha}_{0,\delta}^{\<2>}\eqdef\lim_{\eps\to 0} R_{\delta/\eps}^{\<2>}$ and
$ \bar{\alpha}_{0,\delta}^{\<22>}\eqdef\lim_{\eps\to 0} \bar{\alpha}_{\eps,\delta}^{\<22>}$ exists. Furthermore, 
$u_{0, \delta}$ agrees with the classical solution to
\begin{equ}
\partial_t u_{0, \delta} -\nabla\cdot \bar{A}\nabla u_{0, \delta}= \bar{f} \Big[ -u_{0,\delta}^2 +\frac{3\bar{\alpha}_{0,\delta}}{\sqrt{\det(\bar {A})}} - \frac{9\bar{f}\bar{\alpha}_{0,\delta}}{\det(\bar {A})}\Big]
  u_{0,\delta} +\xi_{0,\delta} \ .
\end{equ}
\item\label{thm:main_item2a}
For $\eps> 0$, the limits $\bar{\alpha}^{\<22>}_{\eps,0}=\lim_{\delta\to 0} \bar{\alpha}^{\<22>}_{\eps,\delta}$  and
  $\beta_{\eps, 0}=\lim_{\delta\to 0} \beta_{\eps,\delta}$ for $\beta\in \big\{c_{\eps,\delta}, \gamma_{\eps, \delta}  \big\}$ exist. Furthermore, there exists a sequence of constants 
$\hat{\alpha}_\eps^{\tau}$ 
such that $u_{\eps,0}$ agrees with the\footnote{Again assuming we choose the same way to affinely parametrise the solution family, e.g.\ by choosing the same cutoff function $\kappa(t)$.
}
 solution, in the sense of \cite{Sin23},
to the equation formally (omitting renormalisation) given by
\begin{align*}
\partial_t u_{\eps, 0} -\nabla \cdot A^\eps \nabla u_{\eps, 0}&=  f^\eps \Big[ -u_{\eps,0}^2 - \frac{3\hat{\alpha}_{\eps}^{\<2>}}{ \sqrt{\det(A_s^\eps)}}+ \frac{3}{\eps} (R^{\<2>}_1 \circ \mathcal{S}^{\eps})  + \frac{9f^\eps \hat{\alpha}_{\epsilon }^{\<22>}}{ \det(A_s^\eps)}  -  \frac{9\bar{f} \bar{\alpha}_{\epsilon, 0 }^{\<22>} }{ \det(\bar {A})} - c_{\eps, 0}\Big] u_{\eps,0} \\
&\qquad + \gamma_{\eps, 0}u_{\eps,0 }
   +\xi \ .
\end{align*}
\item If $f_{ij}(z)$ does not depend on $z$, then $\gamma_{\eps,\delta}=0$. 
\end{enumerate}
\end{theorem}

\begin{remark}
The two constants $\hat{\alpha}_{\eps}^{\<2>}$ and $\hat{\alpha}_{\eps}^{\<22>}$ represent the discrepancy between the renormalisation used here and that used in \cite{Sin23}. Such a discrepancy
is bound to arise since there is no reason in general to expect that the renormalisation used
there behaves well in the oscillatory setting considered here.
\end{remark}

\begin{remark}
Let us note that, since we work only with one specific regularisation, we might as well have written equalities instead of $\sim$ when characterising the divergences. Since this choice would though be regularisation dependent we prefer to state the theorem in a form, which makes it easier for the interested reader to formulate the analogue result for generic homogeneous regularisations. 
The constants $c_{\eps,\delta}, \gamma_{\eps,\delta}$ can be checked to behave as the analogous constants for g-PAM.
\end{remark}

 On the restricted set of parameters  $\triangle_C^{<2}\eqdef \{(\eps,\delta)\in [0,1] \ : \ \eps\leq C \delta^2 \}$, the divergent contributions at mesoscopic scales $\eps\sim \delta$ are sufficiently suppressed and the following holds.
\begin{theorem}\label{thm:phi4 restricted}
Let $\xi, \xi_{\eps,\delta}$ and $v_0$ as well as 
$\alpha^{\<2>}_{\eps,\delta}, \  R^{\<2>}_\lambda,\ \alpha^{\<2>}, \  \alpha^{\<22>}_{\eps,\delta}, \    c_{\eps,\delta}$
be as in Theorem~\ref{thm:phi4} .  
There exist constants\footnote{
Here we write $f\sim g$ (for functions $f, g: (0,1]^2 \to \mathbb{R}$) to mean that 
$f-g$ extends continuously to $[0,1]^2$.
} 
\begin{equ}\label{eq:assymptotics renormalsation constants2}
\bar{\alpha}_{\eps,\delta}^{\<2>} \sim  \frac{\alpha^{\<2>}}{\delta}  
, \qquad 
\bar{\alpha}^{\<22>}_{\eps,\delta}\sim {\alpha^{\<22>}} \abs{\log(\delta)}  \ ,
\end{equ}
and bounded constants $\gamma^{<}_{\eps,\delta}$
such that the constant 
$${c}_{\eps,\delta}^{<}\eqdef c_{\eps,\delta}+ 
3 \int_{[0,1]^4}
\Big[
\frac{\alpha^{\<2>}_{\eps,\delta}}{\sqrt{\det(A_s)}} 
+\frac{1}{\eps} R^{\<2>}_{(\delta/\eps)\vee 1}  - 
\frac{\bar{\alpha}^{\<2>}_{\eps,\delta}}{\sqrt{\det(\bar{A})}} \Big]
+ 
9\int_{[0,1]^4}\frac{\alpha^{\<22>}_{\eps,\delta}}{{\det(A_s)}}  
 $$ 
 is uniformly bounded on $\triangle_C^{<2}\cap \bsquare$
%
and such that if 
 we denote by $u^{<}_{\eps, \delta}$ the solution to
\begin{align*}
\partial_t u^{<}_{\eps, \delta} -\nabla \cdot A^\eps \nabla u^{<}_{\eps, \delta} &=
f^\eps\cdot \Big[ -(u^{<}_{\eps,\delta})^2 + \frac{3\bar{\alpha}_{\eps, \delta}^{\<2>}}{ \sqrt{\det(\bar{A})}}  -  \frac{ 9 \bar{f}\bar{\alpha}_{\epsilon, \delta }^{\<22>} }{ \det(\bar {A})} - c^{<}_{\eps,\delta}\Big] u^{<}_{\eps,\delta} + \gamma^{<}_{\eps,\delta} u^{<}_{\eps,\delta}
 +\xi_{\eps, \delta}\ 
\end{align*}
with initial condition $u^{<}_{\eps, \delta}(0)=\Pi^{\eps,\delta}\<1>(0) + v$ with $v\in \L^\infty$, then there exists a (random) $T>0$ such that the solution map 
\begin{equ}
\triangle^{<2}_C \cap \square\to \L^{0}\left[C([0,T], \mathcal{D}'(\mathbb{T}^2))\right]\;, \qquad (\eps, \delta) \mapsto u^{<}_{\eps, \delta}\;,
\end{equ}
is well defined and has a unique continuous extension to $\triangle^{<2}_C\cap \bsquare$. 
Furthermore,
$u^{<}_{0, \delta}$ agrees with $u_{0, \delta}$ in Theorem~\ref{thm:phi4} for all $\delta\in [0,1]$.
\end{theorem}

  \begin{remark}\label{rem:measo div}
Note that in contrast to the $\Phi^4_2$ equation in \cite{HS23per} and the above theorems on g-PAM, we are here not able to cleanly express the large scale part of all divergences in terms of only objects from homogenisation theory uniformly in $(\eps,\delta)\in (0,1]^2$. To see why, note that
$$\E[(\Pi^{\eps,\delta}\<1>)^2]\sim
\mathbf{1}_{\delta<\eps} \int_{\delta^2}^{\eps^2} \Gamma_\eps + \mathbf{1}_{\delta^2< \eps} \int_{\delta^2\vee \eps^2}^\eps \Gamma_\eps + \int_{\delta^2\vee \eps}^\infty \Gamma_\eps \ .
$$
While the first and last terms in this sum for $ \delta^2<\eps$ satisfy
$$
\mathbf{1}_{\delta<\eps} \int_{\delta^2}^{\eps^2} \Gamma_\eps(z;y,s) dy\,ds= \frac{{\alpha}^{\<2>}_{\eps,\delta}}{\sqrt{\det(A_s(y,s))}} + \mathcal{O}(1) , \quad
\int_{\eps}^\infty\int_{\mathbb{R}} \Gamma^2_\eps(z;y,s) dy\,ds ={ \frac{\bar{\alpha}^{\<2>}\eps^{-1/2}}{\sqrt{\det(A)}}} + \mathcal{O}(1),$$ uniformly in $(\eps,\delta)\in (0,1]$, the middle term is in general unbounded on $\{(\eps,\delta)\in \square  : \eps> \delta^2\}$. We do not know at this point whether its divergent part can in general be expressed more explicitly.
  \end{remark}


 \begin{remark}
We expect that the random $T>0$ in the above theorem could be chosen to be an arbitrary (deterministic) positive time, since at least the main PDE ingredient for the proof of the a priori estimate in \cite{MW20}, the maximum principle, still holds if the Laplace operator is replaced by the uniformly elliptic operator $\nabla \cdot A^\eps \nabla$. 
 \end{remark}


\section{Homogenisation and Regularity Structures}\label{sec:RS}
In this section we first recall some basic background on the theory of regularity structures \cite{Hai14} and fix notations.
Then, the goal of the section is to show how the ansatz in \eqref{eq:formal rewriting} can be implemented. 
To this end, Section~\ref{sec:A slightly weakened topology on kernels} introduces a topology on kernels compatible with the theory of regularity structures in which the (appropriately post processed) two scale expansion error of the kernel $\Gamma_\eps-\bar{\Gamma}  - \eps \psi^{\eps} \nabla\bar{\Gamma}\to 0$ as $\eps \to 0$. Section~\ref{sec:lifting corrector terms} explains how to lift (after appropriate post processing) the corrector terms $\eps \psi^{\eps}\nabla \bar{\Gamma}$ in a way that preserves the property of being $2-\kappa$ regularising and as well as, when applicable, the property that modelled distributions get mapped into function like sectors.
Section~\ref{sec:An abstract fixed point theorem} presents an abstract fixed point theorem\footnote{Up to this point everything has been formulated at a more abstract level than actually used in this article, as these results are also useful in other settings.} and Section~\ref{sec_post_process} shows that the bounds of Section~\ref{sec:Parabolic Operators, Heat Kernels and Homogenisation} can be used in the abstract set-up developed until here. Finally, we perform the rigorous rewriting alluded to in \eqref{eq:formal rewriting} in Section~\ref{sec:Homogenisation in RS}.

\begin{definition}
A regularity structure is a pair $\mathcal{T}=(T,G)$, where $T=\bigoplus_{\alpha\in A} T_\alpha$ is a vector space graded by an index set  $A\subset \mathbb{R}$ bounded from below and without accumulation points, where each component $T_\alpha$ is a Banach space.  Furthermore, $G$ is a group acting on $T$ 
from the left, such that 
$$\Gamma \tau- \tau \in T_{< \alpha}$$
for each $\Gamma\in G$ and $\tau \in T_{\alpha}$
\end{definition}
  
 \begin{definition}
 A model $M=(\Pi,\Gamma)$ for a regularity structure $\mathcal{T}$ consists of a pair of maps 
 $$ \Pi: \mathbb{R}^{d+1} \to L(T, \mathcal{D}'), \qquad \Gamma: \mathbb{R}^{d+1}\times \mathbb{R}^{d+1} \to G$$
 such that, for all $x,y,z\in \mathbb{R}^{d+1}$,
$ \Pi_x \Gamma_{x,y}= \Pi_y$ and $\Gamma_{x,y}\Gamma_{y,z} = \Gamma_{x,z}$.
Furthermore,  for fixed $R>-\abs{\min A}$ and each compact set $\mfK\subset \mathbb{R}^{d+1}$ and $\gamma\in \mathbb{R}$, 
$ \vertiii{M}_{\gamma,\mfK}\eqdef \|\Pi\|_{\gamma, \mfK}+ \|\Gamma\|_{\gamma, \mfK}<+\infty  ,$
where
\begin{equ}\label{eq:model bounds}
\|\Pi\|_{\gamma, \mfK}\eqdef \sup_{\alpha<\gamma} \sup_{\tau\in T_\alpha} \sup_{\lambda\in (0,1]} \sup_{\phi\in \mathfrak{B}_{R}} \frac{ |\Pi_x\tau (\phi^\lambda_x)|}{\lambda^\alpha |\tau| } , \qquad
\|\Gamma\|_{\gamma, \mfK}\eqdef\sup_{\alpha<\gamma} \sup_{\tau\in T_\alpha} \sup_{x,y\in \mfK} \frac{| \Gamma_{x,y}\tau|_{\beta}}{|x-y|^{\alpha-\beta} |\tau|} 
\end{equ}
Similarly, we set
$\vertiii{M;\bar{M}}_{\gamma,\mfK}\eqdef \|\Pi-\bar \Pi \|_{\gamma, \mfK}+ \|\Gamma- \bar{\Gamma}\|_{\gamma, \mfK}\  $  for a second model  $\bar M=(\bar \Pi, \bar \Gamma)$.

\end{definition}  
Recall that for a regularity structure $(T,G)$ a subspace $V= \bigoplus_{\alpha\in A} V_\alpha \subset T$ is called a sector if $\Gamma V\subset V$ for all $\Gamma \in G$. The number $\min \{\alpha \in A: V_\alpha\neq \{0\} \}$ is called the regularity of the sector and $V$ is called function-like if it is non-negative.
\begin{definition}
Given a regularity structure $\mathcal{T}$, a model $M=(T,G)$ on $\mathbb{R}^{d+1}$ and a sector $V\subset T$, recall that the spaces of modelled distributions $\mathcal{D}^{\gamma,\eta}(V)$ consists of functions $f: \{(x,t)\in \mathbb{R}^{d+1}\ : \ t\neq 0 \}\to V$ such that for every compact $\mfK$ the following norms are finite
$$\|f\|_{\gamma,\eta,\mfK} \eqdef \sup_{(x,t)\in \mfK; t\neq 0} \sup_{l<\gamma} \frac{| f(x,t)|_{l}}{(\sqrt{|t|}\wedge 1)^{(\eta-l)\wedge 0} }<\infty$$
and writing $z=(x,t), w=(y,s)$ and $\delta = \sqrt{|t|\wedge |s|}/2$,
$$
\vertiii{f}_{\gamma,\eta,\mfK}\eqdef \|f\|_{\gamma,\eta,\mfK} + \sup_{z,w \in \mfK: |z-w|_\fraks< \delta } \sup_{l<\gamma} 
\frac{| f(w)-\Gamma_{w,z}f(z) |_l}{ |w-z|_\fraks^{\gamma-l}  (\sqrt{|t|}\wedge \sqrt{|s|}\wedge 1)^{\eta-\gamma} }\ .
$$
The distance between two modelled distributions $f \in \mathcal{D}_{M}^{\gamma,\eta}$, $\bar f \in \mathcal{D}_{\bar{M}}^{\gamma,\eta}$ for distinct models $M,\ \bar{M}$ is
$$\vertiii{f;\bar{f}}_{\gamma,\eta,\mfK}\eqdef \|f-\bar{f}\|_{\gamma,\eta,\mfK} 
 + \sup_{z,w \in \mfK: |z-w|_\fraks< \delta} \sup_{l<\gamma} 
\frac{| f(w)-\bar{f}(w)-\Gamma_{w,z}f(z) + \bar{\Gamma}_{w,z}\bar{f}(z)|_l}{ |w-z|_\fraks^{\gamma-l}  (\sqrt{|t|}\wedge \sqrt{|s|}\wedge 1)^{\eta-\gamma} }\ .
$$
For the definition of non-singular modelled distributions $f\in \mathcal{D}^{\gamma}(V)$, see \cite[Def.~3.1]{Hai14}.
\end{definition}

\begin{remark}\label{rem:poly_example}
An important example is the polynomial regularity structure $(\bar{T}, \bar{G})$ given by
$\bar{T}= \bigoplus_{n\in \mathbb{N}} \bar{T}_n$ with  $\bar{T}_n= \text{span} \{\pmb{X}^k \ : \ |k|_\fraks = n\}$
and $\bar{G}\sim (\mathbb{R}^{d+1}, +)$ where the action of $z=(z_1,...,z_{d+1})\in \bar{G}$ is given by 
 $X^k \mapsto (X+z)^k$. Its canonical model is characterised by $\Pi_z X^k= (\cdot-z)^k$.
 
Defining for $f\in C_{\fraks}^R$, the operator $\opbfP^R[f]: z\mapsto \opbfP^R_z[f]\eqdef \sum_{|k|_\fraks<R} \frac{X^k}{k!} D^kf(z)$ one in particular notes that 
$\opP^R_z[f](w)= \Pi_z \opbfP^R_z[f](w)$, which can be used to see for the polynomial model that $\opbfP^R: C_{\fraks}^R\to \mathcal{D}^\gamma(\bar{T})$ is an isomorphism, c.f. \cite[Lem.~2.12]{Hai14} or \cite[Thm.~3.12]{MS23}.
\end{remark}
The following is \cite[Props~6.9~\&~7.2]{Hai14} in our setting.
\begin{prop}
Fix a sector $V$ of regularity $\alpha>-2$ and $R>|\alpha \wedge 0|$ and assume that $-2<\eta<\gamma$. Then, there is a unique (reconstruction) operator 
$\mathcal{R}: \mathcal{D}^{\gamma,\eta}(V) \to \mathcal{C}_{\fraks}^{\alpha\wedge \eta}\ $ satisfying for $z=(w,t)$ with $t\neq 0$
\begin{equ}\label{eq:weighted reconstruction bound}
\big|\big(\mathcal{R}f- \Pi_zf(z)\big)(\phi_z^\lambda)\big|\lesssim \lambda^{\gamma} \sup_{x,y\in \supp \phi_z^\lambda}\frac{|f(x)-\Gamma_{x,y}f(y)|_\zeta}{|x-y|_\fraks^{\gamma-\zeta}}
\end{equ} 
uniformly over $\phi\in \mathfrak{B}_1^{R}$ and $\lambda<\sqrt{t}\wedge 1$, where the implicit constant depends on ${\vertiii{Z} }_{B_1(z)}$.
Denoting for two models $M, \bar{M}$ the associated reconstruction operators by $\mathcal{R}_M, \mathcal{R}_{\bar M}$, the estimate 
$$
\| \mathcal{R}_M f - \mathcal{R}_{\bar{M}} \bar f \|_{\mathcal{C}_{\fraks}^{\alpha\wedge \eta}(\mfK) } \lesssim \vertiii{f; \bar{f}}_{\gamma,\eta, \mfK} + \vertiii{M; \bar{M}}_{\gamma,\mfK}, $$
holds for  $f\in\mathcal{D}^{\gamma,\eta}$ and $\bar{f}\in\mathcal{D}^{\gamma,\eta}$ , where the implicit constant depends on the norms of the involved models and modelled distributions (on $\mfK$, resp.~$\supp(\phi_z^\lambda)$). Finally, an estimate analogous to \eqref{eq:weighted reconstruction bound}
holds for the difference of models and modelled distributions.
\end{prop}
\begin{remark}\label{rem:non_linearity and derivation}
Let us recall some important properties of modelled distributions:
\begin{itemize}
\item If the regularity structure is equipped with a multiplication map on sectors $V, W$, there is an induced (continuous) multiplication map on modelled distributions, see \cite[Prop.~6.12]{Hai14}. If there is a product on a function like sector $V$ such that $V_0$ is generated by the neutral element $\one$ for the 
multiplication, one can furthermore lift composition with a smooth map $G$ to a map
$$\hat{G}:\mathcal{D}^{\gamma,\eta}(V)\to \mathcal{D}^{\gamma,\eta}(V)\ ,$$
for $0\leq \eta \leq \gamma$, see \cite[Prop.~6.13]{Hai14}.
\item If the sector is equipped with abstract differentiation map $\partial_i: V\to T$
and the model is such that $\Pi_x \d_i \tau = \d_i \Pi_x \tau$, it was shown in 
\cite[Prop.~6.15]{Hai14} that this induces a map 
$\partial_i: \mathcal{D}^{\gamma,\eta}\to \mathcal{D}^{\gamma-1,\eta-1} $ such that $\mathcal{R}\partial_i f = \partial_i \mathcal{R}f$.
\end{itemize}
\end{remark}
Lastly, let us recall that one can lift convolutions with singular kernels to regularity structures, \cite{Hai14}. Fix $\beta\in (0,|\fraks|)$,
we shall make the following assumption, which in particular replaces \cite[Ass.~5.4]{Hai14}
\begin{assumption}\label{ass:irrational_beta}
At integer homogeneities the model space is given by polynomials, i.e.\  $T_n= \bar{T}_n$ using the notation of Remark~\ref{rem:poly_example}, and the model acts on it as therein. Furthermore, for all $\alpha\in A$ one has $\alpha+\beta \notin \mathbb{N}$.
\end{assumption}

Finally, recall that a map $I:V\to T$ defined on a sector $V$ is called\footnote{Note that under Assumption~\ref{ass:irrational_beta} one does not need to require $I(\bar{T})\subset \{0\}$ as originally done in \cite{Hai14}.}
 an abstract integration map of order $\beta>0$ if $I(V_\alpha)\subset T_{\alpha+\beta}$ and 
$\big(I \circ \Gamma-\Gamma\circ I\big) V\subset \bar{T}$ for each $\Gamma\in G$. Then, given a model $(\Pi,\Gamma)$ which realises $K$ for $I$ in the sense of \cite[Def.~5.9]{Hai14}, one lifts a singular kernel $K$ to an abstract operator on modelled distributions
\begin{equ}\label{eq:kernel lift}
\mathcal{K}_\gamma f(z)= If(x) +  {J}^K(z)f(z) + \mathcal{N}^K_\gamma f 
\end{equ}
where the three operators 
$${I}: V\to T ,\qquad {J}^K(z): V \to \bar{T}, \qquad \mathcal{N}^K_\gamma: \mathcal{D}^\gamma (V) \to C(\mathbb{R}^{d+1},\bar T)\ ,$$
are defined in \cite[Eq.~5.11, Eq.~5.15 \& Eq.~5.16]{Hai14}.
These results then allow to lift the mild formulation of a singular SPDE to a fixed point problem
in some $\mathcal{D}^{\gamma,\eta}$ space, see \cite[Sec.~7]{Hai14}. As is evident therein and explicitly formulated in \cite[Sec.~2.4]{Sin23} the solution to this fixed point problem is continuous with respect to the kernel if one equips the space of kernels with the topology 
\cite[Def.~1]{Sin23} quantifying \cite[Ass.~5.1]{Hai14} 
and if one works with models compatible with and continuously dependent on those kernels.
The next subsection is dedicated to the observation that one can slightly weaken that topology on kernels and still retain such continuity of the abstract solution map.

%
%

\subsection{A slightly weakened topology on kernels}\label{sec:A slightly weakened topology on kernels}
Let
$K: (\mathbb{R}^{d+1}\times \mathbb{R}^{d+1})\setminus \triangle \to \mathbb{R}$ be a kernel supported on
$ \{(z,z')\in (\mathbb{R}^{d+1})^{2}\ : \ |z-z'|_\fraks  \leq C \} $ for some\footnote{ fixed throughout the article.} $C>0$.
Given a decomposition $K(z,z')= \sum_{n\geq 0} K_n(z,z')$ such that each $K_n$ for $n\geq 1$ is supported on 
$ \{(z,z')\ : \ |z-z'|_\fraks \leq 2^{-n} \} $ we set 
\begin{equ}\label{zwischendefinition}
\vertii{\{K_n\}_n }_{\beta;L,R} \eqdef\sup_{n\in \mathbb{N}} 
 \|\CS_{\beta,n} K_n\|_{C_\fraks^{L,R}}\;,\quad \text{where} \quad
 (\CS_{\beta,n} F)(z) = 2^{(\beta-|\s|)n} F(\mathcal{S}_\fraks^{(2^{n})}z, \mathcal{S}_\fraks^{(2^{n})}\bar z )\;.
\end{equ}
Furthermore, for $\phi\in \mathfrak{B}_{R}$ write $Y^{\lambda}_{z_0,n}(\phi)\eqdef \int_{\mathbb{R}^{d+1}} \phi_{z_0}^{\lambda}(z) K_n(z,\cdot) dz$ 
and set 
$$\llbracket \{K_n\}_{n} \rrbracket_{\beta;R}\eqdef \sup_{z_0\in \mathbb{R}^{d+1}}\sup_{\phi\in \mathfrak{B}_{R}} \sup_{\lambda\in (0,1)} \frac{\sum_{2^{-n}<\lambda} \|Y^{\lambda}_{z_0,n}(\phi)\|_{\mathcal{B}^{2\lambda}_R}}{\lambda^{\beta}} \ , $$
where the notation $\|\ \cdot \ \|_{\mathcal{B}^{2\lambda}_R}$ was introduced in Remark~\ref{rem:useful notation norm}.
%
%

\begin{definition}\label{def:kernel_norm}
Let $\bfK^\beta_{L, R}$ be the space of kernels 
$K: \mathbb{R}^{d+1}\times \mathbb{R}^{d+1}\setminus \triangle \to \mathbb{R}$ supported on
$ \{(z,z')\in (\mathbb{R}^{d+1})^{2}\ : \ |z-z'|_\fraks  \leq C, \ z_{d+1}\geq z'_{d+1}  \} $ 
such that the following norm is finite 
\begin{equ}\label{eq:kernel_norm}
\vertiii{K}_{\beta;L,R} =\inf_{\{K_n\}_{n\geq 0}} \left(   \vertii{\{K_n\}_{n}}_{\beta;L,R} +\llbracket \{K_n\}_{n} \rrbracket_{R,\beta} \right)\;,
\end{equ}
where the infimum is taken over all kernel decompositions $K(z,z')= \sum_{n\geq 0} K_n(z,z')$ such that each $K_n$ for $n\geq 1$ is supported on 
$ \{(z,z')\ : \ |z-z'|_\fraks \leq 2^{-n}, \  z_{d+1}\geq z'_{d+1} \} $. We set $\bfK^\beta_{\infty}\eqdef \bigcap_{L,R>0} \bfK^\beta_{L, R}$.
\end{definition}
\begin{remark}\label{remark, first norm explicit}
Unravelling the definitions one sees that 
\begin{equ}
\vertii{\{K_n\}_n}_{\beta;L,R} =\sup_{n} \, \vertii{K_n}_{\beta;L,R;n}
\end{equ}
where for $G(z,z'): (\mathbb{R}^{d+1})^{\times 2} \to \mathbb{R}$ we write
\begin{equs}
\vertii{G}_{\beta;L,R;n} &= \max_{|l|_\fraks\leq L,|r|_\fraks\leq R}  \frac{\|D_1^{l} D_2^{r} G\|_{L^\infty }}{2^{n(|\fraks|- \beta + |l|_\fraks + |r|_\fraks) }}  
+
\max_{|l|_\fraks\leq L}\frac{\| D^l_1 G\|_{C^{0,R}_\fraks}}{2^{n(|\fraks|- \beta  +l+R)} } \\
&\qquad+
\max_{|r|_\fraks\leq R} \frac{\|D^r_2 G\|_{C^{L,0}_\fraks}}{2^{n(|\fraks|- \beta +L+r)} }
 +
 \frac{\|G\|_{C^{L,R}_\fraks}}{2^{n(|\fraks|- \beta +L+R)} }\ .
\end{equs}
Thus, this norm quantifies \cite[Ass.~5.1, Eq.~5.4]{Hai14}.
\end{remark}

\begin{remark}
The term $\llbracket \{K_n\}_{n} \rrbracket_{\beta;R}$ is a sufficient replacement for \cite[Ass.~5.1, Eq.~5.5]{Hai14}.
Indeed, that assumption  is used in two places therein, namely in the proofs 
of \cite[Lems~5.19 \& 5.15]{Hai14}, where in both cases 
one is given a distribution $F_z$ satisfying $|F_z(\phi^\lambda_z)|\leq C \lambda^{\alpha}$ uniformly over $\phi\in \mathfrak{B}_R$, $\lambda\in (0,1]$ for some $\alpha\in \mathbb{R}$.
Then, the definition of $\llbracket \{K_n\}_{n} \rrbracket_{R}$ guarantees that
$$\sum_{2^{-n}<\lambda} |F_z(Y^{\lambda}_{z,n})|\leq C \llbracket \{K_n\}_{n} \rrbracket_{\beta;R} \lambda^{\alpha+\beta}$$
uniformly over  $\phi\in \mathfrak{B}_R$, $\lambda\in (0,1]$ which is exactly what is required in \cite[Secs~6--7]{Hai14} and \cite[Sec.~2.4]{Sin23}.
\end{remark}
We spell out the following version of the Schauder estimate taking into account continuity with respect to the kernel.

\begin{prop}\label{prop:Schauder_continuity in kernel} 
Let $(T,G)$ be a regularity structure satisfying Assumption~\ref{ass:irrational_beta}, 
 $V$ be a sector of regularity $\alpha>-2$ and $I:V\to T$ an abstract integration map.
 Let  $M$ a model realising the kernel $K\in \bfK^\beta_{L, R}$ for $I$ and
let $f\in\mathcal{D}^{\gamma,\eta}(V)$ with $-2<\eta<\gamma$ and $\gamma>0$. Provided that $\gamma+\beta, \eta+ \beta \notin \mathbb{N}$ and $L>\gamma+\beta$, $R>\alpha\wedge\eta$, it holds that $\mathcal{K}_\gamma f$ defined in \eqref{eq:kernel lift} belongs to $\mathcal{D}^{\gamma+\beta, (\eta \wedge \alpha)+ \beta}$ and satisfies
$$\mathcal{R}\mathcal{K}_\gamma f= K(\mathcal{R}f)\ .$$
Furthermore, considering a second model $\bar{M}$ realising $I$ for $\bar{K}\in \bfK^\beta_{L, R}$ and denoting by $\bar{\mathcal{K} }_\gamma$ its lift
\begin{equ}\label{eq:schauder_eps cont.}
\vertiii{\mathcal{K}{f};  \bar{\mathcal{K}} \bar{f}   }_{\gamma+\beta, (\eta \wedge \alpha)+ \beta, \mfK } \lesssim 
\vertiii{K-\bar{K}}_{\beta;L,R} +
\vertiii{f;\bar f}_{\gamma,\eta, \bar{\mfK}}+\vertiii{ M; \bar{M} }_{\gamma, \bar{\mfK}}  \ ,
\end{equ}
for $f\in \mathcal{D}_{M}^{\gamma,\eta}(V)$, $\bar{f}\in \mathcal{D}_{\bar M}^{\gamma,\eta}(V)$,
where 
the implicit constant depends (continuously) on $\vertiii{K}_{\beta;L,R},$ $\vertiii{\bar K}_{\beta;L,R}$, $\vertiii{f}_{\gamma,\eta, \bar{\mfK}}$, $\vertiii{\bar f}_{\gamma,\eta, \bar{\mfK}}$, $\vertiii{M}_{\gamma,\bar{\mfK}}$ and $\vertiii{\bar M}_{\gamma,\bar{\mfK}}$.  
\end{prop}



\subsubsection{Comparing topologies}
In this section we aim to compare the topology on kernels introduced here with the one introduced in \cite[Def.~1]{Sin23}. First let 
\begin{align*}
\llceil \{K_n\}_n \rrceil_{\beta,R} \eqdef& \max_{|k_1|_\fraks,|k_2|_\fraks<R}  \sup_{n\in \mathbb{N}} \sup_{z'}\frac{ |\int_{\mathbb{R}^d} (z-z')^{k_1} D^{k_2}_2 K_n(z,z') dz|}{2^{-\beta n} }  \\
&\qquad +
\max_{|k_1|_\fraks<R}\sup_{n\in \mathbb{N}} \sup_{|z'-w|_\fraks \leq 1} \frac{ |\int_{\mathbb{R}^d} (z-z')^{k_1}  \big(K_n(z,z') - \opP^R_w[K_n(z,\, \cdot\,)](z') \big) dz|}{2^{-\beta n}|z'-w|^R_\fraks} 
\;, 
\end{align*}
and note that, by \cite[Prop.~A.1]{Hai14}, the quantity $\inf_{\{K_n\}_{n\geq 0}} \bigl( \vertii{ \{K_n\}_n }_{\beta;L,R} + \llceil \{K_n\}_n \rrceil_{\beta;R} \bigr)$ is bounded by the norm introduced
in \cite[Def.~1]{Sin23}. 
Therefore the next lemma shows
that the norm on kernels in Definition~\ref{def:kernel_norm} is indeed weaker than the one introduced therein.
\begin{lemma}\label{lem:weaker_norm}
It holds that for $R,L>0, \beta\in (0,|\fraks|)$
$$ \llbracket \{K_n\}_{n} \rrbracket_{\beta;R} \lesssim \vertii{\{K_n\}_n }_{\beta;0,R} + \llceil \{K_n\}_n \rrceil_{\beta;R} \ . $$
\end{lemma}
\begin{proof}
Let $\phi\in \mathfrak{B}_{R}$, then for $|k|_\fraks<R$
\begin{align*}
D^k Y^{\lambda}_{z_0,n}(\phi)(w)&= \int_{\mathbb{R}^{d+1}} \phi_{z_0}^{\lambda}(z) D^k_w K_n(z,w) dz\\
&= \int_{\mathbb{R}^{d+1}} \opP^{|k|_\fraks}_w [\phi_{z_0}^{\lambda}](z) D^k_w K_n(z,w) dz 
+ \int_{\mathbb{R}^{d+1}} \Big(\phi_{z_0}^{\lambda}- \opP^{|k|_\fraks}_w [\phi_{z_0}^{\lambda}]\Big)(z) D^k_w K_n(z,w) dz\\
&= \sum_{|l|_\fraks<|k|_\fraks} \frac{D^{l}\phi_{z_0}^{\lambda}(w)}{l!} \cdot\int_{\mathbb{R}^{d+1}}  (z-w)^l D^k_w K_n(z,w) dz \\
&\qquad+ \int_{\mathbb{R}^{d+1}} \Big(\phi_{z_0}^{\lambda}- \opP^{|k|_\fraks}_w [\phi_{z_0}^{\lambda}]\Big)(z) D^k_w K_n(z,w) dz\\
\end{align*}
and thus
\begin{align*}
|D^k Y^{\lambda}_{z_0,n}(\phi)(w)|&\leq \sum_{|l|_\fraks<|k|_\fraks} \lambda^{-|\fraks|-|l|_\fraks}2^{-\beta n} \llceil \{K_n\}_n \rrceil_{\beta;R} + \lambda^{-|\fraks|-|k|_\fraks } 2^{-n(|k|_\fraks-|k|_\fraks+\beta)} \vertii{\{K_n\}_n }_{\beta;0,R}.
\end{align*}
Therefore
\begin{equ}\label{eq:loc_bound}
\sum_{2^{-n}<\lambda} |D^k Y^{\lambda}_{z_0,n}(\phi)(w)| \lesssim \big(\vertii{\{K_n\}_n }_{\beta;L,R} + \llceil \{K_n\}_n \rrceil_{\beta,L,R}\big)  \lambda^{-|\fraks|-k+\beta}   \ .
\end{equ}
%
%
%
Next, we turn to estimating
\begin{equs}
 Y^{\lambda}_{z_0,n}(\phi)(w) &  - \opP^{R}_{w_0}[Y^{\lambda}_{z_0,n}(\phi)](w) \\
 &= \int_{\mathbb{R}^{d+1}} \phi_{z_0}^{\lambda}(z)  \Big(K_n(z,w)- \opP^{R}_{w_0}[K_n(z,\cdot)](w)\Big)  dz\\
 &=   \int_{\mathbb{R}^{d+1}} \opP^{R}_w [\phi_{z_0}^{\lambda}](z) \Big(K_n(z,w)- \opP^{R}_{w_0}[K_n(z,\cdot)](w)\Big)  dz \label{eq:local_ll1}\\
 &\qquad + \int_{\mathbb{R}^{d+1}} \Big(\phi_{z_0}^{\lambda}-\opP^{R}_w [\phi_{z_0}^{\lambda}]\Big)(z)  \Big(K_n(z,w)- \opP^{R}_{w_0}[K_n(z,\cdot)](w)\Big)  dz \label{eq:local_ll2}
\end{equs}
One thus sees that 
$$\sum_{2^{-n}<\lambda} \Bigg| \int_{\mathbb{R}^{d+1}} \opP^{R}_w [\phi_{z_0}^{\lambda}](z) \Big(K_n(z,w)- \opP^{R}_{w_0}[K_n(z,\cdot)](w)\Big)  dz\Bigg| \lesssim \lambda^{-|\fraks|-R+\beta} |w-w_0|^R \llceil \{K_n\}_n \rrceil_{\beta;R}\ .$$
For the term \eqref{eq:local_ll2}
we treat separately the cases $|w-w_0|<2^{-n}$ and $|w-w_0|\geq 2^{-n}$. 
Starting with the former, note that 
$$\int_{\mathbb{R}^{d+1}} \Big|\big(\phi_{z_0}^{\lambda}-\opP^{R}_w [\phi_{z_0}^{\lambda}]\big)(z) \Big|\cdot  \Big|K_n(z,w)- \opP^{R}_{w_0}[K_n(z,\cdot)](w)\Big|  dz\lesssim  \lambda^{-|\fraks|-R} 2^{-n\beta} |w-w_0|^R_\fraks \vertii{\{K_n\}_n }_{\beta;0,R}\;,$$
where in the last line we used that the function $z\mapsto K_n(z,w)- \opP^{R}_{w_0}[K_n(z,\cdot)](w)$ is supported on 
$B_{2^{-n}}(w)\cup B_{2^{-n}}(w_0)\subset B_{2^{-n+1}}(w)$.  Summing over $n\in \mathbb{N}$ such that $ |w-w_0|<2^{-n}<\lambda$ completes this case.

Next we turn to the case $|w-w_0|\geq 2^{-n}$. Write 
\begin{align*}
 &\int_{\mathbb{R}^{d+1}} \Big(\phi_{z_0}^{\lambda}-\opP^{R}_w [\phi_{z_0}^{\lambda}]\Big)(z)  \Big(K_n(z,w)- \opP^{R}_{w_0}[K_n(z,\cdot)](w)\Big)  dz\\
 &=\int_{\mathbb{R}^{d+1}} \Big(\phi_{z_0}^{\lambda}-\opP^{R}_w [\phi_{z_0}^{\lambda}]\Big)(z)  K_n(z,w) dz + \int_{\mathbb{R}^{d+1}} \Big(\phi_{z_0}^{\lambda}-\opP^{R}_{w_0} [\phi_{z_0}^{\lambda}]\Big)(z)  \opP^{R}_{w_0}[K_n(z,\cdot)](w)  dz\\
 &\qquad + \int_{\mathbb{R}^{d+1}} \Big(\opP^{R}_{w_0} [\phi_{z_0}^{\lambda}]-\opP^{R}_w [\phi_{z_0}^{\lambda}]\Big)(z)  \opP^{R}_{w_0}[K_n(z,\cdot)](w)  dz
\end{align*}
The first summand we estimate by
\begin{equs}
\int_{\mathbb{R}^{d+1}} \Big|\phi_{z_0}^{\lambda}-\opP^{R}_w [\phi_{z_0}^{\lambda}]\Big|(z)  |K_n(z,w)| dz
&\leq \vertii{\{K_n\}_n}_{\beta;0,0} \lambda^{-|\fraks|-R} 2^{-n(\beta+R)}\\
&\leq \vertii{\{K_n\}_n}_{\beta;0,0} \lambda^{-|\fraks|-R} |w-w_0|^R 2^{-n\beta}\label{eq:loc_e}
\ .
\end{equs}
The second term satisfies similarly
$$\int_{\mathbb{R}^{d+1}} \Big(\phi_{z_0}^{\lambda}-\opP^{R}_{w_0} [\phi_{z_0}^{\lambda}]\Big)(z)  \opP^{R}_{w_0}[K_n(z,\cdot)](w)  dz
=  \sum_{|k|_\fraks<R} \frac{(w-w_0)^k}{k!}  \int_{\mathbb{R}^{d+1}} \Big(\phi_{z_0}^{\lambda}-\opP^{R}_{w_0} [\phi_{z_0}^{\lambda}]\Big)(z) D^k K_n(z,w_0)  dz, $$
where 
$\int_{\mathbb{R}^{d+1}} |(\phi_{z_0}^{\lambda}-\opP^{R}_{w_0} [\phi_{z_0}^{\lambda}])(z)| \cdot |D^k K_n(z,w_0)|  dz \lesssim \vertii{\{K_n\}_n}_{\beta;0,R}   \lambda^{-|\fraks|-R} 2^{-n(\beta+R-k)} $ and therefore
\begin{equs}
\Big|\int_{\mathbb{R}^{d+1}} \Big(\phi_{z_0}^{\lambda}-\opP^{R}_{w_0} [\phi_{z_0}^{\lambda}]\Big)(z)  \opP^{R}_{w_0}[K_n(z,\cdot)](w)  dz \Big|&\lesssim\vertii{\{K_n\}_n}_{\beta;0,R} \sum_{|k|_\fraks<R} |w-w_0|^k \lambda^{-|\fraks|-R} 2^{-n(\beta+R-k)} \\
&\lesssim\vertii{\{K_n\}_n}_{\beta;0,R}  |w-w_0|^R \lambda^{-|\fraks|-R} 2^{-n\beta} \label{eq:loc_ee}
\end{equs}
For the last term, note that using the notation from Remark~\ref{rem:poly_example}
\begin{equs}
\big(\opP^{R}_{w_0} [\phi_{z_0}^{\lambda}]-\opP^{R}_w [\phi_{z_0}^{\lambda}]\big)(z)&= \Pi_{w_0}\left( \opbfP^R_{w_0} [\phi_{z_0}^{\lambda}] - \Gamma_{w_0,w}\opbfP^R_{w} [\phi_{z_0}^{\lambda}]\right)(z) \\
&= \sum_{|k|_\fraks<R} \frac{(z-w_0)^k}{k!} \big(D^k\phi_{z_0}^{\lambda}(w_0)- \opP^{R-k}_w[D^k\phi_{z_0}^{\lambda}](w_0)\big) 
\end{equs}
one finds that 
\begin{align*}
 \int_{\mathbb{R}^{d+1}} \Big(\opP^{R}_{w_0}& [\phi_{z_0}^{\lambda}]-\opP^{R}_w [\phi_{z_0}^{\lambda}]\Big)(z)  \opP^{R}_{w_0}[K_n(z,\cdot)](w)  dz\\
&=\sum_{|l|\fraks<R} \frac{(w-w_0)^l}{l!}  \int_{\mathbb{R}^{d+1}} \Big(\opP^{R}_{w_0} [\phi_{z_0}^{\lambda}]-\opP^{R}_w [\phi_{z_0}^{\lambda}]\Big)(z) D^l K_n(z,w_0)  dz\\
&=\sum_{|l|_\fraks<R} \frac{(w-w_0)^l}{l!}  \int_{\mathbb{R}^{d+1}}  \sum_{|k|\fraks<R} \frac{(z-w_0)^k}{k!} \big(D^k\phi_{z_0}^{\lambda}(w_0)- \opP^{R-k}_w[D^k\phi_{z_0}^{\lambda}](w_0)\big)  D^l K_n(z,w_0)  dz\\
&=\sum_{|l|_\fraks,|k|_\fraks<R} \frac{(w-w_0)^l}{l!} \big(D^k\phi_{z_0}^{\lambda}(w_0)- \opP^{R-k}_w[D^k\phi_{z_0}^{\lambda}](w_0)\big) \int_{\mathbb{R}^{d+1}}  \frac{(z-w_0)^k}{k!}   D^l K_n(z,w_0)  dz\ .
\end{align*}
Noting that 
$|\int_{\mathbb{R}^{d+1}}  \frac{(z-w_0)^k}{k!}   D^l K_n(z,w_0)  dz| \lesssim (\vertii{\{K_n\}_n}_{\beta;0,R}+ \llceil \{K_n\}_n \rrceil_{\beta;R})   2^{-n(\beta + 0\vee (k-l) )}$, 
we thus conclude that since $2^{-n}<|w-w_0|$
\begin{align}
& \int_{\mathbb{R}^{d+1}}  \Big(\opP^{R}_{w_0} [\phi_{z_0}^{\lambda}]-\opP^{R}_w [\phi_{z_0}^{\lambda}]\Big)(z)  \opP^{R}_{w_0}[K_n(z,\cdot)](w)  dz \nonumber \\
&\lesssim \sum_{|l|_\fraks,|k|_\fraks<R} |w-w_0|^{R+l-k} \lambda^{-|\fraks|-R} 2^{-n(\beta + 0\vee (k-l) )}  (\vertii{\{K_n\}_n}_{\beta;0,R}+ \llceil \{K_n\}_n \rrceil_{\beta;R}) \nonumber\\
&\lesssim |w-w_0|^{R}\lambda^{-|\fraks|-R} 2^{-n\beta}  (\vertii{\{K_n\}_n}_{\beta;0,R}+ \llceil \{K_n\}_n \rrceil_{\beta;R})\ . 
\end{align}
%
Finally, combining this with \eqref{eq:loc_e} and \eqref{eq:loc_ee}, and summing over $2^{-n}<\lambda$ completes the proof.
\end{proof}

\subsection{Lifting corrector terms} \label{sec:lifting corrector terms}
For $\beta\in (1,|\fraks|)$ and $\kappa>0$ we shall lift in this section
the maps 
$$ \mathcal{C}^\alpha \to  \mathcal{C}^{\alpha +\beta-\kappa}, \qquad f(\cdot) \mapsto \eps\cdot (\psi\circ \mathcal{S}^{\eps^{-1}}) \cdot K_\eps(f)$$
where $K_\eps$ is a $\beta-1$ regularising kernel at scales larger than $\eps$, to a $\beta-\kappa$ regularising map on 
spaces of modelled distributions.

We first lift multiplication by $\eps$ to an abstract operation, see also \cite{HQ18} for a similar, though distinct situation.
\begin{definition}\label{def:abstract scale}
Given a sector $V$, we say that $\mathcal{E}: {V}\to {T}$ is an abstract multiplication by a scale parameter, if 
\begin{itemize}
\item $\mathcal{E}: {V}_\alpha\to {T}_{\alpha+1}$ for every $\alpha\in A$,
\item $\mathcal{E}|_{\bar{T}}=0$
\item $\mathcal{E}\Gamma \tau - \Gamma \mathcal{E} \tau \in \bar{{T}}$ for every $\tau\in V$ and $\Gamma \in G$.
\end{itemize}
Given $\eps\geq 0$, we say a model $(\Pi, \Gamma)$ realises $\mathcal{E}$ on $V$ at scale $\eps$ if, 
$$J^{\eps}(x)\tau\eqdef\begin{cases}
\eps \sum_{|l|_\fraks<|\tau|+1} \frac{X^l}{l!}  D^l  (\Pi_x \tau)(x) 
& \text{if } \eps>0\\
0 &\text{else.}
\end{cases}
$$
is well defined\footnote{
i.e.  $\lim_{\lambda\to 0}  (\Pi_x \tau)(D^l\phi^\lambda_x) $ exists for every $\phi\in C^\infty_c$ and $|l|_\fraks<|\tau|+1$.
} and 
$\Pi_x \mathcal{E}\tau = \eps \Pi_x \tau - \Pi_x J^{\eps}(x)\tau$ for every $\tau\in V$.

Finally, define $\hat{\mathcal{E}}^\eps_\gamma\eqdef \mathcal{E} + J^{\eps} +\mathcal{N}^\eps_{\gamma} $  whenever 
$\mathcal{N}^\eps_{\gamma} f\eqdef \eps   \sum_{|l|_\fraks<\gamma+1} \frac{X^l}{l!}  D^l  (\mathcal{R}f-\Pi_x f(x))(x)$ is well defined (and with the understanding that 
$\mathcal{N}^0_{\gamma} f\eqdef0$).
\end{definition}
%
%

\begin{assumption}\label{Ass:kernel lift}
For a \textit{scale} $\eps \in [0,1]$, let  
 $K_\eps\in \mfK_{\infty}^{\beta-1}$ be such that $K_\eps(z,z')=0$ whenever $|z-z'|<\eps/2$.
Furthermore, let $\psi\in C(\mathbb{R}^{d+1}, \mathbb{R})$ be Lipschitz continuous and $\mathbb{Z}^{d+1}$- periodic with mean $0$.
\end{assumption}

Throughout this section we make the following assumptions on the regularity structure.

\begin{assumption}\label{Ass:RS1}
We are given a triple of sectors $(V,\, V',\, V'')$ and $\beta\in (1, | \fraks|)$, $\kappa\in (0,1\wedge \beta-1)$ with the property that 
\begin{equ}\label{eq:assumption on intersections of indices}
[\alpha+\beta -1 -\kappa, \alpha+\beta -1]\cap \mathbb{N}= \emptyset \qquad \text{for all } \alpha \in A\;, 
\end{equ} 
such that 
\begin{enumerate}
\item There are abstract integration maps $I: V\to V'$ of regularity $\beta-1$  and  $I^+:V\to T$ of regularity $\beta- {\kappa}$.
\item The sector $V''$ contains the span of an abstract noise symbol $\Psi$ which is of homogeneity $-\kappa$.
\item There is an abstract multiplication by a scale parameter $\mathcal{E}: V''\to T$. 
\item There is an abstract product $V' \otimes \langle \mathbf{1}, \Psi, \mathcal{E}(\Psi)\rangle \to V''$.
\end{enumerate}
\end{assumption}
We also make the following corresponding assumption on models.
\begin{assumption}\label{Ass:RS2}
For a scale $\eps \in [0,1]$ and $(K_\eps, \psi)$ as in Assumption~\ref{Ass:kernel lift} as well as a regularity structure equipped with $(V,\, V',\, V'', I, I^+, \Psi, \mathcal{E})$ as in Assumption~\ref{Ass:RS1}, we shall consider models $M=(\Pi, \Gamma)$ satisfying the following.
\begin{enumerate}
\item It realises $K_\eps$ as a $\beta-1$ regularising kernel for $I$ and, if $\eps>0$ as a $\beta-\kappa$ regularising kernel for $I^+$. (We correspondingly use the notation $J^{K,+}$ and $\mathcal{N}^{K,+}$.)
\item It holds that $\Pi_z \Psi= \psi^\eps$ for all $z\in \mathbb{R}^{d+1}$ for $\eps>0$ and $\Pi_z \Psi= 0$ if $\eps=0$.
\item The model realises $\mathcal{E}$ on $V''$ at scale $\eps\geq 0$.
\item\label{qualitative item} It holds that $\Pi_x (v\cdot \tau)= \Pi_x v \cdot \Pi_x \tau$ for any $v \in V'$ and $\tau\in\{  \mathbf{1}, \Psi, \mathcal{E}(\Psi)\}$.
\end{enumerate}
\end{assumption}

Note that given a sector $V$, it follows from the usual extension theorem, \cite{Hai14}, that one can enlarge the regularity structure to support an abstract 
integration maps $I$ and $I^+$ and
any model $(\Pi,\Gamma)$ can be extended to realise $K_\eps$ for $I$. 
Furthermore, one can always extend the regularity structure to contain $\langle \mathbf{1}, \Psi, \mathcal{E}(\Psi)\rangle$ and carry a product. Since $K_\eps$ is smooth for $\eps>0$ one can furthermore enforce Item~\ref{qualitative item}.  
\begin{remark}
It is actually not necessary to include $I^+$ into the regularity structure, since it mainly serves as a tool in the proof of Proposition~\ref{prop:main_schauder}. We chose to include it, in order not to have to introduce an auxiliary regularity structure in the proof, whose description would have noticeably lengthened the argument.
\end{remark}

The following lemmas will be useful when controlling models in the limit $\eps\to 0$.
\begin{lemma}\label{lem:identities for abstract mult}
In the setting of Assumptions~\ref{Ass:RS1}\&\ref{Ass:RS2} it holds that for all $\tau\in V$,
\begin{equs}\label{eq:JPsiI formula}
J^{\eps}(x)(\Psi I\tau ) &= \eps \psi^\eps(x)\,  (J^{K+}(x)\tau-J^K(x) \tau) \ , \label{eq:JPsiI formula} \\
 \Pi_x \mathcal{E}(\Psi I\tau )&= \eps \psi^\eps(x) \, \Pi_x I^+\tau +  \Pi_x  \mathcal{E}(\Psi) \cdot I\tau \ , \label{eq:Jformula_added_late} \\
\Pi_x \mathcal{E}(\Psi \Gamma_{x,y} I\tau )&= \eps \psi^\eps(x) \, \Pi_x I^+\Gamma_{x,y}\tau +  \Pi_x  \mathcal{E}(\Psi) \cdot\Gamma_{x,y} I \tau \ . \label{eq:Jformula_added_verylate} 
\end{equs}
\end{lemma}
\begin{proof}
Unravelling the definitions we note 
that
\begin{equs}
J^{\eps}(x)(\Psi I\tau ) &=  \eps \sum_{|l|_\fraks<|\tau|+\beta -\kappa} \frac{X^l}{l!}  D^l  (\Pi_x \Psi  I \tau)(x)\\
&= \eps \psi^\eps(x) \sum_{|\tau|+\beta-1-\kappa \leq |l|_\fraks<|\tau|+\beta -\kappa}  \frac{X^l}{l!}  D^l  (\Pi_x  I \tau)(x)\\
\end{equs}
and the first identity follows using \eqref{eq:assumption on intersections of indices}.
Next,
\begin{align*}
 \Pi_x \mathcal{E}(\Psi I\tau )&= \eps \Pi_x (\Psi I\tau )- \Pi_x J^{\eps}(x)(\Psi I\tau )\\
 &= \eps \psi^\eps (\cdot)\, \Pi_x I\tau-  \eps \psi^\eps(x) \sum_{|\tau|+\beta-1-\kappa \leq |l|_\fraks<|\tau|+\beta -\kappa} \frac{\Pi_xX^l (\cdot)}{l!}  D^l  (\Pi_x  I \tau)(x)\\
 &= \eps \psi^\eps (x) \, \Pi_x I\tau-  \eps \psi^\eps(x) \sum_{|\tau|+\beta-1-\kappa \leq |l|_\fraks<|\tau|+\beta -\kappa}  \frac{\Pi_xX^l (\cdot)}{l!}  D^l  (\Pi_x  I \tau)(x)\\
&\qquad+\eps \left(\psi^\eps (\cdot)-\psi^\eps(x)\right) \Pi_x I\tau 
\ ,
\end{align*}
which is indeed the second identity. Next observe the following useful identity.
\begin{equ}\label{eq:J formula}
J^\eps (\Psi\cdot X^k )= \eps \psi^\eps(x) X^k \ .
\end{equ}
To see the last identity, first note that 
\begin{align*}
\Pi_x \mathcal{E}\big( \Psi \Gamma_{x,y} I\tau \big) &= \eps \Pi_x  \Psi \Gamma_{x,y} I\tau  - \Pi_x J^\eps(x) \Psi \Gamma_{x,y} I\tau \big) \\
&= \eps \Pi_x  \Psi \Gamma_{x,y} I\tau  - \Pi_x J^\eps(x) \Psi  I\Gamma_{x,y} \tau \big) - \Pi_x J^\eps(x) \Psi \tilde{p}_{x,y} \big) \ , 
\end{align*}
where $\tilde{p}_{x,y}= \Gamma_{x,y} I\tau- I\Gamma_{x,y} \tau\in \bar{T}$. Thus we find by \eqref{eq:Jformula_added_late} and \eqref{eq:J formula} that
\begin{align*}
\Pi_x \mathcal{E}\big( \Psi \Gamma_{x,y} I\tau \big) 
&=\eps \Pi_x  \Psi \Gamma_{x,y} I\tau -\eps \psi^\eps(x) \Pi_x\big( J^{K+}(x) \Gamma_{x,y} \tau -J^{K}(x) \Gamma_{x,y}\tau + \tilde{p}_{x,y} \big)\\
&=\eps \psi^\eps(\cdot) \Pi_x \Gamma_{x,y} I\tau - \eps \psi^\eps(x) \Pi_x\big( J^{K+}(x) \Gamma_{x,y} \tau -J^{K}(x) \Gamma_{x,y}\tau +\Gamma_{x,y} I\tau- I\Gamma_{x,y} \tau \big)\\
&= \eps \big( \psi^\eps(\cdot) - \eps \psi^\eps(x) \big)\Pi_x \Gamma_{x,y} I\tau 
+ \eps \psi^\eps(x) \Pi_x\big( -J^{K+}(x) \Gamma_{x,y} \tau +J^{K}(x) \Gamma_{x,y}\tau + I\Gamma_{x,y} \tau \big)\\
&= \Pi_x \mathcal{E}(\Psi) \Gamma_{x,y} I\tau 
+ \eps \psi^\eps(x) \Pi_x I^+\Gamma_{x,y} \tau \ .
\end{align*}
\end{proof}
In order to formulate the next lemma, give $\tau\in T$, we write $\langle \tau \rangle\subset T$ 
for the minimal sector containing $\tau$.

\begin{lemma}\label{lem:estimating Gamma with scale} 
In the setting of Assumptions~\ref{Ass:RS1}\&\ref{Ass:RS2}, the value of 
$(\Gamma_{x,y}\mathcal{E} -\mathcal{E}\Gamma_{x,y}) \Psi I\tau\in \bar{T}$
for $\tau\in V$ is uniquely determined by the knowledge of $\psi$, $\Pi_y I\tau$ and $\Pi_y I^+\tau$. 
Similarly, the value of $\Gamma_{x,y}\mathcal{E} \Psi X^k- \mathcal{E}\Gamma_{x,y}\Psi X^k$ is determined by the knowledge of $\psi$.
Furthermore, it holds that for each $\tau \in V_\alpha$
\begin{equ}\label{eq:estimating_Gamma with scale1}
\sup_{\eta\in A} \frac{|(\Gamma_{x,y}\mathcal{E} -\mathcal{E}\Gamma_{x,y}) \Psi I\tau   |_{\alpha+\beta-\kappa-\eta}}{|\tau| |x-y|^{\eta} } \lesssim \eps^{\kappa}, 
\end{equ}
where the implicit constant only depends on $\|M|_{\langle \tau \rangle}\|$,  
the size of the model restricted to $\langle \tau \rangle$.
Furthermore, it holds that
\begin{equ}\label{eq:estimating_Gamma with scale2}
\sup_{\eta\in A}\frac{|\Gamma_{x,y}\mathcal{E} \Psi X^k- \mathcal{E}\Gamma_{x,y}\Psi X^k |_{k+1-\kappa-\eta}}{|x-y|^{\eta}} \lesssim_{k} \eps^{\kappa}\ .
\end{equ}
Finally, for a second model realising $\mathcal{E}$ at the same scale $\eps>0$ and such that $\Pi_x\Psi = \bar{\Pi}_x \Psi $ it holds that 
for any $\tau \in V_\alpha$
\begin{equ}\label{eq:estimating_Gamma with scale3}
\sup_{\eta\in A}  \frac{|\bar{\Gamma}_{x,y}\mathcal{E} \Psi I\tau -\Gamma_{x,y}\mathcal{E} \Psi I\tau|_{\alpha+\beta-\kappa}}{|\tau| |x-y|^{\eta}} \lesssim \eps^{\kappa} \; \|M |_{\langle \tau \rangle}; \bar{M} |_{\langle \tau \rangle}\| 
\end{equ}
where the implicit constants depend only on $\|M |_{\langle \tau \rangle}\|\vee \|\bar{M} |_{\langle \tau \rangle}\|$.

\end{lemma}

\begin{proof}
Observe that $ p_{x,y}= \Gamma_{x,y}\mathcal{E}\Psi I \tau- \mathcal{E} \Psi \Gamma_{x,y} I\tau\in \bar{T}$
satisfies
\begin{align*}
\Pi_x{p_{x,y}}&= \Pi_x\big[\Gamma_{x,y}\mathcal{E}\Psi I \tau- \mathcal{E} \Psi \Gamma_{x,y} I\tau\big]\\
&= \Pi_x\Gamma_{x,y}\mathcal{E}\Psi I \tau -\Pi_x\mathcal{E} \Psi \Gamma_{x,y} I\tau\\
&= \eps \psi^\eps (y) \, \Pi_y I^{+}\tau +  \Pi_y  \mathcal{E}(\Psi) \cdot I\tau 
 - {\eps \psi^\eps(x) \, \Pi_x I^+\Gamma_{x,y}\tau } -  \Pi_x  \mathcal{E}(\Psi) \cdot\Gamma_{x,y} I \tau \\
&=\eps (\psi^\eps(y)-\psi^{\eps}(x))  \Pi_y I^+\tau  
+\eps \psi^{\eps}(x) \Big( \Pi_y I^+\tau  - {\Pi_x I^+\Gamma_{x,y}\tau\big) } 
+\big( \Pi_y  \mathcal{E}(\Psi) -\Pi_x  \mathcal{E}(\Psi) \big) \Pi_y  I\tau \\
&=\eps (\psi^\eps(y)-\psi^{\eps}(x))  \Pi_y I^+\tau  
+\eps \psi^{\eps}(x) \Pi_x \big(  \Gamma_{x,y}I^+  -{  I^+ \Gamma_{x,y}\big)\tau } 
+\eps \big( \psi^\eps(x) -\psi^\eps(y) \big) \Pi_y  I\tau  \ .
\end{align*}
where we used \eqref{eq:Jformula_added_late} and \eqref{eq:Jformula_added_verylate} in the third equality.
Thus, the value of $p_{x,y}=\big(\Gamma_{x,y}\mathcal{E}- \mathcal{E} \Gamma_{x,y}\big) \Psi I\tau$ is determined as claimed in the first line of the lemma.
One similarly finds that
$
\Gamma_{x,y}\mathcal{E} \Psi X^k - \mathcal{E} \Psi   \Gamma_{x,y}X^k =     J_{x}^\eps \Gamma_{x,y}X^k -   \Gamma_{x,y} J_{y}^\eps X^k  = 
 \eps\big( \psi^\eps(x) -\psi^\eps(y)\big) \Gamma_{x,y}X^k \ ,
$
and it only remains to prove the inequalities. Note that \eqref{eq:estimating_Gamma with scale2} can be read off directly.
In order to obtain \eqref{eq:estimating_Gamma with scale1} choose
 $\lambda\sim |x-y|$ and write $\tilde{\phi}\eqdef \phi_{\mathcal{S}^\lambda(
x-y)}$, then
\begin{align*}
|\Pi_x{p_{x,y}}(\phi^\lambda_{x})| &\lesssim \eps |\Pi_y I^+\tau (\tilde{\phi}^\lambda_{y})| 
+\eps \big| \Pi_x \big(  \Gamma_{x,y}I^+  -{  I^+ \Gamma_{x,y}\big)\tau(\phi^\lambda_{x})\big| } 
 +\eps^{\kappa}|x-y|^{1-\kappa}  |\Pi_y I\tau (\tilde{\phi}^\lambda_{y})|\\
&\lesssim \eps  \lambda^{\alpha+\beta}  +\eps \big| \Pi_x \big(  \Gamma_{x,y}I^+  -{  I^+ \Gamma_{x,y}\big)\tau(\phi^\lambda_{x})\big| } +\eps^\kappa |x-y|^{1-\kappa} \lambda^{\alpha+\beta-1}   \\
&\lesssim
\eps^{\kappa}  |x-y|^{\alpha+\beta-\kappa} \ ,
\end{align*}
where used the bounds in the extension theorem \cite[Thm.~5.14]{Hai14} in the second inequality and used \cite[Lem.~5.21]{Hai14} in the third inequality.
Estimate \eqref{eq:estimating_Gamma with scale3} follows very similarly. 
\end{proof}

The following is the analogue of \cite[Prop.~6.16]{Hai14}.
\begin{prop}\label{prop:main_schauder}
In the setting of Assumptions~\ref{Ass:RS1}\&\ref{Ass:RS2}, let $f\in\mathcal{D}^{\gamma,\eta}(V)$ with $0<\eta<\gamma$ and further assume that the sector $V$ has regularity $\alpha$ with $\eta\wedge \alpha>-2$ and $\alpha+ \beta>0$. Then, provided that $\gamma+\beta, \eta+ \beta \notin \mathbb{N}$ and $L>\gamma+\beta$, $R>\alpha\wedge\eta$ one has 
$\hat{\mathcal{E}}^\eps\left( \Psi \cdot \mathcal{K}_\eps (f)\right)\in \mathcal{D}^{\gamma+\beta -\kappa, (\eta \wedge \alpha)+ \beta-\kappa}\ $ uniformly in $\eps\geq 0$ and it holds that 
\begin{equ}\label{eq:reconstructs to the right thing}
\mathcal{R}\hat{\mathcal{E}}^\eps\left( \Psi \cdot \mathcal{K}_\eps (f)\right)
= 
\begin{cases}
\eps\psi^\eps K_\eps(\mathcal{R}f) & \text{if } \eps>0,\\
0 & \text{else.} 
\end{cases}
\end{equ}

Considering two models $M, \bar{M}$ both realising $\mathcal{
E}$ at the same scale $\eps\geq 0$ it holds that 
\begin{equ}\label{eq:schauder_model cont.2}
\vertiii{ \hat{\mathcal{E}}_{\gamma+\beta-1-\kappa}^\eps\left( \Psi \cdot \mathcal{K}_\eps (f)\right); \hat{\mathcal{E}}_{\gamma+\beta-\kappa}^\eps\left( \Psi \cdot \mathcal{K}_\eps (\bar{f})\right)   }_{\gamma+\beta-\kappa, (\eta \wedge \alpha)+ \beta-\kappa, \mfK } \lesssim 
\vertiii{f;\bar f}_{\gamma,\eta, \bar{\mfK}}+ \vertiii{ M, \bar{M} }_{\gamma,\bar{\mfK}} \ ,
\end{equ}
for $f\in \mathcal{D}_{M}^{\gamma,\eta}(V)$, $\bar{f}\in \mathcal{D}_{\bar M}^{\gamma,\eta}(V)$,
where 
the implicit constant depends (continuously) on the size of  $\vertiii{f}_{\gamma,\eta, \bar{\mfK}}$, $\vertiii{\bar f}_{\gamma,\eta, \bar{\mfK}}$, $\vertiii{M}_{\gamma,\bar{\mfK}}$ and $\vertiii{\bar M}_{\gamma,\bar{\mfK}}$.

Lastly, assuming that $M$ realises $\mathcal{
E}$ at scale  $\eps\geq 0$ while $\bar{M}$ realises $\mathcal{
E}$ at scale $\eps=0$, it holds that 
\begin{equs}
&\vertiii{\hat{\mathcal{E}}^\eps_{\gamma+\beta-1-\kappa} \left( \Psi \cdot \mathcal{K}_\eps ({f})\right);  \hat{\mathcal{E}}^0_{\gamma+\beta-1-\kappa}\left( \Psi \cdot \mathcal{K}_0 (\bar{f})\right)   }_{\gamma+\beta-\kappa, (\eta \wedge \alpha)+ \beta-\kappa, \mfK } \\
&\qquad\qquad\qquad\qquad\qquad \lesssim 
\big ( \eps^\kappa\ \vee \ \vertiii{K_\eps-K_0}_{\beta-1;L,R}\big) +
\vertiii{f;\bar f}_{\gamma,\eta, \bar{\mfK}}+\vertiii{ M, \bar{M} }_{\gamma,\bar{\mfK}}  \ . \label{eq:schauder_eps cont.3}
\end{equs}
again for $f\in \mathcal{D}_{M}^{\gamma,\eta}(V)$, $\bar{f}\in \mathcal{D}_{\bar M}^{\gamma,\eta}(V)$ and again with the implicit constant depending on the same quantities as above.

\end{prop}

\begin{proof}
Given space-time points $x=(x_1,...,x_d, x_{d+1}), \  y=(y_1,...,y_d, y_{d+1})\in \mathbb{R}^{1+d}$, we write
$$|x|_P\eqdef 1\wedge |x_{d+1}|, \qquad  |x,y|_P\eqdef 1\wedge |x_{d+1}|\wedge |y_{d+1}|\ .$$
 We shall also write $\tilde{\gamma}\eqdef \gamma+\beta-1-\kappa$ in the proof since this quantity appears so often.
We first check that $\hat{\mathcal{E}}_{\tilde{\gamma}}^\eps\left( \Psi \cdot \mathcal{K}_\eps (f)\right)\in \mathcal{D}^{\tilde{\gamma}+1, (\eta \wedge \alpha)+ \beta-\kappa}$ with a bound uniform in $\eps\in [0,1]$. Without loss of generality, assume that $\vertiii{f}_{\gamma,\eta,\mfK}\leq 1$ and $\vertiii{K}_{\beta-1;L,R}\leq 1$ for some $R,L>(\gamma+\beta) \vee |\alpha|$.

First let $\eps>0$. Since 
\begin{align*}
\hat{\mathcal{E}}^\eps_{\tilde{\gamma}}\left( \Psi \cdot \mathcal{K}_\eps (f)\right)&= \mathcal{E}\left( \Psi \cdot \mathcal{K}_\eps (f)\right) +  J^{\eps}\left( \Psi \cdot \mathcal{K}_\eps (f)\right)
+ \mathcal{N}^\eps_{\tilde{\gamma}}(\Psi \cdot \mathcal{K}_\eps (f)) 
\end{align*}
the estimates on the non-polynomial parts of $\hat{\mathcal{E}}\left( \Psi \cdot \mathcal{K}_\eps (f)\right)$ for 
$\|\hat{\mathcal{E}}_{\tilde{\gamma}}^\eps\left( \Psi \cdot \mathcal{K}_\eps (f)\right) \|_{\gamma,\eta,K}$
follow directly
from the fact that $\mathcal{K}_\eps (f)\in \mathcal{D}^{\gamma+\beta-1, (\eta \wedge \alpha)+ \beta-1}$ by \cite[Prop.~6.16]{Hai14}.
For a scale decomposition as in \eqref{eq:kernel_norm} we shall write $K_\eps= \sum_{n} K_n$
 and $K^{k,\alpha}_{n,xy}(z)\eqdef D_y^k\left(K_n(y,z)- \opP_x^{\alpha+\beta-1} [K_n(\cdot, z)]\right)$ as after Remark~5.17 in \cite{Hai14}.
One observes that
\begin{equs}
\mathcal{N}^\eps_{\tilde{\gamma}}\left( \Psi \cdot \mathcal{K}_\eps (f)\right)(x) &= \eps   
\sum_{\tilde{\gamma}\leq |l|_\fraks<\tilde{\gamma}+1} \frac{X^l}{l!}   \psi^\eps(x) D^l  (\mathcal{R} \mathcal{K}_\eps (f)-\Pi_x  \mathcal{K}_\eps (f)(x) )(x)\\
&= \eps  \sum_{\tilde{\gamma}\leq |l|_\fraks<\tilde{\gamma}+1}  \frac{X^l}{l!}   \psi^\eps(x)  \sum_{n\geq 0} (\mathcal{R}f-\Pi_x f(x) )(K^{l,\gamma}_{n,xy})\\
&= \eps   \sum_{\tilde{\gamma}\leq |l|_\fraks<\tilde{\gamma}+1}  \frac{X^l}{l!}   \psi^\eps(x)   \sum_{n\geq 0}  (\mathcal{R}f-\Pi_x f(x) )(D^l_xK_{n} (x, \cdot)) 
\label{eq:N_formaul raw}
\\
&= \eps    \psi^\eps(x) \left(  \mathcal{N}_{\gamma-\kappa}^{K,+} (f)(x)-  \mathcal{N}_{\gamma-\kappa}^{K} (f)(x)\right)\label{eq:N_formula}
\end{equs} 
where in the second equality we argue exactly
 as for identity just before \cite[Eq.~5.51]{Hai14} and for the forth inequality we use that $K^{l,\gamma}_{n,xy}=D^l_xK_{n} (x, \cdot)$ for $|l|_\fraks\geq \tilde{\gamma}$.
Thus, using \eqref{eq:JPsiI formula} and \eqref{eq:J formula} one finds
\begin{equs}
 J^{\eps}\left( \Psi \cdot \mathcal{K}_\eps (f)\right)
+ \mathcal{N}^\eps_{\tgamma}(\Psi \cdot \mathcal{K}_\eps (f)) &= 
J^{\eps}\left( \Psi \cdot If\right) + J^{\eps}\left( \Psi \cdot (J(x)f(x)+ \mathcal{N}f(x)\right) + \mathcal{N}^\eps_{\tgamma}(\Psi \cdot \mathcal{K}_\eps (f)) \\
&=\eps \psi^\eps(x)\left( J^{K,+}f(x)+  \mathcal{N}_{\tgamma}^{K,+} (f)(x)\right) 
\end{equs}
Therefore, we conclude using the corresponding bounds on  $J^+,\ \mathcal{N}^+_{\tgamma}$ for the $\beta-\kappa$ regularising kernel $\eps^{1-\kappa}K$ in \cite[Prop.~6.16]{Hai14} that 
\begin{equ}\label{eq:poly_part} 
 |\hat{\mathcal{E}}\left( \Psi \cdot \mathcal{K}_\eps (f)\right)|_k\lesssim \eps^{\kappa} \|\psi^\eps\|_{L^\infty} |x|_P^{(\eta\wedge\alpha+\beta-k)\wedge 0}
\end{equ}

To turn to the estimate on the increment, note that
\begin{equs}
&\Gamma_{y,x} \hat{\mathcal{E}}_{\tgamma}^\eps\left( \Psi \cdot \mathcal{K}_\eps (f)\right)(x)- \hat{\mathcal{E}}^\eps_{\tgamma}\left( \Psi \cdot \mathcal{K}_\eps (f)\right)(y)\\
&= \Gamma_{y,x} {\mathcal{E}}\left( \Psi \cdot \mathcal{K}_\eps (f)\right)(x)
+\Gamma_{y,x} J^{\eps}(x)\left( \Psi \cdot \mathcal{K}_\eps (f)\right)(x) 
+\Gamma_{y,x} \mathcal{N}^\eps_{\tgamma} \left(\Psi \cdot \mathcal{K}_\eps (f)(x)\right)\\
&\qquad
- {\mathcal{E}}\left( \Psi \cdot \mathcal{K}_\eps (f)\right)(y)
- J^{\eps}(y)\left( \Psi \cdot \mathcal{K}_\eps (f)\right)(y) 
- \mathcal{N}^\eps_{\tgamma} \left(\Psi \cdot \mathcal{K}_\eps (f)(y)\right)\\
&= {\mathcal{E}} \left(\Gamma_{y,x} \left(\Psi \cdot \mathcal{K}_\eps (f)\right)(x)-  \Psi \cdot \mathcal{K}_\eps (f)(y)\right)    \label{eq:local ref schauder-1}  
\\
&\qquad
+J^\eps(y)\Gamma_{y,x} \left( \Psi \cdot \mathcal{K}_\eps (f)(x)\right)
+\Gamma_{y,x} \mathcal{N}^\eps_{\tgamma} \left(\Psi \cdot \mathcal{K}_\eps (f)(x)\right)\\
&\qquad
- J^{\eps}(y)\left( \Psi \cdot \mathcal{K}_\eps (f)(y) \right)
- \mathcal{N}^\eps_{\tgamma} \left(\Psi \cdot \mathcal{K}_\eps (f)(y)\right)
\end{equs}
where we used that $\Gamma_{y,x} (\mathcal{E}+ J^\eps(x))= (\mathcal{E}+ J^\eps(y))\Gamma_{y,x} $ in the second equality. 
For each $x,y\in \mathbb{R}^{d+1}$ define the bi-linear map $\pi_{x,y}: \mathfrak{K}^{\beta-1}_{\infty}\times \mathcal{D}^{\gamma,\eta}(V)\to \bar{T}\  $
\begin{equs}
\pi_{x,y;\eps}[K,f] \eqdef& J^\eps(y)\Gamma_{y,x} \left( \Psi \cdot \mathcal{K} (f)(x)\right)
+\Gamma_{y,x} \mathcal{N}^\eps_{\tgamma} \left(\Psi \cdot \mathcal{K} (f)(x)\right)\\
&- J^{\eps}(y)\left( \Psi \cdot \mathcal{K} (f)(y) \right)
- \mathcal{N}^\eps_{\tgamma} \left(\Psi \cdot \mathcal{K} (f)(y)\right)\ .
\end{equs}
This allows to write 
\begin{equ}
\Gamma_{y,x} \hat{\mathcal{E}}_{\tgamma}^\eps\left( \Psi \cdot \mathcal{K}_\eps (f)\right)(x)- \hat{\mathcal{E}}^\eps_{\tgamma}\left( \Psi \cdot \mathcal{K}_\eps (f)\right)(y) ={\mathcal{E}} \left(\Gamma_{y,x} \left(\Psi \cdot \mathcal{K}_\eps (f)\right)(x)-  \Psi \cdot \mathcal{K}_\eps (f)(y)\right) + \pi_{x,y;\eps}[K_\eps,f]\;.  
\end{equ}
Thus, the desired estimate on the non-polynomial part now follows directly. 

We shall estimate separately the terms 
$$\pi_{x,y;\eps}[K,f]= \pi_{x,y}^{<,1}  [K,f]+\pi_{x,y;\eps}^{<,2}  [K,f]+ \pi_{x,y;\eps}^{<,3}  [K,f]\ ,$$ where 
\begin{itemize}
\item $\pi_{x,y;\eps}^{<,1}  [K,f]\eqdef J^\eps(y)\left( \Psi \cdot \big(\Gamma_{y,x} \mathcal{K} (f)(x) -  \mathcal{K} (f)(y)\big) \right)$
\item $ \pi_{x,y;\eps}^{<,2}  [K,f]\eqdef \mathcal{N}^\eps_{\tgamma} \left(\Psi \cdot \mathcal{K} (f)(y)\right)$
\item $\pi_{x,y;\eps}^{<,3}  [K,f] \eqdef \Gamma_{y,x} \mathcal{N}^\eps_{\tgamma} \left(\Psi \cdot \mathcal{K} (f)(x)\right)$
\end{itemize}
%
We rewrite
\begin{align*}
&\pi_{x,y;\eps}^{<,1}  [K,f] =J^\eps(y)\left( \Psi \cdot \big(\Gamma_{y,x} \mathcal{K} (f)(x) -  \mathcal{K} (f)(y)\big) \right)\\
&= \sum_{\alpha<\tgamma, \alpha\notin \mathbb{N}}  J^\eps(y)\left( \Psi \cdot Q_\alpha \big(\Gamma_{y,x} \mathcal{K} (f)(x) -  \mathcal{K} (f)(y)\big) \right) + \sum_{\alpha<\tgamma, \alpha\in \mathbb{N}}  J^\eps(y)\left( \Psi \cdot Q_\alpha \big(\Gamma_{y,x} \mathcal{K} (f)(x) -  \mathcal{K} (f)(y)\big) \right)\\
&=\sum_{\alpha<\tgamma, \alpha\notin \mathbb{N}}  J^\eps(y)\left( \Psi \cdot Q_\alpha \big(I (\Gamma_{y,x} f(x) -  f(y))\big) \right) + \sum_{\alpha<\tgamma, \alpha\in \mathbb{N}}  \eps\psi^\eps(y) \cdot Q_\alpha \big(\Gamma_{y,x} \mathcal{K} (f)(x) -  \mathcal{K} (f)(y)\big) \\
&=\underbrace{ \eps\psi^\eps(y)\left( J^{K,+}- J^{K}\right)\big(\Gamma_{y,x} f(x) -f(y))\big)}_{=:\pi_{x,y}^{<,11}  [K,f]}
 + \underbrace{ \sum_{\alpha<\tgamma, \alpha\in \mathbb{N}}  \eps\psi^\eps(y) \cdot Q_\alpha \big(\Gamma_{y,x} \mathcal{K} (f)(x) -  \mathcal{K} (f)(y)\big)}_{=:\pi_{x,y;\eps}^{<,12} [K,f]} \ ,
\end{align*}
where in the second to last line we used \eqref{eq:J formula}
 and in the last line \eqref{eq:JPsiI formula}.

Note that $\pi_{x,y;\eps} [K_n,f]=0$ for $2^{-n}<\eps$, since then  $K_{n}=0$. Thus we shall assume throughout that $2^{-n}\geq\eps$, first considering the case $|x-y|<2^{-n}$.
 Thus, as in  \cite[Eq.~5.46]{Hai14}
$$ |\pi_{x,y;\eps}^{<,11}  [K_n,f]|_k\lesssim \eps^{\kappa}\| \psi\|_{L^\infty}  |x,y|_P^{\eta-\gamma} \sum_{\delta\in B_{k}} |x-y|_{\fraks}^{\gamma-\delta} 2^{(|k|-\beta-\kappa-\delta)n}\ ,$$
 where $B_k\subset \mathbb{N}^{d+1}$ is a finite set satisfying $|k|-(\beta-\kappa)-\delta<0$. 
Arguing as in \cite[Eq.~45\& 5.46]{Hai14} and the third display on \cite[p.~79]{Hai14} one finds for $2^{-n}\geq \eps$
\begin{equs}
 |\pi_{x,y;\eps}^{<,12}  [K_n,f]|_k\lesssim 
 \eps\| \psi\|_{L^\infty}|x,y|_P^{\eta-\gamma}\Big( & 2^{(k-(\beta-1)-\gamma)n } +\sum_{|l|<\gamma+(\beta-1)-k}|x-y|_\fraks^{|l|}2^{(|k+l| -(\beta-1)-\gamma)n} \\
&\qquad +\sum_{\delta\in \tilde{B}_{k}} |x-y|_{\fraks}^{\gamma-\delta} 2^{(|k|-(\beta-1)-\kappa-\delta)n} \Big),
\end{equs} 
 where
 the finite $\tilde{B}_k\subset \mathbb{N}^{d+1}$ satisfies $|k|-(\beta-1)-\delta<0$.
%
Next we rewrite using \eqref{eq:N_formaul raw}
\begin{equs}
 \pi_{x,y;\eps}^{<,2}  [K,f]&=\mathcal{N}^\eps_{\tgamma} \left(\Psi \cdot \mathcal{K} (f)(y)\right)=\eps   \sum_{\tilde{\gamma}\leq |l|_\fraks<\tilde{\gamma}+1}  \frac{X^l}{l!}   \psi^\eps(x)  \sum_{n}  (\mathcal{R}f-\Pi_x f(x) )(D^l_xK_n (x, \cdot)) 
\end{equs}
And thus
$$
 |\pi_{x,y;\eps}^{<,2}  [K_n,f]|_k\lesssim  
 \eps\| \psi\|_{L^\infty} |x,y|_P^{\eta-\gamma}  (2^{-n})^{ \gamma+\beta-1-k} 
$$ 
Finally, $\pi_{x,y;\eps}^{<,3}$ can be treated very similarly to conclude that 
\begin{equ}\label{eq:kernel small scale bound.}
|\pi_{x,y;\eps}  [K_n,f]|_k \lesssim \eps \| \psi\|_{L^\infty}|x,y|_P^{\eta-\gamma} \sum_{\delta>0} |x-y|_\fraks^{\gamma+(\beta-1)-\delta}2^{-\delta n}
\end{equ}

We turn to $2^{-n}>|x-y|_\fraks$
\begin{equs}
&\pi_{x,y;\eps}[K,f]\\
&=J^\eps(y)\left( \Psi \cdot \left(\Gamma_{y,x} \mathcal{K} (f)(x)- \mathcal{K} (f)(y)\right)\right)
+\Gamma_{y,x} \mathcal{N}^\eps_{\tgamma} \left(\Psi \cdot \mathcal{K}_\eps (f)(x)\right)
- \mathcal{N}^\eps_{\tgamma} \left(\Psi \cdot \mathcal{K} (f)(y)\right)\\
&=\sum_{m\in \mathbb{N}, m<\tgamma} J^\eps(y)\left( \Psi \cdot Q_{m}\left(\Gamma_{y,x} \mathcal{K} (f)(x)- \mathcal{K} (f)(y)\right)\right)
+J^\eps(y)\left( \Psi \cdot I\left( \Gamma_{y,x} (f)(x)- (f)(y)\right)\right)\\
&\qquad
+\Gamma_{y,x} \mathcal{N}^\eps_{\tgamma} \left(\Psi \cdot \mathcal{K} (f)(x)\right)
- \mathcal{N}^\eps_{\tgamma} \left(\Psi \cdot \mathcal{K} (f)(y)\right)\\
&=\eps \psi^\eps(y) \Big( \sum_{m\in \mathbb{N}, m<\tgamma}  Q_{m}\left(\Gamma_{y,x} \mathcal{K} (f)(x)- \mathcal{K} (f)(y)\right) \Big) +\eps\psi^\eps(y) \Big( \left( J^{K,+}- J^{K}\right)\big(\Gamma_{y,x} f(x) -f(y))\Big)\\
&\qquad
+ \eps    \psi^\eps(x) \Gamma_{y,x}\left(  \mathcal{N}_{\gamma-\kappa}^{K,+} (f)(x)-  \mathcal{N}_{\gamma-\kappa}^{K} (f)(x)\right)
- \eps    \psi^\eps(y)\left(  \mathcal{N}_{\gamma-\kappa}^{K,+} (f)(y)-  \mathcal{N}_{\gamma-\kappa}^{K} (f)(y)\right)\ , \label{eq:locallastline}
\end{equs}
where in the last line we used \eqref{eq:N_formula}. 
Using $\Gamma_{y,x}(I+J(x))= (I+J(y))\Gamma_{y,x}$ 
\begin{equs}
 \sum_{m\in \mathbb{N}, m<\tgamma}  Q_{m}\left(\Gamma_{y,x} \mathcal{K} (f)(x)- \mathcal{K} (f)(y)\right) 
&=J^K(y)\left(\Gamma_{y,x}  f(x)-f(y) \right) +\Gamma_{y,x} \mathcal{N}_{\gamma}^K (f)(x)
- \mathcal{N}^K_{\gamma}  f(y)\\
&=J^K(y)\left(\Gamma_{y,x}  f(x)-f(y) \right) +\Gamma_{y,x} \mathcal{N}^K_{\gamma-\kappa} (f)(x)
- \mathcal{N}^K_{\gamma-\kappa}  f(y)
\end{equs}
and thus 
\begin{align*}
\pi_{x,y;\eps}[K,f]= 
&\eps \psi^\eps(y) \left( J^K(y)\left(\Gamma_{y,x}  f(x)-f(y) \right)  +\Gamma_{y,x} \mathcal{N}^K_{\gamma-\kappa} (f)(x)
- \mathcal{N}^K_{\gamma-\kappa} f(y)\right)\\
& \quad+\eps\psi^\eps(y) \left( \left( J^{K,+}- J^{K}\right)\big(\Gamma_{y,x} f(x) -f(y))\right)\\
&\quad
+ \eps    \psi^\eps(x) \Gamma_{y,x}\left(  \mathcal{N}^{K,+}_{\gamma-\kappa} (f)(x)-  \mathcal{N}^{K}_{\gamma-\kappa} (f)(x)\right)
- \eps    \psi^\eps(y)\left(  \mathcal{N}^{K,+}_{\gamma-\kappa} (f)(y)-  \mathcal{N}^{K}_{\gamma-\kappa} (f)(y)\right)\\
&=\eps \psi^\eps (y) \left( J^{K,+}(y)\left(\Gamma_{y,x}  f(x)-f(y) \right)  +\Gamma_{y,x} \mathcal{N}^K_{\gamma-\kappa} (f)(x)
- \mathcal{N}^{K,+}_{\gamma-\kappa} f(y)\right)\\
&\qquad
+ \eps    \psi^\eps(x) \Gamma_{y,x}\left(  \mathcal{N}^{K,+}_{\gamma-\kappa} (f)(x)-  \mathcal{N}^{K}_{\gamma-\kappa} (f)(x)\right)\\
&=\underbrace{\eps \psi^\eps(y) \left( J^{K,+}(y)\left(\Gamma_{y,x}  f(x)-f(y) \right)  +\Gamma_{y,x} \mathcal{N}_{\gamma-\kappa}^{K,+} (f)(x)
- \mathcal{N}_{\gamma-\kappa}^{K,+} f(y)\right)}_{=:\pi_{x,y;\eps}^{>,1} [K,f]   } \\
&\quad
-\underbrace{\eps    (\psi^\eps(y)-\psi^\eps(x)) \Gamma_{y,x}\left(  \mathcal{N}^{K,+} _{\gamma-\kappa}(f)(x)-  \mathcal{N}^{K}_{\gamma-\kappa} (f)(x)\right)}_{=:\pi_{x,y;\eps}^{>,2} [K,f]   }
\end{align*}
Note that $\pi^{>,1} [K,f]$ is exactly $\eps^{\kappa} \psi^\eps(y)$ times the polynomial contribution to the increment 
of the modelled distribution obtained by applying the lift of $\eps^{1-\kappa}K\in \mathfrak{K}^{\beta-\kappa}_{\infty}$ to $f$. Thus, by  \cite[Eq.~5.50]{Hai14} and the equation thereafter
\begin{equ}
|\pi^{>,1}_{x,y;\eps} [K_n,f]|_k\lesssim \eps^{\kappa} \|\psi^\eps\|_{L^\infty} |x,y|_P^{\eta-\gamma} \Big( 
\sum_{\delta>0} 2^{\delta n} |x-y|_{\fraks}^{\delta+\gamma+(\beta-\kappa)-k} 
+
\sum_{\zeta<\gamma}
|x-y|_\fraks^{\gamma-\zeta}2^{(k-(\beta-\kappa)-\delta)n }
\Big)
\end{equ}

where the summation over $\delta>0$ is over a finite set.
%
%
Lastly, using \eqref{eq:N_formaul raw} we find
\begin{equs}
|\pi_{x,y;\eps}^{>,2} [K_n,f] |_m &\lesssim \|\psi\|_{C^{1}} (\eps\wedge   |x-y|_\fraks)  \sum_{\tilde{\gamma}\leq |k|_\fraks < \tilde{\gamma}+1}  
\frac{|\Gamma_{y,x}X^k|_m}{k!} | (\Pi_xf(x)-\mathcal{R}f)(D^k K_n(x,\cdot))|
    ,
\end{equs}
where we note that $\gamma+\beta-1 -k<0$ for each term by \eqref{eq:assumption on intersections of indices}.
%
%
Finally, let us complete the proof that $\hat{\mathcal{E}}^\eps\left( \Psi \cdot \mathcal{K}_\eps (f)\right)\in \mathcal{D}^{\gamma+\beta -\kappa, (\eta \wedge \alpha)+ \beta-\kappa}\ $, by checking the required estimate on $\pi_{x,y}[K, f]= \sum_{n} \pi_{x,y}[K_n, f]$. We consider separately the cases
\begin{itemize}
\item We note that for the summands $2^{-n}\leq |x-y|_\fraks$, since only terms with $2^{-n}\geq \eps$ contribute by \eqref{eq:kernel small scale bound.}
\begin{equs}
 \sum_{2^{-n}\leq |x-y|_\fraks}  \pi_{x,y;\eps}[K_n, f]& \lesssim \eps^{\kappa}
\| \psi\|_{L^\infty}|x,y|_P^{\eta-\gamma} \sum_{\delta>0} |x-y|_\fraks^{\gamma+(\beta-1)-\delta-k}2^{-(\delta+1-\kappa) n}\\
&\lesssim \eps^{\kappa }|x,y|_P^{\eta-\gamma} |x-y|_\fraks^{\gamma+\beta-\kappa-k} \label{eq:ppol1}
\end{equs}
\item For $2^{-n}\geq |x-y|_\fraks$ 
 $$ \sum_{2^{-n}> |x-y|_\fraks}  |\pi_{x,y;\eps}^{>,1}[K_n, f]|_k \lesssim \eps^{\kappa} \| \psi\|_{L^\infty} |x,y|_P^{(\eta\wedge \alpha)-\gamma} |x-y|_\fraks^{\gamma+\beta-\kappa-k}$$ follows as direct consequence of \cite[Prop.~6.16]{Hai14}. Finally, it remains to establish
\begin{equ}
  \sum_{2^{-n}> |x-y|_\fraks}  |\pi_{x,y;\eps}^{>,2}[K_n, f]|_k \lesssim \eps^{\kappa} \| \psi\|_{L^\infty} |x,y|_P^{(\eta\wedge \alpha)-\gamma} |x-y|_\fraks^{\gamma+\beta-\kappa-k} \ , \label{eq:ppol2}
  \end{equ} 
 for which one argues exactly as in the last part of the proof of \cite[Prop.~6.16]{Hai14} using that 
for $2^{-n}\in (|x-y|_\fraks, \frac{|x,y|_P}{2})$
$$ 
| (\Pi_xf(x)-\mathcal{R}f)(D^k K_n(x,\cdot))| \lesssim  |x,y|_{P}^{\eta-\gamma} 2^{-n(\gamma-(\beta-1)-k)}
$$
and for $2^{-n}\geq \frac{|x,y|_P}{2}$
\begin{equs}
| (\Pi_xf(x)-\mathcal{R}f)(D^k K_n(x,\cdot))|&\leq   |\mathcal{R}f (D^k K_n(x,\cdot))| + | \Pi_xf(x)(D^k K_n(x,\cdot))|\\
&\lesssim
2^{-n(\alpha\wedge\eta -(\beta-1)-k )}
+\sum_{\alpha\leq \zeta<\gamma} |x,y|^{(\eta-\zeta)\wedge 0}_P  2^{-n(\zeta-(\beta-1)-k)}
\end{equs}
 as well as that only terms with $2^n>\eps$ contribute and that $\gamma+\beta-1 -k<0$. 
\end{itemize}
Finally, we conclude that 
$\vertiii{ \hat{\mathcal{E}}_{\gamma+\beta-1-\kappa}^\eps\left( \Psi \cdot \mathcal{K}_\eps (f)\right)  }_{\gamma+\beta-\kappa, (\eta \wedge \alpha)+ \beta-\kappa, \mfK } \lesssim 
\vertiii{f}_{\gamma,\eta, \bar \mfK}$
where the implicit constant depends only on the size of the model.

Since the proof of \eqref{eq:schauder_model cont.2} adapts directly in the usual way we turn to the proof of {\eqref{eq:schauder_eps cont.3}}.

To shorten notation, we write $F^\eps\eqdef \hat{\mathcal{E}}^\eps_{\gamma+\beta-1-\kappa} \left( \Psi \cdot \mathcal{K}_\eps ({f})\right)$
and $F^0\eqdef\hat{\mathcal{E}}^0_{\gamma+\beta-1-\kappa}\left( \Psi \cdot \mathcal{K}_0 (\bar{f})\right)$. 
We first consider the required estimates for integer homogeneities $k\in \mathbb{N}$, where one directly estimates 
$|F^0(x)-F^\eps(x)|_k = |F^\eps(x)|_k$ using \eqref{eq:poly_part}, while the estimate on 
 $|\Gamma_{y,x} F^0(x)- F^0(y)- (\Gamma_{y,x} F^\eps(x)- F^\eps(y)|_k=|\pi_{x,y}[K,f]|_{k}$
follows directly from \eqref{eq:ppol1} \& \eqref{eq:ppol2}.
  For $\zeta\in A\setminus \mathbb{N}$ one has
$$
|F^0-F^\eps|_\zeta = \big|\mathcal{E}\big( \Psi \cdot (\mathcal{K}_\eps f - \mathcal{K} \bar f )\big)\big|_\zeta \lesssim |\mathcal{K}_\eps f - \mathcal{K} \bar f |_{\zeta-1+\kappa}  
\ 
$$
and
$$|\Gamma_{y,x} F^0(x)- F^0(y)- (\Gamma_{y,x} F^\eps(x)- F^\eps(y))|_\zeta=|\Gamma_{y,x} \mathcal{K} \bar f(x)- \mathcal{K} \bar f(y)- (\Gamma_{y,x} \mathcal{K}_\eps f(x)- \mathcal{K}_\eps f(y))|_{\zeta-1+\kappa} \ ,$$
concluding the by Proposition \ref{prop:Schauder_continuity in kernel}.

Finally, it only remains to check \eqref{eq:reconstructs to the right thing}, which
for $\eps>0$ follows easily by unravelling the definitions and using that $K_\eps$ is smooth, while for $\eps=0$ it is obvious.


%
\end{proof}

\subsection{An abstract fixed point theorem}\label{sec:An abstract fixed point theorem}
Finally, we have collected the novel ingredients in a format that allows to follow the arguments for \cite[Theorem~7.8]{Hai14} to obtain the following abstract fixed point theorem.

\begin{theorem}\label{thm:fixed point}

Let $\beta\in (1,|\fraks|)$ and $\mathcal{T}=(T,G)$ be a regularity structure satisfying Assumption~\ref{ass:irrational_beta}. Let $V$, $\{\bar{V}_i\}_{i=1}^m$ be sectors of  $(T,G)$ of respective regularities $\zeta, \bar{\zeta}_i \in \mathbb{R}$ such that 
$\zeta\leq \min_i \bar{\zeta_i} + \beta$.
Assume that for each $i$ there is an abstract integration map $I^i: V_i \to V$ of regularity $\beta$. Furthermore 
assume that we are given triples $\{(V_j, V'_j,V''_j)\}_{j=1,...,m}$ of sectors of regularities $(\zeta_j, \zeta_j', \zeta_j'')$
and $\kappa>0$ such that $\zeta\leq \min_i {\zeta_j} + \beta-\kappa$ and satisfying 
Assumption~\ref{Ass:RS1}, where we denote by $I^j, I^{j,+}$ and $\Psi_j$ and the respective integration maps and noise symbol and furthermore assume the involved abstract multiplication map is the same for all $j$, i.e.
$$\mathcal{E}: \bigoplus_j V_j'' \to V\subset T .$$

For $L,R>0$, assume that we are given continuous maps
$$G^i: [0,1]\to \bfK_{L,R}^\beta, \quad \eps \mapsto G_{\eps}^i,\qquad \text{and \qquad}
K^j: [0,1]\to \bfK_{L,R}^{\beta-1}, \qquad \eps \mapsto K_{\eps}^i \ ,$$
and functions $\psi_j\in C(\mathbb{R}^{d+1},\mathbb{R})$ such that $\big(K^j_\eps ,  \psi_j\big)$ satisfy Assumption~\ref{Ass:kernel lift} for each $\eps\in [0,1]$.

We denote by $[0,1]\ltimes \mathcal{M}$ all pairs $(\eps, M)\in [0,1]\times \mathcal{M}$ such that $M$ realises $G^i_\eps$ for $I^i$ and
 $(\bar V_j, V'_j,V''_j, I^j, I^{j,+}, \Psi_j, \mathcal{E})$ satisfies Assumption~\ref{Ass:RS2} for the scale $\eps$ with respect to 
$\big(K^j_\eps, \psi_j\big)$.
Furthermore, let $\mathcal{M}^\eps=\{ (\eps, M)\in   [ 0,1]\ltimes \mathcal{M} \}$.

For $\gamma\geq \bar{\gamma}>0, \ \eta<(\bar{\eta}\wedge\bar{\zeta}  )+ \beta-\kappa ,\   \gamma<\bar{\gamma}+ \beta-\kappa ,\ 
 \bar{\eta}\wedge\bar\zeta > (-\beta+\kappa) \vee (- R)  $, let
$$\bar{\opF}_i :\mathcal{D}^{\gamma,\eta}_P(V)\to  \mathcal{D}^{\bar \gamma,\bar\eta}_P(\bar{V}_i),\qquad \opF_j :\mathcal{D}^{\gamma,\eta}_P(V)\to  \mathcal{D}^{\bar \gamma,\bar\eta}_P({V}_j)$$
 be 
  \textit{strongly locally Lipschitz} in the sense of \cite[Sec.~7.3]{Hai14}.

For $\gamma<L$ and $-R<\min(A)$, consider for $T>0$ 
the solution map 
 as a map $$S_T:  \left ([0,1]\ltimes\mathcal{M}\right) \ltimes \mathcal{D}^{\gamma,\eta}_P (T) \to \mathcal{D}^{\gamma,\eta}_P
 , \qquad ((\eps,M^\eps), v)\mapsto U_\eps \ $$
 to the fixed point problem 
\begin{equ}\label{eq:abstract_fixed_point}
U_\eps= \sum_{i=1}^n \mathcal{G}^i \big(\mathbf{R}^+ \bar{\opF}_i(U)\big) +\sum_{j=1}^m \hat{\mathcal{E}}^\eps\Big( \Psi_j\mathcal{K}_\eps^j \big(\mathbf{R}^+{\opF}_j(U)\big) \Big) + v\ .
\end{equ}

Then, for each bounded
  subset  $B\subset \left ([0,1]\ltimes\mathcal{M}\right)\ltimes \mathcal{D}^{\gamma,\eta}_P$
there exists $T>0$ such that the solution map $S_T$ is well defined on $B$. Furthermore, it is locally Lipschitz continuous uniformly in $\eps\in [0,1]$, i.e.\ for  $(Z^\eps, v), (\bar{Z}^\eps, \bar v) \in \big(\mathcal{M}^\eps\ltimes \mathcal{D}^{\gamma,\eta}_P\big) \cap B$, the respective fixed points $U_\eps, \bar{U}_\eps$ satisfy
$$
\| U_\eps, \bar{U}_\eps\|_{\gamma,\eta;T}  \lesssim_B 
   \vertiii{Z^\eps,\bar{Z}^\eps}_{\gamma;O} +\vertiii{v;\bar{v}}_{\gamma,\eta; T} \ ,
$$
uniformly over $\eps\in [0,1]$.
%

  Lastly, for $(Z^\eps, v)\in \big(\mathcal{M}^\eps\ltimes \mathcal{D}^{\gamma,\eta}_P\big) \cap B$ and $(\bar Z, \bar v) \in \big(\mathcal{M}^0\ltimes \mathcal{D}^{\gamma,\eta}_P\big)\cap B$
$$
\| U_\eps, \bar{U}_0\|_{\gamma,\eta;T} \lesssim_B
 \sum_i \|G_\eps^i-G_0^i\|_{\beta;L,R} + \sum_j \eps^\kappa \vee \| K_\eps^j- K_0^j\|_{\beta-1;L,R}
  + \vertiii{Z^\eps,\bar{Z}}_{\gamma;O} 
  +\vertiii{v;\bar{v}}_{\gamma,\eta; T} \ . $$
\end{theorem}
\begin{remark}
Note that, while in the next section we shall work only with first order expansions, the set-up of the above theorem allows to incorporate higher order two scale expansions, c.f. \cite{KMS07}. 
A term in such a higher order expansion would typically be of the form $\eps^k \psi^\eps \nabla^k K_\eps$. This can then be lifted by 
writing it as $\eps \cdot \psi^\eps \cdot \big( \eps^{k-1} K_\eps)$. Alternatively, one could define an abstract multiplication by the scale $\eps^k$ analogously to what was done above.
\end{remark}
 \begin{remark}
Let us point out that lifting multiplication by a scale to an abstract operator is only necessary for equations when it is important that the solution to the abstract fixed point problem takes values in a function-like sector. When working with polynomial non-linearities this is often not necessary, see also Remark~\ref{rem:polynomial-nonlinearities}.
 \end{remark}

\begin{remark}
Note that the smoothness condition in Definition~\ref{def:abstract scale} and the conditions $\{K_\eps\}_{\eps\geq0}\in \mfK_{\infty}^{\beta-1}$ could clearly be relaxed, which we expect to be necessary when studying homogenisation problems involving operators $\mathcal{L}_\eps= \nabla\cdot A(x,t, x/\eps, t/\eps^2)\nabla $.
\end{remark}

\subsection{Post-processing kernel estimates}\label{sec_post_process}

In this section, we shall modify the heat kernel  
and the kernels appearing in the two scale expansion of Section~\ref{sec:periodic_hom} 
by excising the singularity at coinciding time coordinates similarly to \cite[Sec.~2.1]{HS23per}.
We then check that these modified kernels converge in the correct topologies to be able to apply Theorem~\ref{thm:fixed point}.

Fix $\kappa: \mathbb{R}\to [0,1]$ such that 
 \begin{itemize}
 \item $\kappa(t)=0$ for $t<0$ and for $t>2$,
 \item $\kappa(t)=1$ for $t\in (0,1)$,
 \item $\kappa|_{\mathbb{R}_+}$ is smooth.
 \end{itemize}
Write
\begin{equation}\label{eq:rescaled truncations}
\kappa^\eps(t)=\kappa(t/\eps^2) \qquad \text{as well as} \qquad \kappa^\eps_c(t)= \mathbf{1}_{\{t>0\}}\left(1-\kappa^\eps(t) \right) \ .
\end{equation}
%
For later use we also set for $i,j\in \{0,...,d\}$
$$
\bar{\Gamma}_{\eps}^{i,j} (x,t;\zeta, \tau) :
= \begin{cases} 
\kappa^\eps_c(t-\tau) \bar{\Gamma} (x,t;\zeta, \tau) & \text{if } i=j=0,\\
 \kappa^\eps_c(t-\tau)\partial_{x_i}\bar{\Gamma} (x,t;\zeta, \tau) & \text{if }  j=0 , i>0\\
\eps  \kappa^\eps_c(t-\tau) \partial_{\zeta_j}\bar{\Gamma} (x,t;\zeta, \tau) & \text{if } j>0, i=0 ,\\
\eps \kappa^\eps_c(t-\tau) \partial_{x_i}\partial_{\zeta_j} \bar{\Gamma} (x,t;\zeta, \tau) & \text{if } j>0 , i>0.
\end{cases}
$$
as well as 
$$
\bar{\Gamma}_{\eps}^{i,j,+}:
= \begin{cases} 
\bar{\Gamma}_{\eps}^{i,j} & \text{if } i=0,\\
\eps\bar{\Gamma}_{\eps}^{i,j} & \text{if } i>0. 
\end{cases}
$$

We also fix $\chi: \mathbb{R}^d\to [0,1]$ smooth and compactly supported on $[-2/3, 2/3]^{d}\subset \mathbb{R}^d$ such that $\sum_{k\in \mathbb{Z}^d} \chi(x+k)=1$ for all $x\in \mathbb{R}^d$. Finally, for $\Gamma\in \{\bar{\Gamma}, \ {\Gamma}_{\eps}, \  \bar{\Gamma}^{i,j}_\epsilon, \  \bar{\Gamma}^{i,j,+}_\epsilon\}
$ set
\begin{equ}\label{eq:compactly supported kernels}
{K} (t,x;s,y)=  \sum_{k\in \mathbb{Z}^d}\kappa(t-s)\chi(x-y)	{\Gamma}(x,t;y+k,s)\ 
\end{equ}
and denote the resulting kernels by $\bar{K},\  {K}_{\eps}, \ \bar{K}^{i,j}_\epsilon, \ \bar{K}^{i,j,+}_{\epsilon}$ respectively.
In the case $i,j=0$ we shall also sometimes simply write $\bar{K}_\epsilon$ instead of $\bar{K}^{0,0}_\epsilon$.

%

The proof of the following lemma is straightforward, c.f.\ \cite[Lem.~5.5]{Hai14}.
%
\begin{lemma}\label{lem:convergence of kernels easy terms}
For any $L,R>0$ and there exists a constant $C>0$ such that 
$$\|\bar{K}\|_{2 ;L,R} +\max_{j}\|\bar{K}^{0,j}_\eps\|_{2 ;L,R} + \max_{i>0,j\geq 0}\|\bar{K}^{i,j}_\eps\|_{1 ;L,R} + \max_{i,j}\|\bar{K}^{i,j,+}_\eps\|_{2 ;L,R} \leq C$$ 
as well as for any $\beta\in (1,2)$
$$ \|\bar{K} -\bar{K}_\eps \|_{\beta ;L,R} + \max_{j>0}\|\bar{K}^{0,j}_\eps\|_{\beta ;L,R}+\max_{i>0,j\geq 0}\|\bar{K}^{i,j}_\eps\|_{\beta-1 ;L,R} + \max_{i+j>0}\|\bar{K}^{i,j,+}_\eps\|_{\beta ;L,R}  \leq C\eps^{2-\beta}$$ 
uniformly over $\eps\in (0,1]$.
\end{lemma} 
%


%
%

Next, define
 \begin{equ}\label{eq:def_G}
  G_\eps (z,z')\eqdef K_\eps(z,z')-  \sum_{i,j=0}^d \phi_i^\eps(z) \tilde{\phi}_j^\eps(z')\cdot
\bar{K}^{i,j,+}_\eps ( z,z')\ ,
 \end{equ}
 where we use again the convention $\phi_0=\tilde \phi_0 = 1$.

\begin{prop}\label{prop:vanishing of G kernel}
For every $L,R\in (1,2)$, $\kappa>0$ and $\beta\in (1,2)$ it holds that
\begin{equ}\label{eq:van_of_G}
\vertiii{G_{\eps} }_{\beta;L,R} 
\lesssim_{L,R,\kappa, \beta}\left( \eps^{2- {\beta} }\vee \eps^{4-\beta-L-R-\kappa} \vee \eps^{3-\beta-R-\kappa }\right)
\end{equ}
uniformly over $\eps \in (0,1]$.
\end{prop}
Let $\varphi\in \mathcal{C}^{\infty}_{c}(B_2\setminus B_{1/2})$ be such that $\sum_{n=1}^{\infty}\varphi^{2^{-n}}(z)=1$ for every $x\in B_{1/2}\setminus\{0\}$. 
For the proof of this proposition, we shall work with the decomposition $G_{\eps} = \sum_{n=0}^\infty G_{\eps,n}$ defined by
\begin{equ}\label{eq:decomposition_of_G}
{G}_{\eps, n}(z;z')\eqdef   \varphi^{2^{-n}}(z-z')(\Gamma - 	\kappa^{\eps}_c(t-t') \Gamma^{1,1})(z;z')	 
\end{equ}
as well as 
${G}_{\eps, 0}= {G}^{L',R'}_{\eps}- \sum_{n\geq 1} {G}^{L',R'}_{\eps, n}$.

\begin{proof}
We check that $\| G\|_{\beta;L,R}$ is upper bounded by the right hand side of \eqref{eq:van_of_G}, the estimate on $\llbracket \{G_{\eps,n}\}_n \rrbracket_{\beta;R}$ is the content of Lemma~\ref{lem:second_term_in_G_van} below.
For $2^{-n}< \eps/2$ the contribution of ${\Gamma}^{1,1}$ in \eqref{eq:decomposition_of_G} vanishes and it follows from Corollary~\ref{cor:small scale}
that $\sup_{n: 2^{-n}< \eps/2} \|G_{\eps,n}\|_{2;L,R;n} <+\infty.$ 
This then implies that 
$$\sup_{n: 2^{-n}< \eps/2}  \|G_{\eps,n}\|_{\beta;L,R;n} \leq \sup_{n: 2^{-n}< \eps/2}  2^{-n(2-\beta)} \lesssim \eps^{2-\beta}  \ .$$
We turn to $2^{-n}\geq \eps/2$, here it follows that for any $\kappa<1$
$$ \|G_{\eps,n}\|_{2;L,R;n} \lesssim_{\kappa} \eps^{3-L-R-\kappa } 2^{n}$$
by combining Corollary~\ref{cor:pointwise} and Theorem~\ref{thm:further kernel estimate}.
This implies that 
$$\|G_{\eps,n}\|_{\beta;L,R;n} \lesssim \eps^{3-L-R-\kappa} 2^{n(1 -(2-\beta))} \leq \eps^{4-\beta-L-R-\kappa},$$
as claimed.
\end{proof}

%

\begin{lemma}\label{lem:second_term_in_G_van}
For every $R\in (1,2)$ and $\beta<3-R$ and $\kappa>0$, there exists $C>0$ such that
\begin{equ}
\llbracket \{G_{\eps,n}\}_n \rrbracket_{\beta;R} \leq C\eps^{{3-R-\beta-\kappa}} 
\end{equ}
uniformly over $\eps \in (0,1]$. 
\end{lemma}

\begin{proof}
We first consider $2\eps<2^{-n}<\lambda$ and find that
 $$Y^{\lambda}_{n,z_0}(\bar{z})
 =\int \phi_{z:0}^{\lambda}(z) G_{\eps,n}(z,\bar{z}) dz
 = \int \phi^{2^{-n}}(z,\bar{z})  \phi_{z:0}^{\lambda}(z) \left( \Gamma_\eps - \Gamma^{1,1}_\eps\right)(z,\bar{z}) dz\ .$$ 
\begin{itemize}
\item Using Corollary~\ref{cor:pointwise} one finds $|Y^{\lambda}_{n,z_0}|\lesssim\eps \lambda^{-|\fraks|} 2^{-n(2-1)} 
$. 
\item Integrating by parts
\begin{equs}
\nabla Y^{\lambda}_{n, \eps ,z_0}(\bar{z}) 
&=  \int   \phi^{2^{-n}}(z-\bar{z})  \nabla\phi_{z_0}^{\lambda}(z)  \left( \Gamma_\eps - \Gamma^{1,1}_\eps\right)(z,\bar{z})  dz \\
&\qquad+  \int   \phi^{2^{-n}}(z-\bar{z})  \phi_{z_0}^{\lambda}(z)  \nabla_{z}\left( \Gamma_\eps - \Gamma^{1,1}_\eps\right)(z,\bar{z})  dz\\
 &\qquad + \int \phi^{2^{-n}}(z-\bar{z})  \phi_{z_0}^{\lambda}(z)  \nabla_{\bar{z}} \left( \Gamma_\eps - \Gamma^{1,1}_\eps\right)(z,\bar{z}) dz\\
 &\eqdef Y^{(1),\lambda}_{n, \eps ,z_0}(\bar{z})+ Y^{(2),\lambda}_{n, \eps ,z_0}(\bar{z}) +Y^{(3),\lambda}_{n, \eps ,z_0}(\bar{z})\ .
\end{equs} 
As above $|Y^{(1),\lambda}_{n, \eps ,z_0}|\lesssim \eps 2^{-n(2-1)} \lambda^{-1-|\fraks|} $ and by Corollary~\ref{cor:pointwise}  
$$|Y^{(2),\lambda}_{n, \eps ,z_0}|+ |Y^{(3),\lambda}_{n, \eps ,z_0}| \lesssim \eps \log (1+ 2^n\eps^{-1}) 2^{-n(2-2)} \lambda^{-|\fraks|}\lesssim \eps ^{1-\kappa/2} 2^{n\kappa/2}  \lambda^{|\fraks|} \ .$$
\item Similarly, one finds that 
\begin{equs}
\nabla Y^{(1),\lambda}_{n, \eps ,z_0}(\bar{z})
&=  \int   \phi^{2^{-n}}(z-\bar{z})  \nabla^2\phi_{z_0}^{\lambda}(z)  \left( \Gamma_\eps - \Gamma^{1,1}_\eps\right)(z,\bar{z})  dz \\
&\qquad+  \int   \phi^{2^{-n}}(z-\bar{z}) \nabla \phi_{z_0}^{\lambda}(z)  \nabla_{z}\left( \Gamma_\eps - \Gamma^{1,1}_\eps\right)(z,\bar{z})  dz\\
 &\qquad + \int \phi^{2^{-n}}(z-\bar{z}) \nabla \phi_{z_0}^{\lambda}(z)  \nabla_{\bar{z}} \left( \Gamma_\eps - \Gamma^{1,1}_\eps\right)(z,\bar{z}) dz
\end{equs} 
and thus 
$|\nabla Y^{(1),\lambda}_{n, \eps ,z_0}|\lesssim (\lambda^{-1} + 2^{n})  \eps 2^{-n(2-1)} \lambda^{-1-|\fraks|} $ which implies that for $r\in (0,1)$
$$\|  Y^{(1),\lambda}_{n, \eps ,z_0} \|_{C^r}\lesssim (\lambda^{-r} + 2^{nr})  \eps 2^{-n(2-1)} \lambda^{-1-|\fraks|} \lesssim \eps 2^{n(1-r) } \lambda^{-1-|\fraks|} \ . $$
Similarly, one notes that
\begin{equs}
\nabla Y^{(2),\lambda}_{n, \eps ,z_0}(\bar{z})
&=  \int   \phi^{2^{-n}}(z-\bar{z})  \nabla\phi_{z_0}^{\lambda}(z)  \nabla_{z}\left( \Gamma_\eps - \Gamma^{1,1}_\eps\right)(z,\bar{z})  dz \\
&\qquad+  \int   \phi^{2^{-n}}(z-\bar{z})  \phi_{z_0}^{\lambda}(z)  \nabla^2_{z}\left( \Gamma_\eps - \Gamma^{1,1}_\eps\right)(z,\bar{z})  dz\\
 &\qquad + \int \phi^{2^{-n}}(z-\bar{z})  \phi_{z_0}^{\lambda}(z)  \nabla_{z}\nabla_{\bar{z}} \left( \Gamma_\eps - \Gamma^{1,1}_\eps\right)(z,\bar{z}) dz
\end{equs} 
which implies that
$$| \nabla Y^{(2),\lambda}_{n, \eps ,z_0} |\lesssim \eps (\lambda^{-1} + \eps^{-1}2^n + 2^n) \log (1+ 2^n\eps^{-1}) 2^{-n(2-2)} \lambda^{-|\fraks|} \lesssim
 \eps^{-\kappa/2} 2^{n(1+\kappa/2)} \lambda^{-|\fraks|}\ . 
$$ 
Treating the term $\nabla Y^{(3),\lambda}_{n, \eps ,z_0}(\bar{z})$ the same way this implies that
 $$\|  Y^{(2),\lambda}_{n, \eps ,z_0} \|_{C^r}+\|  Y^{(3),\lambda}_{n, \eps ,z_0} \|_{C^r}\lesssim \eps^{1-r-\kappa/2} 2^{n(r+\kappa/2)} \ .$$
\item Finally, one finds similarly to the above that 
$$\sup_{x\in \mathbb{R}^d, t,t'\in \mathbb{R}}\frac{|Y^{\lambda}_{n,z_0}(x,t) -Y^{\lambda}_{n,z_0}(x,t')| }{|t-t'|^{1/2+r/2}\wedge 1 } \lesssim \eps^{1-r-\kappa/2} 2^{n(r+\kappa/2)}\ . $$

\end{itemize}
%
Therefore, 
\begin{align*}
&\sum_{\eps<2^{-n} \leq \lambda} \|Y^{\lambda}_{z_0,n}(\phi)\|_{\mathcal{B}^{2\lambda}_{1+r}}\\
&\lesssim 
\sum_{\eps<2^{-n} \leq \lambda}  \lambda^{|\fraks|}   \|Y^{\lambda}_{n, \eps ,z_0}\|_{L^\infty}
+\sum_{\eps<2^{-n} \leq \lambda} \lambda^{|\fraks|+1} \|\nabla Y^{(1)\lambda}_{n, \eps ,z_0}\|_{L^\infty} 
+\sum_{\eps<2^{-n} \leq \lambda} \lambda^{|\fraks|+1}  \big(\|\nabla Y^{(2)\lambda}_{n, \eps ,z_0}\|_{L^\infty}+\|\nabla Y^{(3)\lambda}_{n, \eps ,z_0}\|_{L^\infty}\big)\\
&\qquad
+ \sum_{\eps<2^{-n} \leq \lambda} \lambda^{|\fraks|+1 +r} \|\nabla Y^{(1)\lambda}_{n, \eps ,z_0}\|_{C^r}
+ \sum_{\eps<2^{-n} \leq \lambda} \lambda^{|\fraks|+1 +r} \big(\|\nabla Y^{(2)\lambda}_{n, \eps ,z_0}\|_{C^r}+ \|\nabla Y^{(3)\lambda}_{n, \eps ,z_0}\|_{C^r}\big) \\
&\qquad +
\sup_{x\in \mathbb{R}^d, t,t'\in \mathbb{R}} \sum_{\eps<2^{-n} \leq \lambda} \lambda^{|\fraks|+1 +r}  \frac{|Y^{\lambda}_{n,z_0}(x,t) -Y^{\lambda}_{n,z_0}(x,t')| }{|t-t'|^{1/2+r/2}\wedge 1}
\\
&= \lambda\eps +\lambda \eps + \lambda \eps^{1-\kappa} +\lambda \eps + \lambda^{1+r}\eps^{1-2r-\kappa}  
\lesssim \lambda\eps^{1-\kappa} + \lambda^{1+r} \eps^{1-2r-\kappa}
 \ .
\end{align*}
Thus for $1+r<\beta$
$$ \frac{\sum_{\eps<2^{-n} \leq \lambda} \|Y^{\lambda}_{z_0,n}(\phi)\|_{\mathcal{B}^{2\lambda}_{1+r}}}{\lambda^{\beta}} \lesssim
\lambda^{1-\beta}\eps^{1-\kappa}+ \lambda^{1+r-\beta}\eps^{1-2r-\kappa} \lesssim \eps^{2-\beta-r-\kappa}
\;.$$
The case  $\eps\sim2^{-n}$ follows as above, since then
 $$Y^{\lambda}_{n,z_0}(\bar{z})
 = \int \phi^{2^{-n}}(x-\bar{x},t-\bar{t}) \kappa^{c}_\eps (t-\bar{t})  \phi_{z:0}^{\lambda}(x,t) \left( \Gamma_\eps - \Gamma^{1,1}_\eps\right)(x,t;\bar{x},\bar{t}) dxdt\ .$$ 
Finally, for
$2^{-n}<\eps/2$, one argues very similarly to above but using the local expansion of Lemma~\ref{lem:close to diagonal} instead of the two scale expansion
to obtain  for $1<\beta\leq 2$
$$ \frac{\sum_{2^{-n}  \leq \eps/2 \wedge  \lambda} \|Y^{\lambda}_{z_0,n}(\phi)\|_{\mathcal{B}^{2\lambda}_{\beta}}}{\lambda^{\beta}} \lesssim \eps^{2-\beta}
\;.$$


%
%

\end{proof}

\subsection{Homogenisation of singular SPDEs}\label{sec:Homogenisation in RS}
In this section we shall rewrite singular SPDEs, using the kernels \eqref{eq:compactly supported kernels}, in such a way that one can apply the results of the previous sections. 
We consider the equation 
$$\partial_t u_\eps -  \nabla\cdot  {A}(x/\eps,t/\eps^2) \nabla u_\eps 
= F(x/\eps , t/\eps^2, u_\eps , 
\nabla u_\eps , \xi ) , \qquad u_\eps(\cdot, 0) = u_{\texttt{in}}\in C^\eta(\mathbb{T}^d) $$
pulled back to $\mathbb{R}^{d+1}$ as the integral equation
\begin{equ}\label{eq:loc_main trick}
u_\eps= K_{\epsilon} \mathbf{R}_+ F(x/\eps , t/\eps^2, u_\eps , 
\nabla u_\eps , \xi  ) + S^{K_\eps}_0 u_{\texttt{in}} , 
\end{equ}
where $S^{K_\eps}_s u_{\texttt{in}}(x,t) \eqdef \int K_{\eps}(x, t, y,s ) u_{\texttt{in}}(y) dy$
(with the usual abuse of notation if $u_{\texttt{in}}$ is a distribution).
Recall the correctors $\phi_i^\eps, \tilde{\phi}_j^\eps$ for $i,j=1,...,d$ and the kernel  $G_\eps$ in \eqref{eq:correctors} resp.~\eqref{eq:def_G}, and 
set for convenience 
 $\tilde\phi_0=1$.
Then, using periodicity of the correctors we rewrite \eqref{eq:loc_main trick}
\begin{equs}\label{eq:lift}
u_\eps
&=  \bar{K}_\eps\big( \mathbf{R}_+ F(x/\eps , t/\eps^2, u_\eps , 
\nabla u_\eps , \xi  )\big) 
+  
G_\epsilon\big( \mathbf{R}_+ F(x/\eps , t/\eps^2, u_\eps , 
\nabla u_\eps , \xi  )\big) \\
&\quad +
\sum_{i,j=1 }^d \eps \phi_i^\eps \cdot
\bar{K}^{i,j}_{\epsilon}
\big(\tilde{\phi}_j^\eps \mathbf{R}_+ F(x/\eps , t/\eps^2, u_\eps , 
\nabla u_\eps , \xi  ) \big)\\
&\quad +
\sum_{j=1 }^d 
\bar{K}^{0,j}_{\epsilon}
\big(\tilde{\phi}_j^\eps \mathbf{R}_+ F(x/\eps , t/\eps^2, u_\eps , 
\nabla u_\eps , \xi  ) \big)
+ S^{\eps}_0 u_{\texttt{in}} \ .
\end{equs}

We assume that the nonlinearity is of the form
$$F_\eps(x/\eps,t/\eps^2, u, Du)= \sum_{\alpha} f_{\alpha}(x/\eps, t/\eps^2) F^{\alpha} (u, D u)\ ,$$
where $\alpha$ runs over some finite index set and the $F^{\alpha}$ are sufficiently regular functions.
We lift $G_\epsilon$, $\bar K_\eps$ and $\bar{K}^{i,j}_{\epsilon}$  to abstract kernels
 $\mathcal{G}_\epsilon$, $\mathcal{\bar K}_\eps$ and $\bar{\mcK}^{i,j}_{\epsilon}$
 , the correctors $\phi_i, \tilde{\phi}_j$ and the functions $f_{\alpha,\beta}$ to abstract noises $\Phi_i, \tilde{\Phi}_j$, resp $\mathbf{f}_{\alpha}\in T$ and multiplication by $\eps$ to abstract multiplication with a scale. Thus, consider the abstract equation
\begin{equs}\label{eq:abstract}
U_\eps&=  
 \mathcal{\bar K}_\eps \big( \mathbf{R}_+ \hat{\opF} (U  , \Xi  ) \big) 
+ 
\sum_{{i,j=1} }^d \hat{\mathcal{E}}^\eps\Big(\Phi_i \cdot
\bar{\mcK}^{i,j}_{\epsilon} 
\big(\tilde{\Phi}_j \cdot \mathbf{R}_+ \hat{\opF}(U  , \Xi  )\big) \Big)\\
&\qquad +\sum_{j=1}^d 
\bar{\mcK}^{0,j}_{\epsilon} 
\big(\tilde{\Phi}_j \cdot \mathbf{R}_+ \hat{\opF}(U  , \Xi  )\big) + 
\mathcal{G}_\epsilon \big( \mathbf{R}_+ \hat{\opF}(U  , \Xi  ) \big) + v^\eps_{\texttt{in}} \ ,
\end{equs}
where $v^\eps_{\texttt{in}}$ is a lift of $P^{\eps} u_{\texttt{in}}$ to  singular modelled distribution, see Section~\ref{sec:convergence initial conditions} below, and 
$\hat{\opF}=  \sum_{\alpha} \mathbf{f}_{\alpha} \hat{F}^{\alpha}(U)
$
with $\hat{F}^{\alpha}$ as in Remark~\ref{rem:non_linearity and derivation}.
\begin{remark}\label{rem:polynomial-nonlinearities}
In the case when the nonlinearity $F_\eps(U, DU) $ is polynomial in the solution, i.e.\ 
$F_\eps(u, Du)= \sum_{\alpha,\beta} f_{\alpha,\beta}(x/\eps, t/\eps^2) u^{\alpha} (\partial_i u)^{\beta_i}$ one can in principle further simplify the abstract equation since one does not need the requirement that the sector on which the nonlinearity acts is function like. Thus one can simply lift \eqref{eq:loc_main trick} as 
$$
U_\eps=  
\sum_{\substack{i,j=0} }^d \Phi_i \cdot
\mcK^{i,j,+}_{\epsilon} 
\big(\tilde{\Phi}_j \cdot \mathbf{R}_+ \opF (U , \Xi  )\big)
+
\mathcal{G}_\epsilon \big( \mathbf{R}_+ \hat{\opF}(U  , \Xi  ) \big) + U_{\texttt{in}} \ ,
$$
%
where $\mcK^{i,j,+}_{\epsilon}$ denotes the abstract lift to a $\beta$ regularising kernels of $K^{i,j,+}_{\epsilon}$.

Let us though mention that this leads to complications when identifying the renormalised equation for the $\Phi^4_3$ equation.
\end{remark}

%

\begin{remark}\label{rem:expanding to different orders}
Note that in \eqref{eq:loc_main trick} we rewrote the integral equation using the two-sided first order two scale expansion of the heat kernel.
Of course one could also rewrite equations by either only expanding in one variable or using higher order expansions, 
which might be necessary depending on the equation under consideration.

The kernel estimates established in Section~\ref{sec_post_process} allows one to consider equations where one can work with $L,R\sim 1+\kappa$, which in particular excludes the KPZ type equations (where $R=1/2+\kappa,\  L=3/2 + \kappa$ for some $\kappa>0$ seems required) for which one would need to improve the estimates of Theorem~\ref{thm:further kernel estimate} by including second order expansions in the first variable. 
In general we expect that working with an $N,N'\in \mathbb{N}$ order expansion in the first resp. second variable allows to treat equations where $L\sim N+ \kappa ,R\sim N'+\kappa$ is needed. The choice of $L,R$ that is required for a specific equation follows from Theorem~\ref{thm:fixed point}.
\end{remark}

%
%

\subsubsection{Convergence of the initial condition}\label{sec:convergence initial conditions}

We first consider the case that $u_{\texttt{in}}\in C^\eta$ for $\eta\in(-1,0]$. For $\gamma\in (0,2)$, let
\begin{equ} \label{eq:initial}
v^\eps_{\texttt{in}}(z)\eqdef\opbfP_{z}^{\gamma} \left[S^{\tilde{G}^{\eps}}  u_{\texttt{in}}\right] + \sum_{i=0}^{d} \hat{\mathcal{E}}^\eps\Phi_i \cdot  \opbfP_{z}^{\gamma}\left[S^{\bar{K}^{i,0}_\eps}u_{\texttt{in}}\right] ,\qquad
v^0_{\texttt{in}}(z)\eqdef  \opbfP_{z}^{\gamma}\left[S^{\bar{K}}u_{\texttt{in}}\right]
\end{equ}
where 
$\tilde{G}^\eps(z;\bar{z}) \eqdef K_\eps(z;\bar{z}) - \eps\sum_{i=0}^{d} \phi_i^\eps(z) \bar{K}^{i,0}_\eps(z;\bar{z})$
and the Taylor lift $\opbfP^\gamma$ was defined in Remark~\ref{rem:poly_example}.

The following two lemmas allow to check that $v^\eps_{\texttt{in}}$ belongs to an appropriate spaces of modelled distributions $\mathcal{D}^{\gamma,\eta}$ and converges in that space to $v^0_{\texttt{in}}$ as $\eps\to 0$.
\begin{lemma}
Consider $-R <\eta <0< \gamma<L$ and $\beta\in (0,2)$ then for every $T>0$ 
$$\|S^{H}  f \|_{C^{\gamma, \eta-2+\beta}_T} \lesssim_T \|H\|_{\beta; L,R} \|f\|_{ C^{\eta}(\mathbb{R}^d) }$$ uniformly over 
$H\in \bfK^\beta_{ L, R}$ and $f\in C^\eta$.
\end{lemma}
\begin{proof}
Let $H= \sum_{n=0}^\infty H_n$ be as in Definition~\ref{def:kernel_norm}. 
For $|k|_\fraks< \gamma$, we see that for $n_t\in \mathbb{N}$ such that $\sqrt{t}\in [2^{-n_t-1}, 2^{n_t})$
\begin{equs}
|D^k S^{H} v (x,t)|&\leq \sum_{n=0}^{n_t} |v_0\left(D^k H_n (x,t;\ \cdot\ ,0)\right) |\lesssim \|H\|_{\beta; L,R} \|v\|_{ C^{\eta}(\mathbb{R}^d) } \sum_{n=0}^{n_t}  2^{-n(\eta -2+\beta-k) } \\
& \lesssim  \|H\|_{\beta; L,R} \|v\|_{ C^{\eta}(\mathbb{R}^d) }  \sqrt{t}^{(\eta -2+\beta -k) \wedge 0}\ . 
\end{equs}
Similarly, we observe that for $t<\tau$ and $(\zeta,\tau)\in \tilde{Q}_{\sqrt{t}}(t,x)\cap  (\mathbb{R}^d \times (0,T])$
\begin{equs}
|S^{H} v (x,t)- \opP^\gamma_{(\zeta, \tau)} [S^{H} v] (x,t) |
+ | \opP^\gamma_{(x,t)} [ S^{H} v] (\zeta, \tau)- S^{H} (\zeta,\tau) | 
 \lesssim_T \|H\|_{\beta; L,R} \|v\|_{ C^{\eta}(\mathbb{R}^d) }  \sqrt{t}^{(\eta -2+\beta -\gamma) \wedge 0}\ . 
\end{equs}
where in the last line we used that $\eta<\gamma$.
\end{proof}
\begin{lemma}
For any $L\in (1,2), R\in [0,1)$, $\beta\in (0,2)$  and $\zeta>0$
$$\|\tilde{G}_\eps\|_{L,R;\beta}\lesssim_\zeta  \eps^{(2-R-\beta)\wedge(3-L-\beta)- \zeta} $$
uniformly over $\eps\in (0,1]$.
\end{lemma}
\begin{proof}
Writing $\tilde{G}_\eps= \sum_{n \in \mathbb{N}} \tilde{G}_{n;\eps}$ the estimate on the terms with $2^{-n}< \eps/2$ follows exactly as in the proof of Proposition~\ref{prop:vanishing of G kernel}.

We turn to $2^{-n}\geq \eps/2$, here it follows by combining Corollary~\ref{cor:pointwise} and Theorem~\ref{thm:further kernel estimate} that for $\zeta<1$
$$ \|G_{\eps,n}\|_{2;L,R;n} \lesssim_{\zeta} \eps^{(1-R) \wedge(2-L) -\zeta } 2^{n}$$
which implies that 
$$\|G_{\eps,n}\|_{\beta;L,R;n} \lesssim \eps^{   (2-L)  \wedge (1-R)-\zeta } 2^{n(1 -(2-\beta))} \leq \eps^{(3-L-\beta)\wedge(2-R-\beta)- \zeta}.$$
\end{proof}

\begin{corollary}\label{cor:initial_cond}
For any $\zeta>0$, and $\gamma \in (0 ,L-\kappa)$, $\beta\in (1,2)$ and $\eta<\beta-2-\kappa$ 
$$\vertiii{v^\eps_{\texttt{in}}- v^0_{\texttt{in}}}_{C^{\gamma, \eta}} \lesssim  \eps^{(2-\beta)\wedge(3-L-\beta)- \zeta}  \|u_{\texttt{in}}\|_{ L^\infty(\mathbb{R}^d) }\ .$$
\end{corollary}
\begin{proof}
By the previous two lemmas
\begin{align*}
\vertiii{v^\eps_{\texttt{in}}- v^0_{\texttt{in}}}_{C^{\gamma, \eta }}
&\lesssim\big( \|\tilde{G}_\eps\|_{L,0;\beta}+ \|\bar{K}-\bar{K}_\eps\|_{L,0;\beta}\big) \|u_{\texttt{in}}\|_{C^{\eta+2-\beta}} 
+ \sum_{i=1}^{d} \vertiii{\opbfP_{z}^{\gamma}\left[S^{\bar{K}^{i,0}_\eps}u_{\texttt{in}}\right]}_{\mathcal{D}^{\gamma +\kappa, \eta+\kappa }}\\
&\lesssim\big( \|\tilde{G}_\eps\|_{L,0;\beta}+ \|\bar{K}-\bar{K}_\eps\|_{L,0;\beta}\big) \|u_{\texttt{in}}\|_{C^{\eta+2-\beta}} 
+ \sum_{i=1}^{d} \|\bar{K}^{i,0}_\eps\|_{L,0;\beta} \|u_{\texttt{in}}\|_{C^{\eta+\kappa+2-\beta }} \\
&\lesssim_\zeta \eps^{(2-\beta)\wedge(3-L-\beta)- \zeta} \|u_{\texttt{in}}\|_{L^\infty} \ .
\end{align*}
\end{proof}

\subsubsection{Convergence of prepared initial conditions}
For the parabolic Anderson model when the nonlinearity involves the derivative of the solution, we shall, as explained in Remark~\ref{rem:initial condition} 
$u_\eps(x,0)= \int_{\mathbb{T}^d} \Gamma_\eps (x,0,y,-\eps) f(y)dy$  as initial condition. That is, we lift the solution 
$v_\eps(x,t)= \int_{\mathbb{T}^d}  \Gamma_\eps (x,t,y,-\eps) f(y)dy$ to the linear equation \eqref{eq:initial condition} 
as 
\begin{equ} \label{eq:initial2}
v^\eps_{\texttt{in}}(z)\eqdef\opbfP_{z}^{\gamma} \left[S_{-\eps}^{\tilde{G}^{\eps}}  u_{\texttt{in}}\right] + \sum_{i=0}^{d}  \hat{\mathcal{E}}^\eps \Phi_i \cdot  \opbfP_{z}^{\gamma}\left[S_{-\eps}^{\bar{K}^{i,0}_\eps}u_{\texttt{in}}\right] ,\qquad
v^0_{\texttt{in}}(z)\eqdef  \opbfP_{z}^{\gamma}\left[S^{\bar{K}}u_{\texttt{in}}\right] \ .
\end{equ}

We observe that \cite[Lem.~7.5]{Hai14} is directly applicable to all terms except the former summand of \eqref{eq:initial2}. This term can be treated using the following lemma.
\begin{lemma}\label{lem:initial2}
Let $\gamma \in (0,3/2)$. Then, for any $\kappa\in (0,1)$ ${\eta}\in [0,(1-2\kappa)\wedge (3/2-\gamma-\kappa)  )$ 
$$ \| S_{-\eps}^{\tilde{G}^\eps} f \|_{{\gamma,\eta} } \lesssim \eps^{\kappa/2} \|f\|_{L^\infty}$$
uniformly over $\eps\in (0,1]$ and $f\in C^\eta(\mathbb{T}^d)$.
\end{lemma}
\begin{proof}
Note that for $ 0\leq t \leq 1$ 
$$ S_{-\eps}^{\tilde{G}^\eps} f(x,t)= \int_{\mathbb{R}^d} \big( \Gamma_\eps(x,t,y,-\eps)-\Gamma^{1,0}_\eps(x,t,y,-\eps)\big) f(y)dy\ .$$
Thus it follows by Theorem~\ref{thm:uniform} that 
$$|S_{-\eps}^{\tilde{G}^\eps} f(x,t)|\lesssim \frac{\eps}{\sqrt{t+\eps}} \lesssim \eps^{1/2}, \qquad |\nabla S_{-\eps}^{\tilde{G}^\eps} f(x,t)|\lesssim \frac{\eps^{1-\kappa/2}}{t+\eps}  \lesssim \frac{\eps^{\kappa/2}  }{t^{(1-\eta)/2}}  \ .$$
Similarly, we observe by Theorem~\ref{thm:further kernel estimate} that for $t<\tau$ and $(\zeta,\tau)\in \tilde{Q}_{\sqrt{t}}(t,x)\cap  (\mathbb{R}^d \times (0,T])$
\begin{equs}
\frac{|S^{H} v (x,t)- \opP^\gamma_{(\zeta, \tau)} [S^{H} v] (x,t) |}{|(x,t)-(\zeta,\tau)|^\gamma_\fraks }
+\frac{ | \opP^\gamma_{(x,t)} [ S^{H} v] (\zeta, \tau)- S^{H} (\zeta,\tau) |  }{|(x,t)-(\zeta,\tau)|^\gamma_\fraks }
 \lesssim \frac{ \eps^{ 2-\gamma-\kappa/2} }{ (t+\eps)^{\frac{1+\gamma}{2}} } \lesssim \frac{\eps^{\kappa/2}}{t^{(\gamma-\eta)/2}}
\end{equs}

%

\end{proof}


%


\section{Application to the g-PAM Equation}\label{sec:application gpam}
At this point we have developed all the necessary tools in order to establish our theorems about the g-PAM and $\Phi^4_3$ equation, following the strategy of the theory of regularity structures.
\begin{enumerate}
\item Constructing a regularity structure and appropriately renormalised models to solve the abstract fixed point problem associated to the equation.
\item Check the form of the renormalised equation.
\item Show convergence to a limiting model when the regularisation is removed and check continuity of 
the map $(\eps,\delta)\mapsto \hat{M}^{\eps,\delta}$
 at $\{\eps=0\}\subset\{(\eps,\delta)\in [0,1]^2\}$.
\end{enumerate}

%
%
%
%
\subsection{Abstract formulation of the equation}\label{sec:abstract form gpam}

We consider the equation
$$\left(\partial_t - \nabla \cdot A(x/\eps,t/\eps^2) \nabla \right) u =  \sum_{i,j=1}^2 f_{i,j}(x/\eps, t/\eps^2) \partial_i u \partial_j u +  \sigma(u) \xi \ .
$$
Recall, that if one were to lift the equation as in \cite{Hai14} this would read
$$U_\eps= 
 \mathcal{K}_\epsilon \big( \mathbf{R}_+ \hat{F}(U, \nabla U)\big) + U_{\texttt{in}}, $$
 where  $$\hat{F}(U, \nabla U)\eqdef \sum_{i,j=1}^2\pmb{f}_{i,j} \partial_i U \partial_j U +  \hat{\sigma}(U) \Xi  $$
and where the $\pmb{f}_{i,j} $ just denote additional noise symbols.
This can be solved using the regularity structure, roughly speaking built using the set of symbols (`trees')
\begin{equ}\label{eq:trees classical g-pam}
\mathfrak{T}_c= \{ \Xi, \   \pmb{f}_{i,j}( \partial_j I\Xi) (\partial_i I\Xi),\  \mathbf{X} \Xi,\  \Xi (I\Xi),\  \mathbf{1},\  I\Xi,\  \mathbf{X}_i \}
\end{equ}
with the homogeneity assignment determined by declaring the noises $\pmb{f}_{i,j}$ to be of homogeneity $-\kappa$ and $\Xi$ to be of homogeneity $-1-\kappa$ and by declaring $I$ to be $2$-regularising. Then the solution has the form
$U= u \mathbf{1} + \sigma(u) I\Xi + u_i \mathbf{X}_i$. This however does not allow to control the limit
$\eps \to 0$ due to the small-scale oscillations in the integral kernel $K_\eps$.

\paragraph{Construction of the RS:}
Instead, we shall lift the equation as in \eqref{eq:abstract}
\begin{equs}\label{eq:abstract_gpam}
U_\eps&=  
 \mathcal{\bar K}_\eps \big( \mathbf{R}_+ \hat{\opF} (U  , \Xi  ) \big) 
+ 
\sum_{\substack{j=1} }^2
\bar{\mcK}^{0,j}_{\epsilon} 
\big(\tilde{\Phi}_j \cdot \mathbf{R}_+ \hat{\opF}(U  , \Xi  )\big) 
+
\sum_{i,j=1 }^2 \hat{\mathcal{E}}^\eps\big(\Phi_i \cdot
\bar{\mcK}^{i,j}_{\epsilon} 
\big(\tilde{\Phi}_j \cdot \mathbf{R}_+ \hat{\opF}(U  , \Xi  )\big)\big) \\
&\qquad +
\mathcal{G}_\epsilon \big( \mathbf{R}_+ \hat{\opF}(U  , \Xi  ) \big) + v^\eps_{\texttt{in}} \ ,
\end{equs}
To construct the regularity structure on which this makes sense,
denote by $\mathfrak{T}$ the minimal normal set such that
\begin{itemize}
\item It contains the symbols $\Psi_i, \tilde{\Psi}_i$.
\item It contains all elements obtained by substituting in a tree $T\in \mathfrak{T}_c$ from \eqref{eq:trees classical g-pam} every occurrence of
$I$ by one of  the following $\{ I, \bar{I},\  I^{0,j}( \tilde{\Psi}_{j}\, \cdot ),\  \mathcal{E}( \Psi_i I^{i,0}( \cdot) ),\  \mathcal{E}( \Psi_i I^{i,j}( \tilde{\Psi}_{j}\, \cdot) )  : i,j>0 \}$.
\item Whenever $\tau_1 \cdot \tau_2\in \mathfrak{T}$ it holds that $\tau_1,\tau_2\in \mathfrak{T}$.
\item Whenever $\tau\in \mathfrak{T}$ is formally obtained by applying an abstract integration $\tau'$, then it holds that $\tau'\in \mathfrak{T}$. 
\end{itemize}
Then, the space $T$ is given by the span of $\mathfrak{T}$. The homogeneity assignment of each tree is determined by declaring the kernels 
$I, \bar{I}, I^{0, j}$ to be $\beta$ regularising and $I^{i,j}$ to be $\beta-1$ regularising for $i>0$ as well as declaring
 $\Xi$ to be of homogeneity $\bar{\zeta}<-1$ and $\pmb{f}_{i,j}, \Psi_i, \tilde{\Psi}_j$ to be of homogeneity $\bar{\zeta}+1$. The numerical value of these constant will be fixed in Section~\ref{sec:proof gpam}.
The structure group $\Gamma$ is then defined be the conditions in the definition of an abstract integration operator \cite[Def.~5.7]{Hai14} and 
and abstract multiplication by a scale in Definition~\ref{def:abstract scale}.

\paragraph{Construction of models}

For $\eps\in[0,1]$, $\delta>0$ and a continuous noise $\xi_{\eps,\delta}\in C$, we shall consider the canonical model $M^{\eps,\delta}=(\Pi^{\eps,\delta}, \Gamma^{\eps,\delta})$ which is characterised by
\begin{enumerate}
\item $\Pi^{\eps,\delta}_x \Xi= \xi_{\eps,\delta} \ , \qquad \Pi^{\eps,\delta}_{x}\Phi_i = \phi^{\eps}_i\ , \qquad  \Pi^{\eps,\delta}_{x}\tilde{\Phi}_j = \tilde{\phi}^{\eps}_j\ , \qquad \Pi^{\eps,\delta}_{x}\pmb{f}_{i,j}= f^\eps_{i,j}$.
\item It acts as the polynomial model on the polynomial sector $\bar{T}\subset T$.
\item The abstract integration map $I$ realises the kernel $G$, $\bar{I}$ realises $K$, $I^{i,j}$ realises $K^{i,j}$ and the model realises the map $\mathcal{E}$ at scale $\eps$.
\item\label{Item:multiplicative} It is multiplicative, i.e. $\Pi_x(\tau \tau')= \Pi_x \tau \Pi_x \tau'$.
\end{enumerate}
We shall construct renormalised models by modifying Item~\ref{Item:multiplicative} above for trees 
trees of negative homogeneity which are not a noise and not planted, i.e.\ trees in 
$\mathfrak{T}_-= \mathfrak{T}_-^{\<Xi2>}\cup \mathfrak{T}_-^{\<2>} \cup \mathfrak{T}_-^{\pmb{f}\<2>}$, where 
\begin{enumerate}
\item $\mathfrak{T}_{\<Xi2>}= \{\Xi I \Xi\}\cup\{ \Xi I^{0,j}\tilde{\Phi} \Xi \ : j\geq 0 \}\cup \{\Xi \mathcal{E}(\Phi_i I^{i,j}\tilde{\Phi}_j \Xi) \ : \ i>0,j\geq 0\}$, \
\item $ \mathfrak{T}_{\<2>,\mu,\nu}= \Big\{ \tau_1\tau_2 \ : \tau_i \in 
\{\partial_{\alpha_i} I \Xi\}\cup\{ \partial_{\alpha_i} I^{0,j}\tilde{\Phi} \Xi \ : j\geq 0 \}\cup \{\partial_{\alpha_i} \mathcal{E}(\Phi_i I^{i,j}\tilde{\Phi}_j \Xi ) \ : \ i>0,j\geq 0\} : \alpha_1=\mu, \alpha_2= \nu \Big\}, $
\item $\mathfrak{T}_{\pmb{f}\<2>{\mu,\nu}} = \Big\{ \pmb{f}_{\mu,\nu} \tau_{\mu,\nu} \ : \tau_{\mu,\nu}\in \mathfrak{T}_{\<2>{\mu,\nu}} 
\} \ .$
\end{enumerate}
That is, we define a renormalised model which agrees with the canonical model, except that 
 for $\tau \in \mathfrak{T}_{\<Xi2>}\cup \mathfrak{T}_{\<2>} $ and for 
$\pmb{f}_{\mu,\nu} \cdot \tau_{\mu,\nu}\in \mathfrak{T}_{\pmb{f}\<2>} $
$$\hat{\Pi}^{\eps,\delta}_x \tau = {\Pi}^{\eps,\delta}_x \tau - g_{\eps,\delta}^{ \tau}\;,
\qquad \hat{\Pi}^{\eps,\delta}_x \pmb{f}_{\mu,\nu}\tau_{\mu,\nu} =f^{\eps}_{\mu,\nu} \hat{\Pi}^{\eps,\delta}_x \tau_{\mu,\nu} - 
\gamma^{\pmb{f}_{\mu,\nu}\tau_{\mu,\nu}}_{\eps,\delta}
$$
where $g_{\eps,\delta}^{ \tau}$ are (fixed) bounded functions and $\gamma^{\mu,\nu}_{\eps,\delta}$ are constants for $\delta>0$.
We denote the resulting model by $\hat{M}^{\eps,\delta}= (\hat{\Pi}^{\eps,\delta}, \hat{\Gamma}^{\eps,\delta})$.

\subsection{The renormalised equation}\label{sec:renorm eq}

In order to identify the effect of renormalisation on the equation, let us define 
$$ 
g_{\eps,\delta}^{ \<Xi2> }= \sum_{\tau\in \mathfrak{T}_{\<Xi2>}}  g_{\eps,\delta}^{ \tau},
 \qquad g_{\eps,\delta}^{ \<2>{\mu, \nu} }= 
\sum_{\tau\in \mathfrak{T}_{\<2>{\mu,\nu}}} 
 g_{\eps,\delta}^{\<2>{\mu,\nu}}, \qquad
\gamma^{\mu,\nu}_{\eps,\delta}\eqdef \sum_{\tau\in \mathfrak{T}_{\<2>{\mu,\nu}}} 
\gamma^{\pmb{f}_{\mu,\nu}\tau_{\mu,\nu}}_{\eps,\delta},
\qquad
\gamma_{\eps,\delta}\eqdef \sum_{\mu,\nu=1}^2 
\gamma^{\mu,\nu}_{\eps,\delta}
$$
\begin{lemma}\label{lem:ren_equ}
Let $U^{\eps,\delta}$ be the solution to \eqref{eq:abstract_gpam} with respect to the model $\hat{M}^{\eps,\delta}$. Then 
$\hat{u}_{\eps,\delta}\eqdef\hat{\mathcal{R}}U^{\eps,\delta}$ satisfies 
\begin{align*}
\partial_t \hat{u}_{\eps,\delta} -\nabla \cdot A^\eps \nabla \hat{u}_{\eps,\delta} &= \sum_{\mu, \nu}
f^\eps_{\mu,\nu}\Big( \partial_\mu \hat{u}_{\eps,\delta} \partial_{\nu}\hat{u}_{\eps,\delta}  - \sigma^2({\hat{u}_{\eps,\delta}}) \cdot 
g_{\eps,\delta}^{ \<2>{\mu, \nu} }
\Big) 
- \sigma^2(\hat{u}_{\eps,\delta}) \gamma_{\eps,\delta}
\\
&\qquad 
+ \sigma(\hat{u}_{\eps,\delta} ) \Big(\xi_{\eps,\delta} - 
\sigma' (\hat{u}_{\eps,\delta}) \cdot
g_{\eps,\delta}^{ \<Xi2> }
\Big)\ .
\end{align*}

\end{lemma}

\begin{proof}
We shall omit $\eps\in [0,1]$ and $\delta\in (0,1]$ in the notation for the sake of readability. 
First note that the solution to \eqref{eq:abstract_gpam}  takes the form 
\begin{equ}\label{eq:abstract sol gpam}
U= u \mathbf{1} + \sigma(u)\left(   I \Xi +\bar{I} \Xi + \sum_{j=1}^2 I^{0,j}( \tilde{\Phi}_{j} \Xi )   +\sum_{i>0,j\geq 0} \mathcal{E}( \Phi_i I^{i,j}( \tilde{\Phi}_{j} \Xi ) ) \right) + \sum_{i=1}^2 u'_i \mathcal{E}(\Phi_i)+
\sum_{i=1}^2 u''_i \mathbf{X}_i\ ,
\end{equ}
for some (continuous) functions $u, u'_{i}, u''_j$ for $i,j=1,2$. 
Therefore, 
$$
\partial_{\nu} U = \sigma(u)\partial_\nu \Big(   I \Xi +\bar{I} \Xi +\sum_{j=1}^2 I^{0,j}( \tilde{\Phi}_{j} \Xi )  +  \sum_{i>0,j\geq 0} \mathcal{E}( \Phi_i I^{i,j}( \tilde{\Phi}_{j} \Xi ) ) \Big) + \sum_{i=1}^2 u'_i \partial_\nu\mathcal{E}(\Phi_i)+ u''_\nu \ ,
$$
which implies that
\begin{align*}
 &\hat{\Pi}_x\big[\pmb{f}_{\mu,\nu} \partial_{\mu} U \partial_{\nu} U \big](x)\\
  &= 
\hat{\Pi}_x\Big[ \pmb{f}_{\mu,\nu}
 \Big( \sigma(u)\partial_\mu \Big(   I \Xi +\bar{I} \Xi +  \sum_{i>0,j\geq 0} \mathcal{E}( \Phi_i I^{i,j}( \tilde{\Phi}_{j} \Xi ) ) \Big) + \sum_{i=1}^2 u'_i \partial_\mu\mathcal{E}(\Phi_i)+ u''_\mu \Big)\\
\qquad &
\qquad\times\Big( \sigma(u)\partial_\nu \Big(   I \Xi +\bar{I} \Xi +  \sum_{i>0,j\geq 0} \mathcal{E}( \Phi_i I^{i,j}( \tilde{\Phi}_{j} \Xi ) ) \Big) + \sum_{i=1}^2 u'_i \partial_\nu\mathcal{E}(\Phi_i)+ u''_\nu \Big)\Big](x) \\
&= f^{\eps}_{\mu,\nu} \cdot \Big(\hat{\Pi}_x\big[ \partial_{\mu} U](x)  \cdot \hat{\Pi}_x \big[\partial_{\nu} U\big](x) -\sigma^2 (u) \cdot 
g^{\<2>{\mu,\nu}}_{\eps,\delta} 
 \Big) - \sigma^2(u) \gamma^{\mu,\nu} .
\end{align*}
Similarly, 
$$
\hat{\sigma} (U) = \sigma(u) \mathbf{1} + \sigma'(u) \Big(\sigma(u)\Big(   I \Xi +\bar{I} \Xi + \sum_{j=1}^2 I^{0,j}( \tilde{\Phi}_{j} \Xi )  + \sum_{i>0,j\geq 0} \mathcal{E}( \Phi_i I^{i,j}( \tilde{\Phi}_{j} \Xi ) ) \Big) + \sum_{i=1}^2 u'_i \mathcal{E}(\Phi_i)+
\sum_{i=1}^2 u''_i \mathbf{X}_i \Big) \ .
$$
Thus
\begin{align*}
\hat{\Pi}_x \big[\Xi \hat{\sigma} (U^\eps)\big](x) &=  \sigma(u)\cdot \hat{\Pi}_x \Xi \\
&\qquad+ \sigma'(u) \sigma(u) \hat{\Pi}_x\Big[\Xi\Big(   I \Xi +\bar{I} \Xi + \sum_{j=1}^2 I^{0,j}( \tilde{\Phi}_{j} \Xi )  + \sum_{i>0,j\geq 0} \mathcal{E}( \Phi_i I^{i,j}( \tilde{\Phi}_{j} \Xi ) )\big](x)  \\
&\qquad+ \hat{\Pi}_x \Big[\Xi\Big(  \sum_{i=1}^2 u'_i \mathcal{E}(\Phi_i)+
\sum_{i=1}^2 u''_i \mathbf{X}_i \Big)\Big](x) \\
&=\hat{\Pi}_x \Xi \cdot  \hat{\Pi}_x \big[{\sigma} (U^{(\eps,\delta)})\big](x) - \sigma(u)\sigma'(u) \cdot \sum_{\tau \in \mathfrak{T}_{\<Xi2>} } g^{\tau}_{\eps,\delta}.
\end{align*}
Therefore, one concludes the proof by evaluating $\hat{\Pi}_x\big[\hat{F}(U)\big](x)$ at $x$. 
%
\end{proof}

%


\subsection{Convergence of the renormalised models} \label{eq: convergence gpam}
Next we shall establish convergence of the renormalised models $\hat{M}^{\eps,\delta}$.
As we shall see this splits into parts of the argument which are rather insensitive to what regularisation of the noise one works with, and other parts where we need to distinguish 
whether one works with 
\begin{equ}
\xi_{\delta}(z)=\xi(\rho_z^\delta) \qquad \text{or}\qquad
\xi_{\epsilon,\delta}(x,t)= \int_{\mathbb{R}^d} \Gamma_\eps(x,t, \zeta, t-\delta^2) \xi(\zeta) \,d\zeta\ .
\end{equ}


For $\tau \in\mathfrak{T}_{\<Xi2>}$
%
 set $g^{\tau}_{\eps,\delta}(z)= \E[\bfPi^{\eps,\delta}\tau(z)]- \mathbf{1}_{\tau=\Xi \bar{I}\Xi  } F_{\eps,\delta}(z)$, where the functions 
$\{F_{\eps,\delta}\}_{\eps\in (0,1], \delta\in [0,1]}$ will be chosen later, depending on the regularisation. 
(But importantly such that $F_{\eps,0}=\lim_{\delta\to 0}F_{\eps,\delta}$ is independent of the regularisation.)

\begin{lemma}\label{lem:conv dumbel1}
Assume that $F_{\eps,\delta}$ is uniformly bounded and $\eps\mathbb{Z}^2\times \eps^2\mathbb{Z}$ periodic with vanishing mean. Then, for $\kappa>0$ sufficiently small 
\begin{enumerate}
\item 
$\E\big[\big| \hat{\Pi}_\star^{(\eps,\delta)} \tau (\varphi^\lambda_{\star})\big|^2\big]\lesssim \lambda^{-\kappa}$ for all $\tau\in \mathfrak{T}_{\<Xi2>}$,
\item
$ \E\big[\big| \big(\hat{\Pi}_\star^{(\eps,\delta)} \Xi \bar{I} \Xi - \hat{\Pi}_\star^{(0,\delta)} \Xi \bar{I} \Xi \big)(\varphi^\lambda_{\star})\big|^2\big]\lesssim \eps^{\kappa} \lambda^{-2\kappa}$,
\item
$ \E\big[\big| \hat{\Pi}_\star^{(\eps,\delta)} \tau (\varphi^\lambda_{\star})\big|^2\big]\lesssim \eps^{\kappa} \lambda^{-2\kappa}$  for $\tau \in \mathfrak{T}_{\<Xi2>}\setminus \{\Xi \bar{I} \Xi\}$,
\end{enumerate}
uniformly over $\eps,\delta\geq 0$, $\varphi\in \mathfrak{B}_0$ and $\star\in \mathbb{R}^{3}$.
\end{lemma}
\begin{proof}
We shall use graphical notation close to the one used in \cite{HP15,HQ18}. For $\tau\in \{\Xi I \Xi, \Xi \bar{I} \Xi \}\cup\{ \Xi I^{0,j}\tilde{\Phi} \Xi \ : \ j=1,2 \}$ we can then represent
	\begin{equ}\label{eq:diagrams}
									\hat{ \Pi}_\star^{(\eps,\delta)} \tau (\varphi^\lambda_{\star}) = \;
						\begin{tikzpicture}[scale=0.35,baseline=0.3cm]
							\node at (0,-1)  [root] (root) {};
							\node at (-2,1)  [dot] (left) {};
							\node at (-2,3)  [bluevar] (left1) {};
							\node at (0,1) [var] (variable1) {};
							\node at (0,3) [var] (variable2) {};
				
							\draw[testfcn] (left) to  (root);
							
							\draw[kernel1] (left1) to (left);
							\draw[arrho] (variable2) to (left1); 
							\draw[arrho] (variable1) to (left); 
						\end{tikzpicture}\;
+
\begin{tikzpicture}[scale=0.35,baseline=0.3cm]
	\node at (0,-1)  [root] (root) {};
	\node at (-2,1)  [dot] (left) {};
	\node at (-2,3)  [bluevar] (left1) {};
	\node at (0,2) [dot] (variable1) {};
	\node at (0,2) [dot] (variable2) {};
	
	\draw[testfcn] (left) to (root);
	
	\draw[kernel1] (left1) to (left);
	\draw[arrho] (variable2) to (left1); 
	\draw[arrho] (variable1) to (left); 
	\end{tikzpicture}\;
	-\; g^{\tau}_{\eps,\delta}\big( \varphi^\lambda_{\star} \big)
 =
\begin{tikzpicture}[scale=0.35,baseline=0.3cm]
							\node at (0,-1)  [root] (root) {};
							\node at (-2,1)  [dot] (left) {};
							\node at (-2,3)  [bluevar] (left1) {};
							\node at (0,1) [var] (variable1) {};
							\node at (0,3) [var] (variable2) {};
				
							\draw[testfcn] (left) to  (root);
							
							\draw[kernel1] (left1) to (left);
							\draw[arrho] (variable2) to (left1); 
							\draw[arrho] (variable1) to (left); 
						\end{tikzpicture}\;
- \;
\begin{tikzpicture}[scale=0.35,baseline=0.3cm]
	\node at (0.0,-1)  [root] (root) {};
	\node at (-1,1)  [dot] (left) {};
	\node at (0,3)  [dot] (top) {};
	\node at (1,1) [bluevar] (right) {};
	
	\draw[testfcn] (left) to  (root);
	
	\draw[kernel] (right) to (root);
	\draw[arrho] (top) to (right); 
	\draw[arrho] (top) to (left); 
	\end{tikzpicture}\;	
	+
\mathbf{1}_{\tau=\Xi \bar{I}\Xi  } F_{\eps,\delta}\big( \varphi^\lambda_{\star} \big) \ ,
					\end{equ}
					where the graphical notation has the following interpretation.
					\begin{itemize}
						\item The node \tikz[baseline=-3] \node [root] {}; represents the point $\star\in \mathbb{R}^{2+1}$, 
						 while the edge \tikz[baseline=-0.1cm] \draw[testfcn] (1,0) to (0,0); represents integration against the rescaled test function $\varphi^{\lambda}_\star$.
						\item The nodes \tikz[baseline=-3] \node [var] {}; represent the kernel variables in the Wiener Chaos representation.
						\item The node \tikz[baseline=-3] \node [dot] {}; represents dummy variables which are to be integrated out.
						\item The node \tikz[baseline=-3] \node [bluevar] {}; represents a dummy variable for $\tau\in \{\Xi I \Xi, \Xi \bar{I} \Xi \}$.  If $\tau= \Xi I^{0,j}\tilde{\Phi} \Xi $ it represents $\tilde{\phi}_j^\eps$ evaluated at a dummy variable. In both cases the dummy variable is integrated out.						
						\item  Edges \tikz[baseline=-0.1cm] \draw[kernel] (0,0) to (1,0); represent integration against a kernel $K(x,t;y,s)$ where
$(s,y)$ and $(t,x)$ are the coordinates of the start and end
						points of the arrow respectively				
						and where $K=G_\eps$ in the case $\tau=\Xi I\Xi								$, $K= \bar{K}_\eps$ if $\tau=\Xi \bar{I}\Xi$ or $K=K_\eps^{0,j}$ if $\tau= \Xi I^{0,j}\tilde{\Phi} \Xi $.
						Similarly, a barred arrow 
						\tikz[baseline=-0.1cm] \draw[kernel1] (0,0) to (1,0);
						represents $K(x,t;y,s) - K(\star;y,s)$.
						\item Edges \tikz[baseline=-0.1cm] \draw[arrho] (0,0) to (1,0); represent the mollifier $\rho^\delta(x-y)$, respectively\footnote{Here in the latter case we commit some abuse of notation, and the kernel $\Gamma_\eps(x,t, y, t-\delta^2) $ should really be split into a compactly supported part on some small enough (equation dependent) ball and a remainder which is very easily treated, see for instance \cite{CS25} where such a step is performed in more detail in a different setting.} 
						$\Gamma_\eps(x,t, y, t-\delta^2) $. 
					\end{itemize}
The desired estimates on contribution to the second Wiener chaos thus follow by arguing directly as in \cite{Hai14}, see also \cite{HQ18, HP15} for notation closer to here, but additionally using that $\|\phi_j^{\eps}\|_{L^\infty}=\|\phi_j\|_{L^\infty}<\infty$ and that the kernels converge in $\bfK^{2-\kappa/2}_{L, 0}$. Similarly, one estimates the contributions of the second term.
Finally, the estimate $|F_{\eps,\delta}\big( \varphi^\lambda_{\star} \big)|\lesssim \eps^{\kappa}\lambda^{-\kappa}$ follows directly from  \cite[Lem.~A1]{Sin23}.
%
%

Similarly, using Lemma~\ref{lem:identities for abstract mult}, we find for $\tau\in \{\Xi \mathcal{E}(\Phi_i I^{i,0} \Xi) \ : \ i>0\}\cup \{\Xi \mathcal{E}(\Phi_i I^{i,j}\tilde{\Phi}_j \Xi) \ : \ i,j>0\}$ that
\begin{equation*}
				\hat{\Pi}_\star^{(\eps,\delta)} \tau (\varphi^\lambda_\star) = \;
						\eps \phi^{\eps}(\star) \Bigg[
						\begin{tikzpicture}[scale=0.35,baseline=0.3cm]
							\node at (0,-1)  [root] (root) {};
							\node at (-2,1)  [dot] (left) {};
							\node at (-2,3)  [bluevar] (left1) {};
							\node at (0,1) [var] (variable1) {};
							\node at (0,3) [var] (variable2) {};
				
							\draw[testfcn] (left) to  (root);
							
							\draw[kernel1] (left1) to (left);
							\draw[arrho] (variable2) to (left1); 
							\draw[arrho] (variable1) to (left); 
						\end{tikzpicture}\;
						- \;
\begin{tikzpicture}[scale=0.35,baseline=0.3cm]
	\node at (0.0,-1)  [root] (root) {};
	\node at (-1,1)  [dot] (left) {};
	\node at (0,3)  [dot] (top) {};
	\node at (1,1) [bluevar] (right) {};
	
	\draw[testfcn] (left) to  (root);
	
	\draw[kernel] (right) to (root);
	\draw[arrho] (top) to (right); 
	\draw[arrho] (top) to (left); 
	\end{tikzpicture}\;	
%
%
\Bigg]
+ 
						\eps\; 
						\begin{tikzpicture}[scale=0.35,baseline=0.3cm]
							\node at (0,-1)  [root] (root) {};
							\node at (-2,1)  [dot] (left) {};
							\node at (-2,3)  [bluevar] (left1) {};
							\node at (0,1) [var] (variable1) {};
							\node at (0,3) [var] (variable2) {};
				
							\draw[testfcn] (left) to  (root);
							
							\draw[kernel] (left1) to (left);
							\draw[arrho] (variable2) to (left1); 
							\draw[arrho] (variable1) to (left); 
							\draw[blue, kernel] (root) to[bend left=70] (left); 
						\end{tikzpicture}\;,
					\end{equation*}
					where we used similar graphical notation as above, but where
					\begin{itemize}
					\item the edge \tikz[baseline=-0.1cm] \draw[blue, kernel] (0,0) to (1,0); represents $\phi^\eps(\cdot) - \phi^{\eps}(\star)$, 
					\item the node \tikz[baseline=-3] \node [bluevar] {}; represents the constant function $1$ for $\tau\in \{\Xi \mathcal{E}(\Phi_i I^{i,0} \Xi) \ : \ i>0\}$ or $\tilde{\phi}_j^\eps$ for $\tau=\Xi \mathcal{E}(\Phi_i I^{i,j}\tilde{\Phi}_j \Xi) \ $, evaluated at a dummy variable which is integrated out.		
					\end{itemize}
					
We first estimate the contributions of the second Wiener Chaos. Here we observe that the former diagram can be treated exactly as above, by using the extra term $\eps$ as part of the kernel. For the latter diagram, note that in the case $\lambda\geq \eps$ this too can be treated as above using the bound $ |\phi^\eps(\cdot) - \phi^{\eps}(\star)|\leq 2\|\phi^\eps\|_{L^\infty}$. For the case $\lambda<\eps$, we estimate $|\phi^\eps(z) - \phi^{\eps}(\star)|\leq \eps^{-1} |z-\star|$ and use that in the relevant domain of integration $|z-\star|\lesssim \lambda$, to obtain
\begin{equ}\label{eq:local diagram bound}
\eps^2\; 
						\begin{tikzpicture}[scale=0.35,baseline=0.3cm]
							\node at (0,-1)  [root] (root) {};
							\node at (-2,1)  [dot] (left) {};
							\node at (-2,3)  [bluevar] (left1) {};
							\node at (0,1) [dot] (variable1) {};
							\node at (0,3) [dot] (variable2) {};
				
							\node at (2,1)  [dot] (right) {};
							\node at (2,3)  [bluevar] (right1) {};
							
							\draw[testfcn] (left) to  (root);
							\draw[kernel] (left1) to (left);
							\draw[arrho] (variable2) to (left1); 
							\draw[arrho] (variable1) to (left); 
							\draw[blue, kernel] (root) to[bend left=70] (left); 
							
							\draw[testfcn] (right) to  (root);
							\draw[kernel] (right1) to (right);
							\draw[arrho] (variable2) to (right1); 
							\draw[arrho] (variable1) to (right); 
							\draw[blue, kernel] (root) to[bend right=70] (right); 
						\end{tikzpicture}\;
						\lesssim \lambda^{2}
						\begin{tikzpicture}[scale=0.35,baseline=0.3cm]
							\node at (0,-1)  [root] (root) {};
							\node at (-2,1)  [dot] (left) {};
							\node at (-2,3)  [bluevar] (left1) {};
							\node at (0,1) [dot] (variable1) {};
							\node at (0,3) [dot] (variable2) {};
				
							\node at (2,1)  [dot] (right) {};
							\node at (2,3)  [bluevar] (right1) {};
							
							\draw[testfcn] (left) to  (root);
							\draw[kernel] (left1) to (left);
							\draw[arrho] (variable2) to (left1); 
							\draw[arrho] (variable1) to (left); 
							
							\draw[testfcn] (right) to  (root);
							\draw[kernel] (right1) to (right);
							\draw[arrho] (variable2) to (right1); 
							\draw[arrho] (variable1) to (right); 
								\end{tikzpicture}\; .
						\end{equ}			
In order to estimate the right hand side of \eqref{eq:local diagram bound} one then uses that the involved kernels represented by  \tikz[baseline=-0.1cm] \draw[kernel] (0,0) to (1,0);  vanish at scales smaller than $\eps$ (which is larger $\lambda$) and that thus there are no integrability problems at those scales.

		\end{proof}


Similarly to above we define for $\tau\in \mathfrak{T}_{\<2>,\mu,\nu}$
$$  
g^{\tau}_{\eps,\delta}(z)=\E[ \bfPi^{\eps,\delta} \tau ](z)-  \mathbf{1}_{\tau=\partial_\mu\bar{I}\Xi \partial_\nu\bar{I}\Xi   } F_{\eps,\delta}^{\mu,\nu} \ .
$$
\begin{lemma}\label{lem:convergene pam cherry}
Let $F^{\mu,\nu}_{\eps,\delta}$ be $\eps\mathbb{Z}^2\times \eps^2\mathbb{Z}$ periodic with vanishing mean and such that 
$\sup_{\eps,\delta} \|F^{\mu,\nu}_{\eps,\delta}\|_{L^p}<\infty $ for every $p<\infty$ .
Then, for $\kappa>0$ sufficiently small
\begin{enumerate}
\item 
$\E\big[\big| \hat{\Pi}_\star^{(\eps,\delta)} \tau (\varphi^\lambda_{\star})\big|^2\big]\lesssim \lambda^{-\kappa}$ for all $\tau\in \mathfrak{T}_{\<2>}$,
\item
$ \E\big[\big| \big(\hat{\Pi}_\star^{(\eps,\delta)} \partial_\mu \bar{I} \Xi\partial_\nu \bar{I} \Xi - \hat{\Pi}_\star^{(0,\delta)} \partial_\mu \bar{I} \Xi\partial_\nu \bar{I} \Xi\big)(\varphi^\lambda_{\star})\big|^2\big]\lesssim \eps^{\kappa} \lambda^{-2\kappa}$ for $\mu,\nu=1,2$,
\item
$ \E\big[\big| \hat{\Pi}_\star^{(\eps,\delta)} \tau (\varphi^\lambda_{\star})\big|^2\big]\lesssim \eps^{\kappa} \lambda^{-2\kappa}$  for $\tau \in \mathfrak{T}_{\<2>}\setminus \{\partial_\mu \bar{I} \Xi\partial_\nu \bar{I} \Xi \ : \mu,\nu=1,2\}$,
\end{enumerate}
uniformly over $\eps>0$, $\delta\geq 0$ and $\varphi\in \mathfrak{B}$, $\star\in \mathbb{R}^{3}$.
\end{lemma}
\begin{proof}
For $\tau_1 \in \{ \partial_{\mu} I \Xi, \partial_{\mu} \bar{I} \Xi,
 \partial_{\mu} I^{0,j}\tilde{\Psi}_j \Xi
 \}$ 
and $\tau_2 \in \{ \partial_{\nu} I \Xi, \partial_{\nu} \bar{I} \Xi,
 \partial_{\nu} I^{0,i}\tilde{\Psi}_i \Xi
 \}$ 
 one finds using similar graphical notation as in the proof of Lemma~\ref{lem:conv dumbel1}
\begin{equ}\label{eq:first cherry term}
			\bigl( \hat{\Pi}_\star^{(\eps,\delta)} \tau_1\tau_2 \bigr)(\varphi^\lambda) = \;
\begin{tikzpicture}[scale=0.35,baseline=0.3cm]
	\node at (0,-1.2)  [root] (root) {};
	\node at (0,0.8)  [dot] (int) {};
	\node at (-1,2.5)  [bluevar] (left) {};
	\node at (1,2.5)  [bluevar] (right) {};
	\node at (-1,3.9)  [var] (topleft) {};
	\node at (1,3.9)  [var] (topright) {};
	\draw[testfcn] (int) to  (root);
	\draw[arrho] (topleft) to (left); 
	\draw[arrho] (topright) to (right); 
	\draw[keps] (left) to (int);
	\draw[keps] (right) to (int);
\end{tikzpicture}\;		+
	 \;
	 \mathbf{1}_{\tau_1=\partial_\mu\bar{I}\Xi   } 
	 \mathbf{1}_{\tau_2=\partial_\nu\bar{I}\Xi   } F_{\eps,\delta}^{\mu,\nu}\big(\varphi^\lambda \big)\;,
		\end{equ}
where 
\begin{itemize}
\item  edges \tikz[baseline=-0.1cm] \draw[keps] (0,0) to (1,0); represent integration against a kernel $K$, where $K\in \{\partial_{\mu} G_\eps, \  \bar{K}^{0,j}_\eps : j=0,1,2\}$,
\item the nodes depicted by \tikz[baseline=-3] \node [bluevar] {};  represent ${\phi}_j^\eps$.
\end{itemize}
Similarly, we find that for $\tau_1\in \{ \partial_{\mu} \mathcal{E} \Psi_i I^{i,j}\tilde{\Psi}_j \Xi \ : \ i>0, j\geq 0\}$, $\tau_2\in \{ \partial_{\mu} I \Xi, \partial_{\mu} \bar{I} \Xi,
 \partial_{\mu} I^{0,j}\tilde{\Psi}_j \Xi
 \ :\   j=1,2 \}$ and 
$$
	\bigl( \hat{\Pi}_\star^{(\eps,\delta)} \tau_1\tau_2 \bigr)(\varphi^\lambda) = \;
\eps \;
\begin{tikzpicture}[scale=0.35,baseline=0.3cm]
	\node at (0,-1.2)  [root] (root) {};
	\node at (0,0.8)  [bluevar] (int) {};
	\node at (-1,2.5)  [bluevar] (left) {};
	\node at (1,2.5)  [bluevar] (right) {};
	\node at (-1,3.9)  [var] (topleft) {};
	\node at (1,3.9)  [var] (topright) {};
	\draw[testfcn] (int) to  (root);
	\draw[arrho] (topleft) to (left); 
	\draw[arrho] (topright) to (right); 
	\draw[keps] (left) to (int);
	\draw[keps] (right) to (int);
\end{tikzpicture}\;		
+ \;
\begin{tikzpicture}[scale=0.35,baseline=0.3cm]
	\node at (0,-1.2)  [root] (root) {};
	\node at (0,0.8)  [bluevar] (int) {};
	\node at (-1,2.5)  [bluevar] (left) {};
	\node at (1,2.5)  [bluevar] (right) {};
	\node at (-1,3.9)  [var] (topleft) {};
	\node at (1,3.9)  [var] (topright) {};
	\draw[testfcn] (int) to  (root);
	\draw[arrho] (topleft) to (left); 
	\draw[arrho] (topright) to (right); 
	\draw[kernel] (left) to (int);
	\draw[keps] (right) to (int);
\end{tikzpicture}\;			,
$$
 where 
 \begin{itemize}
\item  the edge \tikz[baseline=-0.1cm] \draw[kernel] (0,0) to (1,0); represents integration against a kernel $K$, where $K\in \{ \bar{K}^{i,j}_\eps : i=1,2,\ j=0,1,2\}$
\item each \tikz[baseline=-3] \node [bluevar] {};  represents a $\eps\mathbb{Z}^2\times \eps^2\mathbb{Z}$ periodic function (for example$\phi_j^\eps$, $\nabla\phi_j^\eps$, $\phi_i^\eps\phi_j^\eps$ or $\phi_i^\eps \partial_\mu\phi_j^\eps$).
 \end{itemize}
 Finally, similarly for $\tau_1,\tau_2\in \{ \partial_{\mu} \mathcal{E} \Psi_i I^{i,j}\tilde{\Psi}_j \Xi \ : \ i>0, j\geq 0\}$
$$
\big(\hat{ \Pi}_\star^{(\eps,\delta)} \tau_1\tau_2\big) (\varphi^\lambda) = \;
\eps^2 \;
\begin{tikzpicture}[scale=0.35,baseline=0.3cm]
	\node at (0,-1.2)  [root] (root) {};
	\node at (0,0.8)  [bluevar] (int) {};
	\node at (-1,2.5)  [bluevar] (left) {};
	\node at (1,2.5)  [bluevar] (right) {};
	\node at (-1,3.9)  [var] (topleft) {};
	\node at (1,3.9)  [var] (topright) {};
	\draw[testfcn] (int) to  (root);
	\draw[arrho] (topleft) to (left); 
	\draw[arrho] (topright) to (right); 
	\draw[keps] (left) to (int);
	\draw[keps] (right) to (int);
\end{tikzpicture}\;		
+2\eps \;
\begin{tikzpicture}[scale=0.35,baseline=0.3cm]
	\node at (0,-1.2)  [root] (root) {};
	\node at (0,0.8)  [bluevar] (int) {};
	\node at (-1,2.5)  [bluevar] (left) {};
	\node at (1,2.5)  [bluevar] (right) {};
	\node at (-1,3.9)  [var] (topleft) {};
	\node at (1,3.9)  [var] (topright) {};
	\draw[testfcn] (int) to  (root);
	\draw[arrho] (topleft) to (left); 
	\draw[arrho] (topright) to (right); 
	\draw[kernel] (left) to (int);
	\draw[keps] (right) to (int);
\end{tikzpicture}\;				
+\;
\begin{tikzpicture}[scale=0.35,baseline=0.3cm]
	\node at (0,-1.2)  [root] (root) {};
	\node at (0,0.8)  [bluevar] (int) {};
	\node at (-1,2.5)  [bluevar] (left) {};
	\node at (1,2.5)  [bluevar] (right) {};
	\node at (-1,3.9)  [var] (topleft) {};
	\node at (1,3.9)  [var] (topright) {};
	\draw[testfcn] (int) to  (root);
	\draw[arrho] (topleft) to (left); 
	\draw[arrho] (topright) to (right); 
	\draw[kernel] (left) to (int);
	\draw[kernel] (right) to (int);
\end{tikzpicture}\;		
\ .
$$

Observe that all terms with a blue dot vanish in the limit $\eps=0$ since the associated kernels vanish.
Thus the bounds on the terms in the second Wiener Chaos follow
follow again by arguing as in \cite{Hai14}. The contribution to \eqref{eq:first cherry term} in the zeroth Wiener Chaos is estimated using Lemma~\ref{lem:appendix}.
\end{proof}

We define for 
$\tau\in \mathfrak{T}_{\<2>{\mu,\nu}}$
$$\gamma_{\eps,\delta}^{\tau}=
\begin{cases}
\int_{[0,1]^3} f^{\eps}_{\mu,\nu} F^{\tau}_{\eps,\delta} dxdt  & \text{if } \tau \in \mathfrak{T}_{\pmb{f}\<2>{\mu,\nu}} , \\
0 & \text{else. }
\end{cases}
$$
%

\begin{lemma}\label{lem:conv_resonances}
For $\kappa>0$ sufficiently small
\begin{enumerate}
\item \label{item:resonances1}
$\E\big[\big| \hat{\Pi}_\star^{(\eps,\delta)} \tau (\varphi^\lambda_{\star})\big|^2\big]\lesssim \lambda^{-\kappa}$ for all $\tau\in \mathfrak{T}_{\pmb{f}\<2>}$
\item\label{item:resonances2}
$ \E\big[\big| \big(\hat{\Pi}_\star^{(\eps,\delta)} \pmb{f}_{\mu,\nu}\partial_\mu \bar{I} \Xi\partial_\nu \bar{I} \Xi - \hat{\Pi}_\star^{(0,\delta)} \pmb{f}_{\mu,\nu}\partial_\mu \bar{I} \Xi\partial_\nu \bar{I} \Xi\big)(\varphi^\lambda_{\star})\big|^2\big]\lesssim \eps^{\kappa} \lambda^{-\kappa}$ for $\mu,\nu=1,2$.
\item\label{item:resonances3}
$ \E\big[\big| \hat{\Pi}_\star^{(\eps,\delta)} \tau (\varphi^\lambda_{\star})\big|^2\big]\lesssim \eps^{\kappa} \lambda^{-\kappa}$  for $\tau \in \mathfrak{T}_{\pmb{f}\<2>}\setminus \{\partial_\mu \bar{I} \Xi\partial_\nu \bar{I} \Xi \ : \mu,\nu=1,2\}$.
\end{enumerate}
uniformly over $\eps>0$, $\delta\geq 0$ and $\varphi\in \mathfrak{B}$, $\star\in \mathbb{R}^{3}$.
\end{lemma}
\begin{proof}
We first argue how to bound the contributions to the second Wiener Chaos.
For Item~\ref{item:resonances1} and Item~\ref{item:resonances3} this follows from Lemma~\ref{lem:convergene pam cherry} by simply interpreting the oscillatory functions as part of the test function (which is possible since we allow for $\varphi\in \mathfrak{B}_0$). The estimate on Item~\ref{item:resonances2} is slightly more tedious but follows exactly as the bound on the second term of \cite[Eq.~(4.22)]{HS23per} in the proof of Prop.~4.13 therein. 
It remains to consider the contributions to the zeroth chaos, which are of the form $f^{\eps}_{\mu,\nu} F^{\mu,\nu}_{\eps,\delta} - \gamma^{\mu,\nu}$ and  are  bounded by Lemma~\ref{lem:appendix}.

%
\end{proof}

\begin{remark}
Observe that the result actually holds for $f\in L^{p}$ for $p$ large enough depending on $\kappa$.
\end{remark}

\subsection{Identification of divergences and proof of the main results}\label{sec:identification}
In this section we identify appropriate (regularisation dependent) choices for the functions $F_{\eps,\delta}$, $F_{\eps,\delta}^{\mu,\nu}$. 
For this we recall the non-centred analogue of a model, c.f.\  \cite[Sec.~15.5]{frizhairerbook}, on some trees, 
$$\bfPi^{(\eps,\delta)} \Xi I \Xi= \xi_\delta G_\eps(\xi_\delta) \ , \quad 
\bfPi^{(\eps,\delta)}\Xi I^{0,j}\tilde{\Phi}_j \Xi  =
\xi_\delta K_\eps^{0,j} ( \phi^\eps_j \xi_\delta) \ , \quad
\bfPi^{(\eps,\delta)} \Xi \mathcal{E}(\Phi_i I^{i,j}\tilde{\Phi}_j \Xi)= \eps\phi_i^\eps \xi_\delta K_\eps^{0,j} (\phi_j^\eps \xi_\delta) $$
and
$
\bfPi^{(\eps,\delta)}\tau= \Pi^{(\eps,\delta)}_0 \tau
$
for $\tau \in \mathfrak{T}_{\<2>}$.

\subsubsection{Counterterms for Theorem~\ref{thm:g-pam}}
Here we work with the heat kernel regularisation 
\eqref{eq:heat kernel regularisation spatial white noise} of the noise given by 
$$
\xi_{\eps,\delta}(x,t)= \int_{\mathbb{R}^d} \Gamma_\eps(x,t, \zeta, t-\delta^2) \xi(\zeta) \,d\zeta\ .
$$
We define $$F_{\eps,\delta}(z)= \sum_{\tau \in \mathfrak{T}_{\<Xi2>} } \E[\bfPi^{\eps,\delta}\tau](z) - \frac{  \alpha^{\<Xi2>}_{\eps,\delta} }{\sqrt{\det(A_s^\eps(z))}}-   \frac{\bar{\alpha}^{\<Xi2>}_{\eps,\delta}}{\sqrt{\det(\bar{A})}}- c^{\<Xi2>}_{\eps,\delta} \ ,$$
where
\begin{enumerate}
\item ${\alpha}^{\<Xi2>}_{\eps,\delta}\eqdef \frac{1}{4\pi}\int_{\mathbb{R}} ({\tau +2(\delta/\eps)^2} )^{-d/2} \kappa(\eps^2\tau) \kappa(\tau) d\tau$,
\item  $\bar{\alpha}^{\<Xi2>}_{\eps,\delta}\eqdef \frac{1}{4\pi}\int_{\mathbb{R}} ({\tau +2(\delta/\eps)^2} )^{-d/2} \kappa(\eps^2\tau) \kappa_c(\tau) d\tau$,
\item $ c^{\<Xi2>}_{\eps,\delta}\eqdef\int_{[0,1]^3} h_{\eps,\delta}$ for $h_{\eps,\delta}\eqdef \sum_{\tau \in \mathfrak{T}_{\<Xi2>} } \E[\bfPi^{(\eps,\delta)}\tau] - \det(A_s^\eps)^{-1/2}  \alpha^{\<Xi2>}_{\eps,\delta} - \big(\det(\bar{A})\big)^{-1/2} \bar{\alpha}^{\<Xi2>}_{\eps,\delta}$.
\end{enumerate}
\begin{lemma}\label{lem:finite wiener chaos remainder1}
The functions $F_{\eps,\delta}: \mathbb{R}^3\to \mathbb{R}$ satisfy the following properties:
\begin{itemize}
\item $\sup_{\eps,\delta\in (0,1]} \|F_{\eps,\delta}\|_{L^\infty}<\infty$,
\item $F_{\eps,\delta}$ is $\eps\mathbb{Z}^2\times \eps^2 \mathbb{Z}$ periodic and has mean $0$,
\item  the limit $F_{\eps,0}=\lim_{\delta\to 0} F_{\eps,\delta}$ exist as a pointwise limit for any $\eps>0$,
\item $\lim_{\eps\to  0} \|F_{\eps,\delta}\|_{L^\infty}=0$  for each $\delta>0$. 
\end{itemize}
\end{lemma}
\begin{proof}
Observe that 
\begin{align*}
\sum_{\tau \in \mathfrak{T}_{\<Xi2>} } \E[\bfPi^{\eps,\delta}\tau(x,t)] 
&= \int \kappa(t-s) \kappa^\eps(t-s) \Gamma_{\eps}(x,t;y,s-\delta^2) \Gamma_{\eps}(x,t;y,t-\delta^2) dyds \\
&\qquad + \int \kappa(t-s) \kappa_c^\eps(t-s) \Gamma_{\eps}(x,t;y,s-\delta^2) \Gamma_{\eps}(x,t;y,t-\delta^2) dyds +{r^{\<Xi2>}_{\eps,\delta}(x,t)}, 
\end{align*}
where one easily sees that 
$$(0,1]^2 \to L^{\infty}(\mathbb{R}^{1+2}), \qquad  (\eps,\delta)\mapsto r^{\<Xi2>}_{\eps,\delta}= \sum_{k\in \mathbb{Z}^d\setminus\{0\} } \int \kappa(t-s) \Gamma_{\eps}(x,t;y,s-\delta^2) \Gamma_{\eps}(x,t;y+k,t-\delta^2) dyds$$
extends continuously to $(\eps,\delta)\in [0,1]^2$. 
Furthermore, note that
by a direct computation (see \cite[Sec.~3.1.1]{Sin23} for a very similar one)
\begin{align*}
\sqrt{\det(A_s^\eps(x,t))}^{-1}  \alpha^{\<Xi2>}_{\eps,\delta}  &= \int \kappa(t-s) \kappa^\eps(t-s) Z^{*}_{\eps;0}(x,t;y,s-\delta^2) Z^{*}_{\eps;0}(x,t;y,t-\delta^2) dyds\ , \\
\sqrt{\det(\bar{A})}^{-1} \bar{\alpha}^{\<Xi2>}_{\eps,\delta}&=  \int \kappa(t-s) \kappa^\eps(t-s) \bar{\Gamma}(x,t;y,s-\delta^2) \bar{\Gamma}(x,t;y,t-\delta^2) dyds \ .
\end{align*}
We shall estimate the function $h_{\eps,\delta}= h^1_{\eps,\delta}+ h^2_{\eps,\delta}$ where 
%
%
\begin{alignat*}{2}
&h^1_{\eps,\delta}(x,t)= \int \kappa(t-s) \kappa^\eps(t-s)  \Big[&&\Gamma_{\eps}(x,t;y,s-\delta^2) \Gamma_{\eps}(x,t;y,t-\delta^2) \\
&  &&-Z^{*}_{\eps;0}(x,t;y,s-\delta^2) Z^{*}_{\eps;0}(x,t;y,t-\delta^2) \Big]
dyds \ , \\
& h^2_{\eps,\delta}(x,t) = \int \kappa(t-s) \kappa_c^\eps(t-s)\Big[ &&\Gamma_{\eps}(x,t;y,s-\delta^2) \Gamma_{\eps}(x,t;y,t-\delta^2) \\
& &&-\bar{\Gamma}(x,t;y,s-\delta^2) \bar{\Gamma}(x,t;y,t-\delta^2) \Big]
dyds 
\end{alignat*}
term by term.
By a substitution and \eqref{eq:formula for Gamma_epsilon} note that
\begin{align*}
&h^1_{\eps,\delta}(\eps x,\eps^2 t)= \int \kappa(\eps^2(t-s)) \kappa(t-s) \Big[\Gamma_{1}( x, t;y,s-(\delta/\eps)^2) \Gamma_{1}(x,t;y,t-(\delta/\eps)^2)\\
&\qquad\qquad\qquad\qquad\qquad\qquad\qquad\qquad
-Z^{*}_{1;0}(x,t;y,s-(\delta/\eps)^2) Z^{*}_{1;0}(x,t;y,t-(\delta/\eps)^2) \Big]
dyds \ . 
%
\end{align*}
One thus reads off that $\lim_{\eps \to 0} h^1_{\eps,\delta}(\eps x,\eps^2 t)=0$ for $\delta>0$.
Writing
\begin{align*}
&h^1_{\eps,\delta}(\eps x,\eps^2 t)\\
&= \int \kappa(\eps^2(t-s)) \kappa(t-s) \Big[\Gamma_{1}( x, t;y,s-(\delta/\eps)^2) 
-Z^{*}_{1;0}(x,t;y,s-(\delta/\eps)^2) \Big]\Gamma_{1}(x,t;y,t-(\delta/\eps)^2) 
dyds \\
&\quad
+ \int \kappa(\eps^2(t-s)) \kappa(t-s) 
Z^{*}_{1;0}(x,t;y,s-(\delta/\eps)^2) \Big[\Gamma_{1}(x,t;y,t-(\delta/\eps)^2) - Z^{*}_{1;0}(x,t;y,t-(\delta/\eps)^2) \Big]
dyds 
\end{align*} 
 one reads off using the Gaussian heat kernel bound of Theorem~\ref{thm:uniform} and Corollary~\ref{cor:Z*} that 
for $\delta>0$
\begin{equ}\label{eq:h^1 lim delta to 0}
\lim_{\delta \to 0} h^1_{\eps,\delta}(\eps x,\eps^2 t)= \int \kappa(\eps^2(t-s)) \kappa(t-s) \Big[\Gamma_{1}( x, t;x,s) 
-Z^{*}_{1;0}(x,t;x,s) \Big] 
ds \ ,
\end{equ}
as well as that $\sup_{\eps,\delta}\|h^1_{\eps,\delta}\|_{L^\infty}<\infty$. Next write
$  h^2_{\eps,\delta}= h^{2,1}_{\eps,\delta}+h^{2,2}_{\eps,\delta}+ h^{2,3}_{\eps,\delta}$ where

\begin{align*}
 h^{2,1}_{\eps,\delta} (\eps x,\eps^2 t)
 &= \int \kappa(\eps^2(t-s)) \kappa_c(t-s)\Big[ \Gamma_{1}(x,t;y,s-(\delta/\eps)^2) 
-\bar{\Gamma}(x,t;y,s-(\delta/\eps)^2)\Big] \Gamma_{1}(x,t;y,t-(\delta/\eps)^2)
 \\
 h^{2,1}_{\eps,\delta} (\eps x,\eps^2 t)&=
\int \kappa(\eps^2(t-s)) \kappa_c(t-s) \Big[\bar{\Gamma}(x,t;y,s-(\delta/\eps)^2) - \bar{\Gamma}(x,t;x,s-(\delta/\eps)^2)\Big]
 \Gamma_{1}(x,t;y,t-(\delta/\eps)^2)\\
 h^{2,1}_{\eps,\delta} (\eps x,\eps^2 t)&=
 - \int \kappa(\eps^2(t-s)) \kappa_c(t-s) \Big[\bar{\Gamma}(x,t;y,s-(\delta/\eps)^2) - \bar{\Gamma}(x,t;x,s-(\delta/\eps)^2)\Big] \bar{\Gamma}(x,t;y,t-(\delta/\eps)^2) 
\end{align*}
where in all three cases the integral runs over  $y\in \mathbb{R}^2$ and $s\in \mathbb{R}$. 
 By simply using Gaussian heat kernel bounds of Theorem~\ref{thm:uniform} and Corollary~\ref{cor:pointwise} with $N=N'=0$, one reads off 
 that $\lim_{\eps \to 0} h^{2,1}_{\eps,\delta}(\eps x,\eps^2 t)=0$ for $\delta>0$,  
that for $\eps>0$ 
\begin{equ}\label{eq:h^2 lim delta to 0}
\lim_{\delta\to 0} h^{2,1}_{\eps,\delta}(\eps x,\eps^2 t) = \int \kappa(\eps^2(t-s)) \kappa_c(t-s)\Big[ \Gamma_{1}(x,t;y,x) 
-\bar{\Gamma}(x,t;x,s)  \Big]
 dyds \ ,
\end{equ}
and that of $\sup_{\eps,\delta}\|h^{2,1}_{\eps,\delta}\|_{L^\infty}<\infty$ is uniformly bounded.
To control $h^{2,2}_{\eps,\delta}$ note that
\begin{align*}
|h^{2,2}_{\eps,\delta}(\eps x,\eps^2 t)|&\lesssim 
\int \kappa(\eps^2(t-s)) \kappa_c(t-s) \Big[\bar{\Gamma}(x,t;y,s-(\delta/\eps)^2) - \bar{\Gamma}(x,t;x,s-(\delta/\eps)^2)\Big]
 \Gamma_{1}(x,t;y,t-(\delta/\eps)^2)dyds\\
 &\lesssim \int \kappa(\eps^2(t-s)) \kappa_c(t-s) \frac{1}{|t-s+(\delta/\eps)^2|^{3/2}} |x-y| \Gamma_{1}(x,t;y,t-(\delta/\eps)^2)dyds \\
 &\lesssim \delta/\eps \int \kappa(\eps^2(t-s)) \kappa_c(t-s) \frac{1}{|t-s+(\delta/\eps)^2|^{3/2}} ds \\
 &\lesssim \delta/\eps \frac{1}{\sqrt{1+ (\delta/\eps)^2}} \ .
\end{align*}
which implies that $\sup_{\eps,\delta}\|h^{2,2}_{\eps,\delta}\|_{L^\infty}<\infty$ and that  $ \lim_{\delta\to 0} \|h^{2,2}_{\eps,\delta}\|_{L^\infty}=0$
for $\eps>0$. Arguing exactly the same way, one finds that
$\sup_{\eps,\delta}\|h^{2,3}_{\eps,\delta}\|_{L^\infty}<\infty$ and that  $ \lim_{\delta\to 0}\|h^{2,3}_{\eps,\delta}\|_{L^\infty}=0$
for $\eps>0$.
Finally, note that for $\delta>0$ as $\eps\to 0$
\begin{equs}
&|h^{2,2}_{\eps,\delta}(\eps x,\eps^2 t)+h^{2,3}_{\eps,\delta}(\eps x,\eps^2 t)|\\
%
&
 \lesssim \int \kappa(\eps^2(t-s)) \kappa_c(t-s) 
 \frac{1}{|t-s+(\delta/\eps)^2|} \Big|\Gamma_{1}(x,t;y,t-(\delta/\eps)^2)- \bar{\Gamma}(x,t;y,t-(\delta/\eps)^2)  \Big| dyds\\
 &
 \lesssim \frac{\eps}{\delta} \int \kappa(\eps^2(t-s)) \kappa_c(t-s) 
 \frac{1}{|t-s+(\delta/\eps)^2|} ds \to 0 \ . \label{eq:h^3 lim delta to 0}
\end{equs}
\end{proof}
\begin{remark}\label{rem:non-comm1}
Let us observe that we have shown in the previous proof that by dominated convergence $\lim_{\eps \to 0} c_{\eps,\delta}^{\<Xi2>}= 0$ for $\delta>0$,
while combining \eqref{eq:h^1 lim delta to 0},\eqref{eq:h^2 lim delta to 0}\& \eqref{eq:h^3 lim delta to 0} one obtains that
\begin{equs}
\lim_{\eps \to 0}\lim_{\delta\to 0} c_{\eps,\delta}^{\<Xi2>}
&= \int_{[0,1]^3}\int_{\mathbb{R}^{3}}  \kappa(t-s) \Big[\Gamma_{1}( x, t;x,s) 
-Z^{*}_{1;0}(x,t;x,s) \Big] \\
&\qquad\qquad\qquad
+ \kappa_c(t-s)\Big[ \Gamma_{1}(x,t;y,x) 
-\bar{\Gamma}(x,t;x,s)  \Big]
 dyds dxdt \ .
\end{equs}

\end{remark}
Next, we define $$F^{\mu, \nu}_{\eps,\delta}(z)= \sum_{\tau \in \mathfrak{T}_{\<2>{\mu,\nu}} } \E[\bfPi^{\eps,\delta}\tau] -
	\frac{(A_s^\eps)^{-1}_{\mu,\nu}}{\sqrt{\det(A_s^\eps)}} 
\alpha^{\<2>}_{\eps,\delta} -
\sum_{i,j} (\mathbf{1}_{\mu=i}+\partial_{\mu} \phi^\eps_i)	 (\mathbf{1}_{\nu=j}+\partial_{\nu} \phi^\eps_j)	
\frac{\bar{A}^{-1}_{i,j}}{\sqrt{\det(\bar{A})}}
\bar{\alpha}^{\<2>}_{\eps,\delta}- c^{\<2>{\mu, \nu}}_{\eps,\delta} \ ,$$
where
\begin{enumerate}
\item  $\alpha_{\eps,\delta}^{\<2>} \eqdef \frac{1}{8\pi}\int_{\mathbb{R}}\int_{\mathbb{R}} {(2\delta^2 +\sigma+\tau)^{-d/2-1}} \ \kappa^\eps(\tau)\kappa(\tau) d\tau\ \kappa^\eps(\sigma)\kappa(\sigma) d\sigma$,
\item  $\bar{\alpha}_{\eps,\delta}^{\<2>}\eqdef   \frac{1}{8\pi}\int_{\mathbb{R}}\int_{\mathbb{R}} {(2\delta^2 +\sigma+\tau)^{-d/2-1}} \ \kappa_c^\eps(\tau)\kappa(\tau) d\tau\ \kappa_c^\eps(\sigma)\kappa(\sigma) d\sigma     $,
\item $ c^{{\<2>{\mu,\nu}}}_{\eps,\delta}\eqdef \int_{[0,1]^3} h_{\eps,\delta}^{\<2>{\mu,\nu}}$ for 
$$h_{\eps,\delta}^{\<2>{\mu,\nu}}\eqdef
\sum_{\tau \in \mathfrak{T}_{\<2>{\mu,\nu}} } \E[\bfPi^{\eps,\delta}\tau] -
	\frac{(A_s^\eps)^{-1}_{\mu,\nu}}{\sqrt{\det(A_s^\eps)}} 
\alpha^{\<2>}_{\eps,\delta} -
\sum_{i,j} (\mathbf{1}_{\mu=i}+\partial_{\mu} \phi^\eps_i)	 (\mathbf{1}_{\nu=j}+\partial_{\nu} \phi^\eps_j)	
\frac{\bar{A}^{-1}_{i,j}}{\sqrt{\det(\bar{A})}}
\bar{\alpha}^{\<2>}_{\eps,\delta}\ .$$
\end{enumerate}

\begin{lemma}\label{lem:remainder cherry gpam}
The functions $F^{\mu,\nu}_{\eps,\delta}: \mathbb{R}^3\to \mathbb{R}$ satisfy the following properties:
\begin{itemize}
\item $\sup_{\eps,\delta\in (0,1]} \|F^{\mu,\nu}_{\eps,\delta}\|_{L^p}<\infty$ for any $p<\infty$,
\item $F_{\eps,\delta}$ is $\eps\mathbb{Z}^2\times \eps^2 \mathbb{Z}$ periodic and has mean $0$,
\item the limit $F_{\eps,0}= \lim_{\delta\to 0} F_{\eps,\delta}$ exists as a limit in $L^p$ for each $\eps>0$ ,
\item  $\lim_{\eps\to 0} \|F_{\eps,\delta}\|_{L^\infty}=0$ for each $\delta>0$.
\end{itemize}
\end{lemma}

\begin{proof}
Collecting terms one finds that 
\begin{align*}
\sum_{\tau \in \mathfrak{T}_{\<2>{\mu,\nu}} } \E[\bfPi^{\eps,\delta} \tau(x,t) ]&= \int \kappa(t-s)\kappa (t-s') \partial_\mu \Gamma_\eps(x,t;y,s-\delta^2) \partial_{\nu} \Gamma_\eps(x,t;y,s'-\delta^2) 
		dy ds ds' \\
		&\qquad + r^{\<2>{\mu,\nu}}(x,t) \ ,
\end{align*}
where 
$$  (\eps,\delta)\mapsto r^{\<2>{\mu,\nu}}(x,t)\eqdef \sum_{k\in \mathbb{Z}^2\setminus \{0\}} \int \kappa(t-s)\kappa (t-s') \partial_\mu \Gamma_\eps(x,t;y,s-\delta^2) \partial_{\nu} \Gamma_\eps(x,t;y+k,s'-\delta^2) 
		dy ds ds'
	$$
is easily seen to extend continuously to a map $ [0,1]^2 \to L^{\infty}(\mathbb{R}^{1+2})$, 
 and that
by a direct computation 
 \begin{align*}
 &\frac{\big(A_s^\eps(x,t)\big)^{-1}_{i,j}}{\sqrt{\det\big(A_s^\eps(x,t)\big)}} 
\alpha^{\<2>}_{\eps,\delta} +
\sum_{i,j} \big(\mathbf{1}_{\mu=i}+\partial_{\mu} \phi^\eps_i(x,t)\big)	 \big(\mathbf{1}_{\nu=j}+\partial_{\nu} \phi^\eps_j(x,t)\big)	
\frac{\bar{A}^{-1}_{i,j}}{\sqrt{\det(\bar{A})}}
\bar{\alpha}^{\<2>{\mu, \nu}}_{\eps,\delta}\\
&=\int \kappa(t-s)\kappa (t-s') 
	 \Bigg[ 
\kappa^\eps(t-s)\kappa^\eps (t-s') 
{ Z^\eps_{0;\mu} (x,t;y,s-\delta^2)Z^\eps_{0;\nu}(x,t;y,s'-\delta^2)}
\\
		&\quad + \kappa^\eps_c(t-s) \kappa^\eps_c (t-s') \Big[ \sum_{i,j} (\mathbf{1}_{\mu=i}+\partial_{\mu} \phi^\eps_i)	 (\mathbf{1}_{\nu=j}+\partial_{\nu} \phi^\eps_j)	\partial_i \bar{\Gamma} (x,t;y,s-\delta^2)	\partial_j \bar{\Gamma}(x,t;y,s'-\delta^2) \Big]\Bigg]\, dy ds' ds \ .
 \end{align*}
  Therefore, we can write 
$ h_{\eps,\delta}^{\<2>{\mu,\nu}} = h_{\eps,\delta}^{\<2>{\mu,\nu}; 1} + h_{\eps,\delta}^{\<2>{\mu,\nu};2}+  h_{\eps,\delta}^{\<2>{\mu,\nu};3}  
  $
  where
\begin{align*}
& h_{\eps,\delta}^{\<2>{\mu,\nu}; 1}(x,t)= \int \kappa(t-s)\kappa (t-s')  
\kappa^\eps(t-s)\kappa^\eps (t-s')
\Big[
\partial_\mu \Gamma_\eps(x,t;y,s-\delta^2) \partial_{\nu} \Gamma_\eps(x,t;y,s'-\delta^2) 
\\
&\qquad \qquad \qquad \qquad \qquad \qquad \qquad \qquad \qquad 
- {Z^\eps_{0;\mu} (x,t;y,s-\delta^2)Z^\eps_{0;\nu}(x,t;y,s'-\delta^2)}\Big]
\, dy ds' ds \ , 
 \\
& h_{\eps,\delta}^{\<2>{\mu,\nu};2}(x,t)=\int \kappa(t-s)\kappa (t-s') \Big[\kappa^\eps(t-s)\kappa^\eps_c (t-s')  +
\kappa^\eps_c(t-s)\kappa^\eps (t-s')   \Big]\\
&\qquad \qquad \qquad \qquad \qquad \qquad\qquad \qquad \qquad \times  \partial_\mu \Gamma_\eps(x,t;y,s-\delta^2) \partial_{\nu} \Gamma_\eps(x,t;y,s'-\delta^2) 
	\, dy ds' ds \ , \\
& h_{\eps,\delta}^{\<2>{\mu,\nu};3}(x,t)= \int \kappa(t-s)\kappa (t-s') 
 \kappa^\eps_c(t-s) \kappa^\eps_c (t-s') \Big[\partial_\mu \Gamma_\eps(x,t;y,s-\delta^2) \partial_{\nu} \Gamma_\eps(x,t;y,s'-\delta^2) \\
		&\qquad \qquad\qquad- \sum_{i,j} (\mathbf{1}_{\mu=i}+\partial_{\mu} \phi^\eps_i)	 (\mathbf{1}_{\nu=j}+\partial_{\nu} \phi^\eps_j)	\partial_i \bar{\Gamma} (x,t;y,s-\delta^2)	\partial_j \bar{\Gamma}(x,t;y,s'-\delta^2) \Big]\, dy ds' ds \ .
\end{align*}  
%
By the scaling properties \eqref{eq:formula for Gamma_epsilon}\&\eqref{eq:scaling linearisation of derivative} and a substitution  $h_{\eps,\delta}^{\<2>{\mu,\nu};1}(\eps x,\eps^2 t) $ is equal to
\begin{align*}
&\int \kappa(\eps^2(t-s))\kappa  (\eps^2(t-s'))\kappa(t-s)\kappa (t-s')  \Big[\partial_\mu \Gamma_1(x,t;y,s+ (\delta/\eps)^2) \partial_{\nu} \Gamma_1(x,t;y,s'+ (\delta/\eps)^2)  
\\
 &\qquad\qquad \qquad\qquad \qquad\qquad 
-{Z^1_{0;\mu}(x,t;y,s- (\delta/\eps)^2) Z^1_{0;\nu}(x,t;y,s'- (\delta/\eps)^2)} \Big]
		dy ds ds' \ .
\end{align*}
To estimate this term one argues as in \cite[Sec.~3.1.1]{Sin23} using \cite[Lem.~2.8]{Sin23}
and finds that 
$\sup_{\eps,\delta}\|h_{\eps,\delta}^{\<2>{\mu,\nu};1}\|_{L^\infty}<\infty$, that 
$\lim_{\eps\to 0}\|h_{\eps,\delta}^{\<2>{\mu,\nu};1}\|_{L^\infty}=0$ for $\delta>0$  and that for  $\eps>0$
\begin{equs}
\lim_{\eps\to 0} \lim_{\delta\to 0} h_{\eps,\delta}^{\<2>{\mu,\nu};1}(\eps x,\eps^2 t)
&= 
\int \kappa(t-s)\kappa (t-s') \Big[\partial_\mu \Gamma_1(x,t;y,s) \partial_{\nu} \Gamma_1(x,t;y,s') \label{eq:c-const1}
 \\
&\qquad\qquad\qquad\qquad -{Z^1_{0;\mu}(x,t;y,s) Z^1_{0;\nu}(x,t;y,s')} \Big]
		dy ds ds' \ .
\end{equs}
Next, since the argument for the other term is the same, we only bound the first summand in the definition of $h_{\eps,\delta}^{\<2>{\mu,\nu};2}(x,t)$, which is bounded by a multiple of
\begin{align*}
	&
 \int 
\frac{\kappa  (t-s) \kappa (t-s') \kappa^\eps(t-s)\kappa^\eps_c(t-s)}{(t-s+\delta^2)^{(d+1)/2} (t-s'+\delta^2)^{(d+1)/2}}
\exp\Big(-\frac{\kappa |y|^2  }{t-s+\delta^2}   -\frac{\kappa |y|^2  }{t-s'+\delta^2} \Big)
		dy ds ds'	
		\\
	&\lesssim 
 \int 
\frac{\kappa(s)\kappa^\eps(s)}{(s+\delta^2)^{(d+1)/2}} 
\frac{\kappa  (s') \kappa^\eps_c (s') }{(s'+\delta^2)^{(d+1)/2}}
\Big(
{ \frac{ (s+\delta^2)(s'+\delta^2)}{s+s'+2\delta^2} }
 \Big)^{d/2}
		 ds ds'\\	
			 	&\lesssim 
 \int 
\frac{\kappa(s)\kappa^\eps(s)}{(s+\delta^2)^{3/4}}  ds
 \int 
\frac{\kappa  (s') \kappa^\eps_c (s') }{(s'+\delta^2)^{3/4}}
		 ds'
\end{align*}
where we used the elementary inequality $(s+\delta^2) + (s' +\delta^2)\geq (s+\delta^2)^{1/2} (s' +\delta^2)^{1/2}$.
Thus, 
\begin{equ}\label{eq:c-const2}
\sup_{\eps,\delta} \|h_{\eps,\delta}^{\<2>{\mu,\nu};2}\|_{L^\infty}<+\infty\ , \qquad \lim_{\eps\to 0}\sup_{\delta} \|h_{\eps,\delta}^{\<2>{\mu,\nu};2}\|_{L^\infty}.
\end{equ}
Finally, one easily checks pointwise convergence of $h_{\eps,\delta}^{\<2>{\mu,\nu};2}$ as $\delta \to 0$.
We slightly rewrite the remaining term
\begin{align*}
h_{\eps,\delta}^{\<2>{\mu,\nu};3}(x,t)&\eqdef \int \kappa(t-s)\kappa (t-s')  
 \kappa^\eps_c(t-s) \kappa^\eps_c (t-s') 
\text{Err}^\eps_{\mu,\nu}(x,t,y,s- \delta^2,y, s'- \delta^2)
\, dy dy' ds' ds  
\end{align*}
where 
\begin{align*}
\text{Err}^\eps_{\mu,\nu}(x,t;y,s;y',s') &\eqdef  \partial_\mu \Gamma_\eps(x,t;y,s) \partial_{\nu} \Gamma_\eps(x,t;y',s') \\
&-\sum_{i,j} \big(\mathbf{1}_{\mu=i}+(\partial_{\mu} \phi_i)^\eps\big)	 \big(\mathbf{1}_{\nu=j}+(\partial_{\nu} \phi_j)^\eps \big)	\partial_i \bar{\Gamma} (x,t;y,s)	\partial_j \bar{\Gamma}(x,t;y',s') \ .
\end{align*}
Thus	 for $\text{Err}^\eps (x,t,y,s)\eqdef \mathbf{1}_{|t-s|>\eps/2} \max_{\mu=1,2}  \Big|\partial_{\mu} \Gamma_\eps(x,t;y,s) 
		- \sum_{j}  \big(\mathbf{1}_{\mu=j}+(\partial_{\mu} \phi_j)^\eps \big) 
	\partial_j \bar{\Gamma}(x,t;y,s) \Big| \ 	$
\begin{align*}
\Big|\text{Err}^\eps_{\mu,\nu}(x,t;y,s;y',s')\Big| \,
		&\leq \big|
			\partial_\mu \Gamma_\eps(x,t;y,s) \big| \,
	\text{Err}^\eps(x,t,y',s') \\
	& \quad +
	\text{Err}^\eps(x,t,y,s) 
 \Big| \sum_{j}  \big(\mathbf{1}_{\nu=j}+(\partial_{\nu} \phi_j)^\eps \big)
	\partial_j \bar{\Gamma}(x,t;y',s')	\Big| \ .
\end{align*}
Noting that by Theorem~\ref{thm:uniform} both
$$  \int \kappa^\eps_c(t-s) \kappa(t-s) \big| \partial_\mu \Gamma_\eps(x,t;y,s-\delta^2)\big| ds, \qquad 
\int \kappa^\eps_c(t-s) \kappa(t-s) 
 \Big| \sum_{j}  \big(\mathbf{1}_{\nu=j}+(\partial_{\nu} \phi_j)^\eps \big)
	\partial_j \bar{\Gamma}(x,t;y,s-\delta^2)	\Big| ds $$
	are bounded (up to a uniform constant) by
 $
   1+ \sum_{\eps\leq 2^{-n}< 1} 2^n{\mathbf{1}_{|x-y|<2^{-n}}}$ we find that
$$
|h_{\eps,\delta}^{\<2>{\mu,\nu};3}(x,t)|\lesssim \sum_{n: \eps\leq 2^{-n}\leq 1}  \underbrace{2^n \int \kappa^\eps_c(t-s) \kappa(t-s) \big(\mathbf{1}_{|x-y|<2^{-n}} +\mathbf{1}_{n=0} \big)
| \text{Err}^\eps(x,t,y,s'-\delta^2) | dy ds' }_{=:r^n_{\eps,\delta}(x,t)} \ .$$
We bound 
\begin{align*}
&
\|r^n_{\eps,\delta}\|_{\avL^p(Q_{2^{-n}}(x,t))} \\ 
&\leq 2^n \int \Big(\fint_{Q_{2^{-n}}(x,t)} \kappa^\eps_c(t'-s) \kappa(t'-s)  
\big(\mathbf{1}_{|x'-y|<2^{-n}} +\mathbf{1}_{n=0} \big)
 | 
\text{Err}^\eps(x',t',y,s-\delta^2) |^p dx'dt'\Big)^{1/p} dy ds\\
&\leq 2^n \int \big(\mathbf{1}_{|x-y|<2^{-n+1}} +\mathbf{1}_{n=0} \big) \Big(\fint_{Q_{2^{-n}}(x,t)} \mathbf{1}_{[\eps^2,2]} (t'-s)  |  \text{Err}^\eps(x',t',y,s-\delta^2) |^p dx'dt'\Big)^{1/p} dy ds \ .
\end{align*}
Thus one finds using Lemma~\ref{lem:lp kernel bound} that for $n\geq 1$ this is bounded by a multiple of
\begin{equ}\label{eq:Lp-type}
\eps 2^{2n} \int_{[0,1]^3} \frac{\mathbf{1}_{|x-y|<2^{-n+1}}}{|t-s+\delta^2|^{d/2+1/2}} \exp \big( -\kappa\frac{|x-y|^2}{|t-s+\delta^2|} \big)dy ds
 \lesssim \eps 2^{n} \ .
\end{equ}
Similarly one finds that 
\begin{equ}\label{eq:Lp-type2}
\|r^1_{\eps,\delta}\|_{\avL^p(Q_{1}(x,t))} \leq \int \Big(\fint_{Q_{1}(x,t)} \mathbf{1}_{[\eps^2,2]} (t'-s)  |  \text{Err}^\eps(x',t',y,s-\delta^2) |^p dx'dt'\Big)^{1/p} dy ds \lesssim 1
\end{equ}

thus
$$
\|r^n_{\eps,\delta}\|_{L^{p}([0,1]^3)}^{p} \leq  \sum_{z\in \eps\mathbb{Z}^{2}\times \eps^2 \mathbb{Z}, |z|_{\fraks}\lesssim 1 } \|r^n_{\eps,\delta}\|_{L^p(Q_{2^{-n}}(z))}^p \lesssim \eps^p 2^{np}$$
which implies that 
$\|h_{\eps,\delta}^{\<2>{\mu,\nu};2}\|_{L^{p}([0,1]^3)} \lesssim \sum_{n: \eps\leq 2^{-n}\leq 1} \|r^n_{\eps,\delta}\|_{L^{p}([0,1]^3)}  \lesssim 1$.
Next we observe that one can bound the integral in \eqref{eq:Lp-type}\&\eqref{eq:Lp-type2} for $\delta>0$ also by a $\delta$-dependent constant, which implies that
$\|r^n_{\eps,\delta}\|_{\avL^p(Q_{2^{-n}}(x,t)} \lesssim_{\delta} \eps$ and thus
$\lim_{\eps\to 0}\|h_{\eps,\delta}^{\<2>{\mu,\nu};2}\|_{L^{p}[0,1]}=0$ for $\delta>0$.
Lastly, one observes that
\begin{equs}
\lim_{\delta\to 0}h_{\eps,\delta}^{\<2>{\mu,\nu};2}(\eps x, \eps^2 t)
&=\int \kappa(\eps(t-s))\kappa(\eps(t-s')) \kappa_c(t-s) \kappa_c (t-s') \Big[\partial_\mu \Gamma_1(x,t;y,s) \partial_{\nu} \Gamma_1(x,t;y,s') \\
		&\qquad - \sum_{i,j} (\mathbf{1}_{\mu=i}+\partial_{\mu} \phi_i)	 (\mathbf{1}_{\nu=j}+\partial_{\nu} \phi_j)	\partial_i \bar{\Gamma} (x,t;y,s)	\partial_j \bar{\Gamma}(x,t;y,s') \Big]\, dy ds' ds 
\end{equs}
\end{proof}

\begin{remark}
Let us observe that if one were to try to use the (optimal) $L^{\infty}$ estimate in Corollary~\ref{cor:pointwise} instead of Lemma~\ref{lem:lp kernel bound} in the proof of Lemma~\ref{lem:remainder cherry gpam},
 we would necessarily pick up an additional logarithmic factor.
Thus, the use of $L^p$ bound seems to be crucial here.
\end{remark}

\begin{remark}\label{rem:non-comm2}
We observe that it follows from the proof above that 
$\lim_{\eps \to 0} c_{\eps,\delta}^{\<2>}= 0$ for all $\delta>0$, while combining \eqref{eq:c-const1},\eqref{eq:c-const2} with the last equation in the proof of Lemma~\ref{lem:remainder cherry gpam} one finds  that in general 
 $\lim_{\eps \to 0}\lim_{\delta \to 0} c_{\eps,\delta}^{\<2>}$ is non-zero.
\end{remark}

\begin{remark}\label{rem:non-comm3}
Observe that $\gamma_{\eps,\delta}=\int_{[0,1]^3} f^{\eps} F^{\mu,\nu}_{\eps,\delta} = \int_{[0,1]^3} f (F^{\mu,\nu}_{\eps,\delta}\circ\mathcal{S}^{\eps^{-1}}) $. 
Thus, it follows from the forth item of Lemma~\ref{lem:remainder cherry gpam}
that 
 in general
\begin{equ}
\lim_{\delta\downarrow 0} \lim_{\eps\downarrow 0}  \gamma_{\eps,\delta} = 0\neq  \lim_{\eps\downarrow 0}\lim_{\delta\downarrow 0}\gamma_{\eps,\delta} \  .
\end{equ}.
\end{remark}

%
%

\subsubsection{Counterterms for Theorem~\ref{thm:main_translation invariant}}
In this section we identify functions $F_{\eps,\delta}^\flat$, $F_{\eps,\delta}^{\flat,\mu,\nu}$ when working with the regularisation $\xi_\delta= \xi(\rho^\delta_x)$. 
We define similarly to above $$F^\flat_{\eps,\delta}(z)= \sum_{\tau \in \mathfrak{T}_{\<Xi2>} } \E[\bfPi^{\eps,\delta}\tau] -D^{\<Xi2>}_{\delta/\eps}(x,t) - \frac{\bar{\alpha}^{\<Xi2>,\flat}_{\eps,\delta}}{\sqrt{\det (\bar{A})}}- c^{\<Xi2>,\flat}_{\eps,\delta} \ ,$$
where
\begin{equs}
D^{\<Xi2>}_\lambda(x,t) &\eqdef
\int_{(\mathbb{R}^{3})}
\left(\rho^{\lambda}\right)^{*2} (\zeta-x) \kappa (\tau)  
 Z^*_{1,0}(\zeta,\tau; x,t) \, d\zeta d\tau  \\
\frac{\bar{\alpha}^{\<Xi2>,\flat}_{\eps,\delta}}{\sqrt{\det (\bar{A})}}&= \int \big(\kappa(t-s)- \kappa^{\eps}(t-s)) \bar{\Gamma}(y,s,x,t)\big(\rho^{\delta}\big)^{*2}(x-y)dy ds\\
 c^{\<Xi2>,\flat}_{\eps,\delta}&=\int_{[0,1]^3} \int_{\mathbb{R}^{3}}
\big[\kappa^{\eps} (t-s) \big(\Gamma_{\eps}-Z^*_{\eps;0} \big)+ \kappa^{\eps}_c (t-s) \big(\Gamma_{\eps}-\bar{\Gamma} \big)\big](x,t;y,s)\big(\rho^{\delta}\big)^{*2}(x-y) dyds dxdt\;.
\end{equs}
as well as
\begin{align*}
F^{\flat, \mu, \nu}_{\eps,\delta}(z)&= \sum_{\tau \in \mathfrak{T}_{\<2>{\mu,\nu}} } \E[\bfPi^{\eps,\delta}\tau] 
- D^{\<b2>{\mu,\nu}}_{\eps/\delta}(A_s^\eps)
- \sum_{i,j=1}^2 (\mathbf{1}_{\mu=i}+(\partial_{\mu} \phi_i)^\eps)(\mathbf{1}_{\nu=j}+(\partial_{\nu} \phi_j)^\eps) 
\frac{ (\bar{A}^{{-1}})_{i,j}    }{\sqrt{\det(\bar{A})}} \bar{\alpha}^{\<b2>,\flat}_{\eps,\delta} 
 + c^{\<b2>{\mu,\nu},\flat }_{\eps,\delta}
\end{align*}
where
\begin{equs}
D^{\<b2>{\mu,\nu}}_{\lambda}(A_s) &= \int \kappa(t-s)\kappa (t'-s') Z^1_{0;\mu}(x,t;y,s) Z^1_{0;\nu}(x,t;y',s') 
		\big(\rho^{\lambda}\big)^{*2}(y-y')dy dy' ds' ds' \\
\frac{ (\bar{A}^{{-1}})_{i,j}    }{\sqrt{\det(\bar{A})}} \bar{\alpha}^{\<b2>,\flat}_{\eps,\delta} &= \int 
\kappa(t-s)\kappa^\eps_c(t-s)\kappa (t'-s') \kappa^\eps_c (t'-s') \partial_i \bar{\Gamma} (x,t;y,s)	\partial_j \bar{\Gamma}(x,t;y',s')  \big(\rho^{\lambda}\big)^{*2}(y-y')dy\\
c^{{\<2>{\mu,\nu}},\flat }_{\eps,\delta}&= \int_{[0,1]^3} h_{\eps,\delta}^{\<b2>{\mu,\nu},\flat}\\ 
h_{\eps,\delta}^{\<b2>{\mu,\nu},\flat} &= \sum_{\tau \in \mathfrak{T}_{\<2>{\mu,\nu}} } \E[\bfPi^{(\eps,\delta)}\tau] 
- D^{\<b2>{\mu,\nu}}_{\eps/\delta}(A_s^\eps)
- \sum_{i,j=1}^2 (\mathbf{1}_{\mu=i}+(\partial_{\mu} \phi_i)^\eps)(\mathbf{1}_{\nu=j}+(\partial_{\nu} \phi_j)^\eps) 
\frac{ (\bar{A}^{{-1}})_{i,j}    }{\sqrt{\det(\bar{A})}} \bar{\alpha}^{\<b2>,\flat}_{\eps,\delta} \;.
\end{equs}

\begin{lemma}\label{lem:zero chaos flat}
The function $F^{\flat}$ satisfies the properites listed in Lemma~\ref{lem:finite wiener chaos remainder1} and the functions $F^{\flat, \mu, \nu}_{\eps,\delta}$ satisfy the same properties as listed in Lemma~\ref{lem:remainder cherry gpam}.
\end{lemma}
\begin{proof}
The proof follows exactly as the proof of Lemmats~\ref{lem:finite wiener chaos remainder1} and~\ref{lem:remainder cherry gpam}, the only difference being that the $\delta^{2}$ time shift in the heat kernels is replaced everywhere by a convolution with $\rho^{\delta}$.
\end{proof}

\subsubsection{Counterterms for Theorem~\ref{thm:main_translation invariant2}}
Finally, let 
$${c}_{\eps,\delta}^{\<Xi2>,\flat\flat}\eqdef c^{\<Xi2>, \flat}_{\eps,\delta} + \int_{[0,1]^{3}} D^{\<Xi2>}_{\delta/\eps} , \qquad
  {c}_{\eps,\delta}^{{\<b2>{\mu,\nu}, \flat\flat }}=c^{\<b2>{\mu,\nu}, \flat}_{\eps,\delta} + \int_{[0,1]^{3}} D^{\<b2>{\mu,\nu}}_{\delta/\eps}$$
\begin{lemma}
Both, ${c}_{\eps,\delta}^{\<Xi2>,\flat\flat}$ and ${c}_{\eps,\delta}^{{\<b2>{\mu,\nu}, \flat\flat }}$ are bounded on $\triangle^{<}_{C}$ for any $C>0$ and for $\delta>0$ vanish as $\eps\to 0$.
\end{lemma}
\begin{proof}
Boundedness follows directly from the boundedness of $ D^{\tau}_{\lambda}$ for $\lambda>C$ and boundedness of $c^{\tau,\flat}_{\eps,\delta}$.
Similarly, the second claim follows since  $D^{\tau}_{\lambda}\to 0$ as $\lambda\to \infty$ and $c^{\tau, \flat}_{\eps,\delta}\to 0$ as $\eps\to 0$ for $\delta>0$.
\end{proof}

Defining 
$$F^{\flat\flat}_{\eps,\delta}(z)= \sum_{\tau \in \mathfrak{T}_{\<Xi2>} } \E[\bfPi^{(\eps,\delta)}\tau] - \frac{\bar{\alpha}^{\<Xi2>,\flat}_{\eps,\delta}}{\sqrt{\det (\bar{A})}}- c^{\<Xi2>,\flat\flat}_{\eps,\delta} \ ,$$
and 
$$
F^{\flat\flat, \mu, \nu}_{\eps,\delta}(z)= \sum_{\tau \in \mathfrak{T}_{\<2>{\mu,\nu}} } \E[\bfPi^{(\eps,\delta)}\tau] 
- \sum_{i,j=1}^2 (\mathbf{1}_{\mu=i}+(\partial_{\mu} \phi_i)^\eps)(\mathbf{1}_{\nu=j}+(\partial_{\nu} \phi_j)^\eps) 
\frac{ (\bar{A}^{{-1}})_{i,j}    }{\sqrt{\det(\bar{A})}} \bar{\alpha}^{\<b2>,\flat}_{\eps,\delta} 
 + c^{\<b2>{\mu,\nu},\flat\flat }_{\eps,\delta}\ ,
$$
one obtains the following direct corollary of Lemmas~\ref{lem:finite wiener chaos remainder1} and~\ref{lem:remainder cherry gpam}.
\begin{corollary}\label{cor:zero chaos flatflat}
The functions $F^{\flat\flat}$, resp.\ $F^{\flat\flat, \mu, \nu}_{\eps,\delta}$ satisfy the properties listed in Lemma~\ref{lem:finite wiener chaos remainder1} resp.\ Lemma~\ref{lem:remainder cherry gpam} on the restricted domain $(\eps,\delta) \in \triangle^{<}_{C} \cap \bar{\square}$. 
\end{corollary}

\subsubsection{Proof of Theorem~\ref{thm:g-pam}, Theorem~\ref{thm:main_translation invariant} \& Theorem~\ref{thm:main_translation invariant2} }\label{sec:proof gpam}
Finally, we are in the position to prove the main results on the oscillatory g-PAM equation. It remains only to combine what has been done so far which we shall do for all three theorems simultaneously, since the regularisation specific arguments have already been carried out it.
\begin{proof}[Proof of the results on the g-PAM equation]
Observe that the solution map extends continously to $\big((0,1]\times [0,1]\big) \cap \bar{\square}$ by the results on the g-PAM equation in \cite{Sin23}. 
Next we choose the homogeneity assigment for the regularity structure constructed in Section~\ref{sec:abstract form gpam},  imposing that
 that $1<-|\Xi|<\gamma<L,R$ and $\beta<2$ and that
 furthermore $0<\kappa, |\beta-2|, |L-1|, |R-1|$ be sufficiently small (say $1/100$). This puts one in the setting of Theorem~\ref{thm:fixed point}
for the abstract equation \eqref{eq:abstract_gpam}.
To conclude we check the following.
\begin{enumerate}
\item For the initial condition the only distinction to \cite{Hai14} is that the abstract fixed point theorem, Theorem~\ref{thm:fixed point}, requires that $\bar{\eta}> (- R)$. For this reason we take $v_{\mathtt{ in}}^\eps\in C^{\gamma, \eta}$ for $\eta>1-R/2$, which is possible by  Lemma~\ref{lem:initial2} for \eqref{eq:initial2}. 
\item Lemma~\ref{lem:convergence of kernels easy terms} and Proposition~\ref{prop:vanishing of G kernel} guarantee convergence of the involved kernels.
\item Convergence of models follows from the stochastic estimates obtained in Lemma~\ref{lem:conv dumbel1}, Lemma~\ref{lem:convergene pam cherry} and Lemma~\ref{lem:conv_resonances}, which in turn take as input the estimates on the functions $F, F^{\mu,\nu}$ from Lemma~\ref{lem:finite wiener chaos remainder1} \& Lemma~\ref{lem:remainder cherry gpam} when working with kernel regularisation, respectively the same estimates on $F^\flat, F^{\flat, \mu,\nu}$ 
from Lemma~\ref{lem:zero chaos flat}
, resp.~$F^{\flat\flat}, F^{\flat\flat, \mu,\nu}$ from Corollary~\ref{cor:zero chaos flatflat}. Then convergence of the model follows by a Kolmogorov theorem, for trees not containing multiplication with an abstract scale this just follows as in \cite[Thm.~10.7]{Hai14}, while for the remaining trees where the abstract operator $\mathcal{E}$ is present it follows using Lemma~\ref{lem:estimating Gamma with scale}.
\end{enumerate}
Then, we observe that in Lemma~\ref{lem:ren_equ} the counterterms takes the form
$$ \sum_{\tau \in \mathfrak{T}_{\tau} } g_{\eps,\delta}^\tau \ ,
 $$
the exact asymptotic behaviour of which follows by straightforward (but slightly tedious) computations, see \cite[Sec.3.1.1]{Sin23} for very similar calculations. We conclude the proof by Theorem~\ref{thm:fixed point}.
\end{proof}

%
%
%


\section{Application to the \TitleEquation{\Phi^4_3}{Phi43} Equation}\label{sec:application phi}

In this section we prove our main results about the $\Phi^4_3$ equation, Theorem~\ref{thm:phi4} and Theorem~\ref{thm:phi4 restricted}. We follow the same strategy as for g-PAM in Section~\ref{sec:application gpam}.\footnote{Since we only consider one specific regularisation, we could in principle have slightly shortened this section. But the presented structure has the advantage that it makes it clear where the modifications for other regularisations would enter.
} 
Recall that for $\eps>0$, the equation 
$$\left(\partial_t - \nabla \cdot A(x/\eps,t/\eps^2) \nabla \right) u = - f^\eps u^3 + \xi  $$
can be treated in the framework of regularity structures by considering the `lifted' equation
$$
U  = -\mathcal{K}_\eps \big(\mathbf{R}_+( \pmb{f} U^3 +\Xi )\big)  + v_{in} \;,
$$
c.f.\ \cite[Sec.~9]{Hai14}, as well as the expository work~\cite{hai15exp}.
In order to reduce the number of stochastic estimates required, we shall consider the remainder equation formally for  $V=U-\<1>$, where
\begin{equs}
V  
&
= -\mathcal{K}_\eps  \left( \mathbf{R}_+( V^3\cdot \pmb{f} + 3 V^2 \cdot \pmb{f}\<1> + 3 V \pmb{f}\<2>) \right)  -\mathcal{K}_\eps  ( \mathbf{R}_+\pmb{f}\<3> )   + v_{in} \label{eq:abstract_phi^4_3}
\end{equs}
and we interpret the terms $\<1>,\ \<2>,\ \<3>$ as abstract noises.
This rewriting of the equation can thus be solved using a regularity structure built from a normal set of trees containing
%
 $$ \mathfrak{T}_c\eqdef\big\{ 
  \pmb{f}\<2>,\  \pmb{f}\<2>  I(\pmb{f}\<3>)   , \  \<1>, \ \<1> I(\pmb{f}\<3>)  , \ \pmb{f}\<2> I(\pmb{f}\<2>) , \ \mathbf{X}_i \pmb{f}\<2>, \ \mathbf{1},\   I(\pmb{f}\<2>) , \
\pmb{f}\<2> I(\pmb{f}\<2> ) ,
 \   I(\pmb{f}\<2> ),  \  \mathbf{X}_i \big\}\ ,$$
 where the homogeneity assignment is given 
 by $|\pmb{f}|=-\kappa\ , | \<1>|= \zeta \ , |\<2>|= 2\zeta \ , |\<3>|= 3\zeta$, for appropriate $\kappa>0$, $ \zeta<-1/2$
and by declaring $I$ to be $2$-regularising. 
The solution then takes the form
 $$V(z)= v(z) \mathbf{1} + \sum_{i=1}^3 v_i(z) \mathbf{X}_i-  3 v(z)  I(\pmb{f}\<2>)  - I(\pmb{f}\<3>) .$$

\paragraph{Construction of the regularity structure:}
We again lift the equation as in \eqref{eq:abstract}, that is
\begin{equs}
U_\eps&=  
 \mathcal{\bar K}_\eps \big( \mathbf{R}_+ \hat{\opF} (U  , \Xi  ) \big) 
+ 
\sum_{\substack{j=1} }^3 
\bar{\mcK}^{0,j}_{\epsilon} 
\big(\tilde{\Phi}_j \cdot \mathbf{R}_+\hat{\opF} (V ,  \pmb{f},\<1>,\<2>)\big) 
+
\sum_{i, j=1 }^3 \hat{\mathcal{E}}^\eps\big(\Phi_i \cdot
\bar{\mcK}^{i,j}_{\epsilon} 
\big(\tilde{\Phi}_j \cdot \mathbf{R}_+ \hat{\opF} (V ,  \pmb{f},\<1>,\<2>)\big)\big) \\
&\qquad +
\mathcal{G}_\epsilon \big( \mathbf{R}_+\hat{\opF} (V ,  \pmb{f},\<1>,\<2>)\big)
+p_\eps 
 + v^\eps_{\texttt{in}} \ , \label{eq:abstract_phi4}
\end{equs}
for $\hat{\opF} (V ,  \pmb{f},\<1>,\<2>)\eqdef( V^3\cdot \pmb{f} + 3 V^2 \cdot \pmb{f} \<1> + 3 V \cdot \pmb{f}\<2>)$ and where 
\begin{equ}\label{eq:modelled dist}
p_\eps\eqdef
\mathcal{\bar K}_\eps \big( \mathbf{R}_+ \hat{\opF} (U  , \Xi  ) \big) 
+ 
\sum_{\substack{j=1} }^3 
\bar{\mcK}^{0,j}_{\epsilon} 
\big(\tilde{\Phi}_j \cdot \mathbf{R}_+\pmb{f}\<3>\big) 
+
\sum_{i,j=1}^3 \hat{\mathcal{E}}^\eps\big(\Phi_i \cdot
\bar{\mcK}^{i,j}_{\epsilon} 
\big(\tilde{\Phi}_j \cdot \mathbf{R}_+ \pmb{f}\<3>\big)\big) 
+
\mathcal{G}_\epsilon \big( \mathbf{R}_+\pmb{f}\<3>\big) \ ,
\end{equ}
which will be seen in the proof of Theorem~\ref{thm:phi4} to be a modelled distribution belonging to $\mathcal{D}^{\gamma,\eta}$ for any $\gamma>0$ and $\eta<0$ (for the models fixed in the sequel). 

%
Denote by $\mathfrak{T}$ the minimal normal set constructed from $\mathfrak{T}_c$ as for the g-PAM in Section~\ref{sec:abstract form gpam} and by $T$  the span of $\mathfrak{T}$. 
For $\tau\in \mathfrak{T}_c$ we denote by $\mathfrak{T}^\tau\subset \mathfrak{T}$ the set of trees obtained by substituting every occurrence of
$I$ by an 
an element of $\{ I, \bar{I},\  I^{0,j}( \tilde{\Psi}_{j}\, \cdot ),\  \mathcal{E}( \Psi_i I^{i,0}( \cdot) ),\  \mathcal{E}( \Psi_i I^{i,j}( \tilde{\Psi}_{j}\, \cdot) )  : i,j>0 \}$.

The homogeneity assignment of each tree is then determined by declaring 
$I, \bar{I}, I^{0, j}$ to be $\beta$ regularising and $I^{i,j}$ to be $\beta-1$ regularising for $i>0$ as well as declaring
and $ \pmb{f}, \Psi_i, \tilde{\Psi}_j$ to be of homogeneity $-\kappa$ and $| \<1>|= \zeta \ , |\<2>|= 2\zeta \ , |\<3>|= 3\zeta$,
with the values of $\kappa>0$, $ \zeta<-1/2$, $\beta<2$ specified later.
The structure group $\Gamma$ is then implicitly defined by the conditions in the definition of an abstract integration operator \cite[Def.~5.7]{Hai14} and 
abstract multiplication by a scale in Definition~\ref{def:abstract scale}.

\paragraph{Construction of models}
For $\eps\in[0,1]$, $\delta>0$ and a continuous noise $\xi_{\eps,\delta}\in C$, we shall consider, 
for fixed bounded functions $g_{\eps,\delta}^{\<2>}$ and constants $\gamma^{\<2>}_{\eps,\delta}$ 
the (partially renormalised) model $M^{(\eps,\delta)}=(\Pi^{(\eps,\delta)}, \Gamma^{(\eps,\delta)})$ which is characterised by
\begin{enumerate}
\item $ \Pi^{(\eps,\delta)}_{x}\Phi_i = \phi^{\eps}_i\ , \qquad  \Pi^{(\eps,\delta)}_{x}\tilde{\Phi}_j = \tilde{\phi}^{\eps}_j\ , \qquad\Pi^{(\eps,\delta)}_{x}\pmb{f}= f^\eps\ , $
\item $\Pi^{(\eps,\delta)}_x \<1>=\Pi^{(\eps,\delta)} \<1>\eqdef K_\eps(\xi_{\eps,\delta}), \qquad \Pi^{(\eps,\delta)}_x \pmb{f}\<1>=f^\eps K_\eps(\xi_{\eps,\delta}) \ , $
\item $\Pi^{(\eps,\delta)}_x \<2>= \big(\Pi^{(\eps,\delta)} \<1>\big)^2 - g^{\<2>}_{\eps,\delta} \  
,\qquad
 \Pi^{(\eps,\delta)}_x \pmb{f}\<2>=  f^\eps \big(\big(\Pi^{(\eps,\delta)} \<1>\big)^2 - g^{\<2>}_{\eps,\delta}\big)- \gamma^{\<2>}_{\eps,\delta} \ ,
$
\item
$
\Pi^{(\eps,\delta)}_x \<3>= \big(\Pi^{(\eps,\delta)} \<1>\big)^3 - 3g^{\<2>}_{\eps,\delta}\Pi^{(\eps,\delta)} \<1> \ , 
\qquad
 \Pi^{(\eps,\delta)}_x \pmb{f}\<3>= f^\eps K_\eps(\xi_{\eps,\delta})^3 - 3\big(g^{\<2>}_{\eps,\delta}f^\eps + \gamma^{\<2>}_{\eps,\delta}\big)K_\eps(\xi_{\eps,\delta}  ) \ .$
\item It acts as the polynomial model on the polynomial sector $\bar{T}\subset T$.
\item The abstract integration map $I$ realises the kernel $G$, $I^{i,j}$ realises $K^{i,j,+}$.
\item\label{Item:multiplicative2} It is multiplicative, i.e. $\Pi_x(\tau \tau')= \Pi_x \tau \Pi_x \tau'$.
\end{enumerate}

We shall construct (fully) renormalised models by modifying Item~\ref{Item:multiplicative2}  above for certain elements. 
For $I'\in \{ I, \bar{I},  I^{0,j}(\tilde{\Psi}_{j}\, \cdot) ),  \mathcal{E}(\Psi_i I^{i,j}(\tilde{\Psi}_{j}\, \cdot) )) : i,j> 0 \}$, let
\begin{itemize}
\item $\hat{\Pi}^{(\eps,\delta)}_x \<2> I' (\pmb{f}\<2>) \eqdef {\Pi}^{(\eps,\delta)}_x \<2> I' (  \pmb{f}\<2>) -g_{\eps,\delta}^{\<2> I' (\pmb{f}\<2>)} \ ,$ 
\item $\hat{\Pi}^{(\eps,\delta)}_x \pmb{f}\<2> I' ( \pmb{f}\<2>) 
=  f^\eps {\Pi}^{(\eps,\delta)}_x \<2> I' (  \pmb{f}\<2>) -
f^\eps g_{\eps,\delta}^{\<2> I' (\pmb{f}\<2>)}  - \gamma^{ \pmb{f}\<2> I' (\pmb{f}\<2>)} \ ,$
\item $\hat{\Pi}^{(\eps,\delta)}_x \<2> I' ( \pmb{f}\<3>) = {\Pi}^{(\eps,\delta)}_x \<2> I' ( \pmb{f}\<3>) - 3 g_{\eps,\delta}^{\<2> I' (\pmb{f}\<2>)} \Pi \<1>\ , $ 
\item $\hat{\Pi}^{(\eps,\delta)}_x \pmb{f}\<2> I' ( \pmb{f}\<3>) 
=
f^\eps {\Pi}^{(\eps,\delta)}_x \<2> I' ( \pmb{f}\<3>) - 3 \big( f^\eps g_{\eps,\delta}^{\<2> I' (\pmb{f}\<2>)}  + \gamma_{\eps,\delta}^{\<2> I' (\pmb{f}\<2>)} \big) \Pi \<1>  \ ,
$
\end{itemize}
for some bounded functions $g^\tau_{\eps,\delta}$ and constants $\gamma^\tau_{\eps,\delta}$.

\subsection{Renormalised equation}
We define 
$$g^{\<22>}\eqdef\sum_{I'\in \{ I, \bar{I}, I^{0,j}(\tilde{\Psi}_{j}\, \cdot) ),  \mathcal{E}(\Psi_i I^{i,j}(\tilde{\Psi}_{j}\, \cdot) ) : i,j> 0 \}}
g_{\eps,\delta}^{\<2> I' (\pmb{f}\<2>)}$$
and 
$$
\gamma^{\<22>}\eqdef\sum_{I'\in \{ I, \bar{I},  I^{0,j}(\tilde{\Psi}_{j}\, \cdot) ),  \mathcal{E}(\Psi_i I^{i,j}(\tilde{\Psi}_{j}\, \cdot) ) : i,j> 0 \}}
\gamma_{\eps,\delta}^{\<2> I' (\pmb{f}\<2>)}$$
%

%
 
 \begin{lemma}\label{lem:ren_equ2}
Let $\hat{V}^{(\eps,\delta)}$ be the solution to  \eqref{eq:abstract_phi4} with respect to the model $\hat{M}^{(\eps,\delta)}$. Then 
$\hat{u}_{\eps,\delta}\eqdef\Pi_x \<1>+ \hat{\mathcal{R}}V^{(\eps,\delta)}$ satisfies 
\begin{align*}
\partial_t \hat{u}_{\eps,\delta} -\nabla \cdot A^\eps \nabla \hat{u}_{\eps,\delta} &=
- f^\eps \Big[  u^{3}_{\eps,\delta} - 3 g^{\<2>}_{\eps,\delta}  u_{\eps,\delta}  + 9 g^{\<22>}_{\eps,\delta}  u_{\eps,\delta}  \Big] - 3\gamma_{\eps,\delta} ^{\<2>} u_{\eps,\delta}  + 9  \gamma_{\eps,\delta} ^{\<22>} u_{\eps,\delta}  + \xi_{\eps,\delta}\ .
\end{align*}
\end{lemma}

\begin{proof}
We omit ${\eps,\delta}$ in the notation.
First observe that  the solution can be written as
$V= v \mathbf{1} + V_1 + V_2 + V_3$, where $\mathcal{R}V=v$  and we set
\begin{align*}
V_1 &= \sum_{k=1}^3 v'_k\mathcal{E}( \Psi_k)
+ \sum_{k=1}^3 v''_{k}  X_k \\
V_2&=-3v \Bigl(   I( \pmb{f}\<2>)+\bar{I}( \pmb{f}\<2>) + \sum_{j=1}^3 I^{0,j}( \tilde{\Phi}_{j} ( \pmb{f}\<2>) )   +\sum_{i>0,j\geq 0} \mathcal{E}( \Phi_i I^{i,j}( \tilde{\Phi}_{j} ( \pmb{f}\<2>) ) ) \Bigr)\\
V_3 &=- \Bigl(   I( \pmb{f}\<3>)+\bar{I}( \pmb{f}\<3>) + \sum_{j=1}^3 I^{0,j}( \tilde{\Phi}_{j} ( \pmb{f}\<3>) )   +\sum_{i>0,j\geq 0} \mathcal{E}( \Phi_i I^{i,j}( \tilde{\Phi}_{j} ( \pmb{f}\<3>) ) ) \Bigr)\;.
\end{align*}
To derive the renormalised equation, we make the following preliminary observations.
\begin{equs}
&\hat{\Pi}_x  [V^3(x)\cdot \pmb{f} ](x)= v(x) f(x)\ , \qquad \hat{\Pi}_x  [V^2(x) \cdot \pmb{f} \<1> ](x)=v^2(x)   f(x) \Pi_x \<1>(x)\ ,\\
&\hat{\Pi}_x  [V_1(x) \cdot \pmb{f}\<2>](x)=0\ , \qquad \hat{\Pi}_x  [V_2 (x)\cdot \pmb{f}\<2>](x)
=
3v(x) \big( f g^{\<22>}(x) + \gamma^{\<22>}\big)\ ,\\
&\hat{\Pi}_xQ_{\leq 0} [V_3(x) \cdot \pmb{f}\<2>](x)
=
3\Pi^{(\eps,\delta)} \<1>(x) \big( f(x) g^{\<22>}(x) + \gamma^{\<22>}\big) ,
\end{equs}
where for each term we have used that only terms with non-positive homogeneity contribute together with the explicit description of the renormalised map $\hat{\Pi}$ in terms of  ${\Pi}$.
 
%

Thus in particular,
\begin{align*}
\hat{\Pi}_xQ_{\leq 0} V \cdot \pmb{f}\<2>&=\hat{\Pi}_x[ \pmb{f}\<2>](x)\cdot v(x) + 3 \big( v(x) + \Pi \<1>(x)\big) \cdot\big( f g^{\<22>}(x) + \gamma^{\<22>}\big)\\
&=\hat{\Pi}_x[ \pmb{f}\<2>](x)\cdot v(x) + 3 u(x) \cdot\big( f(x) g^{\<22>}(x) + \gamma^{\<22>}\big)\ ,
\end{align*}
and therefore 
\begin{align*}
\hat{\Pi}_xQ_{\leq 0} \hat{\opF} (V ,  \pmb{f},\<1>,\<2>)(x)&= \hat{\Pi}_x [V^3\cdot \pmb{f} + 3 V^2 \cdot \pmb{f} \<1> + 3 V \cdot \pmb{f}\<2>](x)\\
&=f(x)\cdot \mathcal{R}V(x)^3+ 3(\mathcal{R}V(x))^2 f\Pi_x \<1>(x)+ 3
\Pi_x[\pmb{f}\<2>](x) \mathcal{R}V(x) \\
&\qquad+ 9 \mathcal{R}U(x)  \cdot \big( f(x) g^{\<22>}(x) + \gamma^{\<22>}\big)\ .
\end{align*}
We conclude that
\begin{align*}
\big(\partial_t-\nabla A^\eps \nabla \big) v
&= -\Big[ f \big[ v^3 + 3v^2 \Pi_x \<1> + 3 v ( \Pi_x \<1>- g^{\<2>}) \Big]\\
&\qquad +3\gamma^{\<2>} v
-9  \big[f  g^{\<22>}+\gamma^{\<22>} \big] u 
-3f\Big[\big( \Pi_x \<1>\big)^3 
-3 g^{\<2>} \Pi_x \<1> \Big] - 3 \gamma^{\<2>} \Pi_x \<1> 
\\
&= - f \Big[  u^{3}- 3 g^{\<2>} u + 9 g^{\<22>} u \Big] - 3\gamma^{\<2>} u + 9  \gamma^{\<22>} u \ ,
\end{align*}
which is precisely the claim.
\end{proof}
\begin{remark}
Since the form of the nonlinearity does not require us to work in a function like sector, we could in principle lift as in Remark~\ref{rem:polynomial-nonlinearities}, but this seems to lead to complications when defining the renormalised model and identifying the renormalised equation (due to the presence of additional instances of the correctors in the sunset diagram).
\end{remark}
\subsection{Convergence of renormalised models}\label{sec:convergence of ren. models}
Convergence of the the terms 
$$ \Pi^{(\eps,\delta)}_{x}\Phi_i = \phi^{\eps}_i\ , \qquad  \Pi^{(\eps,\delta)}_{x}\tilde{\Phi}_j = \tilde{\phi}^{\eps}_j\ , \qquad\Pi^{(\eps,\delta)}_{x}\pmb{f}= f^\eps\ 
$$
follows directly from Lemma~\ref{lem:appendix}.
Next we turn to
$\Pi^{(\eps,\delta)}_x \<1>=\Pi^{(\eps,\delta)} \<1>\eqdef K_\eps(\xi_{\eps,\delta}) \ $
which unravelling the definition is (formally) equal to 
\begin{equ}\label{eq:linear_guy}
\Pi^{(\eps,\delta)} \<1>(x,t)= \int_{\mathbb{R}^3\times \mathbb{R}} \kappa(t-s) \Gamma_\eps(x,t,y,s-\delta^2) \xi(dy,ds)\ .
\end{equ}
Then \cite[Prop.~5.1]{HS23per} essentially\footnote{Therein, there was an additional time regularisation but the proof adapts ad verbatim.} provides the following bounds.
%
\begin{lemma}\label{lem:stochastic estimates_linear}
Let $\alpha, \alpha' ,\kappa>0$ such that $\alpha'+\kappa<1$. 
Using the shorthand notation $\triangle_{s,t}^{(\eps,\delta)}\<1>_{\eps, \delta} \eqdef \Pi^{(\eps,\delta)} \<1>(\cdot, t)- \Pi^{(\eps,\delta)} \<1>(\cdot, s)$
it holds that
\begin{enumerate}
\item\label{eq:0a} $\E[\langle \Pi^{(\eps,\delta)} \<1>(\ \cdot\ , t) , \psi^\lambda_x \rangle^2]^{\frac{1}{2}} \lesssim \lambda^{-\alpha-d/2+1}$
\item\label{eq:1a} $\E [\langle \triangle_{s,t}^{(\eps,\delta)} \<1> , \psi^\lambda_x \rangle^2]^{\frac{1}{2}} \lesssim |t-s|^{\alpha'/2} \lambda^{-\alpha-\alpha'-/2+1}$
\item\label{eq:1aa} $\E[\langle \triangle_{s,t}^{(\eps,\delta)} \<1>-\triangle_{s,t}\Pi^{{0,\delta}} \<1> , \psi^\lambda_x \rangle^2]^{\frac{1}{2}} \lesssim \eps^{\kappa} |t-s|^{\alpha'/2} \lambda^{-\alpha-\alpha'-\kappa-d/2+1}$
\item\label{eq:1aaa} $\E[\langle \triangle_{s,t}^{(\eps,\delta)} \<1>-\triangle_{s,t}\Pi^{{\eps,0}} \<1> , \psi^\lambda_x \rangle^2]^{\frac{1}{2}} \lesssim \delta^{\kappa} |t-s|^{\alpha'/2} \lambda^{-\alpha-\alpha'-\kappa-d/2+1}\ ,$
\end{enumerate}
uniformly over $|t-s| \vee \lambda \le 1$, $\eps,\delta \in (0,1]$, $x \in \R^d$, 
and $\psi \in \mathfrak{B}^{0}(\mathbb{R}^{d})$.
\end{lemma}
\subsubsection{Renormalisation of Wick powers}

For functions 
$\{F^{\<2>}_{\eps,\delta}\}_{\eps\in (0,1], \delta\in [0,1]}$ that will be chosen later set 
$$g^{\<2>}_{\eps,\delta}(z)= \E[\Pi^{(\eps,\delta)}\tau(z)]- \mathbf{1}_{\tau=\Xi \bar{I}\Xi  } F^{\<2>}_{\eps,\delta}(z) \ , \qquad
\gamma_{\eps,\delta}^{\<2>}= \int_{[0,1]^3} f^{\eps}_{\mu,\nu} F^{\tau}_{\eps,\delta} dxdt\ .$$
Next we set 
$$\|\<1>\|=\|\pmb{f}\<1>\|=-1/2\ ,\qquad \|\<2>\|=\|\pmb{f}\<2>\|=-1\ , \qquad \| \<3>\|=\|\pmb{f}\<3>\|=-3/2. $$
\begin{prop}\label{prop:stochastic estimates}
Assume that $F_{\eps,\delta}$ is uniformly bounded and $\eps\mathbb{Z}^3\times \eps^2\mathbb{Z}$ periodic with vanishing mean.
For every $\alpha>0, \kappa\in [0,1)$ and $T\in \{\<1>,\pmb{f}\<1>,\<2>,\pmb{f}\<2>,\<3>,\pmb{f}\<3>\}$, it holds that 
\begin{enumerate}
\item\label{eq:0b} $\E[\langle  \hat{\Pi}_\star^{(\eps,\delta)}T , \psi^\lambda_{\star} \rangle^2]^{\frac{1}{2}} \lesssim_{\alpha} \lambda^{\|T\|-\alpha}$
\item\label{eq:1b} $\E[\langle  \hat{\Pi}_\star^{(\eps,\delta)}T- \hat{\Pi}_\star^{(0,\delta)}T , \psi^\lambda_{\star} \rangle^2]^{\frac{1}{2}} \lesssim_{\alpha,\kappa}\eps^{\kappa}  \lambda^{\|T\|-\alpha-\kappa}\ $
\end{enumerate} 
uniformly over $\eps, \delta\geq 0$, $\varphi\in \mathfrak{B}_0$ and $\star\in \mathbb{R}^{3}$.
uniformly over $\star\in \mathbb{R}^{3}$, $\lambda\in (0,1]$ and $\psi \in \mathfrak{B}^{2\kappa}(\mathbb{R}^{3})$.
\end{prop}
\begin{proof}
Following the argumentation of \cite[Sec.~4.2]{HS23per} ad verbatim shows the proposition for $T\in \{\<1>,\<2>,\<3>\}$. Then Item~\ref{eq:0b} follows for the remaining trees by simply interpreting the oscillatory function $f^\eps$ as part of the test function. Finally, {Item~\ref{eq:1b}} follows for $T=\pmb{f}\<1>$ and 
$\pmb{f}\<3> $ as well as the the contribution in the Second Wiener Chaos for $
\pmb{f}\<2>$ exactly as the bound on the second term of \cite[Eq.~(4.22)]{HS23per} in the proof of Prop.~4.13 therein. 
The contribution of the mean to Item~\ref{eq:1b} for $T= \pmb{f}\<2>$ is bounded by Lemma~\ref{lem:appendix}.
\end{proof}

\subsubsection{Convergence of Larger trees}
For 
 ${I'\in \{ I, \bar{I},  I^{0,j}(\tilde{\Psi}_{j}\, \cdot) ),  \mathcal{E}(\Psi_i I^{i,j}(\tilde{\Psi}_{j}\, \cdot) ) : i,j> 0 \}}$
 we set
$$
g_{\eps,\delta}^{\<2> I' (\pmb{f}\<2>)}(\star)= 
\E [\bfPi^{(\eps,\delta)} \<2> I' (\pmb{f}\<2>)]- \E [\bfPi^{(\eps,\delta)} \<2> (\star) ] \E [\bfPi^{(\eps,\delta)} I' (\pmb{f}\<2>) (\star) ]
- \mathbf{1}_{\tau=\<2> \bar{I} ( \pmb{f}\<2>)  } F^{\<22>}_{\eps,\delta}(\star) \ ,
$$
where
$\{F^{\<22>}_{\eps,\delta}\}_{\eps\in (0,1], \delta\in (0,1]}$ will be specified later.
\begin{lemma}\label{lem:stochastic estimate sunset}
Assume that $F^{\<22>}_{\eps,\delta}$ is uniformly bounded and $\eps\mathbb{Z}^3\times \eps^2\mathbb{Z}$ periodic with vanishing mean. Then, for $\kappa>0$ sufficiently small 
\begin{enumerate}
\item\label{Item1_sunset}
$\E\big[\big| \hat{\Pi}_\star^{(\eps,\delta)} \tau (\varphi^\lambda_{\star})\big|^2\big]\lesssim \lambda^{-\kappa}$ for all $\tau\in \mathfrak{T}_{\<22>}$,
\item\label{Item2_sunset}
$ \E\big[\big| \big(\hat{\Pi}_\star^{(\eps,\delta)} \<2> \bar{I} ( \pmb{f}\<2>)  - \hat{\Pi}_\star^{(0,\delta)} \<2> \bar{I} ( \pmb{f}\<2>)  \big)(\varphi^\lambda_{\star})\big|^2\big]\lesssim \eps^{\kappa} \lambda^{-2\kappa}$,
\item\label{Item3_sunset}
$ \E\big[\big| \hat{\Pi}_\star^{(\eps,\delta)} \tau (\varphi^\lambda_{\star})\big|^2\big]\lesssim \eps^{\kappa} \lambda^{-2\kappa}$  for $\tau \in \mathfrak{T}_{\<22>}\setminus \{\<2> \bar{I} ( \pmb{f}\<2>) \}$,
\end{enumerate}
uniformly over $\eps,\delta\geq 0$, $\varphi\in \mathfrak{B}_0$ and $\star\in \mathbb{R}^{3}$.
\end{lemma}
\begin{proof}
We use similar graphical notation to \cite{HQ18} find that for the only partially renormalised model $(\Pi^{(\eps,\delta)},\Gamma^{(\eps,\delta)})$ 
we can write
\begin{align*}
{\Pi}^{(\eps,\delta)}_\star \<2> I' ( \pmb{f}\<2>)(\phi_\star^\lambda) &= 
\begin{tikzpicture}[scale=0.25,baseline=0.25cm]
\node at (0,-2)  [root] (root) {};
	\node at (0,0)  [int] (int) {};
	\node at (0,2.5)  [bluevar] (middle) {};
	\node at (-1.6,2)  [var] (left) {};
	\node at (1.6,2)  [var] (right) {};
	\node at (-1.5,4.5)  [var] (tl) {};
	\node at (1.5,4.5)  [var] (tr) {};
	\draw[keps] (left) to  (int);	
	\draw[keps] (right) to (int);
	\draw[keps] (tl) to  (middle);	
	\draw[keps] (tr) to (middle);
	\draw[orange, kernel1] (middle) to (int);
	\draw[testfcn] (int) to  (root);
\end{tikzpicture}\; + 
4\;
\begin{tikzpicture}[scale=0.25,baseline=0.25cm]
\node at (0,-2)  [root] (root) {};
	\node at (0,-0)  [int] (int) {};
	\node at (0,2.5)  [bluevar] (middle) {};
	\node at (-1.6,1.25)  [int] (left) {};
	\node at (1.6,2)  [var] (right) {};
	\node at (1.5,4.5)  [var] (tr) {};
%
	\draw[keps,bend right=40] (left) to  (int);	
	\draw[keps] (right) to (int);
	\draw[keps,bend left=40] (left) to  (middle);	
	\draw[keps] (tr) to (middle);
	\draw[orange, kernel1] (middle) to (int);
	\draw[testfcn] (int) to  (root);
\end{tikzpicture}\; + 
\;
\begin{tikzpicture}[scale=0.25,baseline=0.25cm]
\node at (0,-2)  [root] (root) {};
	\node at (0,0)  [int] (int) {};
	\node at (0,2.5)  [bluevar] (middle) {};
	\node at (-1.6,2)  [var] (left) {};
	\node at (1.6,2)  [var] (right) {};
%
	\draw[keps] (left) to  (int);	
	\draw[keps] (right) to (int);
	\draw[orange, kernel1] (middle) to (int);
	\draw[testfcn] (int) to  (root);
\end{tikzpicture}\; + \;
\begin{tikzpicture}[scale=0.25,baseline=0.25cm]
\node at (0,-2)  [root] (root) {};
	\node at (0,0)  [bluevar] (int) {};
	\node at (1,2.6)  [bluevar] (middle) {};
	\node at (-0.5,4.5)  [var] (tl) {};
	\node at (1,5)  [var] (tm) {};
%
	\draw[orange, kernel1,bend left = 60] (middle) to (int);
	\draw[keps] (tl) to  (middle);	
	\draw[keps] (tm) to  (middle);	
	\draw[testfcn] (int) to  (root);
\end{tikzpicture}\;
 + 
\;
\begin{tikzpicture}[scale=0.25,baseline=0.25cm]
\node at (0,-2)  [root] (root) {};
	\node at (0,0)  [bluevar] (int) {};
	\node at (1,2.6)  [bluevar] (middle) {};
%
	\draw[orange, kernel1,bend left = 60] (middle) to (int);
%
	\draw[testfcn] (int) to  (root);
\end{tikzpicture}\; +  \;
2\; \begin{tikzpicture}[scale=0.25,baseline=0.25cm]
\node at (0,-2)  [root] (root) {};
	\node at (0,0)  [int] (int) {};
	\node at (0,2.5)  [bluevar] (middle) {};
	\node at (-1.6,1.25)  [int] (left) {};
	\node at (1.6,1.25)  [int] (right) {};
%
	\draw[keps,bend right=40] (left) to  (int);	
	\draw[keps,bend left=40] (left) to  (middle);	
	\draw[keps,bend left=40] (right) to  (int);	
	\draw[keps,bend right=40] (right) to  (middle);	
	\draw[orange, kernel1] (middle) to (int);
	\draw[testfcn] (int) to  (root);
\end{tikzpicture}\;,
\end{align*}
where
			\begin{itemize}
						\item the node \tikz[baseline=-3] \node [root] {}; represents the point $\star\in \mathbb{R}^{2+1}$, 
						 while the edge \tikz[baseline=-0.1cm] \draw[testfcn] (1,0) to (0,0); represents integration against the rescaled test function $\varphi^{\lambda}_\star$.
						\item The nodes \tikz[baseline=-3] \node [var] {}; represent the kernel variables in the Wiener Chaos representation.
						\item The node \tikz[baseline=-3] \node [dot] {}; represents dummy variables which are to be integrated out.
						\item The node \tikz[baseline=-3] \node [bluevar] {}; represents either the function $f^\eps$, $F_{\eps,\delta}^{\<2>}$ or $f^\eps F_{\eps,\delta}^{\<2>}$ evaluated at a dummy variable to be integrated out.
						\item  Edges \tikz[baseline=-0.1cm] \draw[keps] (0,0) to (1,0); represent integration against the kernel $\sum_{k\in \mathbb{Z}^3}\int_{y\in \mathbb{R}^{3}}K_{\eps}(x,t;y',s)\Gamma(y',s,y+k,s-\delta)dy$  where
$(s,y)$ and $(t,x)$ are the coordinates of the start and end
						points of the arrow respectively.
						\item Finally, 						
						\tikz[baseline=-0.1cm] \draw[orange, kernel1] (0,0) to (1,0);
						represents $K(x,t;y,s) - K(\star;y,s)$ 
						 where $K=G_\eps$ in the case $I'= I$	,		$K= \bar{K}_\eps$ if $I'=\bar{I}$,  $K(x,t,y,s)=K_\eps^{0,j}(x,t,y,s)\tilde\psi^\eps(y,s)$ if $I'=  I^{0,j}\tilde{\Phi}$, and 
						$K(x,t,y,s)=\psi^\eps_i(y,s)   \eps K_\eps^{i,j}(x,t,y,s)\tilde\psi^\eps_j(y,s) $ if $I'=  \mathcal{E} \big({\Phi}_i I^{i,j}\tilde{\Phi} (\cdot)\big)$.

						
					\end{itemize}
Using that
$$\E [\bfPi^{(\eps,\delta)} \<2> I' (\pmb{f}\<2>)]- \E [\bfPi^{(\eps,\delta)} \<2> (\star) ] \E [\bfPi^{(\eps,\delta)} I' (\pmb{f}\<2>) (\star) ]= 
2\; \begin{tikzpicture}[scale=0.25,baseline=0.1cm]
\node at (0,-2)  [root] (root) {};
	\node at (0,0)  [int] (int) {};
	\node at (0,2.5)  [bluevar] (middle) {};
	\node at (-1.6,1.25)  [int] (left) {};
	\node at (1.6,1.25)  [int] (right) {};
%
	\draw[keps,bend right=40] (left) to  (int);	
	\draw[keps,bend left=40] (left) to  (middle);	
	\draw[keps,bend left=40] (right) to  (int);	
	\draw[keps,bend right=40] (right) to  (middle);	
	\draw[orange, kernel] (middle) to (int);
	\draw[testfcn] (int) to  (root);
\end{tikzpicture}\;.  
$$
we thus can write
\begin{align}
\hat{\Pi}^{(\eps,\delta)}_\star \<2> I' ( \pmb{f}\<2>)(\phi_\star^\lambda) &= 
\begin{tikzpicture}[scale=0.25,baseline=0.25cm]
\node at (0,-2)  [root] (root) {};
	\node at (0,0)  [int] (int) {};
	\node at (0,2.5)  [bluevar] (middle) {};
	\node at (-1.6,2)  [var] (left) {};
	\node at (1.6,2)  [var] (right) {};
	\node at (-1.5,4.5)  [var] (tl) {};
	\node at (1.5,4.5)  [var] (tr) {};
	\draw[keps] (left) to  (int);	
	\draw[keps] (right) to (int);
	\draw[keps] (tl) to  (middle);	
	\draw[keps] (tr) to (middle);
	\draw[orange, kernel1] (middle) to (int);
	\draw[testfcn] (int) to  (root);
\end{tikzpicture}\; + 
4\;
\begin{tikzpicture}[scale=0.25,baseline=0.25cm]
\node at (0,-2)  [root] (root) {};
	\node at (0,-0)  [int] (int) {};
	\node at (0,2.5)  [bluevar] (middle) {};
	\node at (-1.6,1.25)  [int] (left) {};
	\node at (1.6,2)  [var] (right) {};
	\node at (1.5,4.5)  [var] (tr) {};
%
	\draw[keps,bend right=40] (left) to  (int);	
	\draw[keps] (right) to (int);
	\draw[keps,bend left=40] (left) to  (middle);	
	\draw[keps] (tr) to (middle);
	\draw[orange, kernel1] (middle) to (int);
	\draw[testfcn] (int) to  (root);
\end{tikzpicture}\; + 
\;
\begin{tikzpicture}[scale=0.25,baseline=0.25cm]
\node at (0,-2)  [root] (root) {};
	\node at (0,0)  [int] (int) {};
	\node at (0,2.5)  [bluevar] (middle) {};
	\node at (-1.6,2)  [var] (left) {};
	\node at (1.6,2)  [var] (right) {};
%
	\draw[keps] (left) to  (int);	
	\draw[keps] (right) to (int);
	\draw[orange, kernel1] (middle) to (int);
	\draw[testfcn] (int) to  (root);
\end{tikzpicture}\; + \;
\begin{tikzpicture}[scale=0.25,baseline=0.25cm]
\node at (0,-2)  [root] (root) {};
	\node at (0,0)  [bluevar] (int) {};
	\node at (1,2.6)  [bluevar] (middle) {};
	\node at (-0.5,4.5)  [var] (tl) {};
	\node at (1,5)  [var] (tm) {};
%
	\draw[orange, kernel1,bend left = 60] (middle) to (int);
	\draw[keps] (tl) to  (middle);	
	\draw[keps] (tm) to  (middle);	
	\draw[testfcn] (int) to  (root);
\end{tikzpicture}\;
 + 
\;
\begin{tikzpicture}[scale=0.25,baseline=0.25cm]
\node at (0,-2)  [root] (root) {};
	\node at (0,0)  [bluevar] (int) {};
	\node at (1,2.6)  [bluevar] (middle) {};
%
	\draw[orange, kernel1,bend left = 60] (middle) to (int);
%
	\draw[testfcn] (int) to  (root);
\end{tikzpicture}\; 
%
%
+2\; \begin{tikzpicture}[scale=0.25,baseline=0.25cm]
\node at (0,-2)  [root] (root) {};
	\node at (0,0)  [int] (int) {};
	\node at (0,2.5)  [bluevar] (middle) {};
	\node at (-1.6,1.25)  [int] (left) {};
	\node at (0,1.25)  [int] (right) {};
	\draw[keps,bend right=40] (left) to  (int);	
	\draw[keps,bend left=40] (left) to  (middle);	
	\draw[keps] (right) to  (int);	
	\draw[keps] (right) to  (middle);	
	\draw[orange, kernel, bend left= 90] (middle) to (root);
	\draw[testfcn] (int) to  (root);
\end{tikzpicture}\; 
+ \mathbf{1}_{I'= \bar{I}  }  \langle F^{\<22>}_{\eps,\delta}, \phi_\star^\lambda \rangle \ .
\label{eq:sunset_scematic}
\end{align}
We first consider the case when
 ${I'\in \{ \bar{I},  I^{0,j}(\tilde{\Psi}_{j}\, \cdot) ): j> 0 \}}$.
we obtain Item~\ref{Item1_sunset} and Item~\ref{Item3_sunset} by and using that all periodic functions appearing are all uniformly bounded and then arguing exactly as in \cite{Hai14} (see \cite{HS23m} for the corresponding lemmas involving singular kernels in the non-translation invariant setting). 

For Item~\ref{Item2_sunset} one argues similarly, but uses Lemma~\ref{lem:appendix} to gain the factor $\eps^\kappa$ when performing the integral corresponding to the top blue variable in the scematic representation \eqref{eq:sunset_scematic}.
Next we observe that we can argue exactly the same way in the case $I'=I$ 
since $\| G_{\eps}\|_{\beta,L,0}\lesssim \eps^{\kappa}$ for some $\beta<2$ and $L>1$,
which is sufficient for the argument just explained. 

Finally, the case
 ${I'\in \{ \mathcal{E}(\Psi_i I^{i,j}(\tilde{\Psi}_{j}\, \cdot) ) : i,j> 0 \}}$
 is an adaptation of the analogue case for g-PAM in the proof of {Lemma~\ref{lem:conv dumbel1}}.
\end{proof}
Since the proof of the next lemma is only simpler, we leave it to the reader.
\begin{lemma}
Under the assumption of Lemma~\ref{lem:stochastic estimate sunset} the analogue estimates hold for the trees belonging to $\mathfrak{T}^{\<31>}$.
\end{lemma}

\begin{lemma}
Assume that $F^{\<22>}_{\eps,\delta}$ is uniformly bounded and $\eps\mathbb{Z}^3\times \eps^2\mathbb{Z}$ periodic with vanishing mean. Then, for $\kappa>0$ sufficiently small 
\begin{enumerate}
\item 
$\E\big[\big| \hat{\Pi}_\star^{(\eps,\delta)} \tau (\varphi^\lambda_{\star})\big|^2\big]\lesssim \lambda^{-1/2-\kappa}$ for all $\tau\in \mathfrak{T}_{\<32>}$,
\item
$ \E\big[\big| \big(\hat{\Pi}_\star^{(\eps,\delta)} \<2> \bar{I} ( \pmb{f}\<3>)  - \hat{\Pi}_\star^{(0,\delta)} \<2> \bar{I} ( \pmb{f}\<3>)  \big)(\varphi^\lambda_{\star})\big|^2\big]\lesssim \eps^{\kappa} \lambda^{-1/2 -2\kappa}$,
\item
$ \E\big[\big| \hat{\Pi}_\star^{(\eps,\delta)} \tau (\varphi^\lambda_{\star})\big|^2\big]\lesssim \eps^{\kappa} \lambda^{-1/2-2\kappa}$  for $\tau \in \mathfrak{T}_{\<23>}\setminus \{\<2> \bar{I} ( \pmb{f}\<3>) \}$,
\end{enumerate}
uniformly over $\eps,\delta\geq 0$, $\varphi\in \mathfrak{B}_0$ and $\star\in \mathbb{R}^{3}$.
\end{lemma}
\begin{proof}
From the definition one finds that using the graphical notation of Lemma~\ref{lem:stochastic estimate sunset} we can write
\begin{align}
\hat{\Pi}^{(\eps,\delta)}_\star \<2> I' ( \pmb{f}\<3>)(\phi_\star^\lambda) &= 
\begin{tikzpicture}[scale=0.25,baseline=0.25cm]
\node at (0,-2)  [root] (root) {};
	\node at (0,0)  [int] (int) {};
	\node at (0,2.5)  [bluevar] (middle) {};
	\node at (-1.6,2)  [var] (left) {};
	\node at (1.6,2)  [var] (right) {};
	\node at (-1.5,4.5)  [var] (tl) {};
	\node at (1.5,4.5)  [var] (tr) {};
	\node at (0,5)  [var] (tm) {};	
	\draw[keps] (left) to  (int);	
	\draw[keps] (right) to (int);
	\draw[keps] (tl) to  (middle);	
	\draw[keps] (tr) to (middle);
	\draw[keps] (tm) to  (middle);	
	\draw[orange,  kernel1] (middle) to (int);
	\draw[testfcn] (int) to  (root);
\end{tikzpicture}\; + 
6\;
\begin{tikzpicture}[scale=0.25,baseline=0.25cm]
\node at (0,-2)  [root] (root) {};
	\node at (0,-0)  [int] (int) {};
	\node at (0,2.5)  [bluevar] (middle) {};
	\node at (-1.6,1.25)  [int] (left) {};
	\node at (1.6,2)  [var] (right) {};
	\node at (1.5,4.5)  [var] (tr) {};
	\node at (0,6)  [var] (tm) {};
	\draw[keps,bend right=40] (left) to  (int);	
	\draw[keps] (right) to (int);
	\draw[keps,bend left=40] (left) to  (middle);	
	\draw[keps] (tr) to (middle);
	\draw[keps] (tm) to  (middle);	
	\draw[orange,  kernel1] (middle) to (int);
	\draw[testfcn] (int) to  (root);
\end{tikzpicture}\; + 
3\;
\begin{tikzpicture}[scale=0.25,baseline=0.25cm]
\node at (0,-2)  [root] (root) {};
	\node at (0,0)  [int] (int) {};
	\node at (0,2.5)  [bluevar] (middle) {};
	\node at (-1.6,2)  [var] (left) {};
	\node at (1.6,2)  [var] (right) {};
	\node at (1,4.6)  [var] (tr) {};
	\draw[keps] (left) to  (int);	
	\draw[keps] (right) to (int);
	\draw[keps,bend left = 60] (tr) to (middle);
	\draw[orange,  kernel1] (middle) to (int);
	\draw[testfcn] (int) to  (root);
\end{tikzpicture}\; + \;
\begin{tikzpicture}[scale=0.25,baseline=0.25cm]
\node at (0,-2)  [root] (root) {};
	\node at (0,0)  [int] (int) {};
	\node at (1,2.6)  [bluevar] (middle) {};
	\node at (-0.5,4.5)  [var] (tl) {};
	\node at (2.5,4.5)  [var] (tr) {};
	\node at (1,5)  [var] (tm) {};
%
	\draw[orange, kernel1,bend left = 60] (middle) to (int);
	\draw[keps] (tl) to  (middle);	
	\draw[keps] (tr) to (middle);
	\draw[keps] (tm) to  (middle);	
	\draw[testfcn] (int) to  (root);
\end{tikzpicture}\; + 
3\;
\begin{tikzpicture}[scale=0.25,baseline=0.25cm]
\node at (0,-2)  [root] (root) {};
	\node at (0,0)  [int] (int) {};
	\node at (1,2.6)  [bluevar] (middle) {};
	\node at (2,4.6)  [var] (tr) {};
%
	\draw[orange, kernel1,bend left = 60] (middle) to (int);
%
	\draw[keps,bend left = 60] (tr) to (middle);
	\draw[testfcn] (int) to  (root);
\end{tikzpicture}\; + 6\;
\begin{tikzpicture}[scale=0.25,baseline=0.25cm]
\node at (0,-2)  [root] (root) {};
	\node at (0,0)  [int] (int) {};
	\node at (0,2.5)  [bluevar] (middle) {};
	\node at (-1.6,1.25)  [int] (left) {};
	\node at (1.6,2)  [var] (right) {};
%
	\draw[keps,bend right=40] (left) to  (int);	
	\draw[keps] (right) to (int);
	\draw[keps,bend left=40] (left) to  (middle);	
	\draw[orange,  kernel1] (middle) to (int);
	\draw[testfcn] (int) to  (root);
\end{tikzpicture}\;\\
& + 3\; \left( 2\; 
\begin{tikzpicture}[scale=0.25,baseline=0.25cm]
\node at (0,-2)  [root] (root) {};
	\node at (0,0)  [int] (int) {};
	\node at (0,2.5)  [bluevar] (middle) {};
	\node at (-1.6,1.25)  [int] (left) {};
	\node at (1.6,1.25)  [int] (right) {};
	\node at (0,5)  [var] (tm) {};
	\draw[keps,bend right=40] (left) to  (int);	
	\draw[keps,bend left=40] (left) to  (middle);	
	\draw[keps,bend left=40] (right) to  (int);	
	\draw[keps,bend right=40] (right) to  (middle);	
	\draw[keps] (tm) to  (middle);	
	\draw[orange,  kernel] (middle) to (int);
	\draw[testfcn] (int) to  (root);
\end{tikzpicture}\;
- 
\langle g^{\<2> I' ( \pmb{f} \<2>)}_{\eps,\delta} \Pi^{(\eps,\delta)}
 \<1> , \phi^\lambda_{\star} \rangle
\right)
+ 6\;
\begin{tikzpicture}[scale=0.25,baseline=0.25cm]
\node at (0,-2)  [root] (root) {};
	\node at (0,0)  [int] (int) {};
	\node at (0,2.5)  [bluevar] (middle) {};
	\node at (-1.6,1.25)  [int] (left) {};
	\node at (0,1.25)  [int] (right) {};
	\node at (0,5)  [var] (tm) {};
%
	\draw[keps,bend right=40] (left) to  (int);	
	\draw[keps,bend left=40] (left) to  (middle);	
	\draw[keps] (right) to  (int);	
	\draw[keps] (right) to  (middle);	
	\draw[keps] (tm) to  (middle);	
	\draw[orange, kernel, bend left= 90] (middle) to (root);
	\draw[testfcn] (int) to  (root);
\end{tikzpicture}\;. \nonumber
\end{align}	
The only term for which it may not be obvious that the arguments in \cite{Hai14} (see again \cite[Sec.~17]{HS23m} for non-translation invariant variants of the estimates),
adapt directly in the analogous manner to the estimates in Section~\ref{eq: convergence gpam} is the one in brackets above. 
That term can be written as 
$$
\left( 2\; 
\begin{tikzpicture}[scale=0.25,baseline=0.25cm]
\node at (0,-2)  [root] (root) {};
	\node at (0,0)  [int] (int) {};
	\node at (0,2.5)  [bluevar] (middle) {};
	\node at (-1.6,1.25)  [int] (left) {};
	\node at (1.6,1.25)  [int] (right) {};
	\node at (0,5)  [var] (tm) {};
	\draw[keps,bend right=40] (left) to  (int);	
	\draw[keps,bend left=40] (left) to  (middle);	
	\draw[keps,bend left=40] (right) to  (int);	
	\draw[keps,bend right=40] (right) to  (middle);	
	\draw[keps] (tm) to  (middle);	
	\draw[orange,  kernel] (middle) to (int);
	\draw[testfcn] (int) to  (root);
\end{tikzpicture}\;
-
\langle g^{\<2> I' ( \pmb{f} \<2>)}_{\eps,\delta} \Pi^{(\eps,\delta)}
 \<1> , \phi^\lambda_{\star} \rangle 
\right) 
=
 2 \begin{tikzpicture}[scale=0.25,baseline=0.25cm]
\node at (0,-2)  [root] (root) {};
	\node at (0,0)  [int] (int) {};
	\node at (0,2.5)  [dot] (middle) {};
	\node at (0,5)  [var] (tm) {};
%
	\draw[keps] (tm) to  (middle);	
	\draw[orange,  kernelBig] (middle) to (int);
	\draw[testfcn] (int) to  (root);
\end{tikzpicture}\;
+ \mathbf{1}_{I'=\bar{I}} F^{\<22>},
$$
where \tikz[baseline=-0.1cm] \draw[orange, kernelBig] (1,0) to (0,0); represents the renormalised kernel $\mathcal{R}(K\cdot f^\eps)$, see \eqref{eq:app_ren},
which acts on Hölder functions as $\mathcal{R}(Kf^\eps)(F)\eqdef\int K(x,y)f^\eps(y) [F(x)-F(y)]dy$
with $K=G_\eps$ in the case $I'= I$,	$K= \bar{K}_\eps$ if $I'=\bar{I}$,  $K(x,t,y,s)=K_\eps^{0,j}(x,t,y,s)\tilde\psi^\eps(y,s)$ if $I'=  I^{0,j}\tilde{\Phi}$, and 
$K(x,t,y,s)=\psi^\eps_i(y,s)   \eps K_\eps^{i,j}(x,t,y,s)\tilde\psi^\eps_j(y,s) $ if $I'=  \mathcal{E} \big({\Phi}_i I^{i,j}\tilde{\Phi} (\cdot)\big)$.
Thus the desired bounds follow using Lemma~\ref{lem:ap_osc_kernel}.
%
\end{proof}
For ${I'\in \{ I, \bar{I},  I^{0,j}(\tilde{\Psi}_{j}\, \cdot) ),  \mathcal{E}(\Psi_i I^{i,j}(\tilde{\Psi}_{j}\, \cdot) ) : i,j> 0 \}}$ define 
$$
\gamma_{\eps,\delta}^{\<2> I' (\pmb{f}\<2>)}=
\begin{cases}
\int_{[0,1]^3} f^{\eps}_{\mu,\nu} F^{\tau}_{\eps,\delta} dxdt  & \text{if } I'= \bar{I}  , \\
0 & \text{else. }
\end{cases}
$$
and 
$$\|\pmb{f}\<2> I' ( \pmb{f}\<2>)\|=0\ ,\qquad \|\pmb{f}\<2> I' ( \pmb{f}\<3>)\|=-1/2\ . $$
Finally, one further obtains the following estimates along similar lines as the proof of Lemma~\ref{lem:conv_resonances}.
\begin{lemma}\label{lem:conv dumbel2}
Assume that $F_{\eps,\delta}$ is uniformly bounded and $\eps\mathbb{Z}^2\times \eps^2\mathbb{Z}$ periodic with vanishing mean. Then, for $\kappa>0$ sufficiently small 
\begin{enumerate}
\item 
$\E\big[\big| \hat{\Pi}_\star^{(\eps,\delta)} \tau (\varphi^\lambda_{\star})\big|^2\big]\lesssim \lambda^{\|\tau\|-\kappa}$ for all $\tau\in \mathfrak{T}^{\pmb{f}\<2>  I(\pmb{f}\<2>)}\cup \mathfrak{T}^{\pmb{f}\<2>  I(\pmb{f}\<3>)}$,
\item
$ \E\big[\big| \big(\hat{\Pi}_\star^{(\eps,\delta)} \pmb{f}\<2>  \bar{I} ( \pmb{f}\<2>) - \hat{\Pi}_\star^{(0,\delta)} \pmb{f}\<2>  \bar{I} ( \pmb{f}\<2>) \big)(\varphi^\lambda_{\star})\big|^2\big]\lesssim \eps^{\kappa} \lambda^{-2\kappa}$,
\item
$ \E\big[\big| \big(\hat{\Pi}_\star^{(\eps,\delta)} \pmb{f}\<2>  \bar{I} ( \pmb{f}\<3>) - \hat{\Pi}_\star^{(0,\delta)} \pmb{f}\<3>  \bar{I} ( \pmb{f}\<2>) \big)(\varphi^\lambda_{\star})\big|^2\big]\lesssim \eps^{\kappa} \lambda^{-1/2-2\kappa}$,
\item
$ \E\big[\big| \hat{\Pi}_\star^{(\eps,\delta)} \tau (\varphi^\lambda_{\star})\big|^2\big]\lesssim \eps^{\kappa} \lambda^{\|\tau\| -2\kappa}$  for all $\tau\in \big( \mathfrak{T}^{\pmb{f}\<2>  I(\pmb{f}\<2>)} \cup \mathfrak{T}^{\pmb{f}\<2>  I(\pmb{f}\<3>)}\big) \setminus \{\pmb{f}\<2>  \bar{I} ( \pmb{f}\<3>) ,\ \pmb{f}\<2>  \bar{I} ( \pmb{f}\<2>) \} $,
\end{enumerate}
uniformly over $\eps,\delta\geq 0$, $\varphi\in \mathfrak{B}_0$ and $\star\in \mathbb{R}^{3}$.
\end{lemma}

%

\subsection{Identification of divergences and proof of the main results}\label{sec:identificationPhi4}
\subsubsection{Counterterms for Theorem~\ref{thm:phi4}}\label{sec:count_phi_full}

Next, we define $$F^{\<2>}_{\eps,\delta}(z)=  \E[(\Pi^{\eps,\delta}\<1>)^2] -
	\frac{\alpha^{\<2>}_{\eps,\delta}}{\sqrt{\det(A_s^\eps)}} 
-\frac{1}{\eps} R^{\<2>}_{(\delta/\eps)\vee 1} \circ \mathcal{S}^\eps  
- c^{\<2>}_{\eps,\delta} \ ,$$
where (for $d=3$)
\begin{enumerate}
\item ${\alpha}^{\<2>}_{\eps,\delta}\eqdef(4\pi)^{-d/2}\frac{1}{\eps}\int_{(\delta/\eps)^2}^1 s^{-d/2} ds $
\item  $R^{\<2>}_\lambda(z)\eqdef \int_{\mathbb{R}^3\times[\lambda^2,\infty)}\Gamma_1^2(z;y,s) dyds $
for $\lambda\geq 1$.
\item $ c^{\<2>}_{\eps,\delta}\eqdef\int_{[0,1]^3} h^{\<2>}_{\eps,\delta}$ for $h^{\<2>}_{\eps,\delta}\eqdef  \E[(\Pi^{\eps,\delta}\<1>)^2] - \det(A_s^\eps)^{-1/2}  \alpha^{\<2>}_{\eps,\delta}  -\frac{1}{\eps} R^{\<2>}_{(\delta/\eps)\vee 1} \circ \mathcal{S}^\eps  $.
\end{enumerate}

\begin{lemma}
The functions $F^{\<2>}_{\eps,\delta}: \mathbb{R}^3\to \mathbb{R}$ satisfy the following properties:
\begin{itemize}
\item $\sup_{\eps,\delta\in (0,1]} \|F^{\<2>}_{\eps,\delta}\|_{L^\infty}<\infty$,
\item $F^{\<2>}_{\eps,\delta}$ is $\eps\mathbb{Z}^2\times \eps^2 \mathbb{Z}$ periodic and has mean $0$,
\item  the limits $F^{\<2>}_{\eps,0}=\lim_{\delta\to 0} F^{\<2>}_{\eps,\delta}$ exist as pointwise limits for any $\eps>0$,
\item $\lim_{\delta\to 0} \|F^{\<2>}_{\eps,\delta}\|_{L^\infty}=0$ for each $\delta>0$.
\end{itemize}
\end{lemma}
\begin{proof}
By \eqref{eq:linear_guy}
\begin{equs}
\E[(\Pi^{\eps,\delta}\<1>(x,t))^2]&=  \int_{\mathbb{R}^3\times \mathbb{R}} \kappa^2(t-s) \Gamma^2_\eps(x,t,y,s-\delta^2) dyds + r^{\<2>}_{\eps,\delta}(x,t)
\end{equs} 
where 
$$r^{\<2>}_{\eps,\delta}(x,t)\eqdef \sum_{k\in \mathbb{Z}^3\setminus 0} \int_{\mathbb{R}^3\times \mathbb{R}} \kappa^2(t-s) \Gamma_\eps(x,t,y,s-\delta^2) \Gamma_\eps(x,t,y+k,s-\delta^2) dyds\ ,$$
which extends to a continuous map $[0,1]^{2}\to L^{\infty}\mathbb{R}^{1+3}$.
Next, using that $\kappa(t)=1$ for $t\in (0,1) $
\begin{align*}
&\int_{\mathbb{R}^3\times \mathbb{R}} \kappa^2(t-s) \Gamma^2_\eps(x,t,y,s-\delta^2) dyds =
\int_{\mathbb{R}^3\times \mathbb{R}} \kappa^2(s-\delta^2) \Gamma^2_\eps(x,t,y,t-s) dyds &\\
&=\mathbf{1}_{\delta<\eps} \int_{\delta^2}^{\eps^2} \Gamma^2_\eps(x,t,y,t-s) dyds 
+ \int_{\eps^2\vee \delta^2}^\infty \kappa^2(s-\delta^2) \Gamma^2_\eps(x,t,y,t-s) dyds 
\end{align*}
Thus, we write $h^{\<2>}_{\eps,\delta}= h^1_{\eps,\delta}+h^2_{\eps,\delta} + r^{\<2>}_{\eps,\delta}$ where
\begin{align*}
h_{\eps,\delta}^1(x,t)&\eqdef \mathbf{1}_{\delta<\eps} \int_{\mathbb{R}^3}\int_{\delta^2}^{\eps^2} \Gamma^2_\eps(x,t,y,t-s) dyds - \det(A_s^\eps)^{-1/2}  \alpha^{\<2>}_{\eps,\delta} \\
h_{\eps,\delta}^2(x,t)&\eqdef  \int_{\eps^2\vee \delta^2}^\infty \int_{\mathbb{R}^3} \kappa^2(s-\delta^2) \Gamma^2_\eps(x,t,y,t-s) dyds -\frac{1}{\eps} R_{(\delta/\eps)\vee 1} \circ \mathcal{S}^\eps  \ .
\end{align*}
We find that 
\begin{equs}
h^1_{\eps,\delta}(\eps x,\eps^2 t)
&=\mathbf{1}_{\delta<\eps} \frac{1}{\eps} \int_{(\delta/\eps)^2}^{1} \int_{\mathbb{R}^3} \big[\Gamma^2_1 (x, t,y,t-s) - (Z^*)^2(x,t; y,t-s)\big] dyds \ , \label{eq:cherry rhs local to bound.}
\end{equs}
where we used that 
$$
\det(A_s)^{-1/2}  \alpha^{\<2>}_{\eps,\delta}=\mathbf{1}_{\delta<\eps} \frac{1}{\eps} \int_{(\delta/\eps)^2}^{1} (Z^*)^2(x,t; y,t-s)dyds \ .
$$
Then \eqref{eq:cherry rhs local to bound.} is bounded, and is seen to converge for $\eps>0$ as $\delta\to 0$ exactly as in \cite[Sec.~3.1.3]{Sin23}. 
We also directly read off the expression that it vanishes for $\eps<\delta$.
Next, for
\begin{align*}
h^2(\eps x,\eps t)
&= \frac{1}{\eps} \int_{1\vee (\delta/\eps)^2}^\infty \int_{\mathbb{R}^3} \kappa^2(\eps^2 s-\delta^2) \Gamma^2( x, t,y,t-s) dyds -\frac{1}{\eps} 
\int_{\mathbb{R}^3\times[(\delta/\eps)^2\vee 1,\infty)}\Gamma_1^2(x,t,y,s) dyds
 \\
 &= \frac{1}{\eps} \int_{\mathbb{R}^3\times[(\delta/\eps)^2\vee 1,\infty)} \big[ \kappa^2(\eps^2 s-\delta^2) -1] \Gamma^2( x, t,y,t-s) dyds 
\end{align*}
note that $\kappa^2(s/\eps^2-\delta^2) -1= 0$ unless $\eps^2 s -\delta^2>1$ which is contained in $s> \eps^{-2}$ and thus 
$$|h^2(\eps x,\eps t)|\leq \frac{1}{\eps} \int_{\mathbb{R}^3\times[\eps^{-2}, \infty)}  \Gamma^2( x, t,y,t-s) dyds 
\leq  \int_{\mathbb{R}^3\times[1, \infty)}  \Gamma_\eps^2( x, t,y,t-s) dyds 
$$
from which we read off uniform boundedness and convergence as $\delta\to 0$ for $\eps>0$. Finally, since 
$$h_{\eps,\delta}^2(x,t)\eqdef  \int_{\eps^2\vee \delta^2}^\infty \int_{\mathbb{R}^3} \Big[\kappa^2(s-\delta^2) \Gamma^2_\eps(x,t,y,t-s)  - \bar{\Gamma}^2(x,t,y,t-s)\Big]  dyds \ ,$$
we see that $ \|h_{\eps,\delta}^2\|_{L^\infty}\to 0$ for fixed $\delta>0$.

\end{proof}
Next define $$F^{\<22>}_{\eps,\delta}(z)= \sum_{
I'}
\E [\bfPi^{(\eps,\delta)} \<2> I' (\pmb{f}\<2>)(z)]- \E [\bfPi^{(\eps,\delta)} \<2> (z) ] \E [\bfPi^{(\eps,\delta)} I' (\pmb{f}\<2>) (z) ] - \frac{  2\alpha^{\<22>}_{\eps,\delta}f^\eps }{\det(A_s^\eps(z))}-   \frac{2\bar{\alpha}^{\<22>}_{\eps,\delta} \bar{f} }{\det(\bar{A})}- c^{\<22>}_{\eps,\delta} \ ,$$
where the sum over runs over $I'\in \{ I, \bar{I},  I^{0,j}(\tilde{\Psi}_{j}\, \cdot) ),  \mathcal{E}(\Psi_i I^{i,j}(\tilde{\Psi}_{j}\, \cdot) ) : i,j> 0 \}$ and the constants are defined as follows, denoting by $H_{t}(x)$ the fundamental solution of the operator $\partial_{t} - \triangle$.
\begin{enumerate}
\item 
Set 
$ \alpha^{\<22>}_{\eps,\delta} = \int_{\mathbb{R}^{3+1}}  \kappa(2 s) H_\tau (y) \big(\tilde{H}_{1; \delta/\eps}( s, y)\big)^2 dy  ds \ $
where
$$
\tilde{H}_{1;\delta/\epsilon} (t;x)
\eqdef \int_{ 
\substack{
(\delta/\eps)^2<t-s<1\\
(\delta/\eps)^2<s<1
 }}
 \int_{\mathbb{R}^3} H_{t-s} (x-y) 
 H_{s} (y) 
\,dy  ds \ .
$$

\item\label{itemmmm2}  
$\bar{\alpha}^{\<22>}_{\eps,\delta}\eqdef\int_{\mathbb{R}^{3+1}} \kappa(s)
(1-\kappa(2\eps^{-2}s)) H_{s}(y)\left(\tilde{{H}}^{>_\eps}_{0;\delta}(y,s) \right)^2
 \ dy  ds $
where 
$$\tilde{{H}}^{>_\eps}_{0;\delta}(x,t)= 
 \int_{\mathbb{R}^3}\int_{
\substack{ s>\eps^2\\
 t-s>\eps^2}} H_{t-s} (x-y) 
 H_s (y)  \kappa(t-\tau-\delta^2) \kappa(-\tau-\delta^2) dyds$$

\item $ c^{\<22>}_{\eps,\delta}\eqdef\int_{[0,1]^3} h^{\<22>}_{\eps,\delta}$ for 
$$h^{\<22>}_{\eps,\delta}\eqdef\sum_{I'} \E [\bfPi^{(\eps,\delta)} \<2> I' (\pmb{f}\<2>)]- \E [\bfPi^{(\eps,\delta)} \<2> ] \E [\bfPi^{(\eps,\delta)} I' (\pmb{f}\<2>) 
](z) - \frac{ 2 \alpha^{\<22>}_{\eps,\delta}f^\eps }{\det(A_s^\eps(z))}-   \frac{2 \bar{\alpha}^{\<22>}_{\eps,\delta} \bar{f} }{\det(\bar{A})} \ . $$
\end{enumerate}

\begin{lemma}\label{lem:finite wiener chaos remainder12}
The functions $F^{\<22>}_{\eps,\delta}: \mathbb{R}^3\to \mathbb{R}$ satisfy the following properties:
\begin{itemize}
\item $\sup_{\eps,\delta\in (0,1]} \|F^{\<22>}_{\eps,\delta}\|_{L^\infty}<\infty$,
\item $F^{\<22>}_{\eps,\delta}$ is $\eps\mathbb{Z}^2\times \eps^2 \mathbb{Z}$ periodic and has mean $0$,
\item  the limits $F^{\<22>}_{\eps,0}=\lim_{\delta\to 0} F^{\<22>}_{\eps,\delta}$ exist as pointwise limits for any $\eps>0$,
\item $\lim_{\eps\to  0} \|F^{\<22>}_{\eps,\delta}\|_{L^\infty}=0$  for each $\delta>0$. 
\end{itemize}
\end{lemma}
\begin{proof}
Unraveling the definition one finds that 
\begin{equs}\label{eq:decomp of sunset div}
 \sum_{
I'}
 &\E [\bfPi^{(\eps,\delta)} \<2> I' (\pmb{f}\<2>)(x,t)]- \E [\bfPi^{(\eps,\delta)} \<2> (x,t) ]\E [\bfPi^{(\eps,\delta)} I' (\pmb{f}\<2>)  
(x,t)]\\
&= 2
\int_{\mathbb{R}^{3+1}} \kappa(t-s) \Gamma_\eps (x,t; y,s) f^\eps(y,s)\left(\tilde{\Gamma}_{\eps;\delta}(x,t; y,s)+\tilde{r}_{\eps,\delta}(x,t; y,s)\right)^2 \ dy  ds
\end{equs}
where
\begin{align}
\tilde{\Gamma}_{\epsilon;\delta}(x,t, \bar{x},\bar{t})&= 
 \int_{\mathbb{R}^3\times \mathbb{R}} \Gamma_\eps (x,t;\eta, \tau-\delta^2) 
 \Gamma_\eps ( \bar{x},\bar{t};\eta, \tau-\delta^2) 
\,d\eta \ \kappa(t-\tau) \kappa(\bar{t}-\tau) d\tau \label{eq:tilde_kernel}
\\
\tilde{r}_{\eps,\delta}(x,t, \bar{x},\bar{t})&=
 \sum_{k\in \mathbb{Z}^3\setminus\{0\}}
\int_{\mathbb{R}^3\times \mathbb{R}} \Gamma_\eps (x,t;\eta, \tau-\delta^2) 
 \Gamma_\eps ( \bar{x},\bar{t};\eta+k, \tau-\delta^2) 
\,d\eta \ \kappa^2(t-\tau) d\tau \ .\nonumber 
\end{align}
Next we decompose into $\Gamma_\eps^{<}(x,t,y,s)\eqdef\mathbf{1}_{t-s<\eps^2}\Gamma_\eps(x,t,y,s)$ and $\Gamma_\eps^{>}(x,t,y,s)\eqdef\mathbf{1}_{t-s>\eps^2}\Gamma_\eps(x,t,y,s)$ and write
\begin{equ}\label{eq:rhs_tilde_kernels}
\tilde{\Gamma}_{\epsilon;\delta}=\tilde{\Gamma}_{\epsilon;\delta}^{>}+ \tilde{\Gamma}_{\epsilon;\delta}^{<} + 2\tilde{\Gamma}_{\epsilon;\delta}^{\neq}
\end{equ}
where the kernels on the right hand side of \eqref{eq:rhs_tilde_kernels} are given by replacing in \eqref{eq:tilde_kernel} the instances of $\Gamma_\eps$ by $\tilde{\Gamma}_{\epsilon;\delta}^{<}$ or $\tilde{\Gamma}_{\epsilon;\delta}^{>}$ or one of each respectively.
Similarly, set 
$\Gamma_\eps^{\lesssim}= \kappa(2\eps^{-2}(t-s) ) \Gamma_\eps(x,t,y,s)$ and $\Gamma^{\gtrsim}=\big(1 -\kappa(2\eps^{-2}(t-s) ) \big)\Gamma_\eps(x,t,y,s) $.
Accordingly we write 
$$
h^{\<22>}_{\eps,\delta}=  2 \sum_{ \substack{  a\in \{\lesssim,\gtrsim\},\\
 b,c\in \{<,>,\neq\}}} h_{\eps,\delta}^{abc} + \text{`terms involving an occurrence of $\tilde{r}_{\eps,\delta}$'}
 $$
where
\begin{align*}
h_{\eps,\delta}^{\lesssim <<}(x,t) &\eqdef
\int_{\mathbb{R}^{3+1}} \kappa(t-s) \Gamma^{\lesssim}_\eps (x,t; y,s) f^\eps(y,s)\left(\tilde{\Gamma}^{<}_{\eps;\delta}(x,t; y,s)\right)^2 \ dy  ds 
- \frac{  \alpha^{\<22>}_{\eps,\delta}f^\eps }{\det(A_s^\eps(z))}\\
h_{\eps,\delta}^{\gtrsim >>}(x,t)&\eqdef\int_{\mathbb{R}^{3+1}} \kappa(t-s) \Gamma^{\gtrsim}_\eps (x,t; y,s) f^\eps(y,s)\left(\tilde{\Gamma}^>_{\eps;\delta}(x,t; y,s) \right)^2 \ dy  ds-\frac{\bar{\alpha}^{\<22>}_{\eps,\delta} f^\eps }{{\det(\bar{A})}}
\end{align*}
and in the remaining cases
$$
h_{\eps,\delta}^{abc}(x,t) \eqdef
\int_{\mathbb{R}^{3+1}} \kappa(t-s) \Gamma^{a}_\eps (x,t; y,s) f^\eps(y,s)\tilde{\Gamma}^{b}_{\eps;\delta}(x,t; y,s)\tilde{\Gamma}^{c}_{\eps;\delta}(x,t; y,s) \ dy  ds \ .
$$
One straightforwardly checks that the terms involving $\tilde{r}_{\eps,\delta}$ contributing to
\eqref{eq:decomp of sunset div} extend to a continuous map $(0,1]^2 \to L^{\infty}(\mathbb{R}^{1+3})$. Therefore we focus on the terms in the left sum therein.

\paragraph{Estimate on \TitleEquation{h_{\eps,\delta}^{\lesssim <<}}{h}:}
 We first note using that $\kappa(s)=1$ for $s\in (0,2)$  
$$\tilde{\Gamma}_{\epsilon;\delta}^{<}(x,t; \bar{x},\bar{t})= \int_{
\substack{ 
{\delta^2<t-s<\eps^2}
\\
{\delta^2<\bar{t}-s<\eps^2}
}} \int_{\mathbb{R}^3} \Gamma^<_\eps (x,t;\eta, s) 
 \Gamma^<_\eps ( \bar{x},\bar{t};\eta, s) 
\, d\eta ds  $$
Therefore 
\begin{align*}
\tilde{\Gamma}_{\epsilon;\delta}^{<}(\eps x,\eps^2 t; \eps \bar{x},\eps^2 \bar{t})
&=\eps^{d+2- 2d}
 \int_{
\substack{ 
 {(\delta/\eps)^2<t-s<1}
\\
 {(\delta/\eps)^2<\bar{t}-s<1}
}} \int_{\mathbb{R}^3} \Gamma_1 (x,t;\eta, s) 
 \Gamma_1 ( \bar{x},\bar{t};\eta, s) 
\, d\eta ds  \\
&=: \eps^{d+2- 2d}\tilde{\Gamma}_{\epsilon/\delta ; 1}( x, t;  \bar{x}, \bar{t})= \eps^{-1}\tilde{\Gamma}_{\epsilon/\delta ; 1}( x, t;  \bar{x}, \bar{t})
\end{align*}
Thus we write
\begin{align*}
h_{\eps,\delta}^{\lesssim < <}(\eps x,\eps^2 t) &=
\int_{\mathbb{R}^{3+1}} \Gamma^{\lesssim}_\eps (\eps x,\eps^2 t; y,s) f^\eps(y,s)\left(\tilde{\Gamma}^{<}_{\eps;\delta}(\eps x,\eps^2 t; y,s)\right)^2 \ dy  ds 
- \frac{  \alpha^{\<22>}_{\eps,\delta}f}{{\det(A_s(z))}}\\
&= 
\int_{\mathbb{R}^{3+1}}  \Gamma^{\lesssim}_1 ( x, t; y, s) f(y,s)\left(\tilde{\Gamma}^<_{1; \delta/\eps}( x, t; y, s)\right)^2 \ dy  ds 
- \frac{  \alpha^{\<22>}_{\eps,\delta}f}{\det(A_s(z))}\\
\end{align*}
Thus recalling that 
$$\frac{  \alpha^{\<22>}_{\eps,\delta}}{\det(A_s(z))} = \int_{\mathbb{R}^{3+1}}  \kappa(2 (t- s)) Z^* ( x, t; y, s) f(y,s)\left(\tilde{Z}_{1; \delta/\eps}( x, t; y, s)\right)^2 \ dy  ds \ ,$$
where 
\begin{align*}
&\tilde{Z}^*_{1;\delta/\epsilon} (x,t, \bar{x},\bar{t})
\eqdef \int_{ 
\substack{
(\delta/\eps)^2<t-s<1\\
(\delta/\eps)^2<\bar t -s<1
 }}
 \int_{\mathbb{R}^d} Z^{*} (\eta, \tau;x,t;) 
 Z^* ( \eta, \tau';\bar{x},\bar{t}) 
\,d\eta  d\tau \ d\tau'\ .
\end{align*}
We thus find that 
\begin{align*}
&h_{\eps,\delta}^{\lesssim<<}(\eps x,\eps^2 t) \\
&= \int_{\mathbb{R}^{3+1}} f(y,s)\Big[ \Gamma^{\lesssim}_1 ( x, t; y, s) \left(\tilde{\Gamma}_{1; \delta/\eps}( x, t; y, s)\right)^2 
- \kappa(2(t-s))Z^* ( x, t; y, s) f(y,s)\left(\tilde{Z}_{1; \delta/\eps}( x, t; y, s)\right)^2\Big]
 dy  ds \\
& +  \int_{\mathbb{R}^{3+1}} [f(y,s) - f(x,t)] \kappa(2(t-s)) Z^* ( x, t; y, s) f(y,s)\left(\tilde{Z}_{1; \delta/\eps}( x, t; y, s)\right)^2  dy  ds \ ,
 \end{align*}
where we can estimate the first term as in \cite{Sin23} and the latter term by using the regularity of $f$.
\paragraph{Estimate on \TitleEquation{h_{\eps,\delta}^{>>>}}{h}:} Writing 
 ${\Gamma}^{\gtrsim_\eps}_0(x,t,y,s)\eqdef (1-\kappa(2\eps^{-2}(t-s)))\bar \Gamma(x,t,y,s)$ and $\tilde{{\Gamma}}^{>_\eps}_{0;\delta}$ for the kernel obtained as in 
\eqref{eq:tilde_kernel} with $\Gamma^\eps$ replaced by $\mathbf{1}_{t-s>\eps^2}\bar{\Gamma}(x,t,y,s)$ we find that
\begin{align*}
&h_{\eps,\delta}^{\gtrsim>>}(x,t ) \\
&= \int_{\mathbb{R}^{3+1}} f^\eps(y,s) \kappa(t-s) \Big[\Gamma^{\gtrsim}_\eps (x,t; y,s) \left(\tilde{\Gamma}^{>}_{\eps;\delta}(x,t; y,s) \right)^2
-{\Gamma}^{{\gtrsim}_\eps}_0 (x,t; y,s) \left(\tilde{{\Gamma}}^{>_\eps}_{0;\delta}(x,t; y,s) \right)^2\Big]
 \ dy  ds 
 \\
 &\qquad+\int_{\mathbb{R}^{3+1}} [f^\eps(y,s)- \bar{f}] \kappa(t-s)
{\Gamma}^{\gtrsim_\eps}_0 (x,t; y,s) \left(\tilde{{\Gamma}}^{>_\eps}_{0;\delta}(x,t; y,s) \right)^2
 \ dy  ds \ .
\end{align*}
The first term is thus bounded using Corollary~\ref{cor:pointwise}, while we use Lemma~\ref{lem:appendix} for the second summand.

\paragraph{Estimate on \TitleEquation{\sum_{a}h^{ab\neq}}{sum hab} and \TitleEquation{\sum_{a}h^{a\neq c}}{sum hac}.}
Note that for $\bar t = t + r$
\begin{equs}
0 \leq \tilde{K}^{\neq}_{\eps,\delta} (\bar{x},\bar{t}, x,t)&=
 \int_{\substack{\delta^2<t-s<\eps^2\\
\eps^2<\bar{t}-s
 }  
 }  \int_{\mathbb{R}^3} \Gamma_\eps (x,t;\eta, s) 
 \Gamma_\eps ( \bar{x},\bar{t};\eta, s) 
\,d\eta \kappa(\bar{t}-s-\delta) ds
\\
&\lesssim
 \int_{\substack{0<\tau<\eps^2\\
\eps^2-r<\tau<3-r
 }   }   \frac{1}{(2\tau+r)^{d/2}} \exp \Big(-\kappa\frac{|x-\bar{x}|^2}{2\tau+r} \Big) ds  \label{eq:veryloc}
\end{equs}
and that the integral is empty if $r<0$.
%
 We estimate
\begin{align*}
\sum_{a}h^{ab\neq}_{\eps,\delta}(x,t)&=\int_{\mathbb{R}^3\times[0,\eps^2/2]} \kappa(t-s) \Gamma (x,t; y,s) f^\eps(y,s)
\tilde{\Gamma}^{b}_{\eps;\delta}(x,t; y,s) \tilde{\Gamma}^{\neq}_{\eps;\delta}(x,t; y,s) 
 \ dy  ds\\
&+ \int_{\mathbb{R}^3\times(\eps^2/2,3)} \kappa(t-s) \Gamma (x,t; y,s) f^\eps(y,s)\tilde{\Gamma}^{b}_{\eps;\delta}(x,t; y,s) \tilde{\Gamma}^{\neq}_{\eps;\delta}(x,t; y,s)  \ dy  ds \ .
\end{align*} 
For the first term note that
when  $r<\eps^2/2$ \eqref{eq:veryloc} is bounded by $(\eps^2+r)^{-1/2}$
and thus the first summand is bounded by 
$$\int_0^{\eps^2/2} \int_{\mathbb{R}^{3}} \frac{1}{(\eps^2+r)^{1/2}r^{d/2+1/2}} \exp \Big(-\kappa\frac{|x-\bar{x}|^2}{r} \Big) dx dr \lesssim  
\int_0^{\eps^2/2}\frac{1}{(\eps^2+r)^{1/2}r^{1/2}}  dr
 \lesssim 1\ .
$$
The second integral, note that 
for $r>\eps^2/2$ \eqref{eq:veryloc} can be crudely bounded by $\eps^2/r^{d/2}$.
Therefore,
$$\int_{\eps^2/2}^3 \int_{\mathbb{R}^{3}} \frac{\eps^4}{r^{3d/2}}  \exp \Big(-\kappa\frac{|x-\bar{x}|^2}{r} \Big) dx dr \lesssim  
\int_{\eps^2/2}^3\frac{\eps^4}{r^{d}}  dr\lesssim 1 \ .
$$

\paragraph{Estimate on \TitleEquation{\sum_{a} h^{a ><}}{sum h}:} Observe the following inequalities 
\begin{align*}
0\leq& \tilde{K}^{<}_{\eps,\delta} (x,t,\bar{x},\bar{t})=
 \int_{\substack{\delta^2<t-s<\eps^2\\
\delta^2<\bar{t}-s<\eps^2
 }  }  
 \int_{\mathbb{R}^3} \Gamma_\eps (x,t;\eta, s) 
 \Gamma_\eps ( \bar{x},\bar{t};\eta, s) 
\,d\eta  ds \lesssim {(t-\bar{t})^{-1/2}}\\
0 \leq & \tilde{K}^{>}_{\eps,\delta} (x,t,\bar{x},\bar{t})=
 \int_{\substack{\eps^2<t-s<3\\
\eps^2<\bar{t}-s<3
 }  }  
 \int_{\mathbb{R}^3} \Gamma_\eps (x,t;\eta, s) 
 \Gamma_\eps ( \bar{x},\bar{t};\eta, s)  
\,d\eta  ds \lesssim \eps^{-1}\ ,
\end{align*}
%
thus
$$\sum_{a} h^{a ><} \leq \int_{|t-\bar{t}|<\eps} \frac{\eps^{-1}}{(t-\bar{t})^{1/2+d/2}} \Big(-\kappa\frac{|x-\bar{x}|^2}{t-\bar{t}} \Big) dx dt \lesssim \int_{|t-\bar{t}|<\eps} \frac{\eps^{-1}}{(t-\bar{t})^{1/2}}dt \lesssim 1$$ 
\end{proof}

\subsubsection{Counterterms for Theorem~\ref{thm:phi4 restricted}}\label{sec:count_phi_rest}
Define 
\begin{align*}
&F^{\<2>,<}_{\eps,\delta}(z)=  \E[(\Pi^{\eps,\delta}\<1>)^2] -
	\frac{\bar{\alpha}^{\<2>,<}_{\eps,\delta}}{\det(\bar{A})}
- c^{\<2>,<}_{\eps,\delta} \ ,\\
&F^{\<22>,<}_{\eps,\delta}(z)= \sum_{
I'}
\E [\bfPi^{(\eps,\delta)} \<2> I' (\pmb{f}\<2>)(x,t)]- \E [\bfPi^{(\eps,\delta)} \<2> (x,t) ]\E [\bfPi^{(\eps,\delta)} I' (\pmb{f}\<2>)  
(x,t)] -   \frac{2\bar{\alpha}^{\<22>}_{\eps,\delta} \bar{f} }{\det(\bar{A})}- c^{\<22>,<}_{\eps,\delta} \ ,
\end{align*}
where
\begin{enumerate}
\item $\bar{\alpha}^{\<2>,<}_{\eps,\delta}\eqdef(4\pi)^{-d/2} \eps^{-1}\int_{(\delta/\eps)^2}^\infty s^{-d/2} ds $,
\item $ c^{\<2>,<}_{\eps,\delta}\eqdef\int_{[0,1]^3} h^{\<2>,<}_{\eps,\delta}$ for $h^{\<2>,<}_{\eps,\delta}\eqdef   \E[(\Pi^{\eps,\delta}\<1>)^2] -
	\frac{\bar{\alpha}^{\<2>,<}_{\eps,\delta}}{\det(\bar{A})}$,
\item let $\bar{\alpha}^{\<22>,<}_{\eps,\delta}$ be
  as in Item~\ref{itemmmm2} above {Lemma~\ref{lem:finite wiener chaos remainder12}},
\item $ c^{\<22>,<}_{\eps,\delta}\eqdef\int_{[0,1]^3} h^{\<22>,<}_{\eps,\delta}$ for $h^{\<22>,<}_{\eps,\delta}\eqdef   \sum_{
I'}
 \E[\Pi^{\eps,\delta}\<2> I' (\pmb{f}\<2>)
](z) -   \frac{\bar{\alpha}^{\<22>}_{\eps,\delta} \bar{f} }{\det(\bar{A})} $.
\end{enumerate}
\begin{lemma}
For any $C>0$ and $\tau \in \{\<2>,\<22>\}$ the functions $F^{\tau,<}_{\eps,\delta}: \mathbb{R}^3\to \mathbb{R}$ satisfy the following
\begin{itemize}
\item $\sup_{\eps,\delta\in (0,1]:\eps\leq C\delta^2 } \|F^{\tau,<}_{\eps,\delta}\|_{L^\infty}<\infty$,
\item $F^{\tau,<}_{\eps,\delta}$ is $\eps\mathbb{Z}^2\times \eps^2 \mathbb{Z}$ periodic and has mean $0$,
\item $\lim_{\eps\to  0} \|F^{\tau,<}_{\eps,\delta}\|_{L^\infty}=0$  for each $\delta>0$. 
\end{itemize}
\end{lemma}
\begin{proof}
We first check the claim for $\tau= \<2>$. Note that 
$$h^{\<2>,<}_{\eps,\delta}= h^{\<2>}_{\eps,\delta} +
\frac{\alpha^{\<2>}_{\eps,\delta}}{\sqrt{\det(A_s^\eps)}} 
+\frac{1}{\eps} R^{\<2>}_{(\delta/\eps)\vee 1} \circ \mathcal{S}^\eps  - 
\frac{\bar{\alpha}^{\<2>,<}_{\eps,\delta}}{\sqrt{\det(\bar{A})}}.$$
Since $\sup_{\eps,\delta\in (0,1]:\eps\leq C\delta^2 } {\alpha^{\<2>}_{\eps,\delta}}<\infty$ and $\lim_{\eps\to 0} \alpha^{\<2>}_{\eps,\delta}=0$ it remains to note that
\begin{align*}
\Big|\frac{1}{\eps} R^{\<2>}_{(\delta/\eps)\vee 1} \circ \mathcal{S}^\eps  - 
\frac{\bar{\alpha}^{\<2>,<}_{\eps,\delta}}{\det(\bar{A})}\Big| 
&\leq \frac{1}{\eps} \int_{\mathbb{R}^3\times[(\delta/\eps)^2\vee 1,\infty)} |\Gamma_1^2(z;y,s)-\bar{\Gamma}^2(z;y,s)|  dyds\\
&\lesssim\frac{1}{\eps} \int_{(\delta/\eps)^2\vee 1}^\infty {s^{-d/2-1/2} } \lesssim \frac{1}{\eps}(\eps/\delta)^2\lesssim \frac{\eps}{\delta^2} \ .
\end{align*}
For $\tau= \<22>$ the claim is simpler since 
$h^{\<22>,<}_{\eps,\delta}= h^{\<22>}_{\eps,\delta} + \frac{\bar{\alpha}^{\<22>}_{\eps,\delta}}{{\det(\bar{A})}} \ .$
\end{proof}

\subsubsection{Proof of Theorems~\ref{thm:phi4} and ~\ref{thm:phi4 restricted}}
Finally, we conclude the main results about the $\Phi^4$ equation by combining what has been done so far. 

\begin{proof}[Proof of the results on the $\Phi^4_3$ equation]
Continuity of the solution map on $\big((0,1]\times [0,1]\big) \cap \bar{\square}$ follows directly from the results in \cite{Sin23}. 
We choose the homogeneity assignment for the regularity structure constructed in the beginning of Section~\ref{sec:application phi} by imposing that
 $\kappa>0$, $\beta<2$ $ \zeta<-1/2$ are such that
$|\kappa| +|\zeta-1/2|\leq 1/1000$ and  $5/100<|\beta-2|<6/100$ .

Thus we want to apply Theorem~\ref{thm:fixed point} to \eqref{eq:abstract_phi4}, where we shall solve the equation in $\mathcal{D}^{\gamma,\eta}$ for the choice of $\eta= 1/100$ and 
$|2\zeta -\kappa|<\gamma<L,R <1+ 1/100$. Thus, we check the remaining assumptions.
\begin{enumerate}
\item Since we work with uniformly bounded initial conditions for the remainder equation, we can lift it by Corollary~\ref{cor:initial_cond}.
\item Convergence of the involved kernels follows from Lemma~\ref{lem:convergence of kernels easy terms} and Proposition~\ref{prop:vanishing of G kernel}.
\item Convergence of models follows from the stochastic estimates in Section~\ref{sec:convergence of ren. models} and the bounds in Sections~\ref{sec:count_phi_full} and \ref{sec:count_phi_rest} as in the proof of Theorem~\ref{thm:g-pam}.
\item Next we need to check that $p_\eps$ defined in \eqref{eq:modelled dist} almost surely belongs to $\mathcal{D}^{\gamma, \eta}$.
We observe that this follows for all terms except $\mathcal{G}_\epsilon \big( \mathbf{R}_+\pmb{f}\<3>\big)$ directly from the Schauder estimates for singular modelled distributions together with the the model bounds on $\tilde{\Phi}_j \<3>$. To conclude the same for
$\mathcal{G}_\epsilon \big( \mathbf{R}_+  \<3>    \big) $, we additionally observe that 
 defining the model on $I^{G_\eps} (\<3> )$ by stochastic estimates and Kolmogorov only requires bounds on $\vertiii{G_\eps}_{\beta; 1/2, 0}$.
 Lastly, convergence in $\eps\to 0$ follows similarly.
%
%
\end{enumerate}
Finally, the form and asymptotic behaviour of counterterms is checked by straightforward (but slightly tedious) computations, which concludes the proof. 
\end{proof}

\begin{appendix}

\section{Periodic Functions and Renormalised Kernels}\label{ap:A}

\begin{lemma}\label{lem:appendix}
Let $p\in [1,\infty]$ and $\kappa\in (0,1)$. It holds that for any $f\in \L^p(\mathbb{R}^{d+1})$ that is $(\eps\mathbb{Z})^{d}\times (\eps^2\mathbb{Z})$ periodic 
\begin{equ}
\left\|f- \int_{[0,1]^{d+1}} f \right\|_{C_\fraks^{-\kappa-1/p}}\leq \eps^{\kappa} \|f\|_{L^p([0,1]^{d+1})}\ .
\end{equ}
\end{lemma}
\begin{proof}
Without loss of generality assume that $f$ has vanishing mean, i.e.\
$\int_{[0,1]^{d+1}} f$. We shall show that
\begin{equation}\label{eq:periodic_improved convergence}
\left|\int_{\mathbb{R}^{d+1}} f(z) \phi^{\lambda}_\star (z) \,dz \right|\lesssim \|f\|_{L^p} \eps^{\kappa} \lambda^{-\kappa-1/p}
\end{equation}
uniformly over $\phi\in C_c(B_1)$ satisfying $\|\phi\|_{C^\kappa_\fraks }<1$ and $\star\in \mathbb{R}^{d+1}$. 
%
%
By translation, it suffices to consider only the case $\star=0$. 
Note that for $\eps>\lambda$ 
we simply use H\"older's inequality with $1/p+1/q=1$
$$
\left|\int_{\mathbb{R}^{d+1}} f(z) \phi^{\lambda}_\star (z) \,dz \right|\lesssim \|f\|_{L^p} \|\phi^\lambda\|_{L^q}   
\leq \|f\|_{L^p} \eps^{\kappa} \lambda^{-\kappa-1/p}\ .
$$
In the case $\eps<\lambda$, let $F(z)= f(\mathcal{S}^{\eps}z)$. Then,
\begin{align*}
&\int_{\mathbb{R}^{d+1}} f(z)  \phi^{\lambda} (z)\,dz = \eps^{|\fraks|} \int_{\mathbb{R}^{d+1}} F(\mathcal{S}^{\lambda^{-1}} z) \phi (\mathcal{S}^{\eps^{-1}} z) \\
&= \eps^{|\fraks|} \sum_{h\in \mathbb{Z}^{d+1}} \int_{\mathcal{S}^\lambda ([0,1]^{d+1}+ h)} F(\mathcal{S}^{\lambda^{-1}}z)\left( \phi (\mathcal{S}^{\eps^{-1}} z) - \phi (\mathcal{S}^{\eps^{-1}} h)\right)\ .
\end{align*}
Thus, $\left|\int_{\mathbb{R}^{d+1}} f(z) \phi^{\lambda}_x (z) \right|$ is bounded by
\begin{align*}
\eps^{|\fraks|} \sum_{h\in \mathbb{Z}^{d+1}} \| F(\mathcal{S}^{\lambda^{-1}} \cdot)\|_{L^p(\mathcal{S}^\lambda ([0,1]^{d+1}+ h))} \|  \phi (\mathcal{S}^{\eps^{-1}} \cdot ) - \phi (\mathcal{S}^{\eps^{-1}} h)\|_{L^q( \phi (\mathcal{S}^{\eps^{-1}} z) - \phi (\mathcal{S}^{\eps^{-1}} h))}\ .
\end{align*}
Observing that $\| F(\mathcal{S}^{\lambda^{-1}} \cdot)\|_{L^p(\mathcal{S}^\lambda ([0,1]^{d+1}+ h))} = \lambda^{- |\fraks|/p}
\| F\|_{L^p ([0,1]^{d+1})}= \lambda^{- |\fraks|/p}
\| f\|_{L^p ([0,1]^{d+1})}$
and that 
$\| \phi (\mathcal{S}^{\eps^{-1}} z) - \phi (\mathcal{S}^{\eps^{-1}} h)\|_{L^{q} (\mathcal{S}^\lambda ([0,1]^{d+1}+ h))   } \leq
\frac{\eps^\kappa}{\lambda^\kappa} \|\phi\|_{C^\kappa}\| 1 \|_{L^{q} (\mathcal{S}^\lambda ([0,1]^{d+1}+ h))   } 
\leq \frac{\eps^\kappa}{\lambda^\kappa} \|\phi\|_{C^\kappa} \lambda^{-|\fraks|/q}$, we conclude that
\begin{align*}
\left|\int_{\mathbb{R}^{d+1}} f(z) \phi^{\lambda}_x (z) \right| 
& \leq \eps^{|\fraks|} 
\lambda^{- |\fraks|/p}
\| f\|_{L^p ([0,1]^{d+1})}
\frac{\eps^\kappa}{\lambda^\kappa} \|\phi\|_{C^\kappa} \lambda^{-|\fraks|/q}
N_{\lambda, \eps}\\
& \leq 
\frac{\eps^\kappa}{\lambda^\kappa} 
\| f\|_{L^p ([0,1]^{d+1})}
\|\phi\|_{C^\kappa} \ ,
\end{align*}
where
$N_{\lambda, \eps}:= |\{h\in \mathbb{Z}^d \ : \ \supp(\phi(\eps \ \cdot \ )) \cap \mathcal{S}^\lambda([0,1]^{d+1}+h)\neq \emptyset   \}|\lesssim \eps^{-|\fraks|}\lambda^{|\fraks|}$.
%
%
\end{proof}

Let the kernel
$K: \mathbb{R}^{d+1}\times \mathbb{R}^{d+1}\setminus \triangle \to \mathbb{R}$ be supported on a set of finite distance from the diagonal $\triangle$. 
Write using the notation of \eqref{zwischendefinition} 
\begin{equ}
 \|{K}\|_{\beta;L,R}:= \inf \vertii{\{K_n\}_n }_{\beta;L,R}\ ,
\end{equ}
where the infimum is taken over all decomposition $K(z,z')= \sum_{n\geq 0} K_n(z,z')$ such that each $K_n$ for $n\geq 1$ is supported on 
$ \{(z,z')\in (\mathbb{R}^{d+1})^{\times 2}\ : \ |z-z'|_\fraks \leq 2^{-n} \} $.
For such a kernel $K$ with $\|{K}\|_{\beta;L,0}<\infty$ for 
 $\beta\in (-1,0]$, set 
\begin{equ}\label{eq:app_ren}
\big(\mathcal{R}(K)\phi\big)(z):= \int_{\mathbb{R}^{d+1}}K(z,w) [ \phi(z)-\phi(w)] dw \ .
\end{equ}
 Accordingly, for a second kernel $G=\sum_n G_n$ decomposed as described above, we set  
$$\mathcal{R}(K)\star G(z,w):=\sum_{n}\big(\mathcal{R}(K) G_n(\cdot, w)\big)(z)\ .$$
Given a bounded function $f\in L^\infty$ we write $\big(K\cdot f\big)(z,w):= K(z,w) \cdot f(w)$.

\begin{lemma}\label{lem:ap_osc_kernel}
Let $\beta_1\in (-1,0)$, $\beta_2\in (0,|\fraks|)$ be such that $\beta_1+\beta_2\in (0,|\fraks|)$,  let $f\in L^\infty$ and let 
 $K_1$ and $K_2$ be kernels as above. 
Then,  whenever $1>L'>|\beta_1|+L$ it holds that 
\begin{equ}\label{eq:ap_unif}
\|\mathcal{R}({K}_1\cdot f^\eps)\star K_2\|_{\beta_1+\beta_2,L,0}\lesssim \|f\|_{L^\infty} \|{K_1}\|_{\beta_1;L,0} \|{K_2}\|_{\beta_2;L',0}\ .
\end{equ}
uniformly over $\eps\in (0,1]$.
Furthermore, assume that $f$ is $\mathbb{Z}^{d+1}$ periodic and has vanishing mean. Then, for $\kappa<R_1\wedge L_2\wedge 1$ and $ 0\leq L<L_2-|\beta_1|$ satisfying
$\beta_1+\beta_2>0$ and
 $\beta_1+\beta_2+\kappa+ L<|\fraks|$
\begin{equ}\label{eq:ap_conv}
\|\mathcal{R}({K}_1\cdot f^\eps)\star K_2\|_{\beta,L,0}\lesssim \eps^\kappa \|f\|_{L^\infty} \|{K_1}\|_{\beta_1;L_1,R_1} \|{K_2}\|_{\beta_2;L_2,0}\ .
\end{equ}
uniformly over $\eps\in (0,1]$.
\end{lemma}
\begin{proof}
Observing that 
$\|{K}_1\cdot f^\eps\|_{\beta_1,L,0} \leq \|f\|_{L^\infty} \|{K}_1\|_{\beta_1,L,0} $, we assume without loss of generality that $f=1$ when checking \eqref{eq:ap_unif}.
Indeed,
$$\big(\mathcal{R}(K_1)\star K_{2,n}\big)(x,z)= \int K_1(x,y)[K_{2,n}(y,z)-K_{2,n}(x,z)]dy \ .$$
We observe that for $|x-z|>2^{-n+1}$, the the renormalisation has no effect and we bound
$$
\int K_1(x,y)[K_{2,n}(y,z)-K_{2,n}(x,z)]dy= \int_{|y-z|<2^{-n}} K_1(x,y)K_{2,n}(y,z)dy \lesssim|x-z|^{-|\fraks|+\beta_1} 2^{-n\beta_2} \ .
$$
For  $|x-z|\leq 2^{-n+1}$, we bound using that $L'>|\beta_1|$
\begin{equs}
& \int K_1(x,y)[K_{2,n}(y,z)-K_{2,n}(x,z)]dy \\
& \lesssim \int_{|y-z|<2^{-n+1}} K_1(x,y) |x-y|^{L'} dy\cdot  2^{n(|\fraks|-\beta_2 + L')}
+K_{2,n}(x,z) \int_{|y-z|\geq 2^{-n+1}} K_1(x,y) dy\\
&\lesssim 2^{n(|\fraks|-\beta_1-\beta_2 )} \ .\label{bound_dif_eps ..}
\end{equs}
Thus summing over $n\in \mathbb{N}$ we conclude that 
$\big|\big(\mathcal{R}({K}_1)\star K_2\big)(x,z)\big|\lesssim |x-z|^{-|\fraks|+\beta_1+\beta_2}$. 
To bound the H\"older norms, we proceed exactly the same way. The main difference is that we require $L'>|\beta_1|+L$ in the case $|x-z|\leq 2^{-n+1}$.

Next, we check \eqref{eq:ap_conv}. Proceeding similarly to above we note that 
 for $|x-z|>2^{-n+1}$, 
$$
\int K_1(x,y)f^{\eps}(y)[K_{2,n}(y,z)-K_{2,n}(x,z)]dy=
  \int \rho^{2^{-n}}_z(y) f^\eps(y) dy\ ,
$$
where we have set $\rho^{2^{-n}}_z(y):=K_1(x,y) K_{2,n}(y,z)$. Using the notation of Remark~\ref{rem:useful notation norm}, we note that 
$$\|\rho^{2^{-n}}(y)\|_{\mathfrak{B}_{R_1\wedge L_2}^{2^{-n}}}\lesssim |x-z|^{-|\fraks|+\beta_1}2^{-n\beta_2}$$
uniformly in $n\in \mathbb{N}, x\in \mathbb{R}^{d}$ such that $|x-z|> 2^{-n+1}$.
Thus, by Lemma~\ref{ap:A}  $\kappa<R_1\wedge L_2 \wedge 1$ 
$$
\big| \int \rho^{2^{-n}}_z(y) f^\eps(y) dy\big| \lesssim 
\eps^{\kappa} |x-y|^{-|\fraks|+\beta_1} 2^{n(\kappa-\beta_2)} \ .
$$
We turn to the case 
  $|x-z|\leq 2^{-n+1}<\eps$. Then, it follows from the bound \eqref{bound_dif_eps ..} that
\begin{align*}
& \int K_1(x,y)f^\eps(y) [K_{2,n}(y,z)-K_{2,n}(x,z)]dy \lesssim 2^{n(|\fraks|-\beta_1-\beta_2 )} \lesssim \eps^{\kappa} 2^{n(|\fraks|-\beta_1-\beta_2 +\kappa )}\ .
\end{align*}
Finally, consider the case $|x-z|\vee \eps\leq 2^{-n}$. Fix $\phi\in C_c^\infty(B_5)$ such that $\phi|_{B_4}=1$ and write $\phi^n=\phi\circ \mathcal{S}^{2^{-n}}$. Then, we can write
\begin{equs}
 &\int K_1(x,y)f^\eps(y) [K_{2,n}(y,z)-K_{2,n}(x,z)]dy\\
 &=  \int \phi^n(x-y) K_1(x,y)f^\eps(y) [K_{2,n}(y,z)-K_{2,n}(x,z)]dy
 -K_{2,n}(x,z)  \int \big(1-\phi^n(x-y)\big) K_1(x,y)f^\eps(y) dy\\
 &= \int\rho^{2^{-n}}_x f^\eps(y)dy - \sum_{m\leq n} K_{2,n}(x,z)  \int \big(1-\phi^n(x-y)\big) K_{1,m}(x,y)f^\eps(y) dy \ , \label{eq:lasteqref}
\end{equs}
where 
$\rho^{2^{-n}}_x(y)= \phi^n(x-y) K_1(x,y) [K_{2,n}(y,z)-K_{2,n}(x,z)]$
is seen to satisfy 
$\|\rho^{2^{-n}}\|_{\mathfrak{B}_{R_1\wedge L_2}^{2^{-n}}}\lesssim 2^{n(|\fraks|-\beta_1-\beta_2)}\  $.
Thus, for $\kappa<R_1\wedge L_2\wedge 1$ we find that 
$$\big| \int\rho^{2^{-n}}_x f^\eps(y)dy\big| \lesssim \eps^{\kappa} 2^{n(|\fraks|-\beta_1-\beta_2+\kappa)}$$
The estimate the second term of \eqref{eq:lasteqref} note that
$\rho^{2^{-m}}_x(y):=\big(1-\phi^n(x-y)\big) K_{1,m}(x,y)$ satisfies
$\|\rho^{2^{-m}}\|_{\mathfrak{B}_{R_1}^{2^{-m}}}\lesssim 2^{-m\beta_1} $ and therefore, for $\kappa< R_1\wedge 1$
$$\Big| \sum_{m\leq n} K_{2,n}(x,z)  \int \big(1-\phi^n(x-y)\big) K_{1,m}(x,y)f^\eps(y) dy\Big|\lesssim \eps^{\kappa} 2^{n(|\fraks|-\beta_1-\beta_2+\kappa)} \ .$$
Thus the bound \eqref{eq:ap_conv} for $L=0$ follows by summing over $n\in \mathbb{N}$. 
To obtain the bound for $L>0$ one proceeds analogously.
\end{proof}
%

\end{appendix}

\bibliographystyle{Martin}
\bibliography{./Periodic.bib}{}

\end{document}